\documentstyle{amsart}

\newtheorem{Theorem}{Theorem}[section]

\newtheorem{Lemma}[Theorem]{Lemma}

\newtheorem{Remark}[Theorem]{Remark}

\newtheorem{Definition}[Theorem]{Definition}

\begin{document}

\title{Strong Toroidalization of dominant morphisms of 3-folds}

\author{Steven Dale Cutkosky}
\thanks{Research    partially supported by NSF}

\maketitle

\section{Introduction}

Suppose that $f:X\rightarrow Y$ is a dominant morphism of algebraic varieties, over a field ${\bf k}$ of characteristic zero.
If $X$ and $Y$ are nonsingular,
$f:X\rightarrow Y$ is toroidal if there are simple normal crossing divisors $D_X$ on $X$ and $D_Y$ on $Y$ such that
$f^{-1}(D_Y)=D_X$, and $f$ is  locally given by monomials in appropriate etale local parameters on $X$.
The precise definition of this concept is in \cite{AK}, \cite{KKMS} and Definition \ref{Def247} of this paper. The problem of toroidalization is to determine, given a dominant morphism $f:X\rightarrow Y$,
 if there exists a commutative diagram 
 \begin{equation}\label{eq393}
 \begin{array}{rcl}
 X_1&\stackrel{f_1}{\rightarrow}&Y_1\\
 \Phi\downarrow&&\downarrow\Psi\\
 X&\stackrel{f}{\rightarrow}&Y
 \end{array}
 \end{equation}
 such that $\Phi$ and $\Psi$ are products of blow ups of nonsingular subvarieties, $X_1$ and $Y_1$ are nonsingular,
 and there exist simple normal crossing divisors $D_{Y_1}$ on $Y_1$ and $D_{X_1}=f^{-1}(D_{Y_1})$ on $X_1$ such that $f_1$ is toroidal
 (with respect to $D_{X_1}$ and $D_{Y_1}$).     This  is stated in Problem 6.2.1. of \cite{AKMW}.
 
 A stronger form of toroidalization is also asked for in  \cite{AKMW}, which we will call strong toroidalization.
 Suppose that $f:X\rightarrow Y$ is a dominant morphism of nonsingular projective varieties, $D_Y$ is a SNC divisor on $Y$ and $D_X=f^{-1}(D_Y)$ is a 
 SNC divisor on $X$ such that the locus $\text{sing}(f)$ where the morphism $f$ is not smooth is contained in $D_X$.
 The problem of strong toroidalization is to determine if there exists a commutative diagram (\ref{eq393})
 such that $\Phi$ and $\Psi$ are products of blow ups of nonsingular centers which are supported in the preimages of
 $D_X$ and $D_Y$ respectively, and make SNCs with the respective preimages of $D_X$ and $D_Y$, and $f_1$ 
 is toroidal with respect to $D_{Y_1}=\Psi^{-1}(D_Y)$ and $D_{X_1}=\Phi^{-1}(D_X)$.

 Toroidalization, and related concepts, have been considered earlier in different contexts,
 mostly for morphisms of surfaces. Strong torodialization is the strongest structure theorem which could be true for
 general morphisms.  The concept of torodialization fails completely in positive characteristic. A simple example in characterisitc $p>0$ is obtained from the map of curves $s=t^p(t^p+1)$.
 
 In the case when $Y$ is a curve, toroidalization follows from embedded resolution of singularities (\cite{H}).
 When $X$ and $Y$ are surfaces, there are several proofs in print (\cite{AkK}, Corollary 6.2.3 \cite{AKMW}, \cite{CP}, \cite{Mat}).  They all make use of special properties of the
 birational geometry of surfaces.  An outline of proofs of the above cases can be found in the introduction to
 \cite{C3}.

 In \cite{C3}, strong toroidalization  is solved in the case when $X$ is a 3-fold and $Y$ is a surface,
 In Theorem 0.1 of \cite{C5} we prove toroidalization of birational morphisms of 3-folds. 
 In Theorem 1.1 of \cite{C7}, we prove strong toroidalization of birational morphisms of 3-folds.
 In this paper, we prove  strong  toroidalization for dominant morphisms of 3-folds.

\begin{Theorem}\label{Theorem3}  Suppose that $f:X\rightarrow Y$ is a dominant morphism of nonsingular projective 3-folds  over an
algebraically closed field ${\bf k}$ of characteristic 0. Further suppose that there is a simple normal crossings (SNC) divisor $D_Y$ on $Y$
such that  $D_X=f^{-1}(D_Y)$ is a SNC divisor which contains the non smooth locus of the map $f$.
 Then there exists a commutative diagram of morphisms
$$
\begin{array}{rll}
X_1&\stackrel{f_1}{\rightarrow}&Y_1\\
\Phi\downarrow&&\downarrow\Psi\\
X&\stackrel{f}{\rightarrow}&Y
\end{array}
$$
where $\Phi,\Psi$ are products of possible blow ups for the preimages of $D_X$ and $D_Y$ respectively,
  and $f_1$ is toroidal with respect to $D_{Y_1}=\Psi^{-1}(D_Y)$ and $D_{X_1}=\Phi^{-1}(D_X)$.
\end{Theorem}

A 3-fold over a field ${\bf k}$ is an abstract variety over ${\bf k}$ of dimension 3.

A possible blow up on a nonsingular 3-fold with toroidal structure is the blow up of a point or a nonsingular curve contained in the toroidal structure which makes SNCs with the toroidal structure.

We deduce  Theorem \ref{Theorem3} from the following strong toroidalization theorem for morphisms of (possibly singular) varieties.

\begin{Theorem}\label{Theorem1} Suppose that $f:X\rightarrow Y$ is a dominant morphism of 3-folds over an
algebraically closed field ${\bf k}$ of characteristic 0. Further suppose that there is an equidimensional codimension 1 reduced subscheme  $D_Y$ of $Y$
such that $D_Y$ contains the singular locus of $Y$, and $D_X=f^{-1}(D_Y)$ contains the non smooth locus of the map $f$.
 Then there exists a commutative diagram of morphisms
$$
\begin{array}{rll}
X_1&\stackrel{f_1}{\rightarrow}&Y_1\\
\Phi\downarrow&&\downarrow\Psi\\
X&\stackrel{f}{\rightarrow}&Y
\end{array}
$$
where $\Phi,\Psi$ are products of blow ups of nonsingular curves and points supported above $D_X$ and $D_Y$ respectively,
 $D_{Y_1}=\Psi^{-1}(D_Y)$ is a simple normal crossings  divisor on $Y_1$,
$D_{X_1}=f_1^{-1}(D_{Y_1})$ is a simple normal crossings  divisor on $X_1$ and $f_1$ is toroidal with respect to $D_{Y_1}$ and $D_{X_1}$.
\end{Theorem}

If we relax some of the restrictions in the definition of toroidalization, there are other constructions
producing a toroidal morphism $f_1$, which 
are valid for arbitrary dimensions of $X$ and $Y$.
In \cite{AK} it is shown that a diagram (\ref{eq393}) can be constructed where $\Phi$ is weakened to being a
modification (an arbitrary birational morphism).  In \cite{C1}, \cite{C2} and \cite{C4}, it is shown that a diagram (\ref{eq393})
can be constructed where $\Phi$ and $\Psi$ are locally products of blow ups of nonsingular centers and $f_1$ is locally toroidal, but the morphisms $\Phi$, $\Psi$ and $f_1$
may not be separated.  This construction is obtained by patching local solutions, at least one of which contains the center of any given valuation.

\section{Notation}\label{SectionNotation}
Throughout this paper, ${\bf k}$ will be an algebraically closed field of characteristic zero. A curve, surface or 3-fold is
an abstract variety over ${\bf k}$ of respective dimension 1, 2 or 3.
If $X$ is a variety, and $p\in X$ is a nonsingular point, then regular parameters at $p$ are regular parameters in ${\cal O}_{X,p}$.
Formal regular parameters at $p$ are regular parameters in $\hat{\cal O}_{X,p}$.
 If $X$ is a variety and $V\subset X$ is a subvariety, then
${\cal I}_V\subset {\cal O}_X$ will denote the ideal sheaf of $V$.
If $V$ and $W$ are subvarieties of a variety $X$, we denote the scheme theoretic intersection 
$Y=\text{spec}({\cal O}_X/{\cal I}_V+{\cal I}_W)$ by $Y=V\cdot W$.

Let $f:X\rightarrow Y$ be a morphism of varieties. We will denote the singular locus of $f$ by $\text{sing}(f)$ the closed set of points $p\in X$ such that $f$ is not smooth at $p$. If $D$ is a Cartier divisor on $Y$, then $f^{-1}(D)$ will denote
the reduced divisor $f^*(D)_{red}$.

 Suppose that $a,b,c,d\in{\bf Q}$. Then we will write $(a,b)\le (c,d)$ if
$a\le b$ and $c\le d$.

A toroidal structure on a nonsingular variety $X$ is a simple normal crossing divisor (SNC divisor) $D_X$ on $X$.

We will say that a nonsingular curve $C$ which is a subvariety of a nonsingular 3-fold $X$ with toroidal structure
$D_X$ makes simple normal crossings (SNCs) with $D_X$ if for all $p\in C$, there exist regular parameters
$x,y,z$ at $p$ such that $x=y=0$ are local equations of $C$, and $xyz=0$ contains the support of $D_X$ at $p$.

Suppose that $X$ is a nonsingular 3-fold with toroidal structure $D_X$. If $p\in D_X$ is on the intersection of three components of $D_X$ then $p$ is called a 3-point. If $p\in D_X$
is on the intersection of two components of $D_X$ (and is not a 3-point) then $p$ is called a 2-point. If $p\in D_X$
is not a 2-point or a 3-point, then $p$ is called a 1-point. If $C$ is an irreducible component of the intersection of two
components of $D_X$, then $C$ is called a 2-curve. $\Sigma(X)$ will denote the closed locus of 2-curves on $X$.

A possible center on a nonsingular 3-fold $X$ with toroidal structure defined by a SNC divisor $D_X$, is a point
on $D_X$ or a nonsingular curve in $D_X$ which makes SNCs with $D_X$. A possible center on a nonsingular surface $S$
with toroidal structure defined by a SNC divisor $D_S$  is a point on $D_S$. We will also call the blow up of a possible center a possible blow up.

Observe that if $\Phi:X_1\rightarrow X$ is the blow up of a possible center, then $D_{X_1}=\Phi^{-1}(D_X)$ is a SNC
divisor on $X_1$. Thus $D_{X_1}$ defines a toroidal structure on $X_1$. All blow ups $\Phi:X_1\rightarrow X$
considered in this paper will be of possible centers, and we will impose the toroidal structure on $X_1$ defined by
$D_{X_1}=\Phi^{-1}(D_X)$.

By a general point $q$ of a variety $V$, we will mean a point $q$ which 
satisfies conditions which hold on some nontrivial open subset of $V$.
The exact open condition which we require will generally be clear from context.
By a general section of a coherent sheaf ${\cal F}$ on a projective variety $X$, we mean the section corresponding 
to a general point of the $k$-linear space $\Gamma(X,{\cal F})$.

If $X$ is a variety, ${\bf k}(X)$ will denote the function field of $X$. A 0-dimensional valuation $\nu$ of ${\bf k}(X)$ is
a valuation of ${\bf k}(X)$ such that $\bold k$ is contained in the valuation ring $V_{\nu}$ of $\nu$ and the residue field of
$V_{\nu}$ is $\bf k$.  If $X$ is a projective variety which is birationally equivalent to $X$, then there exists a
unique (closed) point $p_1\in X_1$ such that $V_{\nu}$ dominates ${\cal O}_{X_1,p_1}$. $p_1$ is called the center of $\nu$
on $X_1$. If $p\in X$ is a (closed) point, then there exists a 0-dimensional valuation $\nu$ of ${\bf k}(X)$ such that $V_{\nu}$ dominates ${\cal O}_{X,p}$
(Theorem 37, Section 16, Chapter VI \cite{ZS}). For $a_1,\ldots, a_n\in {\bf k}(X)$, $\nu(a_1),\ldots, \nu(a_n)$ are rationally dependent if there exist
$\alpha_1,\ldots, \alpha_n\in{\bf Z}$ which are not all zero, such that $\alpha_1\nu(a_1)+\cdots+\alpha_n\nu(a_n)=0$ (in the value group of $\nu$).
Otherwise, $\nu(a_1),\ldots, \nu(a_n)$ are rationally independent.

If $x_1,\ldots,x_n$ are indeterminates, and $M_i=\prod_{j=1}^n x_j^{a_{ij}}$ are monomials
for $1\le i\le m$, then we will denote $\text{rank}(M_1,\ldots,M_m)=\text{rank}(a_{ij})$.

\section{toroidal morphisms and prepared morphisms}

Suppose that $X$ is a nonsingular variety with toroidal structure $D_X$. We will say that an ideal sheaf ${\cal I}\subset {\cal O}_X$ is toroidal
if ${\cal I}$ is locally generated by monomials in local equations of components of $D_X$.

Suppose that $q\in X$. We say that $u,v,w$ are (formal) permissible
parameters at $q$ (for $D_X$) if $u,v,w$ are regular parameters in  $\hat {\cal O}_{X,q}$
such that $u=0$ is a local equation of $D_X$ at $q$ if $q$ is a 1-point, $uv=0$ is a local equation of $D_X$ at $q$ if $q$ is a 2-point and $uvw=0$ is a local equation of $D_X$ at $q$ if $q$ is a 3-point.
$u,v,w$ are algebraic permissible parameters if we further have that $u,v,w\in{\cal O}_{X,q}$.

\begin{Definition}\label{torf} Let $f:X\rightarrow Y$ be a dominant morphism of nonsingular 3-folds with toroidal structures $D_Y$ on $Y$ and $D_X=f^{-1}(D_Y)$ on $X$ such that $\text{sing}(f)\subset D_X$. Suppose that $u,v,w$ are (possibly formal) permissible parameters
at $q\in Y$. Then
$u,v$ are {\bf toroidal forms} at $p\in f^{-1}(q)$ if there exist permissible parameters $x,y,z$
in $\hat{\cal O}_{X,p}$ such that
\begin{enumerate}
\item[1.]  $q$ is a 2-point or a 3-point, $p$ is a 1-point and 
\begin{equation}\label{eqTF1}
u=x^a,
v=x^b(\alpha+y)
\end{equation}
where $0\ne \alpha\in  k$.
\item[2.] $q$ is 2-point or a 3-point, $p$ is a 2-point and 
\begin{equation}\label{eqTF21}
u=x^ay^b,
v=x^cy^d
\end{equation}
with $ad-bc\ne 0$.
\item[3.] $q$ is a 2-point or a 3-point, $p$ is a 2-point and 
\begin{equation}\label{eqTF22}
u=(x^ay^b)^k,
v=(x^ay^b)^t(\alpha+z)
\end{equation}
where $0\ne\alpha\in  k$, $a,b,k,t>0$ and $\text{gcd}(a,b)=1$.
\item[4.] $q$ is a 2-point or a 3-point, $p$ is a 3-point and 
\begin{equation}\label{eqTF3}
u=x^ay^bz^c,
v=x^dy^ez^f
\end{equation}
where 
$$
\text{rank}\left(\begin{matrix} a&b&c\\ d&e&f\end{matrix}\right)=2.
$$
\item[5.] $q$ is a 1-point, $p$ is a 1-point and 
\begin{equation}\label{eqTF01}
u=x^a,
v=y
\end{equation}
\item[6.] $q$ is a 1-point, $p$ is a 2-point and 
\begin{equation}\label{eqTF02}
u=(x^ay^b)^k,
v=z
\end{equation}
with $a,b,k>0$ and $\text{gcd}(a,b)=1$
\end{enumerate}
\end{Definition}

Regular parameters $x,y,z$ as in Definition \ref{torf} will be called permissible
parameters for $u,v,w$ at $p$.

\begin{Lemma}\label{Lemmatorf}
Suppose that $q\in Y$, $p\in f^{-1}(q)$ and $u,v,w$ are permissible
parameters at $q$ such that $u,v$ are toroidal forms at $p$. Then there exist
permissible parameters $x,y,z$ for $u,v,w$ at $p$ such that
an expression of
Definition \ref{torf} holds for $u$ and $v$, and one of the following respective forms for $w$ holds at $p$.
\begin{enumerate}
\item[1.]  $q$ is a 2-point or a 3-point, $p$ is a 1-point, $u,v$ satisfy (\ref{eqTF1})  and 
\begin{equation}\label{eqTF1w}
w=g(x,y)+ x^cz
\end{equation}
where $g$ is a series.
\item[2.] $q$ is 2-point or a 3-point, $p$ is a 2-point, $u,v$ satisfy (\ref{eqTF21})  
and 
\begin{equation}\label{eqTF21w}
w=g(x,y)+x^ey^fz
\end{equation}
where $g$ is a series.
\item[3.] $q$ is a 2-point or a 3-point, $p$ is a 2-point, $u,v$ satisfy (\ref{eqTF22}) and 
\begin{equation}\label{eqTF22w}
w=g(x^ay^b,z)+x^cy^d
\end{equation}
where $g$ is a series and $\text{rank}(u,x^cy^d)=2$.
\item[4.] $q$ is a 2-point or a 3-point, $p$ is a 3-point, $u,v$ satisfy (\ref{eqTF3}) and 
\begin{equation}\label{eqTF3w}
w=g(x,y,z)+N
\end{equation}
where $g$ is a series in monomials $M$ in $x,y,z$ such that $\text{rank}(u,v,M)=2$,
and $N$ is a monomial in $x,y,z$ such that $\text{rank}(u,v,N)=3$.
\item[5.] $q$ is a 1-point, $u,v$ satisfy (\ref{eqTF01})  and 
\begin{equation}\label{eqTF01w}
w=g(x,y)+ x^cz
\end{equation}
where $g$ is a series.
\item[6.] $q$ is a 1-point, $p$ is a 2-point, $u,v$ satisfy (\ref{eqTF02})  and 
\begin{equation}\label{eqTF02w}
w=g(x^ay^b,z)+x^cy^d
\end{equation}
where $g$ is a series and $\text{rank}(u,x^cy^d)=2$.
\end{enumerate}
\end{Lemma}

Lemma \ref{Lemmatorf} and its proof are Lemma 3.2 \cite{C5}.

\begin{Definition}\label{Def247} (\cite{KKMS}, \cite{AK})
A normal variety $\overline X$ with a SNC divisor $D_{\overline X}$ on $\overline X$ is called toroidal if for every point $p\in \overline X$ there exists an affine
toric variety $X_{\sigma}$, a point $p'\in X_{\sigma}$ and an isomorphism of $k$-algebras
$$
\hat{\cal O}_{\overline X,p}\cong \hat{\cal O}_{X_{\sigma},p'}
$$
such that the ideal of $D_{\overline X}$ corresponds to the ideal of $X_{\sigma}-T$ (where $T$ is the  torus in $X_{\sigma}$). Such a pair
$(X_{\sigma},p')$ is called a local model at $p\in \overline X$. $D_{\overline X}$ is called a toroidal structure on $\overline X$.

A dominant morphism $\Phi:\overline X\rightarrow \overline Y$ of toroidal varieties with SNC divisors $D_{\overline Y}$
on $\overline Y$ and  $D_{\overline X}=\Phi^{-1}(D_{\overline Y})$ on $\overline X$, 
is called toroidal at $p\in\overline X$, and we will say that $p$ is a toroidal point of $\Phi$ if with $q=\Phi(p)$, there exist local models
$(X_{\sigma},p')$ at $p$, $(Y_{\tau},q')$ at $q$ and a toric morphism $\Psi:X_{\sigma}\rightarrow Y_{\tau}$ such that the following
diagram commutes:
$$
\begin{array}{rll}
\hat{\cal O}_{\overline X,p}&\leftarrow &\hat{\cal O}_{X_{\sigma},p'}\\
\hat\Phi^*\uparrow&&\hat\Psi^*\uparrow\\
\hat{\cal O}_{\overline Y,q}&\leftarrow&\hat{\cal O}_{Y_{\tau},q'}.
\end{array}
$$
$\Phi:\overline X\rightarrow \overline Y$ is called toroidal (with respect to $D_{\overline Y}$ and $D_{\overline X}$) if $\Phi$ is toroidal at all $p\in \overline X$.
\end{Definition}
The following is the list of toroidal forms for a dominant morphism $f:X\rightarrow Y$
of nonsingular 3-folds with toroidal structure $D_Y$ and $D_X=f^{-1}(D_X)$. Suppose that $p\in D_X$, $q=f(p)\in D_Y$,
and $f$ is toroidal at $p$.
Then there exist permissible parameters $u,v,w$ at $q$ and permissible parameters $x,y,z$ for $u,v,w$ at $p$ such that
one of the following forms hold:
\begin{enumerate}
\item[1.] $p$ is a 3-point and $q$ is a  3-point,
$$
\begin{array}{ll}
u&=x^ay^bz^c\\
v&=x^dy^ez^f\\
w&=x^gy^hz^i,
\end{array}
$$
where $a,b,d,e,f,g,h,i\in{\bf N}$ and
$$
\text{Det}\left(\begin{array}{lll}
a&b&c\\
d&e&f\\
g&h&i
\end{array}\right)\ne 0.
$$
\item[2.] $p$ is a 2-point and $q$ is a 3-point,
$$
\begin{array}{ll}
u&=x^ay^b\\
v&=x^dy^e\\
w&=x^gy^h(z+\alpha)
\end{array}
$$
with $0\ne \alpha\in  k$ and $a,b,d,e,g,h\in{\bf N}$ satisfy $ae-bd\ne 0$.
\item[3.] $p$ is a 1-point and $q$ is a 3-point,
$$
\begin{array}{ll}
u&=x^a\\
v&=x^d(y+\alpha)\\
w&=x^g(z+\beta)
\end{array}
$$
with $0\ne \alpha,\beta\in  k$, $a,d,g>0$.
\item[4.] $p$ is a 2-point and $q$ is a 2-point,
$$
\begin{array}{ll}
u&=x^ay^b\\
v&=x^dy^e\\
w&=z
\end{array}
$$
with $ae-bd\ne 0$
\item[5.] $p$ is a 1-point and $q$ is a 2-point,
$$
\begin{array}{ll}
u&=x^a\\
v&=x^d(y+\alpha)\\
w&=z
\end{array}
$$
with $0\ne \alpha\in  k$, $a,d>0$.
\item[6.] $p$ is a 1-point and $q$ is a 1-point,
$$
\begin{array}{ll}
u&=x^a\\
v&=y\\
w&=z
\end{array}
$$
with $a>0$.
\end{enumerate}

\begin{Definition}\label{Def125} Let notation be as in Definition \ref{torf}. $u,v,w$ have a (non toroidal) {\bf monomial form} at $p\in f^{-1}(q)$ if
\begin{enumerate}
\item[1.] $u,v$ have a form (\ref{eqTF1}) at $p$, $q$ is a 2-point and
$w=x^c(z+\alpha)$ for some $c\in{\bf N}$ with $c>0$, $\alpha\in{\bf k}$.
\item[2.] $u,v$ have a form (\ref{eqTF21}) at $p$, $q$ is a 2-point and
$w=x^ey^f(z+\alpha)$ for some $e,f\in{\bf N}$, $e+f>0$, $\alpha\in{\bf k}$.
\item[3.] $u,v$ have a form (\ref{eqTF22}) at $p$ and $w=x^cy^d$ with $ad-bc\ne 0$.
\item[4.] $u,v$ have a form (\ref{eqTF3}) at $p$, $q$ is a 2-point and $w=x^gy^hz^i$ with
$$
\text{det}\left(\begin{array}{lll}a&b&c\\d&e&f\\g&h&i\end{array}\right)\ne 0.
$$
\item[5.] $u,v$ have a form (\ref{eqTF01}) at $p$ and $w=x^b(z+\alpha)$ with $b\in{\bf N}$, $b>0$, $\alpha\in{\bf k}$.
\item[6.] $u,v$ have a form (\ref{eqTF02}) at $p$ and $w=x^cy^d$ with $ad-bc\ne 0$.
\end{enumerate}
\end{Definition}

A prepared morphism $\Phi_X:X\rightarrow S$ from a nonsingular 3-fold $X$ to a nonsingular  surface $S$ (with respect to toroidal structures $D_S$ and $D_X=\Phi_X^{-1}(D_S)$) is defined in Definition 6.5 \cite{C3}.

\begin{Remark}\label{Remark1}
Suppose that $f:X\rightarrow Y$ is a dominant proper morphism of nonsingular 3-folds with toroidal structures determined by SNC divisors
$D_Y$ on $Y$ and $D_X=f^{-1}(D_Y)$ on $X$, and $D_X$ contains the singular locus of the morphism $f$.
With our assumptions on $f$, $f$ is generically finite.
Recall that the fundamental locus of a generically finite morphism $f:X\rightarrow Y$ of nonsingular proper varieties is
$\{p\in Y\mid \text{dim }f^{-1}(p)>0\}$. The fundamental locus is a closed set of codimension $\ge 2$ in $Y$.
Let $\overline X$ be the normalization of $Y$ in the function field of $X$, with induced finite morphism $\lambda:\overline X\rightarrow Y$. The branch
locus of $\lambda$ is contained in the SNC divisor $D_Y$.  Let $E$ be  an irreducible  component of $D_Y$. By Abhyankar's lemma (c.f. XIII 5.3 \cite{G}), the irreducible
components of $\lambda^{-1}(E)$ are disjoint. Thus the irreducible components of $D_X$ which dominate $E$ are disjoint.
\end{Remark}

\begin{Definition}\label{Def31}
A dominant morphism $f:X\rightarrow Y$ of nonsingular 3-folds with toroidal structures determined by SNC divisors $D_Y$ on $Y$, and $D_X=f^{-1}(D_Y)$ on $X$ such that the singular locus of $f$ is contained in $D_X$ is {\bf prepared} for $D_Y$ and $D_X$ if:
\begin{enumerate}
\item[1.] If $q\in Y$ is a 3-point,  $u,v,w$ are  permissible parameters at $q$
and $p\in f^{-1}(q)$, then $u,v$ and $w$ are each a unit (in $\hat{\cal O}_{X,p}$) times a monomial in local equations of the toroidal
structure $D_X$ at $p$ . Furthermore, there exists a permutation of $u,v,w$ such that
$u,v$ are toroidal forms  at $p$.
\item[2.] If $q\in Y$ is a 2-point, $u,v,w$ are permissible parameters at $q$ and $p\in f^{-1}(q)$, then either
\begin{enumerate}
\item $u,v$ are toroidal forms  at $p$ or
\item $p$ is a 1-point and there exist regular parameters $x,y,z\in\hat{\cal O}_{X,p}$ such that there is
 an expression 
$$
\begin{array}{ll}
u&=x^a\\
v&=x^c(\gamma(x,y)+x^dz)\\
w&=y
\end{array}
$$
where $\gamma$ is a unit series and $x=0$ is a local equation of $D_X$, or
\item $p$ is a 2-point and there exist regular parameters $x,y,z$ in $\hat{\cal O}_{X,p}$ such that there is
 an expression 
$$
\begin{array}{ll}
u&=(x^ay^b)^k\\
v&=(x^ay^b)^l(\gamma(x^ay^b,z)+x^cy^d)\\
w&=z
\end{array}
$$
where $a,b>0$, $\text{gcd}(a,b)=1$, $ad-bc\ne 0$, $\gamma$ is a unit series and $xy=0$ is a local equation of $D_X$.
\end{enumerate}
\item[3.] If $q\in Y$ is a 1-point,  and $p\in f^{-1}(q)$,
 then there exist  permissible parameters $u,v,w$ at $q$ such that 
  $u,v$ are
 toroidal forms  at $p$.
\end{enumerate}
\end{Definition}

\begin{Lemma}\label{LemmaT1}
Suppose that $X,Y$ are projective, $f:X\rightarrow Y$ is a prepared morphism of 3-folds, and $q\in D_Y$ is a 1-point. Then there exist algebraic permissible parameters $u,v,w$ at $q$ such that a form (\ref{eqTF01}) or (\ref{eqTF02}) of Definition \ref{torf} (for this fixed choice of $u,v,w$) holds at $p$ for all $p\in f^{-1}(q)$.
\end{Lemma}

\begin{pf} Let $q\in Y$ be a 1-point, $u,v,w$ be permissible parameters at $q$, $\overline X$ be the normalization of $Y$ in ${\bf k}(X)$.  We have a natural commutative diagram
$$
\begin{array}{lll}
X&\stackrel{f}{\rightarrow}&Y\\
g&\searrow&\uparrow\lambda\\
&&\overline X\end{array}.
$$
Let $\lambda^{-1}(q)=\{q_1,\ldots,q_n\}$. Abhyankar's Lemma implies that $\overline X$ is nonsingular above a neighborhood of $q$, and there exist permissible parameters $u_1,v_i,w_i$ at $q_i$ and $r_i\in{\bf N}$ such that
$$
u=u_i^{r_i}, v=v_i, w=w_i.
$$
Lemma 3.5 \cite{C5} (applied to $g:X\rightarrow \overline X$) implies $u_i, \alpha v_i+\beta w_i$ are toroidal forms at all points of $g^{-1}(q_i)$ for a dense open set of $(\alpha,\beta)\in {\bf k}^2$. The condition is that $u_i=\alpha v_i+\beta w_i=0$ intersects the (nonsingular at $q_i$) fundamental locus of $g$ transversally at $q_i$. Thus there exist $\alpha,\beta\in {\bf k}$ such that $u,\alpha v+\beta w$ are toroidal forms at all $p\in f^{-1}(q)$.
\end{pf}

\begin{Lemma}\label{LemmaT2} Suppose that $f:X\rightarrow Y$ is a prepared morphism of 3-folds and $E$ is a component of $D_X$ such that $f(E)=C$ is a curve which is not a 2-curve. Suppose that $q\in C$ is a 1-point. Then all formal branches of $C$ are smooth at $q$, and there exist permissible parameters $u,v,w$ at $q$ such that $u=w=0$ are (formal) local equations of $C$ at $q$.
\end{Lemma}

\begin{pf} Suppose that $p\in f^{-1}(q)\cap E$. There exist permissible parameters $u,v,w$ at $q$ such that one of the following forms hold.

$p$ is a 1-point: 
\begin{equation}\label{eqT3}
u=x^a, v=y, w=g(x,y)+x^nz
\end{equation}
where  $x=0$ is a local equation of $E$, or

$p$ is a 2-point: 
\begin{equation}\label{eqT4}
u=(x^ay^b)^i,v=z, w=g(x^ay^b,z)+x^cy^d
\end{equation}
where $a,b>0$, $ad-bc\ne0$ and $x=0$ is a local equation of $E$.

In case (\ref{eqT3}), we have an expression $w=\phi(y)+x\Omega$ in $\hat{\cal O}_{X,p}$ ($n>0$ since $f(E)=C$). $u=w-\phi(y)=0$ is a local equation of a branch of $C$. 

In case (\ref{eqT4}), one of the following must hold in $\hat{\cal O}_{X,p}$:
\begin{equation}\label{eqT5} w=\phi(z)+xy\Omega,
\end{equation}

\begin{equation}\label{eqT6} w=\phi(z)+xy\Omega+x^c
\end{equation}
where $c>0$,  or 

\begin{equation}\label{eqT7}
w=\phi(z)+xy\Omega+y^d
\end{equation}
where $d>0$.

In cases (\ref{eqT5}) or (\ref{eqT6}), we have
$$
0=u=w-\phi(z)
$$
is a local equation of a branch of $C$ at $q$.

Suppose that case (\ref{eqT7}) holds. $f(E)=C$ a curve implies that there exists an irreducible series
$\Psi(\overline x,\overline y)$ in a power series ring ${\bf k}[[\overline x,\overline y]]$ such that $x$ divides $\Psi(v,w)$ in $\hat{\cal O}_{X,p}$. 
$$
\Psi(v,w)=\Psi(z,\phi(z)+xy\Omega+y^d)
$$
implies $\Psi(0,y^d)=0$ (set $x=z=0$), which implies that $\overline x\mid \Psi(\overline x,\overline y)$. Since $\Psi(\overline x,\overline y)$ is irreducible, we have that $\Psi=\overline x\lambda $ where $\lambda$ is a unit series. This is a contradiction.
\end{pf}

We recall Theorem 1.3 of \cite{C7}.

\begin{Theorem}\label{Theorem2}(Theorem 1.3 \cite{C7})
Suppose that $f:X\rightarrow Y$ is a dominant morphism of nonsingular projective 3-folds
over an algebraically closed field $k$ of characteristic zero, with toroidal structures determined by SNC divisors $D_Y$ on $Y$ and $D_X=f^{-1}(D_Y)$ on $X$
such that  $D_X$ contains the singular locus of $f$.
Then there exists a commutative diagram
$$
\begin{array}{rll}
X_1&\stackrel{f_1}{\rightarrow}&Y_1\\
\Phi\downarrow&&\downarrow\Psi\\
X&\stackrel{f}{\rightarrow}&Y
\end{array}
$$
such that $\Psi$ and $\Phi$ are products  of possible blow ups  for the preimages of $D_Y$, $D_X$ respectively, such that $f_1$ is prepared
for  $D_{Y_1}=\Psi^{-1}(D_Y)$ and $D_{X_1}=\Phi^{-1}(D_X)$, and $D_{X_1}$. 
\end{Theorem}

\section{The $\tau$ invariant}

Throughout this section, we assume that $f:X\rightarrow Y$ is a dominant morphism of nonsingular 3-folds.

\begin{Definition}\label{Def221}
 Suppose that $f:X\rightarrow Y$ is prepared.

Suppose that $p\in X$. Define $\tau_f(p)=\tau(p)=-\infty$ if there exist permissible parameters $u,v,w$ at $q$ such that $u,v,w$ are toroidal forms at $p$. (The existence of such formal parameters $u,v,w$ implies the existence of algebraic permissible parameters).

 Suppose that $p\in X$ is a 3-point, and $\tau(p)\ne-\infty$.
Suppose that $u,v,w$ are permissible parameters at $q=f(p)$.   Then   there is an expression (after possibly permuting
$u,v,w$ if $q$ is a 3-point) 
\begin{equation}\label{eq16}
\begin{array}{ll}
u&=x^ay^bz^c\\
v&=x^dy^ez^f\\
w&=\sum_{i\ge 0} \alpha_iM_i +N
\end{array}
\end{equation}
where $x,y,z$ are permissible parameters at $p$ for $u,v,w$, $\text{rank}(u,v)=2$,
 the sum in $w$ is over (possibly infinitely many) monomials $M_i=x^{a_i}y^{b_i}z^{c_i}$ in $x,y,z$ such that
$\text{rank}(u,v,M_i)=2$,
$\text{deg}(M_i)\le\text{deg}(M_j)$ if $i<j$, $N$ is a monomial in $x,y,z$ such that 
$\text{rank}(u,v,N)=3$
and $N\not\,\mid M_i$ for any $M_i$ in the series $\sum\alpha_i M_i$.

If $q$ is a 3-point and $u,v,w$ is not a toroidal form
(at $p$), we necessarily have (since $f$ is prepared) that 
\begin{equation}\label{eq224}
\sum \alpha_iM_i=M_0\gamma
\end{equation}
 where
 $\gamma$
is a unit series in the monomials $\frac{M_i}{M_0}$ (in $x,y,z$)  such that $$
\text{rank}(u,v,M_0)=\text{rank}(u,v,\frac{M_i}{M_0})=2
$$
for all $i$, and $M_0\mid N$.

Define a group $H_p=H_{f,p}$ as follows.
The Laurent monomials  in $x,y,z$ form a group under multiplication.
We define $H_p=H_{f,p}$ to be the subgroup generated by $u,v$ and the terms $M_i$
appearing in the expansion (\ref{eq16}). We will write the group $H_p$
additively as:
$$
H_p=H_{f,p}=(a,b,c){\bf Z}+(d,e,f){\bf Z}+\sum(a_i,b_i,c_i){\bf Z}.
$$
Define a subgroup $A_p$ of $H_p$ by:
$$
A_p=A_{f,p}=\left\{\begin{array}{ll} (a,b,c){\bf Z}+(d,e,f){\bf Z}+(a_0,b_0,c_0){\bf Z}&\text{ if $q$ is 3-point}\\
(a,b,c){\bf Z}+(d,e,f){\bf Z}&\text{if $q$ is a 2-point}.
\end{array}\right.
$$
Define
$$
L_p=L_{f,p}=H_p/A_p,
$$
$$
\tau(p)=\tau_f(p)=\mid L_p\mid.
$$

 Suppose that $p\in X$ is a 2-point, and $\tau(p)\ne-\infty$.
Suppose that $u,v,w$ are permissible parameters at $q=f(p)$ (which satisfy the conclusions of Lemma \ref{Lemma1} if $f(p)=q$ is a 1-point). Then there is an expression (after possibly permuting
$u,v,w$ if $q=f(p)$ is a 3-point) of one the following forms: 

$q$ is 2-point or a 3-point: 
\begin{equation}\label{eqT10}
\begin{array}{ll}
u&=x^ay^b\\
v&=x^cy^d\\
w&=\sum_{i\ge 0} \alpha_{ij}x^{i}y^{j} +x^ey^f(z+\beta)
\end{array}
\end{equation}
where $x,y,z$ are permissible parameters at $p$ for $u,v,w$, 
and $x^ey^f\not\,\mid x^{a_i}y^{b_i}$ for all $i$, or 

$q$ is a 2-point or a 3-point: 

\begin{equation}\label{eqT11}
\begin{array}{ll}
u&=(x^ay^b)^k\\
v&=(x^ay^b)^t(z+\beta)\\
w&=\sum_{i\ge 0} \alpha_i(z)(x^ay^b)^i +x^cy^d
\end{array}
\end{equation}
where $\text{gcd}(a,b)=1$, $x,y,z$ are permissible parameters at $p$ for $u,v,w$, 
and $x^cy^d\not\,\mid (x^ay^b)^i$ for all $i$, or 

$q$ is a 2-point: 

\begin{equation}\label{eqT12}
\begin{array}{ll}
u&=(x^ay^b)^k\\
v&=\sum\alpha_i(z)(z^ay^b)^i+x^cy^d\\
w&=z
\end{array}
\end{equation}
where $\text{gcd}(a,b)=1$, $x,y,z$ are permissible parameters at $p$ for $u,v,w$, 
and $x^cy^d\not\,\mid (x^ay^b)^i$ for all $i$, or 

$q$ is a 1-point: 

\begin{equation}\label{eqT13}
\begin{array}{ll}
u&=(x^ay^b)^k\\
v&=z\\
w&=\sum_{i\ge 0} \alpha_i(z)(x^ay^b)^i +x^cy^d
\end{array}
\end{equation}
where $\text{gcd}(a,b)=1$, $x,y,z$ are permissible parameters at $p$ for $u,v,w$, 
and $x^cy^d\not\,\mid (x^ay^b)^i$ for all $i$.

If $q$ is a 3-point, we necessarily have (since $f$ is prepared) that 
\begin{equation}\label{eqT14}
w=M_0\gamma
\end{equation}
 where
 $\gamma$
is a unit series.

In equation (\ref{eqT12}), we have 
\begin{equation}\label{eqT70}
v=(x^ay^b)^{i_0}\gamma
\end{equation}
where $\gamma$ is a unit series.

Suppose that (\ref{eqT10}) holds. Then define
$$
\tau(p)=\ell\left(((a,b){\bf Z}+(c,d){\bf Z}+\sum_{\alpha_{ij}\ne 0} (i,j){\bf Z})/((a,b){\bf Z}+(c,d){\bf Z})\right)
$$
if $q$ is a 2-point,

$$
\tau(p)=\ell\left(((a,b){\bf Z}+(c,d){\bf Z}+\sum_{\alpha_{ij}\ne 0} (a_i,b_i){\bf Z})/((a,b){\bf Z}+(c,d){\bf Z}+(a_0,b_0){\bf Z})\right)
$$
if $q$ is a 3-point, and $M_0=x^{a_0}y^{b_0}$ in (\ref{eqT14}).

Suppose that (\ref{eqT11}) holds. Then define
$$
\tau(p)=\ell\left((k{\bf Z}+t{\bf Z}+\sum_{\alpha_i\ne 0} i{\bf Z})/(k{\bf Z}+t{\bf Z})\right)
$$
if $q$ is a 2-point,

$$
\tau(p)=\ell\left((k{\bf Z}+t{\bf Z}+\sum_{\alpha_i\ne 0} i{\bf Z})/(k{\bf Z}+t{\bf Z}+i_0{\bf Z})\right)
$$
if $q$ is a 3-point, and $M_0=x^{i_0}\gamma$ in (\ref{eqT14}).

Suppose that (\ref{eqT12}) (and (\ref{eqT70})) hold. Then define
$$
\tau(p)=\ell\left((k{\bf Z}+\sum_{\alpha_i\ne 0} i{\bf Z})/(k{\bf Z}+i_0{\bf Z})\right).
$$

Suppose that (\ref{eqT13})  holds. Then define
$$
\tau(p)=\ell\left((k{\bf Z}+\sum_{\alpha_i\ne 0} i{\bf Z})/k{\bf Z}\right).
$$

 Suppose that $p\in X$ is a 1-point, and $\tau(p)\ne-\infty$.
Suppose that $u,v,w$ are permissible parameters at $q=f(p)$ which satisfy the conclusions of Lemma \ref{Lemma1} if $f(p)=q$ is a 1-point. Then there is an expression (after possibly permuting
$u,v,w$ if $q=f(p)$ is a 3-point) of one the following forms: 

$q$ is 2-point or a 3-point 
\begin{equation}\label{eqT15}
\begin{array}{ll}
u&=x^a\\
v&=x^b(\alpha+y)\\
w&=\sum_{i<c} \alpha_i(y)x^i +x^c(z+\beta)
\end{array}
\end{equation}
where $x,y,z$ are permissible parameters at $p$ for $u,v,w$, 
 or 
 
 $q$ is 2-point 

\begin{equation}\label{eqT16}
\begin{array}{ll}
u&=x^a\\
v&=\sum_{i<c} \alpha_i(y)x^i +x^c(z+\beta)\\
w&=y
\end{array}
\end{equation}
where $x,y,z$ are permissible parameters at $p$ for $u,v,w$, 
 or 
 
 $q$ is a 1-point 

\begin{equation}\label{eqT17}
\begin{array}{ll}
u&=x^a\\
v&=y\\
w&=\sum_{i<c} \alpha_i(y)x^i +x^c(z+\beta)\\
\end{array}
\end{equation}
where $x,y,z$ are permissible parameters at $p$ for $u,v,w$. 

In case (\ref{eqT15}) we have 
\begin{equation}\label{eqT71}
w=x^{i_0}\gamma
\end{equation}
where $\gamma$ is a unit series if $q$ is a 3-point.
In case (\ref{eqT16}) we have 
\begin{equation}\label{eqT72}
v=x^{i_0}\gamma
\end{equation}
where $\gamma$ is a unit series.

Suppose that (\ref{eqT15}) holds. Then define
$$
\tau(p)=\ell\left((a{\bf Z}+b{\bf Z}+\sum_{\alpha_{i}\ne 0} i{\bf Z})/(a{\bf Z}+b{\bf Z})\right)
$$
if $q$ is a 2-point,

$$
\tau(p)=\ell\left((a{\bf Z}+b{\bf Z}+\sum_{\alpha_{i}\ne 0} i{\bf Z})/(a{\bf Z}+b{\bf Z}+i_0{\bf Z})\right)
$$
if $q$ is a 3-point.

Suppose that (\ref{eqT16}) holds. Then define
$$
\tau(p)=\ell\left((a{\bf Z}+\sum_{\alpha_{i}\ne 0} i{\bf Z})/(a{\bf Z}+i_0{\bf Z})\right).
$$

Suppose that (\ref{eqT17}) holds. Then define
$$
\tau(p)=\ell\left((a{\bf Z}+\sum_{\alpha_{i}\ne 0} i{\bf Z})/a{\bf Z}\right).
$$

\end{Definition}

Observe that $\tau(p)<\infty$ in Definition \ref{Def221}, since $H_p$ is a finitely generated
group, and $H_p/A_p$ is a torsion group.

We define 
$$
\tau(X)=\tau_f(X)=\text{max}\{\tau_f(p)\mid p\in X\}.
$$
\begin{Lemma}\label{Lemma0} Suppose that $f:X\rightarrow Y$ is prepared and $p\in D_X$. Then
$\tau_f(p)$ is independent of choice of permissible  parameters $u,v,w$ at $q=f(p)$ and permissible parameters $x,y,z$ at $p$ for $u,v,w$. 
\end{Lemma}

\begin{pf}
This is proved in Lemma 3.10 of \cite{C5} when $p$ is a 3-point.  We give here a proof when $p$ is a 1-point and $q=f(p)$ is a 2-point or a 3-point. The proof in the remaining cases is similar.

This is not difficult to prove when there exists a computation of $\tau_f(p)$ which gives $-\infty$, so we will assume that $\tau_f(p)\ge 0$ for all computations of $\tau_f(p)$.

Assume that $p$ is a 1-point, and that $q=f(p)$ is a 2-point or a 3-point.  Suppose that $u,v,w$ are permissible parameters at $q$, such that $u,v,w$ have an expression  of the form (\ref{eqTF1}) and (\ref{eqT15}), or of the form  2 (b) of Definition \ref{Def31} and (\ref{eqT16}).

We will first show that for fixed permissible parameters $u,v,w$ at $q$, $\tau_f(p)$ is independent of choice of permissible parameters $x,y,z$ for $u,v,w$ at $p$ which have an expression of type (\ref{eqTF1}) (and (\ref{eqT15}).

We have  that $u,v,w$ have an expression 
\begin{equation}\label{eqT211}
\begin{array}{ll}
u&=x^a\\
v&=x^b(\alpha+y)\\
w&=\sum_{j=0}^n\alpha_{i_j}(y)x^{i_j}+x^c(z+\beta)
\end{array}
\end{equation}
where $i_j<c$ for all $j$, and $\alpha_{i_j}(y)\ne 0$ for all $j$.

If $\overline x,\overline y, \overline z$ are other permissible parameters for $u,v,w$ at $q$, which also give an expansion of type 2 (and (\ref{eqT15})), then we have
$x=\omega \overline x$ where $\omega^a=1$, and thus $y=\omega^{-b}\overline y$, $z=\omega^{-c}\overline z$.  Thus the computation of $\tau_f(p)$ for $\overline x,\overline y,\overline z$ is the same as for $x,y,z$.

By a similar calculation, for fixed permissible parameters $u,v,w$ at $q$, we  may  show that $\tau$ is independent of choice of permissible parameters $x,y,z$ for $u,v,w$ at $p$ which have an expression of type 2 (b) of Definition \ref{Def31} and (\ref{eqT16}) at $p$.

Now suppose that $q$ is a 2-point, and $u,v,w$ admit expressions of both the form 2 (b) of Definition \ref{Def31} and of the form (\ref{eqTF1}) at $p$. Let $\tau_1$ be the computation of $\tau_f(p)$ from a representation of $u,v,w$  in the form 2 (b) of Definition \ref{Def31}, and let $\tau_2$ be the computation of $\tau_f(p)$ from a representation of $u,v,w$ in the form of (\ref{eqTF1}).

Since  $u,v,w$ have an expression of type 2 (b) at $p$, there are permissible parameters $x,y,z$ at $p$ such that there is an expression 
\begin{equation}\label{eqT212}
\begin{array}{ll}
u&=x^a\\
v&=\sum_{j=0}^n\alpha_{i_j}(y)x^{i_j}+x^c(z+\beta)\\
w&=y
\end{array}
\end{equation}
of the form (\ref{eqT16}), where $\alpha_{ij}$ are all  nonzero and $i_j<c$ for all $j$.
$\tau_1$ is the computation of $\tau_f(p)$ for this expression, since 
 $u,v,w$ also have an expression of type (\ref{eqTF1}) at $p$, as $u,v$ are toroidal forms at $p$. Since $\tau_1\ge 0$,  we have
$$
v=(\overline\alpha+y\lambda(y))x^{i_0}+\alpha_{i_1}(y)x^{i_1}+\cdots+\alpha_{i_n}(y)x^{i_n}+x^c(z+\beta)
$$
where $0\ne\overline\alpha=\alpha_{i_0}(0)\in {\bf k}$ and $\lambda(y)$ is a unit series. Set
$$
\overline y=y\lambda(y)+\alpha_{i_1}(y)x^{i_1-i_0}+\cdots+\alpha_{i_n}(y)x^{i_n-i_0}+x^{c-i_0}(z+\beta).
$$
Then 
\begin{equation}\label{eqT206}
y=\lambda(y)^{-1}\overline y-\lambda(y)^{-1}\alpha_{i_1}(y)x^{i_1-i_0}-\cdots-\lambda(y)^{-1}\alpha_{i_n}(y)x^{i_n-i_0}-x^{c-i_0}\lambda(y)^{-1}(z+\beta).
\end{equation}
Iterate, substituting (\ref{eqT206}) into successive iterations of (\ref{eqT206}), to get a series $\overline P$ such that
$$
y=\overline P(\overline y,x^{i_1-i_0},\ldots,x^{i_n-i_0})+x^{c-i_0}\Omega
$$
for some $\Omega\in\hat{\cal O}_{X,p}$.

We compute the Jacobian determinant 
\begin{equation}\label{eqT205}
\text{Det}\left(\begin{array}{l}\partial(u,v,w)\\ \partial(x,y,z)\end{array}\right)=-ax^{a+c-1}.
\end{equation}
We have 
$$
\begin{array}{ll}
u&=x^a\\
v&=x^{i_0}(\overline y+\overline\alpha)\\
w&=\overline P(\overline y,x^{i_1-i_0},\ldots, x^{i_n-i_0})+x^{c-i_0}\Omega,
\end{array}
$$
and 
\begin{equation}\label{eqT207}
\text{Det}\left(\begin{array}{l} \partial(u,v,w)\\ \partial(x,\overline y,z)\end{array} \right)=ax^{a+c-1}\frac{\partial \Omega}{\partial z}.
\end{equation}
Since (\ref{eqT205}) and (\ref{eqT207}) differ by multiplication of a unit series, 
$\frac{\partial \Omega}{\partial z}(0,0,0)\ne 0$. Set $\overline\beta=\Omega(0,0,0)$, $\overline z=\Omega-\overline\beta$. We have
$$
\begin{array}{ll}
u&=x^a\\
v&=x^{i_0}(\overline y+\alpha)\\
w&=\overline P(\overline y,x^{i_1-i_0},\ldots, x^{i_n-i_0})+x^{c-i_0}(\overline \beta+\overline z).
\end{array}
$$
$\tau_2$ is the value of $\tau_f(p)$ for this expression, computed from (\ref{eqT15}).  We have 
\begin{equation}\label{eqT210}
\begin{array}{ll}
\tau_2&\le \ell\left((a{\bf Z}+i_0{\bf Z}+\sum_{j=1}^n(i_j-i_0){\bf Z})/(a{\bf Z}+i_0{\bf Z})\right)\\
&=\ell\left((a{\bf Z}+\sum_{j=0}^ni_j{\bf Z})/(a{\bf Z}+i_0{\bf Z})\right)\\
&=\tau_1.
\end{array}
\end{equation}

We now show that $\tau_1\le\tau_2$. Since $u,v,w$ have an expression of type (\ref{eqTF1}) at $p$, we have an expression 
\begin{equation}\label{eq251}
\begin{array}{ll}
u&=x^a\\
v&=x^b(\alpha+y)\\
w&=P(x,y)+x^c(z+\beta)
\end{array}
\end{equation}
of the form (\ref{eqT15}), where 
$$
P=\alpha_{i_0}(y)+\alpha_{i_1}(y)x^{i_1}+\cdots+\alpha_{i_n}(y)x^{i_n},
$$
$\alpha_{i_j}(y)\ne 0$ for all $j$, and $i_n<c$.

$\tau_2$ is the computation of $\tau_f(p)$ for this expression. 
Since $\tau_2\ge0$, and 
$u,v,w$ also have an expression of type 2 (b) of Definition \ref{Def31}, we have   $\frac{\partial P}{\partial y}(0,0)\ne 0$. We then have
$$
\alpha_{i_0}(y)=y\lambda(y)
$$
where $\lambda$ is a unit series. Set $\overline z=P+x^c(z+\beta)$. 
 We have 
\begin{equation}\label{eqT208}
\begin{array}{ll}
y&=\lambda(y)^{-1}[\overline z-\alpha_{i_1}(y)x^{i_1}-\cdots-\alpha_{i_n}(y)x^{i_n}-x^c(z+\beta)]\\
&=\lambda(y)^{-1}\overline z-\lambda(y)^{-1}\alpha_{i_1}(y)x^{i_1}-\cdots-\lambda(y)^{-1}\alpha_{i_n}(y)x^{i_n}-\lambda(y)^{-1}x^c(z+\beta).
\end{array}
\end{equation}
Iterate, substituting (\ref{eqT208}) into successive iterations of (\ref{eqT208}), to get series $\overline P$ and $\Omega$ such that
$$
y=\overline P(\overline z,x^{i_1},\ldots,x^{i_n})+x^c\Omega(\overline z,x,z).
$$
Thus we have an expression
$$
\begin{array}{ll}
u&=x^a\\
v&=x^b(\alpha+\overline P(\overline z,x^{i_1},\ldots, x^{i_n}))+x^{b+c}\Omega\\
w&=\overline z.
\end{array}
$$

We compute the Jacobian determinant from (\ref{eq251}),
$$
\text{Det}\left(\begin{array}{l} \partial(u,v,w)\\ \partial(x,y,z)\end{array}\right)=ax^{a+b+c-1}.
$$
We compare with 
$$
\text{Det}\left(\begin{array}{l} \partial(u,v,w)\\ \partial(x, z,\overline z)\end{array}\right)=ax^{a+b+c-1}\frac{\partial \Omega}{\partial z},
$$
to see that
$\frac{\partial \Omega}{\partial z}(0,0,0)\ne 0$. Set $\overline\beta=\Omega(0,0,0)$, $\overline y=\Omega-\overline\beta$. We have
$$
\begin{array}{ll}
u&=x^a\\
v&=x^b(\alpha+\overline P(\overline z,x^{i_1},\ldots, x^{i_n}))+x^{b+c}(\overline\beta+\overline y)\\
w&=\overline z
\end{array}
$$
in the form of 2 (b) of Definition \ref{Def31}. 
 $\tau_1$ is the value of $\tau_f(p)$ for this expression.  We compute from (\ref{eqT16}), 
\begin{equation}\label{eqT209}
\begin{array}{ll}
\tau_1&\le \ell\left((a{\bf Z}+b{\bf Z}+\sum_{j=0}^ni_j{\bf Z})/(a{\bf Z}+b{\bf Z})\right)\\
&=\tau_2.
\end{array}
\end{equation}
By (\ref{eqT210}) and (\ref{eqT209}) we see that $\tau_1=\tau_2$.

From our calculations so far, we conclude that for fixed $u,v,w$, $\tau_f(p)$ is independent of choice of permissible parameters $x,y,z$ for $u,v,w$ at $p$.

Now we will show that if $u_1,v_1,w_1$ are a permutation of $u,v,w$ such that $u_1,v_1,w_1$ have one of the forms (\ref{eqTF1}) and (\ref{eqT15}) or 2 (b) of Definition \ref{Def31} and (\ref{eqT16}) at $p$, then we obtain the same value of $\tau_f(p)$.

First suppose that $u_1=v, v_1=u$, $w_1=w$, and $u,v$ are toroidal forms at $p$ (so that (\ref{eqTF1}) holds). Then $u,v,w$ have an expression (\ref{eqT211}) at $p$. We  have
$$
\begin{array}{ll}
u_1=v&=\overline x^b\\
v_1=u&=\overline x^a(\overline\alpha+\overline y)
\end{array}
$$
where $x=\overline x(\alpha+y)^{-\frac{1}{b}}$, $\overline\alpha=\alpha^{-\frac{a}{b}}$, and $y=\lambda(\overline y)\overline y$
for an appropriate  unit series $\lambda(\overline y)$.

Set
$$
\overline\beta=\beta\alpha^{-\frac{c}{b}},
\overline z=(z+\beta)(\alpha+\lambda(\overline y)\overline y)^{-\frac{c}{b}}-\overline\beta.
$$
We have
$$
\begin{array}{ll}
w &=\sum_{j=0}^n\alpha_{i_j}(\lambda(\overline y)\overline y)(\alpha+\lambda(\overline y)\overline y)^{-\frac{i_j}{b}}\overline x^{i_j}+\overline x^c(z+\beta)(\alpha+\lambda(\overline y)\overline y)^{-\frac{c}{b}}\\
&=\sum_{j=0}^n \overline \alpha_{ij}(\overline y)\overline x^{i_j}+\overline x^c(\overline\beta+\overline z)
\end{array}
$$
where $\overline\alpha_{ij}(\overline y)\ne 0$ for all $j$.

Let $\tau_1$ be the calculation of $\tau_f(p)$ for the variables $u,v,w$, and let $\tau_2$ be the calculation of $\tau$ for the variables $u_1,v_1,w_1$.
 We see that $\tau_1=\tau_2$.

Now suppose that  $u_1=v, v_1=u, w_1=w$, and $u,v,w$  have an expression of the form 2 (b) of Definition \ref{Def31} and (\ref{eqT16}). $u,v,w$ have an expression (\ref{eqT212}) at $p$. We have
$$
\begin{array}{ll}
u&=x^a\\
v&=x^{i_0}(\sum_{j=0}^n\alpha_{i_j}(y)x^{i_j-i_0}+x^{c-i_0}(z+\beta))\\
w&=y
\end{array}
$$
where $\alpha_{i_0}(y)$ is a unit series.
Set 
$$
\gamma=\sum_{j=0}^n\alpha_{i_j}(y)x^{i_j-i_0}+x^{c-i_0}(z+\beta).
$$
Define $\overline x$ by 
\begin{equation}\label{eqT213}
x=\overline x\gamma^{-\frac{1}{i_0}}.
\end{equation}
We have a series $P$ and $\Omega\in\hat{\cal O}_{X,p}$ such that 
$$
\begin{array}{ll}
u_1&=v=\overline x^{i_0}\\
v_1&=u=\overline x^a\gamma^{-\frac{a}{i_0}}=\overline x^a(P(y,x^{i_1-i_0},\ldots,x^{i_n-i_0})+\overline x^{c-i_0}\Omega)\\
w&=y
\end{array}
$$
By iterating substitution of (\ref{eqT213}) in $P(y,x^{i_1-i_0},\ldots, x^{i_n-i_0})$, we see that there is a  series $\overline P(y,\overline x^{i_1-i_0},\ldots, \overline x^{i_n-i_0})$ and $\overline\Omega\in\hat{\cal O}_{X,p}$ such that
$$
v_1=\overline x^a(\overline P(y,\overline x^{i_1-i_0},\ldots, \overline x^{i_n-i_0})+\overline x^{c-i_0}\overline\Omega).
$$

Since the Jacobians
$$
\text{Det}\left(\begin{array}{l}\partial(u,v,w) \\ \partial(x,y,z)\end{array}\right)=-ax^{a+c-1}
$$
and
$$
\text{Det}\left(\begin{array}{l} \partial(u_1,v_1,w_1)\\ \partial(\overline x,y,z)\end{array}\right)=-i_0\overline x^{a+c-1}\frac{\partial \overline\Omega}{\partial z}
$$
differ by multiplication by a unit, $\frac{\partial\Omega}{\partial z}(0,0,0)\ne 0$. 

Set 
$$
\overline\beta=\overline\Omega(0,0,0),\,\, \overline z=\Omega-\overline\beta.
$$
Then $\overline x, y,\overline z$ are regular parameters in $\hat{\cal O}_{X,p}$, and 
$$
\begin{array}{ll}
u_1&=\overline x^{i_0}\\
v_1&=\overline x^a\overline P(y,\overline x^{i_1-i_0},\ldots,\overline x^{i_n-i_0})+\overline x^{a+c-i_0}(\overline\beta+\overline z)\\
w&=y.
\end{array}
$$

Let $\tau_1$ be the calculation of $\tau_f(p)$ for the variables $u,v,w$, and let $\tau_2$ be the calculation of $\tau_f(p)$ for the variables $u_1,v_1,w_1$ (from (\ref{eqT16}). We see that $\tau_1\le\tau_2$.
From this argument for the change of variables $u_1,v_1,w_1$ to $u,v,w$, we see that $\tau_2\le\tau_1$, so that $\tau_1=\tau_2$.

Now suppose that $q$ is a 3-point, $u,v,w$ have an expression (\ref{eqTF1}) and (\ref{eqT15}), and $u_1=w,v_1=v,w_1=u$ also have an expression (\ref{eqTF1}) and (\ref{eqT15}).

We have 
\begin{equation}\label{eqT330}
\begin{array}{ll}
u&=x^a\\
v&=x^b(\alpha+y)\\
w&=x^{i_0}(\sum_{j=0}^n\alpha_{i_j}(y)x^{i_j-i_0}+x^{c-i_0})(z+\beta))
\end{array}
\end{equation}
where $n<c$, $\alpha_{i_j}(y)\ne 0$ for all $j$, and $\alpha_{i_0}(y)$ is a unit series. Set
$$
\gamma=\sum_{j=0}^n\alpha_{i_j}(y)x^{i_j-i_0}+x^{c-i_0}(z+\beta).
$$
Set $x=\overline x\gamma^{-\frac{1}{i_0}}$.
We have
$$
\begin{array}{ll}
u_1=w&=\overline x^{i_0}\\
v_1=v&=\overline x^b\gamma^{-\frac{b}{i_0}}(\alpha+y)\\
w_1=u&=\overline x^a\gamma^{-\frac{a}{i_0}}.
\end{array}
$$

As we are assuming that all calculations of $\tau_f(p)$ are $\ge 0$, we have that $c>0$. Since $w,v$ are assumed to be toroidal forms at $p$, we then have
$$
\frac{\partial}{\partial y}(\gamma^{-\frac{b}{i_0}}(\alpha+y))(0,0,0)\ne 0.
$$
Set 
$$
\overline\alpha=\gamma^{-\frac{b}{i_0}}(0,0,0)\alpha,\,\,
\overline y=\gamma^{-\frac{b}{i_0}}(\alpha+y)-\overline\alpha.
$$

We have 
\begin{equation}\label{eqT214}
x=\overline x[(\alpha+y)^{-1}(\overline y+\overline\alpha)]^{\frac{1}{b}}.
\end{equation}

There exists a unit series $\lambda(y)$, series $\overline\alpha_{r_1,\ldots r_n}(y)$, and $\overline\Omega\in\hat{\cal O}_{X,p}$, such that 
$$
\overline y=y\lambda(y)+\sum_{r_1,\ldots,r_n>0}\overline\alpha_{r_1,\ldots,r_n}(y)
x^{(i_1-i_0)r_1}\cdots x^{(i_n-i_0)r_n}+x^{c-i_0}\overline\Omega.
$$
Thus 

\begin{equation}\label{eqT215}
y=\overline y\lambda(y)^{-1}-\sum_{r_1,\ldots,r_n>0}\lambda(y)^{-1}\overline \alpha_{r_1,\ldots,r_n}(y)x^{(i_1-i_0)r_1}\cdots x^{(i_n-i_0)r_n}
-x^{c-i_0}\lambda(y)^{-1}\overline\Omega.
\end{equation}

There exists a series 
\begin{equation}\label{eqT216}
P(y,x^{i_1-i_0},\ldots,x^{i_n-i_0})
\end{equation}
such that $u=\overline x^a(P+\overline x^{c-i_0}\Omega)$
for some $\Omega\in\hat{\cal O}_{X,p}$.

Substituting (\ref{eqT214}) and (\ref{eqT215}) into successive iterations of (\ref{eqT216}), we see that there is a series
$\overline P(\overline y,\overline x^{i_1-i_0},\ldots,\overline x^{i_n-i_0})$ and $\Omega'\in\hat{\cal O}_{X,p}$ such that 
$$
w_1=u=\overline x^a(\overline P(\overline y,\overline x^{i_1-i_0},\ldots,\overline x^{i_n-i_0})+\overline x^{c-i_0}\Omega').
$$

Thus there is an expression 
\begin{equation}\label{eqT187}
\begin{array}{ll}
u_1=w&=\overline x^{i_0}\\
v_1=v&=\overline x^b(\overline\alpha+\overline y)\\
u&=\overline x^a(\overline P(\overline y,\overline x^{i_1-i_0},\ldots,\overline w_1=x^{i_n-i_0})+\overline x^{c-i_0}\Omega').
\end{array}
\end{equation}

We compute the Jacobian determinants
$$
\text{Det}\left(\begin{array}{l} \partial(u,v,w)\\ \partial(x,y,z)\end{array}\right)=a x^{a+b+c-1}
$$
and   
$$
\text{Det}\left(\begin{array}{l} \partial(u_1,v_1,w_1)\\ \partial(\overline x,\overline y,z)\end{array}\right)=i_0\overline x^{a+b+c-1}\frac{\partial \Omega'}{\partial z}.
$$
Thus $\frac{\partial \Omega'}{\partial z}(0,0,0)
\ne 0$.
Set
$$
\overline\beta=\Omega'(0,0,0),\,\, \overline z=\Omega'-\overline\beta.
$$

 $\overline x,\overline y,\overline z$ are regular parameters in $\hat{\cal O}_{X,p}$.
 We have an expression 
 \begin{equation}\label{eqT217}
 \begin{array}{ll}
 u_1=w&=\overline x^{i_0}\\
 v_1=v&=\overline x^b(\overline\alpha+\overline y)\\
 w_1=u&=\overline x^a\overline P(\overline y,\overline x^{i_1-i_0},\ldots,\overline x^{i_n-i_0})+\overline x^{a+c-i_0}(\overline\beta+\overline z)
 \end{array}
 \end{equation}
 of the form (\ref{eqT15}).

Let $\tau_1$ be the computation of $\tau_f(p)$ from (\ref{eqT217}), $\tau_2$ be the computation of $\tau_f(p)$ from (\ref{eqT217}).  We see that $\tau_1\le \tau_2$.

By this argument applied to the change of variables $u_1=w,v_1=v,w_1=u$ to $u,v,w$, we see that $\tau_2\le\tau_1$, so that $\tau_1=\tau_2$.

There is a similar calculation if $u_1=u, v_1=w, w_1=v$.

Now all other permutations $u_1,v_1,w_1$ of $u,v,w$ such that $u_1,v_1,w_1$ have a form (\ref{eqTF1}) and (\ref{eqT15}) at $p$ (so that $u_1,v_1$ are toroidal forms at $p$) are a composition of change of variables of the three types analyzed above.  We have shown that $\tau_f(p)$ is invariant under such change of variables. Hence we get the same evaluation of $\tau_f(p)$ for all possible permuations of the variables $u,v,w$ and $u_1,v_1,w_1$.

Now assume that $u_1=u, v_1=v$ and  $w_1$ are permissible parameters at $q$.

Suppose that $u,v,w$ have an expression of the form (\ref{eqTF1}) at $p$. We then have an expression (\ref{eqT211}) of $u,v,w$. 
Let $\tau_1$ be the computation of $\tau_f(p)$ for this expression. Set
$$
P=\sum_{j=0}^n\alpha_{ij}(y)x^{i_j}.
$$
We have an expansion 
\begin{equation}\label{eqT219}
w_1=\sum_{i=0}^{\infty}g_i(u,v)w^i
\end{equation}
and $w_1=w\gamma$ where $\gamma(u,v,w)$ is a unit series if $q=f(p)$ is a 3-point.

Substitute $w=P+x^c(z+\beta)$ into (\ref{eqT219}) to get
$$
w_1=\sum_{i=0}^{\infty}\sum_{j=0}^ig_i(u,v)\binom{i}{j}P^{i-j}x^{jc}(z+\beta)^j.
$$
We see that there exists a series $\overline P(y,x^a,x^b,x^{i_0},\ldots,x^{i_n})$ such that
$$
w_1=\overline P(y,x^a,x^b,x^{i_0},\ldots,x^{i_n})+x^c\Omega
$$
for some series $\Omega$. Comparing the Jacobian determinants
$$
\text{Det}\left(\begin{array}{l}\partial(u,v,w)\\ \partial(x,y,z)\end{array}\right)=ax^{a+b+c-1}
$$
and  
$$
\text{Det}\left(\begin{array}{l}\partial(u,v,w_1)\\ \partial(x,y,z)\end{array}\right)=ax^{a+b+c-1}\frac{\partial\Omega}{\partial z},
$$
we see that 
$\frac{\partial\Omega}{\partial z}(0,0,0)\ne 0$. Set
$$
\overline\beta=\Omega(0,0,0),\,\, \overline z=\Omega-\overline\beta.
$$
We have
$$
\begin{array}{ll}
u&=x^a\\
v&=x^b(\alpha+y)\\
w_1&=\overline P(y,x^a,x^b,x^{i_0},\ldots,x^{i_n})+x^c(\overline z+\overline\beta)
\end{array}
$$
of the form of (\ref{eqT15}). Further, we have $w_1=x^{i_0}\overline\gamma$ where $\overline\gamma(\overline x,\overline y,\overline z)$ is a unit series if $q=f(p)$ is a 3-point. Let $\tau_2$ be the computation of $\tau_f(p)$ for this expression. We see that $\tau_2\le\tau_1$. Now applying this argument to the change of variables from $u,v,w_1$ to $u,v,w$, we see that $\tau_1=\tau_2$.

Now suppose that $u,v,w$ have an expression of the form 2 (b) of Definition \ref{Def31} and (\ref{eqT16}) (so that $q=f(p)$ is a 2-point). We have
$$
\begin{array}{ll}
u&=x^a\\
v&=\sum_{j=0}^n\alpha_{i_j}(y)x^{i_j}+x^c(z+\beta)\\
w&=y
\end{array}
$$
where $\alpha_{i_j}(y)$ are all nonzero and $i_j<c$ for all $j$.

Let $\tau_1$ be the computation of $\tau_f(p)$ for $u,v,w$.

By the formal implicit function theorem, there exists a unit series $\gamma(u,v,w)$ and a series $\phi(u,v)$ with $\phi(0,0)=0$ such that $w_1=\gamma(w-\phi)$. We have
$$
w_1=\gamma(0,0,y)y+x\Omega
$$
for some series $\Omega(x,y,z)$. We may thus set $\overline y=w_1$, and have that $x,\overline y,z$ are regular parameters in $\hat{\cal O}_{X,p}$.

There exist a unit series $\overline\alpha_0(y)$ and series $\overline\alpha_{r_0,\ldots,r_n}(y)$ such that there is an expression
$$
\overline y=\overline\alpha_0(y)y+\sum_{r_0+\cdots+r_n>0}\overline\alpha_{r_0,\ldots,r_n}(y)x^{i_0r_0+\cdots+i_nr_n}+x^c\overline\Omega
$$
with $\overline\alpha_0(0)\ne 0$.

We thus have an expression 
\begin{equation}\label{eqT220}
y=\overline \alpha_0(y)^{-1}\overline y-\sum_{r_0+\cdots+r_n>0}\overline\alpha_0(y)^{-1}\overline\alpha_{r_0,\ldots,r_n}(y)x^{i_0r_0+\cdots+i_nr_n}+x^c\overline\alpha_0(y)^{-1}\overline\Omega.
\end{equation}

We substitute (\ref{eqT220}) into successive iterations of
$$
v=\sum_{j=0}^n\alpha_{i_j}(y)x^{i_j}+x^c(z+\beta)
$$
to obtain an expression
$$
v=P(\overline y, x^{i_0},\ldots, x^{i_n})+x^c\Omega.
$$
Comparing the Jacobian determinants
$$
\text{Det}\left(\begin{array}{l}\partial(u,v,w)\\ \partial(x,y,z)\end{array}\right)=-ax^{a+c-1}
$$
and  
$$
\text{Det}\left(\begin{array}{l}\partial(u,v,w_1)\\ \partial(x,\overline y,z)\end{array}\right)=-ax^{a+c-1}\frac{\partial\Omega}{\partial z},
$$
we see that 
$\frac{\partial\Omega}{\partial z}(0,0,0)\ne 0$. Set
$$
\overline\beta=\Omega(0,0,0),\,\, \overline z=\Omega-\overline\beta.
$$
$x,\overline y,\overline z$ are regular parameters in $\hat{\cal O}_{X,p}$, and we have an expression
$$
\begin{array}{ll}
u&=x^a\\
v&=\overline P(x^{i_0},\ldots,x^{i_n})+x^c(\overline z+\overline\beta)\\
w_1&=\overline y
\end{array}
$$
of the form of (\ref{eqT16}).  Let $\tau_2$ be the computation of $\tau_f(p)$ for $u,v,w_1$. We have
 $\tau_2\le\tau_1$. Now applying this argument to the change of variables from $u,v,w_1$ to $u,v,w$, we see that $\tau_1=\tau_2$.

Using the above techniques, we can show that we have the same computation of $\tau_f(p)$ for a change of variables 
$$
u_1=\lambda_1u,
v_1=\lambda_2v,
w_1=w
$$
where $\lambda_1,\lambda_2$ are unit series.

Finally, suppose that $u_1,v_1,w_1$ are arbitrary permissible parameters at $q$. Then $u_1,v_1,w_1$ may be obtained from $u,v,w$ be a series of changes of variables of the types considered above. It follows that the computation of $\tau_f(p)$ for the variables $u,v,w$ and $u_1,v_1,w_1$ are the same.

\end{pf}

\begin{Lemma}\label{LemmaT111} Suppose that $f:X\rightarrow Y$ is prepared. Then
$\tau_f$ is  upper semi continuous on $D_X$.
\end{Lemma}

\begin{pf} This follows from a local calculation, computing a deformation in etale local coordinates of the possible prepared forms.

Suppose that $p^*\in X$. We must find a Zariski open neighborhood $U$ of $p^*$ in $X$ such that $p\in U\cap D_X$ implies
$\tau_f(p)\le\tau_f(p^*)$.

We work out a few cases.  Suppose that $p*\in X$ is a 3-point, $\tau_f(p^*)\ge 0$ and there exist (algebraic) permissible parameters $u,v,w$ at $q*=f(p*)$ such that $u,v$ are toroidal forms at $p*$. Then there exist regular parameters $\tilde x,\tilde y,\tilde z\in{\cal O}_{X,p*}$, and units $\gamma_1,\gamma_2\in{\cal O}_{X,p*}$  such that 
$$
\begin{array}{ll}
u&=\tilde x^a\tilde y^b\tilde z^c\gamma_1\\
v&=\tilde x^d\tilde y^e\tilde z^f\gamma_2
\end{array}
$$
with
$$
\text{rank}\left(\begin{array}{lll}
a&b&c\\
d&e&f
\end{array}\right)=2.
$$

There exist rational numbers $\overline a_{ij}$ such that if we set
$$
x=\tilde x\gamma_1^{\overline a_{11}}\gamma_2^{\overline a_{12}},
y=\tilde x\gamma_1^{\overline a_{21}}\gamma_2^{\overline a_{22}},
z=\tilde x\gamma_1^{\overline a_{31}}\gamma_2^{\overline a_{32}},
$$
then 
$$
u=x^ay^bz^c,\,\,
v=x^dy^ez^f.
$$
From the Jacobian determinant
$$
\text{Det}\left(\begin{array}{l} \partial(u,v,w)\\ \partial(x,y,z)\end{array}\right)
$$
we see that 
we have an expansion 
\begin{equation}\label{eqT191}
\begin{array}{ll}
u&=x^ay^bz^c\\
v&=x^dy^ez^f\\
w&=\sum_{i\ge 0} \alpha_ix^{a_i}y^{b_i}z^{c_i}+x^gy^hz^i\gamma
\end{array}
\end{equation}
where $\gamma$ is a unit series,
$$
\text{Det}\left(\begin{array}{lll}
a&b&c\\
d&e&f\\
g&h&i
\end{array}
\right)\ne 0,
$$

and for all $i$,
$\alpha_i\ne 0$,
$$
\text{Det}\left(\begin{array}{lll}
a&b&c\\
d&e&f\\
a_i&b_i&c_i
\end{array}
\right)= 0,
$$
and $(a_i,b_i,c_i)\not\ge(g,h,i)$.

We can compute $\tau_f(p^*)$ by changing variables, multiplying $x,y,z$ by appropriate rational powers of $\gamma$. We obtain that
$$
\begin{array}{ll}
u&=\hat x^a\hat y^b\hat z^c\\
v&=\hat x^d\hat y^e\hat z^f\\
w&=\sum_{i\ge0} \alpha_i\hat x^{a_i}\hat y^{b_i}\hat z^{c_i}+\hat x^g\hat y^h\hat z^i.
\end{array}
$$
Thus
$$
\tau_f(p^*)=\ell\left(((a,b,c){\bf Z}+(d,e,f){\bf Z}+\sum_{i\ge 0} (a_i,b_i,c_i){\bf Z})/((a,b,c){\bf Z}+(d,e,f){\bf Z})\right)
$$
if $q^*$ is a 2-point,
$$
\tau_f(p^*)=\ell\left(((a,b,c){\bf Z}+(d,e,f){\bf Z}+\sum (a_i,b_i,c_i){\bf Z})/((a,b,c){\bf Z}+(d,e,f){\bf Z}+(a_0,b_0,c_0){\bf Z})\right)
$$
if $q^*$ is a 3-point.

We compute the Jacobian 
\begin{equation}\label{eqT197}
\text{Det}\left(\begin{array}{l}
\partial(u,v,w)\\
\partial(\hat x,\hat y,\hat z)
\end{array}
\right)=\text{Det}\left(\begin{array}{lll}
a&b&c\\
d&e&f\\
g&h&i
\end{array}
\right) \hat x^{a+d+g-1}\hat y^{b+e+h-1}\hat z^{c+f+i-1}.
\end{equation}

Let $\overline e$ be a common denominator of the $\overline a_{ij}$. There exists an affine neighborhood $U=\text{spec}(A)$ of $p^*$ such that $\tilde x,\tilde y,\tilde z$ are uniformizing parameters on $U$, and $\gamma_1,\gamma_2$ are units on $U$. Let
$$
B=A[\gamma_1^{\frac{1}{\overline e}},\gamma_2^{\frac{1}{\overline e}}],
$$
$V=\text{spec}(B)$. After possibly shrinking $U$ to a smaller neighorhood of $p^*$, we may assume that $V$ is an etale cover of $U$ and $x,y,z$ are uniformizing parameters on $V$.
Let $\Lambda:V\rightarrow U$ be the natural morphism.

Suppose that $p\in D_X\cap U$ is a 2-point. We will show that $\tau_f(\Lambda(p))\le\tau_f(p^*)$. After possibly interchanging $x,y,z$, we may assume that ${\cal O}_{V,p}$  has regular parameters $x,y,\tilde z=z-\alpha$ for some $0\ne\alpha\in{\bf k}$. After interchanging $u,v$ if necessary, we have three possible cases:
\begin{enumerate}
\item[1.]  $ae-bd\ne0$, 
\item[2.] $ae-bd\ne0$, $(d,e)\ne(0,0)$,
\item[3.]  $(d,e)=0$. 
\end{enumerate}

We will analyze these three cases in turn.

\vskip .2truein
\noindent {\bf Suppose that $ae-bd\ne 0$}, 
Define  (formal) regular parameters $\overline x,\overline y, \overline z$ at $p$ by choosing $\lambda_1,\lambda_2\in{\bf Q}$ so that 
$$
x=\overline x z^{\lambda_1},
y=\overline yz^{\lambda_2}, z=\overline z+\alpha
$$
satisfy
$$
u=\overline x^a\overline y^b,
v=\overline x^d\overline y^e.
$$

We thus have an expression 
\begin{equation}\label{eqT190}
u=\overline x^a\overline y^b,
v=\overline x^d\overline y^e,
w=P(\overline x,\overline y)+\overline x^m\overline y^n\Omega
\end{equation}
for some series $P(\overline x,\overline y)$ and $\Omega\in\hat{\cal O}_{X,p}$,
where $\frac{\partial\Omega}{\partial z}(0,0,0)\ne 0$.

We compute the Jacobian determinant

$$
\text{Det}\left(\begin{array}{l}
\partial(u,v,w)\\
\partial(\overline x,\overline y,\overline z)
\end{array}
\right)=(ae-bd)\frac{\partial\Omega}{\partial \overline z}\overline x^{a+d+m-1}\overline y^{b+e+n-1}
$$
so that $\frac{\partial \Omega}{\partial \overline z}(0,0,0)\ne 0$.
Comparing with (\ref{eqT197}), we see that  $m=g$ and $n=h$.

 We compute that
 $$
 \frac{\partial}{\partial \overline x}=z^{\lambda_1}\frac{\partial}{\partial  x},
\frac{\partial}{\partial \overline y}=z^{\lambda_2}\frac{\partial}{\partial  y},
\frac{\partial}{\partial \overline z}=\lambda_1\frac{x}{z}\frac{\partial}{\partial x}+\lambda_2\frac{y}{z}\frac{\partial}{\partial y}+\frac{\partial}{\partial  z}.
$$
Thus
$$
\frac{\partial^{m+n} w}{\partial\overline x^m\partial \overline y^n}=z^{m\lambda_1+n\lambda_2}
\frac{\partial^{m+n}w}{\partial x^m\partial y^n}.
$$

From (\ref{eqT191}) we see that $\frac{\partial^{m+n} w}{\partial x^m\partial y^n}\in (x,y)$ if $(m,n)\ne (a_i,b_i)$ for any $i$, and $(m,n)\not\ge (g,h)$.  We have an expansion
$$
P(\overline x,\overline y)=\sum_{(m,n)\not\ge(g,h)}\frac{1}{m!n!}\frac{\partial^{m+n} w}{\partial \overline x^m\partial \overline y^n}(p)\overline x^m\overline y^n.
$$
Thus there is an expansion
$$
P=\sum\beta_i\overline x^{a_i}\overline y^{b_i}
$$
with 
$$
\beta_i=\alpha^{m\lambda_1+n\lambda_2}\frac{\partial^{m+n}w}{\partial x^m\partial y^n}(p).
$$

We conclude that we may make a formal change of variables in (\ref{eqT190}), setting $\beta=\Omega(0,0,0)$ and $z^*=\Omega-\beta$, to get 
\begin{equation}\label{eqT192}
\begin{array}{ll}
u&=\overline x^a\overline y^b\\
v&=\overline x^d\overline y^e\\
w&=\sum\beta_i\overline x^{a_i}\overline y^{b_i}+\overline x^g\overline y^h(z^*+\beta).
\end{array}
\end{equation}

We now compare $\tau_f(p)$ to $\tau_f(p^*)$. Let $q=f(p)$.

Suppose that $q^*$ is a 3-point. If $q$ is a 2-point, then either $(a_0,b_0)=(0,0)$, or
$g=h=0$ in (\ref{eqT192}). $(a_0,b_0)=(0,0)$ implies $ae-bd=0$ which is not possible. Thus $g=h=0$, and we see that $\tau_f(p)=-\infty\le\tau_f(p^*)$. Suppose that $q$ is a 3-point. Then we compute $\tau_f(p)$ from (\ref{eqT192}) and (\ref{eqT10}), to get
$$
\begin{array}{ll}
\tau_f(p)&= \ell\left(((a,b){\bf Z}+(d,e){\bf Z}+\sum_{\beta_i\ne 0, (a_i,b_i)\not\ge (g,h)}(a_i,b_i){\bf Z})/((a,b){\bf Z}+(d,e){\bf Z}+(a_0,b_0){\bf Z})\right)\\
&\le\tau_f(p^*).
\end{array}
$$

Suppose that $q^*$ is a 2-point. Then $q$ is a 2-point, and we have from (\ref{eqT192}) and (\ref{eqT10}) that
$$
\begin{array}{ll}
\tau_f(p)&= \ell\left(((a,b){\bf Z}+(d,e){\bf Z}+\sum_{\beta_i\ne 0, (a_i,b_i)\not\ge (g,h)}(a_i,b_i){\bf Z})/((a,b){\bf Z}+(d,e){\bf Z})\right)\\
&\le\tau_f(p^*).
\end{array}
$$

\vskip .2truein
\noindent {\bf Suppose that $ae-bd\ne0$ and $(d,e)\ne(0,0)$.}

There exist $\overline a,\overline b,t,k$ such that 
$(a,b)=k(\overline a,\overline b)$, $(d,e)=t(\overline a,\overline b)$ with $\text{gcd}(\overline a,\overline b)=1$, $k,t\ne 0$.

Define  (formal) regular parameters $\overline x,\overline y, \overline z$ at $p$ by choosing $\lambda_1,\lambda_2,\lambda_3\in{\bf Q}$ so that 
$$
x=\overline x z^{\lambda_1},
y=\overline yz^{\lambda_2}, z^{\lambda_3}=\overline z+\overline\alpha
$$
satisfy
$$
u=(\overline x^{\overline a}\overline y^{\overline b})^k,
v=(\overline x^{\overline a}\overline y^{\overline b})^t(\overline z+\overline\alpha).
$$

We thus have an expression 
\begin{equation}\label{eqT193}
u=(\overline x^{\overline a}\overline y^{\overline b})^k,
v=(\overline x^{\overline a}\overline y^{\overline b})^t(\overline z+\overline\alpha),
w=P(\overline x^{\overline a}\overline y^{\overline b},\overline z)+\overline x^m\overline y^n\Omega
\end{equation}
where $\Omega$ is a unit series and $\overline a n-\overline b m\ne 0$.

We compute the Jacobian determinant

$$
\text{Det}\left(\begin{array}{l}
\partial(u,v,w)\\
\partial(\overline x,\overline y,\overline z)
\end{array}
\right)=(\overline x^{\overline a}\overline y^{\overline b})^{t+k-1}\overline x^{m+\overline a-1}\overline y^{n+\overline b-1}\gamma,
$$
where 
$$
\gamma =k(\overline b m-\overline a n)\Omega+k\overline b\overline x\frac{\partial\Omega}{\partial\overline x}-k\overline a\overline y\frac{\partial\Omega}{\partial \overline y}
$$
 is a unit series. Comparing with (\ref{eqT197}), we see that
 $m=g$ and $n=h$.

 We compute that
 $$
 \frac{\partial}{\partial \overline x}=z^{\lambda_1}\frac{\partial}{\partial  x},
\frac{\partial}{\partial \overline y}=z^{\lambda_2}\frac{\partial}{\partial  y},
\frac{\partial}{\partial \overline z}=\frac{\lambda_1}{\lambda_3}xz^{-\lambda_3}\frac{\partial}{\partial x}+\frac{\lambda_2}{\lambda_3}yz^{-\lambda_3}\frac{\partial}{\partial y}+\frac{1}{\lambda_3}z^{1-\lambda_3}\frac{\partial}{\partial z}.
$$
Thus
$$
\frac{\partial^{m+n} w}{\partial\overline x^m\partial \overline y^n}=z^{m\lambda_1+n\lambda_2}
\frac{\partial^{m+n}w}{\partial x^m\partial y^n}.
$$

From (\ref{eqT191}) we see that $\frac{\partial^{m+n} w}{\partial x^m\partial y^n}\in (x,y)$ if $(m,n)\ne (a_i,b_i)$ for any $i$, and $(m,n)\not\ge (g,h)$.  
There is thus an expansion
$$
P(\overline x^{\overline a}\overline y^{\overline b},z)=\sum_{j,k}\beta_{j,k}(\overline x^{\overline a}\overline y^{\overline b})^j\overline z^k
$$
with 
$$
\beta_{j,k}=\frac{1}{k!(j\overline a)!(j\overline b)!}\frac{\partial^k}{\partial \overline z^k}\left(\frac{\partial^{j(\overline a+\overline b)} w}{\partial\overline x^{j\overline a}\partial\overline y^{j\overline b}}\right)(p).
$$
Thus $\beta_{j,k}=0$ if $(j\overline a,j\overline b)\ne(a_i,b_i)$ for some $i$, and $(j\overline a,j\overline b)\not\ge(g,h)$.

We conclude that we may make a formal change of variables in (\ref{eqT193}), setting 
$$
\overline x=x^*\gamma_1\text{ and }\overline y= y^*\gamma_2,
$$
with appropriate unit series $\gamma_1$ and $\gamma_2$, to get 
\begin{equation}\label{eqT194}
\begin{array}{ll}
u&=((x^*)^{\overline a}(y^*)^{\overline b})^k\\
v&=(x^*)^{\overline a}(y^*)^{\overline b})^t(\overline z+\overline\alpha)\\
w&=P((x^*)^{\overline a}(y^*)^{\overline b},\overline z)+(x^*)^g(y^*)^h.
\end{array}
\end{equation}

We now compare $\tau_f(p)$ to $\tau_f(p^*)$. Let $q=f(p)$.

Suppose that $q^*$ is a 3-point.  Suppose that $q$ is a 3-point. We compute $\tau_f(p)$ from (\ref{eqT194}) and (\ref{eqT11}), to get
$$
\tau_f(p)= \ell\left((k{\bf Z}+t{\bf Z}+\sum_{\beta_{j,k}\ne 0, (j\overline a,j\overline b)\not\ge (g,h)}j{\bf Z})/(k{\bf Z}+t{\bf Z}+j_0{\bf Z})\right),
$$
where $j_0(\overline a,\overline b)=(a_0,b_0)$.
We have a surjection
$$
\begin{array}{l}
((a,b,c){\bf Z}+(d,e,f){\bf Z}+\sum (a_i,b_i,c_i){\bf Z})/((a,b,c){\bf Z}+(d,e,f){\bf Z}+(a_0,b_0,c_0){\bf Z})\\
\rightarrow (k{\bf Z}+t{\bf Z}+\sum_{j(\overline a,\overline b)=(a_i,b_i)}j{\bf Z})/(k{\bf Z}+t{\bf Z}+j_0{\bf Z})
\end{array}
$$
defined by $j(\overline a,\overline b,0)+l(0,0,1)\mapsto j$.
Thus 
$\tau_f(p)\le\tau_f(p^*)$.

If $q$ is a 2-point, then we must have that $(a_0,b_0)=(0,0)$ or $g=h=0$. But we cannot have $g=h=0$ since
$$
\text{Det}\left(\begin{array}{lll}
a&b&c\\
d&e&f\\
g&h&i
\end{array}\right)\ne 0.
$$

   If $(a_0,b_0)=(0,0)$, then we have $(a_0,b_0,c_0)=(0,0,c_0)$, from which we conclude that $\tau_f(p)\le \tau_f(p^*)$.

Suppose that $q^*$ is a 2-point. Then $q$ is a 2-point.

We compute $\tau_f(p)$ from (\ref{eqT194}) and (\ref{eqT11}), to get
$$
\tau_f(p)= \ell\left((k{\bf Z}+t{\bf Z}+\sum_{\beta_{j,k}\ne 0, (j\overline a,j\overline b)\not\ge (g,h)}j{\bf Z})/(k{\bf Z}+t{\bf Z})\right).
$$
As in the case when $q^*$ is a 3-point, we conclude that $\tau_f(p)\le\tau_f(p^*)$.

\vskip .2truein
\noindent {\bf Suppose that $d=e=0$.} We have $(a,b)=k(\overline a,\overline b)$ with $\overline a,\overline b>0$ and $\text{gcd}(\overline a,\overline b)=1$.

Define   regular parameters $\overline x,\overline y, \overline z$ in $\hat{\cal O}_{X,p}$ by choosing $\lambda_1,\lambda_2\in{\bf Q}$ and $\overline\alpha=\alpha^f$, so that 
$$
x=\overline x z^{\lambda_1},
y=\overline yz^{\lambda_2}, z^{f}=\overline z+\overline\alpha
$$
satisfy
$$
u=(\overline x^{\overline a}\overline y^{\overline b})^k,
v=\overline z+\overline\alpha.
$$

We thus have an expression 
\begin{equation}\label{eqT195}
u=(\overline x^{\overline a}\overline y^{\overline b})^k,
\overline v=v-\overline\alpha=\overline z,
w=P(\overline x^{\overline a}\overline y^{\overline b},\overline z)+\overline x^m\overline y^n\Omega
\end{equation}
where $u,\overline v$ are toroidal forms at $q=f(p)$, $P$ is a series, $\Omega$ is a unit series and $\overline a n-\overline b m\ne 0$.

We compute the Jacobian

$$
\text{Det}\left(\begin{array}{l}
\partial(u,v,w)\\
\partial(\overline x,\overline y,\overline z)
\end{array}
\right)=\overline x^{a+m-1}\overline y^{b+n-1}\gamma.
$$
where 
$$
\gamma=(mb-an)\Omega+b\overline x\frac{\partial \Omega}{\partial \overline x}-a\overline y\frac{\partial \Omega}{\partial \overline y}
$$
 is a unit series. Comparing with (\ref{eqT197}), we see that
 $m=g$ and $n=h$.

 We compute that
 $$
 \frac{\partial}{\partial \overline x}=z^{\lambda_1}\frac{\partial}{\partial  x},
\frac{\partial}{\partial \overline y}=z^{\lambda_2}\frac{\partial}{\partial  y},
\frac{\partial}{\partial \overline z}=\frac{\lambda_1}{f}\frac{x}{z^f}\frac{\partial}{\partial x}+\frac{\lambda_2}{f}\frac{y}{z^f}\frac{\partial}{\partial y}+\frac{1}{f}z^{1-f}\frac{\partial}{\partial z}.
$$
Thus
$$
\frac{\partial^{m+n} w}{\partial\overline x^m\partial \overline y^n}=z^{m\lambda_1+n\lambda_2}
\frac{\partial^{m+n}w}{\partial x^m\partial y^n}.
$$

From (\ref{eqT191}) we see that $\frac{\partial^{m+n} w}{\partial x^m\partial y^n}\in (x,y)$ if $(m,n)\ne (a_i,b_i)$ and $(m,n)\not\ge (g,h)$.  
There is an expansion
$$
P=\sum_{j,k}\beta_{j,k}(\overline x^{\overline a}\overline y^{\overline b})^j\overline z^k
$$
with 
$$
\beta_{j,k}=\frac{1}{k!(j\overline a)!(j\overline b)!}\frac{\partial^k}{\partial \overline z^k}\left(\frac{\partial^{j(\overline a+\overline b)} w}{\partial\overline x^{j\overline a}\partial\overline y^{j\overline b}}\right)(p).
$$
Thus $\beta_{j,k}=0$ if $(j\overline a,j\overline b)\ne(a_i,b_i)$ and $(j\overline a,j\overline b)\not\ge(g,h)$.

We conclude that we may make a formal change of variables in (\ref{eqT195}), setting 
$$
\overline x=x^*\gamma_1\text{ and }\overline y=y^*\gamma_2,
$$
with appropriate unit series $\gamma_1$ and $\gamma_2$, to get 
\begin{equation}\label{eqT196}
\begin{array}{ll}
u&=((x^*)^{\overline a}(y^*)^{\overline b})^k\\
\overline v&=v-\overline\alpha=\overline z\\
w&=P((x^*)^{\overline a}(y^*)^{\overline b},\overline z)+(x^*)^g(y^*)^h.
\end{array}
\end{equation}

We now compare $\tau_f(p)$ to $\tau_f(p^*)$. Let $q=f(p)$.

Suppose that $q^*$ is a 3-point. Then 
$q$ is a 1-point or a 2-point. If $q$ is a 2-point, then $u,w,\overline v$ are permissible parameters at $q$ such that there is  an expression
$$
\begin{array}{ll}
u&=((x^*)^{\overline a}(y^*)^{\overline b})^k\\
w&=P((x^*)^{\overline a}(y^*)^{\overline b},\overline z)+(x^*)^g(y^*)^h\\
\overline v&=\overline z
\end{array}
$$ 
of the form of 2 (c) of Definition \ref{Def31}, and (\ref{eqT12}).

We have
$$
\tau_f(p)=\ell\left(
(k{\bf Z}+\sum_{\beta_{j,k}\ne 0, (j\overline a,j\overline b)\not\ge (g,h)}j{\bf Z})/(k{\bf Z}+j_0{\bf Z})\right),
$$
where $j_0(\overline a,\overline b)=(a_{i_0},b_{i_0})$.

We thus have a surjection

$$
\begin{array}{l}
((a,b,c){\bf Z}+(d,e,f){\bf Z}+\sum (a_i,b_i,c_i){\bf Z})/((a,b,c){\bf Z}+(d,e,f){\bf Z}+(a_0,b_0,c_0){\bf Z})\\
\rightarrow 
(k{\bf Z}+\sum_{j(\overline a,\overline b)=(a_i,b_i)}j{\bf Z})/(k{\bf Z}+j_0{\bf Z})
\end{array}
$$
defined by $j(\overline a,\overline b,0)+l(0,0,1)\mapsto j$.
Thus 
$\tau_f(p)\le\tau_f(p^*)$.

If $q^*$ is a 3-point and $q$ is a 1-point, then  we have permissible parameters $u,\overline v,\overline w$ at $q$, defined by
$$
\begin{array}{ll}
u&=((x^*)^{\overline a}(y^*)^{\overline b})^k\\
\overline v&=\overline z\\
\overline w&=w-\beta=P((x^*)^{\overline a}(y^*)^{\overline b},\overline z)-\beta+(x^*)^g(y^*)^h
\end{array}
$$ 
of the form (\ref{eqTF02}) and (\ref{eqT13}), where $0\ne \beta=P(0,0)$ so that $(a_{j_0},b_{j_0})=(0,0)$.

We have
$$
\tau_f(p)=\ell\left(
k{\bf Z}+\sum_{\beta_{j,k}\ne 0, (j\overline a,j\overline b)\not\ge (g,h)}j{\bf Z})/k{\bf Z}\right).
$$

We  have a surjection

$$
\begin{array}{l}
((a,b,c){\bf Z}+(d,e,f){\bf Z}+\sum (a_i,b_i,c_i){\bf Z})/((a,b,c){\bf Z}+(d,e,f){\bf Z}+(a_0,b_0,c_0){\bf Z})\\
\rightarrow 
(k{\bf Z}+\sum_{j(\overline a,\overline b)=(a_i,b_i)}j{\bf Z})/(k{\bf Z}+j_0{\bf Z})
\end{array}
$$
defined by $j(\overline a,\overline b,0)+l(0,0,1)\mapsto j$.
Thus 
$\tau_f(p)\le\tau_f(p^*)$.

There is a similar analysis if $q^*$ is a 2-point.

We have completed the analysis for a 2-point $p\in V$. 

Now suppose that $p\in D_X\cap V$ is a 1-point. Recall that we are assuming that  $p^*$ is a 3-point.
We will show that $\tau_f(\Lambda(p))\le \tau_f(p^*)$.
After possibly interchanging  $x,y,z$, we may suppose that ${\cal O}_{V,p}$ has regular parameters $x,y-\alpha, z-\beta$ for some $0\ne\alpha\in{\bf k}$, $0\ne\beta\in{\bf k}$. After  interchanging $u,v$ and $y,z$ if necessary, we have two possible cases: $a,d\ne 0$, $e\ne 0$  and $a\ne 0$, $d=0$, $e\ne0$. We will analyze these two cases in turn.

\vskip .2truein
\noindent {\bf Suppose that $a,d\ne 0$ and $e\ne 0$.}

Define  regular parameters $\overline x,\overline y,\overline z$ in $\hat{\cal O}_{X,p}$ by choosing $\lambda_{ij}\in {\bf Q}$, $\overline\alpha=\alpha^{\lambda_{21}}\beta^{\lambda_{22}}$, so that
$$
x=\overline xy^{\lambda_{11}}z^{\lambda_{12}},
\overline y +\overline \alpha=y^{\lambda_{21}}z^{\lambda_{22}},
\overline z + \beta = z
$$
satisfy
$$
u=\overline x^a,
v=\overline x^d(\overline y+\overline \alpha).
$$
We thus have an expression 
\begin{equation}\label{eqT198}
\begin{array}{ll}
u&=\overline x^a\\
v&=\overline x^d(\overline y+\overline \alpha)\\
w&=P(\overline x,\overline y)+\overline x^m\Omega
\end{array}
\end{equation}
where $\frac{\partial \Omega}{\partial \overline z}(0,0,0)\ne 0$.

We compute the Jacobian
$$
\text{Det}\left(\begin{array}{l}
\partial(u,v,w)\\
\partial(\overline x,\overline y,\overline z)
\end{array}
\right)=\overline x^{ a+d+m-1}\frac{\partial\Omega}{\partial \overline z}.
$$
 Comparing with (\ref{eqT197}), we see that
 $m=g$.

We have
$$
\begin{array}{ll}
x&=\overline x (\overline y+\overline\alpha)^{\frac{\lambda_{11}}{\lambda_{21}}}
(\overline z+\beta)^{\lambda_{12}-\frac{\lambda_{11}\lambda_{22}}{\lambda_{21}}}\\
y&=(\overline y+\overline\alpha)^{\frac{1}{\lambda_{21}}}
(\overline z+\beta)^{-\frac{\lambda_{22}}{\lambda_{21}}}\\
z&=\overline z+\beta,
\end{array}
$$
from which we see that

 $$
 \begin{array}{ll}
 \frac{\partial}{\partial \overline x}&=y^{\lambda_{11}}z^{\lambda_{12}}\frac{\partial}{\partial x}\\
\frac{\partial}{\partial \overline y}&=
\frac{\lambda_{11}}{\lambda_{21}}xy^{-\lambda_{21}}z^{-\lambda_{22}}\frac{\partial}{\partial x}+\frac{1}{\lambda_{21}}y^{1-\lambda_{21}}z^{-\lambda_{22}}\frac{\partial}{\partial y}\\
\frac{\partial}{\partial \overline z}&=
\left(\frac{\lambda_{12}\lambda_{21}-\lambda_{11}\lambda_{22}}{\lambda_{21}}\right)xz^{-1}\frac{\partial}{\partial x}-\frac{\lambda_{22}}{\lambda_{21}}yz^{-1}\frac{\partial}{\partial y}+\frac{\partial}{\partial z}.
\end{array}
$$
Thus

$$
\frac{\partial^{m} w}{\partial\overline x^m}=y^{m\lambda_{11}}z^{m\lambda_{12}}
\frac{\partial^{m}w}{\partial x^m}.
$$

From (\ref{eqT191}) we see that $\frac{\partial^m w}{\partial x^m}\in (x)$ if $m\ne a_i$ for some $i$ and $m\not\ge g$.  
There is an expansion
$$
P=\sum_{j,k}\beta_{j,k}\overline x^{j}\overline y^{k}
$$
with 
$$
\beta_{j,k}=\frac{1}{k!j!}\frac{\partial^k}{\partial \overline y^k}\left(\frac{\partial^j w}{\partial\overline x^{j}}\right)(0,0,0).
$$
Thus $\beta_{j,k}=0$ if $j\ne a_i$ and $j\not\ge g$.

We may make a formal change of variables in (\ref{eqT198}), setting $\overline \beta=\Omega(0,0,0)$ and $z^*=\Omega-\overline\beta$, to get 
\begin{equation}\label{eqT199}
\begin{array}{ll}
u&=\overline x^a\\
v&=\overline x^d(\overline y+\overline \alpha)\\
w&=\sum \beta_{j,k}\overline x^j\overline y^k+\overline x^g(z^*+\overline \beta).
\end{array}
\end{equation}

We now compare $\tau_f(p)$ to $\tau_f(p^*)$. Let $q=f(p)$.

Suppose that $q^*$ is a 3-point.  If $q$ is a 2-point then we have that $a_0=0$, and
$$
\tau_f(p)=
\ell\left( (a{\bf Z}+d{\bf Z}+\sum_{j<g} j{\bf Z})/(a{\bf Z}+d{\bf Z})\right)\le\tau_f(p^*).
$$
If $q^*$ is a 3-point and $q$ is a 3-point, then 
$$
\tau_f(p)=
\ell\left( (a{\bf Z}+d{\bf Z}+\sum_{\beta_{jk}\ne 0,j<g} j{\bf Z})/(a{\bf Z}+d{\bf Z}+a_0{\bf Z})\right)\le\tau_f(p^*).
$$

Suppose that $q^*$ is a 2-point. Then $q$ is a 2-point, and by a similar calculation, $\tau_f(p)\le\tau_f(p^*)$.

\vskip .2truein
\noindent {\bf Suppose that $a\ne0$, $d= 0$, $e\ne 0$.}

Define (formal) regular parameters $\overline x,\overline y,\overline z$ at $p$ by choosing $\lambda_{1j}\in {\bf Q}$ and $\overline\alpha=\alpha^e\beta^f$ so that
$$
x=\overline xy^{\lambda_{11}}z^{\lambda_{12}},
\overline y +\overline \alpha=y^{e}z^{f},
\overline z + \beta = z
$$
satisfy
$$
u=\overline x^a,
\overline v=v-\overline \alpha=\overline y.
$$
We thus have an expression 
\begin{equation}\label{eqT200}
\begin{array}{ll}
u&=\overline x^a\\
\overline v&=v-\overline\alpha=\overline  y\\
w&=P(\overline x,\overline y)+\overline x^m\Omega
\end{array}
\end{equation}
where $\frac{\partial \Omega}{\partial \overline z}(0,0,0)\ne 0$.

We compute the Jacobian
$$
\text{Det}\left(\begin{array}{l}
\partial(u,v,w)\\
\partial(\overline x,\overline y,\overline z)
\end{array}
\right)=\overline x^{ a+m-1}\frac{\partial\Omega}{\partial\overline z}.
$$
 Comparing with (\ref{eqT197}), we see that
 $m=g$.

We have
$$
\begin{array}{ll}
x&=\overline x (\overline y+\overline\alpha)^{\frac{\lambda_{11}}{e}}
(\overline z+\beta)^{\lambda_{12}-\frac{\lambda_{11}f}{e}}\\
y&=(\overline y+\overline\alpha)^{\frac{1}{e}}
(\overline z+\beta)^{-\frac{f}{e}}\\
z&=\overline z+\beta,
\end{array}
$$
from which we see that 

 $$
 \begin{array}{ll}
 \frac{\partial}{\partial \overline x}&=y^{\lambda_{11}}z^{\lambda_{12}}\frac{\partial}{\partial x}\\
\frac{\partial}{\partial \overline y}&=
\frac{\lambda_{11}}{e}xy^{-e}z^{-f}\frac{\partial}{\partial x}+\frac{1}{e}y^{1-e}z^{-f}\frac{\partial}{\partial y}\\
\frac{\partial}{\partial \overline z}&=
\frac{(\lambda_{12}e-\lambda_{11}f)}{e}xz^{-1}\frac{\partial}{\partial x}-\frac{f}{e}yz^{-1}\frac{\partial}{\partial y}+\frac{\partial}{\partial z}
.\end{array}
$$
Thus
$$
\frac{\partial^{m} w}{\partial\overline x^m}=y^{m\lambda_{11}}z^{m\lambda_{12}}
\frac{\partial^{m}w}{\partial x^m}.
$$

From (\ref{eqT191}) we see that $\frac{\partial^m w}{\partial x^m}\in (x)$ if $m\ne a_i$ for some $i$ and $m\not\ge g$.  
There is an expansion
$$
P=\sum_{j,k}\beta_{j,k}\overline x^{j}\overline y^{k}
$$
with 
$$
\beta_{j,k}=\frac{1}{k!j!}\frac{\partial^k}{\partial \overline y^k}\left(\frac{\partial^j w}{\partial\overline x^{j}}\right)(0,0,0).
$$
Thus $\beta_{j,k}=0$ if $j\ne a_i$ and $j\not\ge g$.

We may make a formal change of variables in (\ref{eqT200}), setting $\overline \beta=\Omega(0,0,0)$ and $z^*=\Omega-\overline\beta$, to get 
\begin{equation}\label{eqT201}
\begin{array}{ll}
u&=\overline x^a\\
\overline v&=\overline y\\
w&=\sum \beta_{j,k}\overline x^j\overline y^k+\overline x^g(z^*+\overline \beta).
\end{array}
\end{equation}

We now compare $\tau_f(p)$ to $\tau_f(p^*)$. Let $q=f(p)$.

Suppose that $q^*$ is a 3-point.  If $q$ is a 1-point then we have that $a_0=0$, and
$u,\overline v,\overline w=w-w(q)$ are permissible parameters at $q$ which have a form (\ref{eqTF01}) and (\ref{eqT17}) at $p$.
$$
\tau_f(p)=
\ell\left( (a{\bf Z}+\sum_{\beta_{jk}\ne 0, j<g} j{\bf Z})/a{\bf Z}\right)\le\tau_f(p^*).
$$
If $q$ is a 2-point, then $u,w,\overline v$ are premissible parameters at $q$ which have a form 2 (b) of Definition \ref{Def31} and (\ref{eqT16}), and
$$
\tau_f(p)=
\ell\left( (a{\bf Z}+\sum_{\beta_{jk}\ne 0, j<g} j{\bf Z})/(a{\bf Z}+a_0{\bf Z})\right)\le\tau_f(p^*).
$$

Suppose that $q^*$ is a 2-point. Then $q$ is a 1-point, and by a similar calculation, $\tau_f(p)\le\tau_f(p^*)$.

This completes the analysis that $\tau_f(p)\le \tau_f(p^*)$ if $p^*$ is a 3-point and $p\in U$.

The analysis when $p^*$ is a 2-point or a 1-point is simpler, and we leave it to the reader.

We conclude that $\tau_f$ is upper semi-continuous.

\end{pf}

\begin{Definition}\label{Def326} Suppose that $f:X\rightarrow Y$ is a prepared, proper morphism, and $\tau\in{\bf N}$ is such that $\tau_f(X)\le\tau$. Let $G_X(f,\tau)=\{p\in X\mid \tau_f(X)=\tau\}$,
$G_Y(f,\tau)=f(G_X(f,\tau))$. We will say that $f$ is $\tau$-prepared if $G_Y(f,\tau)$ contains no 2-curves and no 3-points.
\end{Definition}

By Lemma \ref{LemmaT111}, $G_X(f,\tau)$ is a closed subset of $X$, and since $f$ is proper, $G_Y(f,\tau)$ is a closed subset of $Y$.

\section{Super Parameters}

Throughout this section, we assume that $f:X\rightarrow Y$ is a dominant, proper morphism of nonsingular 3-folds.

\begin{Definition}\label{Def357} Suppose that $f:X\rightarrow Y$ is prepared, and $q\in Y$ is a
2-point. Permissible parameters $u,v,w$ at $q$ are {\bf super parameters} for $f$ at $q$ if 
at all
$p\in f^{-1}(q)$, there exist permissible parameters $x,y,z$ for $u,v,w$ at $p$ such that we have one of the forms: 
\begin{enumerate}
\item[1.] $p$ is a 1-point 
\begin{equation}\label{eqT221}
\begin{array}{ll}
u&=x^a\\
v&=x^b(\alpha+y)\\
w&=x^c\gamma(x,y)+x^d(z+\beta)
\end{array}
\end{equation}
where $\gamma$ is a unit series (or zero), $0\ne\alpha\in \bold k$ and $\beta\in \bold k$,
\item[2.] $p$ is  a 2-point of the form of  (\ref{eqTF21}) of Definition \ref{torf} 
\begin{equation}\label{eqT222}
\begin{array}{ll}
u&=x^ay^b\\
v&=x^cy^d\\
w&=x^ey^f\gamma(x,y)+x^gy^h(z+\beta)
\end{array}
\end{equation}
where $ad-bc\ne 0$, $\gamma$ is a unit series (or zero), and $\beta\in \bold k$.
\item[3.] $p$ is a 2-point of the form of  (\ref{eqTF22}) of Definition \ref{torf} 
\begin{equation}\label{eqT223}
\begin{array}{ll}
u&=(x^ay^b)^k\\
v&=(x^ay^b)^t(\alpha+z)\\
w&=(x^ay^b)^l\gamma(x^ay^b,z)+x^cy^d
\end{array}
\end{equation}
where $0\ne\alpha\in \bold k$, $ad-bc\ne 0$ and $\gamma$ is a unit series (or zero).
\item[4.] $p$ is a 3-point 
\begin{equation}\label{eqT224}
\begin{array}{ll}
u&=x^ay^bz^c\\
v&=x^dy^ez^f\\
w&=x^gy^hz^i\gamma+x^jy^kz^l
\end{array}
\end{equation}
where $\text{rank}(u,v,x^jy^kz^l)=3$, $\text{rank}(u,v,x^gy^hz^i)=2$
and
$\gamma$ is a unit series in monomials $M$ such that $\text{rank}(u,v,M)=2$ (or $\gamma$ is zero).
\end{enumerate}

Suppose that  $q\in Y$ is a 
1-point. Permissible parameters $u,v,w$ at $q$ are {\bf super parameters} for $f$ at $q$ if 
at all
$p\in f^{-1}(q)$, there exist permissible parameters $x,y,z$ for $u,v,w$ at $p$ such that we have one of the forms:
\begin{enumerate}
\item[5.] $p$ is a 1-point 
\begin{equation}\label{eqT280}
\begin{array}{ll}
u&=x^a\\
v&=y\\
w&=x^c\gamma(x,y)+x^d(z+\beta)
\end{array}
\end{equation}
where $\gamma$ is a unit series (or zero) and $\beta\in \bold k$,
\item[6.] $p$ is a 2-point 
\begin{equation}\label{eqT279}
\begin{array}{ll}
u&=(x^ay^b)^k\\
v&=z\\
w&=(x^ay^b)^l\gamma(x^ay^b,z)+x^cy^d
\end{array}
\end{equation}
where  $ad-bc\ne 0$ and $\gamma$ is a unit series (or zero).
\end{enumerate}
\end{Definition}

\begin{Lemma}\label{LemmaT204} Suppose that $f:X\rightarrow Y$ is prepared, and $\Phi:X_1\rightarrow X$ is the blow up of a 2-curve or a 3-point. Let $f_1=f\circ\Phi:X_1\rightarrow Y$. Then $f_1$ is prepared, and $\tau_{f_1}(p)\le\tau_f(\Phi(p))$ for all $p\in D_{X_1}$.

If $q\in Y$ and $u,v,w$ are super parameters for $f$ at $q$, then $u,v,w$ are super parameters for $f_1$ at $q$.
\end{Lemma}

\begin{pf} We prove this in the case when $\Phi:X_1\rightarrow X$ is the blow up of a 2-curve $C$. The case when $\Phi$ is the blow up of a 3-point is similar.

Suppose that $p\in C$. Let $q=f(p)$. Then there are permissible parameters $u,v,w$ at $q$ and $x,y,z$  for $u,v,w$ at $p$ such that either $u,v$ are toroidal forms at $p$, or a form 2 (c) of Definition \ref{Def31} holds at $q$. Further,  $x=y=0$ are local equations of $C$ at $p$.

The most difficult case is when $p$ is a 3-point, $q=f(p)$ is a 3-point and $\tau_f(p)\ge 0$. The other cases are similar. Assume that this case holds.

We have an expansion of the form of (\ref{eq16}) 
\begin{equation}\label{eqT331}
\begin{array}{ll}
u&=x^ay^bz^c\\
v&=x^dy^ez^f\\
w&=\sum_{(a_i,b_i,c_i)\not\ge (g,h,i)}\alpha_ix^{a_i}y^{b_i}z^{c_i}+x^gy^hz^i
\end{array}
\end{equation}
with $\alpha_i\ne 0$ for all $i$. There are $\phi_i,\psi_i\in{\bf Q}$ such that
$$
(a_i,b_i,c_i)=\phi_i(a,b,c)+\psi_i(d,e,f)
$$
for all $i$.

 Suppose that $p_1\in\Phi^{-1}(p)$. Then (after possibly interchanging $x$ and $y$) $\hat{\cal O}_{X_1,p_1}$ has regular parameters $x_1,y_1,z$ where
$$
x=x_1, y=x_1(y_1+\alpha)
$$
for some $\alpha\in{\bf k}$.
We have
$$
\begin{array}{ll}
u&=x_1^{a+b}(y_1+\alpha)^bz^c\\
v&=x_1^{d+e}(y_1+\alpha)^ez^f\\
w&=\sum_{(a_i,b_i,c_i)\not\ge (g,h,i)}\alpha_ix_1^{a_i+b_i}(y_1+\alpha)^{b_i}z^{c_i}+x_1^{g+h}(y_1+\alpha)^hz^i.
\end{array}
$$
We may assume that $\tau_{f_1}(p_1)\ge 0$.

\vskip .2truein
\noindent{\bf Case 1. Assume that $0\ne\alpha$ and $(a+b)f-c(d+e)\ne 0$.}  There exist regular parameters $\overline x_1,\overline y_1,\overline z_1$ in $\hat{\cal O}_{X_1,p_1}$ and $\overline\beta\in {\bf k}$ such that
$$
\begin{array}{ll}
u&=\overline x_1^{a+b}\overline z_1^{c}\\
v&=\overline x_1^{d+e}\overline z_1^f\\
w&=\sum_{(a_i+b_i,c_i)\not\ge (g+h,i)}\alpha_i\overline x_1^{a_i+b_i}\overline z_1^{c_i}+\overline x_1^{g+h}\overline z_1^i(\overline y_1+\overline\beta)
\end{array}
$$
of the form of  (\ref{eqT10}). The homomorphism $\Lambda:{\bf Z}^3\rightarrow {\bf Z}^2$ defined by
$(x,y,z)\mapsto (x+y,z)$ induces a surjection
$$
\begin{array}{l}
((a,b,c){\bf Z}+(d,e,f){\bf Z}+\sum_{i\ge  0}(a_i,b_i,c_i){\bf Z})/
((a,b,c){\bf Z}+(d,e,f){\bf Z}+(a_0,b_0,c_0){\bf Z})\\
\rightarrow
((a+b,c){\bf Z}+(d+e,f){\bf Z}+\sum_{i\ge 0}(a_i+b_i,c_i){\bf Z})/
((a+b,c){\bf Z}+(d+e,f){\bf Z}+(a_0+b_0,c_0){\bf Z}).
\end{array}
$$
Thus
$$
\begin{array}{ll}
\tau_f(p)&\ge\ell(
((a+b,c){\bf Z}+(d+e,f){\bf Z}+\sum_{i\ge 0}(a_i+b_i,c_i){\bf Z})\\
&\,\,\,/((a+b,c){\bf Z}+(d+e,f){\bf Z}+(a_0+b_0,c_0){\bf Z})
)\\
&\ge \ell(
((a+b,c){\bf Z}+(d+e,f){\bf Z}+\sum_{(a_i+b_i,c_i)\not\ge(g+h,i)}(a_i+b_i,c_i){\bf Z})\\
&\,\,\,/
((a+b,c){\bf Z}+(d+e,f){\bf Z}+(a_0+b_0,c_0){\bf Z}))\\
&=\tau_{f_1}(p_1).
\end{array}
$$
\vskip .2truein
\noindent{\bf Case 2. Assume $0\ne\alpha$, and $(a+b)f-c(d=e)=0$.} Then there exist $\overline a,\overline b\in{\bf N}$ such that  $(a+b,c)=k(\overline a,\overline b)$, $(d+e,f)=t(\overline a,\overline b)$ with $k,t\ne 0$, $\text{gcd}(\overline a,\overline b)=1$.

There exist regular parameters $\overline x_1,\overline y_1,\overline z_1$ in $\hat{\cal O}_{X_1,p_1}$ and $0\ne\overline\alpha\in{\bf k}$ such that
$$
\begin{array}{ll}
u&=\overline x_1^{a+b}\overline z_1^c=(\overline x_1^{\overline a}\overline z_1^{\overline b})^k\\
v&=\overline x_1^{d+e}\overline z_1^f(\overline y_1+\overline\alpha)=
(\overline x_1^{\overline a}\overline z_1^{\overline b})^t(\overline y_1+\overline\alpha)\\
w&=\sum_{(a_i+b_i,c_i)\not\ge(g+h,i)}\alpha_i(\overline x_1^a\overline z_1^b)^{\phi_ik+\psi_it}(\overline y_1+\overline\alpha)^{\psi_i}
+\overline x_1^{g+h}\overline z_1^i
\end{array}
$$
of the form of (\ref{eqT11}).
We have
$$
\tau_{f_1}(p_1)=\ell\left(
k{\bf Z}+t{\bf Z}+\sum_{(a_i+b_i,c_i)\not\ge(g+h,i)}(\phi_ik+\psi_it){\bf Z})/
(k{\bf Z}+t{\bf Z}+(\phi_0k+\psi_0t){\bf Z})\right).
$$
As in the argument of Case 1, we see that
$$
\begin{array}{ll}
\tau_f(p)&\ge \ell( ((a+b,c){\bf Z}+(d+e,f){\bf Z}+\sum_{i\ge 0}(a_i+b_i,c_i){\bf Z})\\
&\,\,\, /((a+b,c){\bf Z}+(d+e,f){\bf Z}+(a_0+b_0,c_0){\bf Z}))\\
&=\ell((k(\overline a,\overline b){\bf Z}+t(\overline a,\overline b){\bf Z}
+\sum_{i\ge 0}(\phi_ik(\overline a,\overline b)+\psi_it(\overline a,\overline b){\bf Z})\\
&\,\,\, /(k(\overline a,\overline b){\bf Z}+t(\overline a,\overline b){\bf Z}+(\phi_0k(\overline a,\overline b)+\psi_0t(\overline a,\overline b)){\bf Z}))\\
&=\ell( (k{\bf Z}+t{\bf Z}+\sum_{i\ge 0}(\phi_ik+\psi_it){\bf Z})\\
&\,\,\, /(k{\bf Z}+t{\bf Z}+(\phi_0k+\psi_0t){\bf Z}))\\
&\ge\tau_{f_1}(p_1).
\end{array}
$$

\vskip .2truein
\noindent{\bf Case 3. $0=\alpha$.} There exist regular parameters $\overline x_1,\overline y_1,\overline z_1$ in $\hat{\cal O}_{X_1,p_1}$ such that 

We have
$$
\begin{array}{ll}
u&=\overline x_1^{a+b}\overline y_1^b \overline z^c\\
v&= \overline x_1^{d+e} \overline y_1^e \overline z^f\\
w&=\sum_{(a_i+b_i,b_i,c_i)\not\ge (g+h,h,i)} \alpha_i \overline x_1^{a_i+b_i} \overline y_1^{b_i} \overline z^{c_i}+ \overline x_1^{g+h} \overline y_1^h \overline z_1^i
\end{array}
$$
of the form of (\ref{eq16}). Thus
$$
\begin{array}{ll}
\tau_{f_1}(p_1)&=\ell(
((a+b,b,c){\bf Z}+(d+e,e,f){\bf Z}+\sum_{(a_i+b_i,b_i,c_i)\not\ge(g+h,h,i)}
(a_i+b_i,b_i,c_i){\bf Z})\\
&/
((a+b,b,c){\bf Z}+(d+e,e,f){\bf Z}+(a_0+b_0,b_0,c_0){\bf Z}))\\
&\le\tau_f(p).
\end{array}
$$

The fact that super parameters for $f$ at $q\in Y$ are super parameters for $f_1$ at $q$ follows from substitution of local equations for $\Phi$ into the forms of Definition \ref{Def357} defining super parameters, and making an appropriate change of variables.

\end{pf}

\begin{Lemma}\label{LemmaT127}
Suppose that $f:X\rightarrow Y$ is prepared, $q\in Y$ and $u,v,w$ are super parameters at $q$. Then there exists a sequence of blow ups of 2-curves $\Phi:X_1\rightarrow X$ such that $f_1=f\circ\Phi:X_1\rightarrow Y$ is prepared, $\tau_{f_1}(p)\le \tau_f(\Phi(p))$ for $p\in D_{X_1}$, $u,v,w$ are super parameters at $q$ for $f_1$ and if $p\in f_1^{-1}(q)$ with $\tau_{f_1}(p)>0$, then $w=0$ is a divisor supported on $D_X$ at $p$.
\end{Lemma}

\begin{pf}
The fact that for a sequence of blow ups of 2-curves $\Phi:X_1\rightarrow Y$, $f_1=f\circ\Phi:X_1\rightarrow Y$ is prepared and $\tau_{f_1}(p)\le\tau_f(\Phi(p))$ for $p\in D_{X_1}$ follows from Lemma \ref{LemmaT204}.

The condition that $u,v,w$ are super parameters at $q$ is preserved by blowup of 2-curves above $X$.

If $p\in f_1^{-1}(q)$ is a 1-point, then 1 or 5 of Definition \ref{Def357} hold, and if $\tau_{f_1}(p)>0$, then  $w=0$ is a divisor supported on $D_X$ at $p$.

  We may construct $\Phi$ so that if $p\in f_1^{-1}(q)$, then
$(x^ey^f,x^gy^h)\hat{\cal O}_{X_1,p}$ is principal if 2 of Definition \ref{Def357} holds,
$((x^ay^b)^l,x^cy^d)\hat{\cal O}_{X,p}$ is principal if 3 holds, $(x^gy^hz^i, x^jy^kz^l)\hat{\cal O}_{X_1,p}$ is principal if 4 holds, $((x^ay^b)^l,x^cy^d)\hat{\cal O}_{X_1,p}$ is principal if 6 holds.

We see that in all these cases that $w=0$ is a divisor supported on $D_X$ at $p$ if $\tau_{f_1}(p)>0$,
 so that the conclusions of the lemma hold.
\end{pf}

\begin{Lemma}\label{Lemma1} Suppose that $f:X\rightarrow Y$ is prepared and
$C\subset Y$ is a 2-curve. Then there exists a commutative diagram
$$
\begin{array}{rll}
X_1&\stackrel{f_1}{\rightarrow}&Y_1\\
\Phi_1\downarrow&&\downarrow\Psi_1\\
X&\stackrel{f}{\rightarrow}&Y
\end{array}
$$
where $\Psi_1:Y_1\rightarrow Y$ is the blow up of $C$, $\Phi_1:X_1\rightarrow X$ is a product of blow ups of 2-curves, 
$\Phi_1$ is an isomorphism above $f^{-1}(Y-C)$ and $f_1$ is prepared. If $p_1\in X_1$ and $p=\Phi_1(p_1)$, then
$$
\tau_{f_1}(p_1)\le \tau_f(p).
$$
If $f$ is $\tau$-prepared then $f_1$ is $\tau$-prepared.
\end{Lemma}

\begin{pf} Most of this is proven in Lemma 5.2 \cite{C5}. It only remains to  check that
$\tau_{f_1}(p_1)\le\tau_f(p)$ when $p_1$ is a 1 or 2-point. 

By Lemma \ref{LemmaT127}, we know that $\tau_{f\circ\Phi_1}(p_1)\le\tau_f(p)$. We verify that $\tau_{f_1}(p_1)\le\tau_{f\circ\Phi_1}(p_1)$.

 This can be seen from consideration of local equations of $f\circ\Phi_1$, $\Phi_1$ and $\Psi_1$, using Lemma \ref{Lemma0} and the methods of Lemma \ref{Lemma0} and Lemma \ref{LemmaT111}.
\end{pf}

\begin{Lemma}\label{LemmaT47} 
Suppose that $f:X\rightarrow Y$ is prepared and $q\in Y$ is a 2-point. Then there exists a commutative diagram 
\begin{equation}\label{eqT48}
\begin{array}{rll}
X_1&\stackrel{f_1}{\rightarrow}&Y_1\\
\Phi\downarrow&&\downarrow\Psi\\
X&\stackrel{f}{\rightarrow}&Y
\end{array}
\end{equation}
where $\Phi$ and $\Psi$ are products of blow ups of 2-curves, such that $f_1$ is prepared, 
$\tau_{f_1}(p)\le\tau_f(\Phi(p))$ for $p\in D_{X_1}$, and there exist no points of form 2 (b) or 2 (c) of Definition \ref{Def31} for $f_1$ above points of $\Psi^{-1}(q)$.
\end{Lemma}

\begin{pf} There exist sequences of blow ups of 2-curves $Y_1\rightarrow Y$ such that the rational map $X\rightarrow Y_1$ is defined at all points $p\in f^{-1}(q)$ such that $f$ has at $p$ an expression of form 2 (b) or 2 (c) of Definition \ref{Def31}, and $p$  maps to 1-point. By Lemma \ref{Lemma1}, by blowing up 2-curves above $X$, we can construct $f_1$ which has the desired property.
\end{pf}

\begin{Definition}\label{DefT40} Suppose that $f:X\rightarrow Y$ is $\tau$-prepared, $p\in D_X$ and $q=f(p)$. Let $u,v,w$ be permissible parameters at $q$. We say that $w$ is  good  at $p$ for $f$ if one of the following expressions holds:

$p$ a 1-point, $q$ a 1-point 
\begin{equation}\label{eqT41}
u=x^a,v=y,w=\sum_{a\not\,\mid i}a_{ij}x^iy^j+x^n(z+\alpha)
\end{equation}

$p$ a 2-point, $q$ a 1-point 
\begin{equation}\label{eqT42}
u=(x^ay^b)^k, v=z, 
w=\sum_{k\not\,\mid i}a_{ij}(x^ay^b)^iz^j+x^cy^d
\end{equation}

$p$ a 1-point, $q$ a 2-point 
\begin{equation}\label{eqT43}
u=x^a, v=x^b(\alpha+y), w=\sum_{d\not\,\mid i}a_{ij}x^iy^j+x^n(z+\alpha)
\end{equation}
where $d=\text{gcd}(a,b)$.

$p$ a 2-point, $q$ a 2-point 
\begin{equation}\label{eqT44}
u=x^ay^b, v=x^cy^d,
w=\sum_{(i,j)\not\in {\bf Z}(a,b)+{\bf Z}(c,d)} a_{ij}x^iy^j+x^ey^fz
\end{equation}

$p$ a 2-point, $q$ a 2-point 
\begin{equation}\label{eqT46}
u=(x^ay^b)^k, v=(x^ay^b)^t(\alpha+z),
w=\sum _{d\not\,\mid i}a_{ij}(x^ay^b)^iz^j+x^cy^d
\end{equation}
where $d=\text{gcd}(k,t)$.

$p$ a 3-point, $q$ a 2-point 
\begin{equation}\label{eqT45}
u=x^ay^bz^c,
v=x^dy^ez^f,
w=\sum a_{ijk}x^iy^jz^k+x^gy^hz^i
\end{equation}

where the sum is over $i,j,k$ such that
$$
\text{Det}\left(\begin{array}{lll}
a&b&c\\
d&e&f\\
i&j&k
\end{array}\right)=0, (i,j,k)\not\in {\bf Z}(a,b,c)+{\bf Z}(d,e,f).
$$

\end{Definition}

\begin{Definition}\label{DefT327} Suppose that $f:X\rightarrow Y$ is $\tau$-prepared,
$p\in D_X$ and $q=f(p)$ is a 1-point. Let $u,v,w$ be permissible parameters at $q$. We say that $w$ is weakly good at $p$ for $f$ if one of the following forms hold:
\begin{enumerate}
\item[1.] $p$ is a 1-point,
$$
u=x^a, v=y,
w=\sum_{j=0}^m a_{\sigma_j}(y)x^{\sigma_j}+x^n(z+\alpha)
$$
where $\alpha\in {\bf k}$, $\sigma_0<\sigma_1<\cdots<\sigma_m<n$, $\sigma_i$ are all nonzero, and $a\not\,\mid\sigma_0$.

\item[2.] $p$ is a 2-point,
$$
u=(x^ay^b)^k, v=z,
w=\sum_{j=0}^m a_{\sigma_j}(z)(x^ay^b)^{\sigma_j}+x^cy^d
$$
where $\text{gcd}(a,b)=1$, $ad-bc\ne 0$, $\sigma_m(a,b)\not >(c,d)$,
$\sigma_0<\sigma_1<\cdots<\sigma_m<n$, $\sigma_i$ are all nonzero, and $k\not\,\mid\sigma_0$.
\end{enumerate}
\end{Definition}

\begin{Remark}\label{RemarkT284} Suppose that $f;X\rightarrow Y$ is $\tau$-prepared and $\tau_f(p)=0$, $u,v,w$ are permissible parameters at $q=f(p)$, and $w$ is good (weakly good) at $p$ for $f$. Then $u,v,w$ are monomial forms (Definition \ref{Def125}) at $p$.
\end{Remark}

\begin{Remark}\label{RemarkT174} Suppose that $f;X\rightarrow Y$ is $\tau$-prepared.
Observe that if $q$ is a 1-point, and $u,v,w$ are permissible parameters at $q$ satisfying the conclusions of Lemma \ref{LemmaT1}, then for all $p\in f^{-1}(q)$, there exists a series $\phi_p(u,v)$ such that $w-\phi_p(u,v)$ is good (weakly good) at $p$ for $f$.
\end{Remark}

\begin{Lemma}\label{LemmaT128}
Suppose that $f:X\rightarrow Y$ is prepared, $q\in Y$ and $u,v,w\in{\cal O}_{Y,q}$ are permissible parameters at $q$. Suppose that $p\in f^{-1}(q)$ and there exists $\phi(u,v)\in\hat{\cal O}_{Y,q}$ such that $w-\phi(u,v)$ is good (weakly good) at $p$ for $f$.

Suppose that $p$ is an $n$-point. Then there exists an affine neighborhood $V=\text{Spec}(S)$ of $p$   such that $w-\phi(u,v)$ is good (weakly good) at $p'$ for $f$ for all $n$-points $p'\in f^{-1}(q)\cap V$.
\end{Lemma}

\begin{pf} We will prove this in the case that  $p$ and $q$ are 1-points, and a form (\ref{eqT41}) of Definition \ref{DefT40} holds for $u,v,w-\phi$ in $\hat{\cal O}_{X,p}$.
The proof in the other cases is similar.

There exists an affine neighborhood $V=\text{spec}(S)$ of $p$, regular parameters $x,y,z\in\hat{\cal O}_{X,p}$ and a finite etale morphism $\pi:V_1=\text{spec}(S_1)\rightarrow S$ such that $x,y,z$ are uniformizing parameters on $V_1$, and regular parameters in ${\cal O}_{V_1,p'}$ for $p'\in \pi^{-1}(p)$ such that 
$$
\begin{array}{ll}
u&=x^a\\
v&=y\\
w&=\sum_{i<n}a_{ij}x^iy^j+x^n(\gamma(x,y,z)z+\Omega(x,y))
\end{array}
$$
in $\hat{\cal O}_{V_1,p'}=\hat{\cal O}_{X,p}$ where $\gamma$ is a unit series, $\Omega(x,y)$ is  a series.

Let $U=\text{Spec}(R)$ be an affine neighborhood of $q$ such that $f(V)\subset U$.

In $\hat{\cal O}_{X,p}$, 
$$
w-\phi(u,v)=\sum_{a\not\,\mid i, i<n}a_{ij}x^iy^j+x^n(\gamma z+\tilde\Omega)
$$
where $\tilde\Omega(x,y)$ is a series.
We see that 
$$
x^n\text{ divides }\frac{\partial(w-\phi)}{\partial z}\text{ in }\hat{\cal O}_{X,p},
$$
so that 
$$
x^n\text{ divides }\frac{\partial(w-\phi)}{\partial z}\text{ in }{\cal O}_{X,p}\otimes_{R_q}\hat{R_q}
$$
and thus
$$
x^n\text{ divides }\frac{\partial(w-\phi)}{\partial z}\text{ in }{\cal O}_{S_1,p'}\otimes_{R_q}\hat{R_q}
$$
at all points $p'\in\pi^{-1}(p)$.

We also have that 
$$
\frac{\partial^{n+1}(w-\phi)}{\partial z\partial x^n}(p)\ne 0
$$
which implies 
$$
\frac{\partial^{n+1}(w-\phi)}{\partial z\partial x^n}(p')\ne 0
$$
at all $p'\in \pi^{-1}(p)$.

Finally, we see that 
$$
x\text{ divides }\frac{\partial^i(w-\phi)}{\partial x^i}
$$
in ${\cal O}_{X,p}\otimes_{R_q}\hat R_q$ if $a$ divides $i$ and $i<n$, and thus
$$
x\text{ divides }\frac{\partial^i(w-\phi)}{\partial x^i}
$$
in ${\cal O}_{S_1,p'}\otimes_{R_q}\hat{R_q}$ at all points $p'$ of $\pi^{-1}(p)$.

Thus there exists a Zariski closed subset $C$ of $V_1$ which is disjoint from $\pi^{-1}(p)$ such that if $\overline p\in (f\circ\pi)^{-1}(q)\cap(V_1-C)$, then 
\begin{equation}\label{eqT129}
x^n\text{ divides } \frac{\partial(w-\phi)}{\partial z}\text{ in }{\cal O}_{S_1,\overline p}\otimes_{R_q}\hat R_q
\end{equation}

\begin{equation}\label{eqT130}
 \frac{\partial^{n+1}(w-\phi)}{\partial z\partial x^n}
 \text{ is a unit in }{\cal O}_{S_1,\overline p}\otimes_{R_q}\hat R_q
\end{equation}
and 
\begin{equation}\label{eqT131}
x\text{ divides } \frac{\partial^i(w-\phi)}{\partial x^i}\text{ in }{\cal O}_{S_1,\overline p}\otimes_{R_q}\hat R_q
\end{equation}
if $i<n$ and $a$ divides $i$.

Let $\overline C=\pi(C)$. $\pi:\text{Spec}(S_1)-\pi^{-1}(\overline C)\rightarrow V-\overline C$ is finite etale. Let $\overline V=\text{spec}(\overline S)$ be an affine neighborhood of $p$ in $V-\overline C$, and let $\pi^{-1}(\overline V)=\text{Spec}(\overline S_1)$. After replacing $V$ with $\overline V$, $S$ with $\overline S$, $V_1$ with $\pi^{-1}(\overline V)$ and $S_1$ with $\overline S_1$, we have that (\ref{eqT129}), (\ref{eqT130}) and (\ref{eqT131}) hold at all $\overline p\in (f\circ\pi)^{-1}(q)$.

Now consider the expression of $u,v,w-\phi(u,v)$ at $\overline p\in (f\circ\pi)^{-1})(q)$.
There exists $\alpha\in {\bf k}$ such that $x,y,z-\alpha$ are regular parameters at $\overline p$.
We have
$$
\begin{array}{ll}
u&=x^a\\
v&=y\\
w-\phi&=\sum \frac{1}{i!j!k!}\frac{\partial^{i+j+k}(w-\phi)}{\partial x^i\partial y^j\partial z^k}(0,0,\alpha)x^iy^j(z-\alpha)^k.
\end{array}
$$

(\ref{eqT129}) implies
$$
\frac{\partial^{i+j+k}(w-\phi)}{\partial x^i\partial y^j\partial z^k}(0,0,\alpha)=0
$$
if $i<n$ and $k\ge 1$,
(\ref{eqT131}) implies
$$
\frac{\partial^{i+j}(w-\phi)}{\partial x^i\partial y^j}(0,0,\alpha)=0
$$
if $a$ divides $i$ and $i<n$, and 
(\ref{eqT130}) implies
$$
\frac{\partial^{n+1}(w-\phi)}{\partial x^n\partial z}(0,0,\alpha)\ne0.
$$
Thus
$$
w-\phi=\sum_{i<n, a\not\,\mid i}\frac{1}{i!j!}\frac{\partial^{i+j}(w-\phi)}{\partial x^i\partial y^j}(0,0,\alpha)x^iy^j+x^n(\gamma'(z-\alpha)+\Omega'(x,y))
$$
where $\beta\in{\bf k}$ and $\gamma'$ is a unit series, so that $w-\phi$
is good at $\overline p$.
\end{pf}

\begin{Lemma}\label{LemmaT135} Suppose that $f:X\rightarrow Y$ is prepared, and $C\subset Y$ is an irreducible curve in the fundamental locus of $f$ such that $C$ contains a 1-point.

Suppose that $U\subset Y$ is an affine open subset, with uniformizing parameters $u,v,w$ such that $u,v,w$ are regular parameters in ${\cal O}_{Y,q}$ for a 1-point $q\in C\cap U$ such that $u=w=0$ are local equations of $C$.
Then for a general point $\overline q$ of $C\cap U$, and appropriate $\alpha\in {\bf k}$, $u,\overline v=v-\alpha, w$ are permissible parameters  at $\overline q$ and for $p\in f^{-1}(\overline q)$, either $p$ is a 1-point and we have a form at $p$ 
\begin{equation}\label{eqT160}
\begin{array}{ll}
u&=x^a\\
\overline v&=y\\
w&=\sum_{i<n} \phi_i(y)x^i+x^n(z+\delta)
\end{array}
\end{equation}
with $\delta\in{\bf k}$ and $\phi_i(0)\ne 0$ whenever $\phi_i\ne 0$,  or $p$ is a 2-point with a form at $p$ 
\begin{equation}\label{eqT332}
\begin{array}{ll}
u&=(x^ay^b)^t\\
\overline v&=z\\
w&=\sum \phi_i(z)(x^ay^b)^i+x^cy^d
\end{array}
\end{equation}
with $ad-bc\ne 0$ and $\phi_i(0)\ne 0$ whenever $\phi_i\ne 0$.
\end{Lemma}

\begin{pf} $u,v,w$ are permissible parameters at the 1-point $q\in C$. Suppose $p\in f^{-1}(q)$. 
Then there exists a Zariski open neighborhood $V=V^p=\text{Spec}(S)$ of $p$ in $X$, and an etale neighborhood $
W=W^p=\text{Spec}(S_1)$ of $V^p$ with uniformizing parameters $x,y,z$ in $S_1$, with induced morphism 
$$
\pi:\text{Spec}(S_1)\rightarrow \text{Spec}({\bf k}[x,y,z]),
$$
 such that $x,y,z$ are regular parameters in ${\cal O}_{W^p,p_1}$ for $p_1\in\pi^{-1}(p)$, and
by Lemma \ref{LemmaT1} and its proof,
we have one of the following cases:

\vskip .2 truein
\noindent {\bf Case 1.}  Suppose that $p$ is a 1-point. Then we have in $\hat{\cal O}_{X,p}=\hat S_p$: 
\begin{equation}\label{eqT252}
u=x^a, v=y,
w=\sum_{i<n} \phi_i(y)x^i+x^n(\gamma z+\psi(x,y))
\end{equation}
where $\gamma$ is a unit series.
In (\ref{eqT252}),  for $i<n$ we have 
\begin{equation}\label{eqT256}
\frac{1}{i!}\frac{\partial^iw}{\partial x^i}=\phi_i(y)+x\Omega\in\hat S_p
\end{equation}
for some $\Omega\in\hat S_p$ and 
$$
\frac{\partial^{n+1}w}{\partial x^n\partial z}(0,0,0)\ne 0.
$$

We can choose $V^p$ so that for $p_1\in W^p\cap D_X$, with regular parameters 
$$
x,\overline y=y-\alpha, \overline z=z-\beta
$$
 in $\hat{\cal O}_{W^p,p_1}$, for $i<n$ we have 
\begin{equation}\label{eqT254}
\frac{\partial^{i+1} w}{\partial z\partial x^i}(0,\alpha,\beta)=0
\end{equation}
and 
\begin{equation}\label{eqT255}
\frac{\partial^{n+1}w}{\partial z\partial x^n}(0,\alpha,\beta)\ne 0.
\end{equation}

We can choose $V^p$ so that for $i<n$, all irreducible components of 
$x=\frac{\partial^iw}{\partial x^i}=0$ in $W^p$ contain (a preimage of) $p$.

\vskip .2 truein
\noindent {\bf Case 2.}  Suppose that $p$ is a 2-point. Then we have in $\hat{\cal O}_{X,p}=\hat S_p$: 
\begin{equation}\label{eqT260}
u=(x^ay^b)^t, v=z,
w=\sum \phi_i(z)(x^ay^b)^i+x^cy^d\gamma
\end{equation}
with $\text{gcd}(a,b)=1$, $ad-bc\ne 0$ and $c>ai$ or $d>bi$ for all $i$ in the series. Further, $\gamma$ is a unit series.

We have $\Omega_1,\Omega_2\in\hat S_p=\hat{\cal O}_{X,p}$ such that 
\begin{equation}\label{eqT261}
\frac{1}{j!k!}\frac{\partial^{j+k}w}{\partial x^j\partial y^k}=
\left\{\begin{array}{ll}
x\Omega_1+y\Omega_2&\text{ if }jb-ka\ne 0 \text{ and }j<c\text{ or }k<d\\
\phi_i(z)+x\Omega_1+y\Omega_2&\text{ if there exists $i$ such that }(j,k)=i(a,b)\\
& \text{ and }j<c\text{ or }k<d
\end{array}
\right.
\end{equation}
There exists $\Omega_1\in\hat{\cal O}_{X,p}$ such that 
\begin{equation}\label{eqT257}
\frac{1}{j!}\frac{\partial^jw}{\partial x^j}=\left\{\begin{array}{ll}
x\Omega_1&\text{ if $j<c$ and there do not exist $k,i$ such that $(j,k)=i(a,b)$}\\
y^{ib}\phi_i(z)+x\Omega_1&\text{ if $j<c$ and there exist $k,i$ such that $(j,k)=i(a,b)$}.
\end{array}\right.
\end{equation}
There exists $\Omega_1\in\hat{\cal O}_{X,p}$ such that
$$
\frac{1}{k!}\frac{\partial^kw}{\partial y^k}=\left\{\begin{array}{ll}
y\Omega_1&\text{ if $k<d$ and there do not exist $j,i$ such that $(j,k)=i(a,b)$}\\
x^{ia}\phi_i(z)+y\Omega_1&\text{ if $k<d$ and there exist $j,i$ such that $(j,k)=i(a,b)$}.
\end{array}\right.
$$
Furthermore,
$$
\frac{\partial^{c+d}w}{\partial x^c\partial y^d}(0,0,0)\ne 0.
$$

We can choose $V^p$ so that 
\begin{enumerate}
\item for $ai<c$ or $bi<d$, all irreducible components of 
$$
x=y=\frac{\partial^{i(a+b)}w}{\partial x^{ia}\partial y^{ib}}=0
$$
in $W^p$ contain (a preimage of) $p$,
\item  for $j<c$, all irreducible components of
$$
x=\frac{\partial^jw}{\partial x^j}=0
$$
in $W^p$ contain (a preimage of) $p$, and
\item for $k<d$, all irreducible components of 
$$
y=\frac{\partial^kw}{\partial y^k}=0
$$
in $W^p$ contain (a preimage of) $p$.
\end{enumerate}

There exist $V_1=V^{p_1},\ldots, V_n=V^{p_n}$ such that $\{V_1,\ldots, V_n\}$ is an affine cover of $f^{-1}(q)$.
We may assume that $V_1,\ldots, V_n$ is an affine cover of $f^{-1}(U)$.

Suppose that $\overline q\in C\cap U$ is a general point. Then ${\cal O}_{Y,q}$ has regular parameters $(u,\overline v=v-\alpha,w)$ for a general $\alpha\in{\bf k}$.
$u,\overline v,w$ are permissible parameters at $\overline q$. Suppose that $\overline p\in f^{-1}(\overline q)$. $\overline p\in V^p=V_i$ for some $i$.
We identify $p$ and $\overline p$ with points in $\pi^{-1}(p)$, $\pi^{-1}(\overline p)$.

Suppose that $p$ is a 1-point. There exists $\beta\in{\bf k}$ such that $x,\overline y=y-\alpha, z-\beta$ are regular parameters in $\hat{\cal O}_{X,\overline p}$. We have an expression 
\begin{equation}\label{eqT169}
\begin{array}{ll}
u&=x^a\\
\overline v&=\overline y=y-\alpha\\
w&=\sum \frac{1}{i!j!k!}\frac{\partial^{i+j+k}w}{\partial x^i\partial y^j\partial z^k}(0,\alpha,\beta)x^i(y-\alpha)^j(z-\beta)^k
\end{array}
\end{equation}
in $\hat{\cal O}_{X,\overline p}$.

By (\ref{eqT254}) and (\ref{eqT255}),  in (\ref{eqT169}), we have
$$
w=\sum_{i<n}\overline\phi_i(\overline y)x^i+x^n(\overline z\overline\gamma+\overline\Psi(x,\overline y))
$$
where $\overline\gamma$ is a unit series, and 
$$
\overline\phi_i(\overline y)=\sum_{j=0}^{\infty}\frac{1}{i!j!}\frac{\partial^{i+j}w}{\partial x^i\partial y^j}(0,\alpha,\beta)\overline y^j.
$$
If $\phi_i(y)=0$ we have $\overline\phi_i(\overline y)=0$.

 Suppose that $\phi_i(y)\ne 0$. 
We will show that $\overline\phi_i(0)\ne 0$. If $\phi_i(0)\ne 0$, then $\frac{\partial^iw}{\partial x^i}$ does not vanish on $W^p\cap D_X$, so that 
$$
\overline\phi_i(0)=\frac{\partial^iw}{\partial x^i}(\overline p)\ne 0.
$$

Suppose that $\phi_i(0)=0$. Further suppose that 
$$
\overline\phi_i(0)=\frac{\partial^iw}{\partial x^i}(\overline p)= 0.
$$
Then there exists an irreducible curve $\Lambda$ which is a component of $x=\frac{\partial^iw}{\partial x^i}=0$ in $W^p$ containing $\overline p$.

By our construction of $W^p$, we may assume that our choice of preimage of $p$ in $W^p$ satisfies $p\in\Lambda$. Let $I_{\Lambda}$ be the prime ideal of $\Lambda$ in $S_1$. $\frac{\partial^iw}{\partial x^i}, x\in I_{\Lambda}\hat{(S_1)}_p$ implies $\phi_i(y)\in I_{\Lambda}\hat{(S_1)}_p$.
Since $I_{\Lambda}\hat{(S_1)}_p$ is reduced, we have $y\in I_{\Lambda}\hat{(S_1)}_p$.
As $(S_1)_p\rightarrow\hat{S_1}=\hat{\cal O}_{X,p}$ is faithfully flat, we have $y\in I_{\Lambda}(S_1)_p$.

Since $I_{\Lambda}$ is a prime ideal and $(I_{\Lambda})_p\ne (S_1)_p$, we have that $y\in I_{\Lambda}$.

The ideal of $\overline p$ in $(S_1)_{\overline p}$ is $(x,y-\alpha,z-\beta)$.
$$
I_{\Lambda}(S_1)_{\overline p}\subset (x,y-\alpha,z-\beta)
$$
implies $y\in (x,y-\alpha,z-\beta)$ which is impossible since $0\ne\alpha$. Thus we have $\overline\phi_i(0)\ne 0$.

Suppose that $p$ is a 2-point and $\overline p$ is a 2-point. Then there is $\alpha\in{\bf k}$ such that $x,y,z-\alpha$ are regular parameters in $\hat{\cal O}_{X,\overline p}$, and we have an expression 
\begin{equation}\label{eqT170}
\begin{array}{ll}
u&=(x^ay^b)^t\\
\overline v&=\overline z=z-\alpha\\
w&=\sum \frac{1}{i!j!k!}\frac{\partial^{i+j+k}w}{\partial x^i\partial y^j\partial z^k}(0,0,\alpha)x^iy^j(z-\alpha)^k
\end{array}
\end{equation}
in $\hat{\cal O}_{X,\overline p}$.

We have 
$$
\frac{\partial^{c+d}w}{\partial x^c\partial y^d}(0,0,\alpha)\ne 0.
$$
 Further if $j<c$ or $k<d$, $jb-ka\ne 0$ and $l\ge 0$, we have
$$
\frac{\partial^{j+k+l}w}{\partial x^j\partial y^k\partial z^l}(0,0,\alpha)=0.
$$
Thus in (\ref{eqT170}) we have
$$
w=\sum\overline\phi_i(z)(x^ay^b)^i+x^cy^d\overline\gamma, and
$$
where $\overline \gamma$ is a unit series, the sum is over $i$ such that $ai<c$ or $bi<d$, and
$$
\overline\phi_i(\overline z)=\sum_{k=0}^{\infty}\frac{1}{(ia)!(jb)!k!}\frac{\partial^{ia+jb+k}w}{\partial x^{ia}\partial y^{ib}\partial z^k}(0,0,\alpha)\overline z^k.
$$
By a similar analysis as for Case 1, we see that
if $\phi_i(z)=0$ then $\overline\phi_i(\overline z)=0$ and if $\phi_i(z)\ne 0$, we have $\overline\phi_i(0)\ne 0$.

Suppose that $p$ is a 2-point and $\overline p$ is a 1-point. Then after possibly interchanging $x$ and $y$, there exist $\alpha,\beta\in{\bf k}$ 
and regular parameters $x,y-\beta,z-\alpha$ in $\hat{\cal O}_{X,\overline p}$ with $0\ne\alpha,\beta$.

Set 
$$
\overline x=xy^{\frac{b}{a}},
\overline y=y-\beta,
\overline z=z-\alpha
$$
Then $\overline x,\overline y,\overline z$ are regular parameters in $\hat{\cal O}_{X,\overline p}$. From the Jacobian of $u,\overline v,w$ we see that we have an expression
$$
\begin{array}{ll}
u&=\overline x^{at}\\
\overline v&=\overline z\\
w&=P(\overline x,\overline z)+\overline x^c\Omega
\end{array}
$$
where $P$ and $\Omega$ are series, and $P$ has degree $<c$ in $\overline x$.

We have 
$$
\begin{array}{ll}
x&=\overline x(\overline y+\beta)^{-\frac{b}{a}}\\
y&=\overline y+\beta\\
z&=\overline z+\alpha
\end{array}
$$
Further,
$$
\begin{array}{ll}
\frac{\partial}{\partial \overline x}&=y^{-\frac{b}{a}}\frac{\partial}{\partial x}\\
\frac{\partial}{\partial \overline y}&=-\frac{b}{a}xy^{-1}\frac{\partial}{\partial x}+\frac{\partial}{\partial y}\\
\frac{\partial}{\partial \overline z}&=\frac{\partial}{\partial z}.
\end{array}
$$
We have an expansion
$$
P=\sum_{j<c}\overline\phi_j(\overline z)\overline x^j
=\sum_{j<c}\left(\sum_{k=0}^{\infty}\frac{\partial^{j+k}w}{\partial \overline z^k\partial \overline x^j}(\overline p)\overline z^k\right)\overline x^j
$$
where
$$
\frac{\partial^{j+k}w}{\partial \overline z^k\overline x^j}(\overline p)
=\beta^{-j\frac{b}{a}}\frac{\partial^{j+k}w}{\partial z^k\partial x^j}(\overline p).
$$
By (\ref{eqT257}) we have that $\overline \phi_j(\overline z)=0$ if there do not exist $k,i$ such that $(j,k)=i(a,b)$.

Suppose that $j<c$ and there exists $k,i$ such that $(j,k)=i(a,b)$ and $\phi_j(z)\ne0$. Suppose that $\overline\phi_j(0)=0$. 
Then
$$
\overline\phi_j(0)=0=\frac{\partial^{j}w}{\partial \overline x^j}(\overline p)
=\beta^{-j\frac{b}{a}}\frac{\partial^{j}w}{\partial x^j}(\overline p)
$$
implies there exists an irreducible curve $\Lambda$ which is a component of
$$
x=\frac{\partial^j w}{\partial x^j}=0
$$
in $W^p$ which contains $p$ and $\overline p$. Let $I_{\Lambda}$ be the prime ideal of $\Lambda$ in $S_1$. $\frac{\partial^jw}{\partial x^j}, x\in I_{\Lambda}\widehat{(S_1)}_p$ implies 
$$
y^{ib}\phi_i(z)\in I_{\Lambda}\widehat{(S_1)_p}.
$$
 Since $y\not\in I_{\Lambda}\widehat{(S_1)_p}$, as $\Lambda$ is not a 2-curve, we have $z\in I_{\Lambda}\widehat{(S_1)}_p$. Thus $z\in I_{\Lambda}(S_1)_p$.
 As $(I_{\Lambda})_p\ne (S_1)_p$, we have that $z\in I_{\Lambda}$, which is a contradiction, since $\alpha\ne 0$.
  Thus $\overline\phi_i(0)\ne 0$.

\end{pf}

\begin{Lemma}\label{LemmaT107} Suppose that $f:X\rightarrow Y$ is prepared.  Then there exists an open subset $V$ of $Y$ such that $V\cap C\ne \emptyset$ for every integral curve $C\subset D_Y$ contained in the fundamental locus of $f$ which contains a 1-point, and if $U=f^{-1}(V)$, $\overline f=f\mid U$, then there exists a commutative diagram
$$
\begin{array}{rll}
U_1&\stackrel{\overline f_1}{\rightarrow}&V_1\\
\overline \Phi_1\downarrow&&\downarrow\overline \Psi_1\\
U&\stackrel{\overline f}{\rightarrow}&V
\end{array}
$$
such that $\overline\Phi_1$ and $\overline\Psi_1$ are products of blow ups of curves which dominate a curve  $C$ contained in the fundamental locus of $\overline f$ which are possible centers (for the preimage of $D_V=D_Y\cap V$) and $\overline f_1$ is toroidal.
\end{Lemma}

\begin{pf} 
By Lemma \ref{LemmaT135}, there exists an open set $V\subset Y$ such that 
$V\cap C=\emptyset$ if $C$ is a 2-curve or an isolated point contained in the fundamental locus of $f$,
$V\cap C\ne\emptyset$ for all curves $C$ contained in the fundamental locus of $f$ which contain a 1-point, and if $\overline q\in C\cap V$ is a 1-point, then there exist permissible parameters $u,v,w$ at $\overline q$ such that $u=w=0$ are local equations of $C$, and if 
 $p_1\in f^{-1}(q)$, then
 we have permissible parameters $x,y,z$ in $\hat{\cal O}_{X,p_1}$ such that 
$p_1$ is a 1-point: 
\begin{equation}\label{eqT132}
\begin{array}{ll}
u&=x^a\\
v&=y\\
w&=\sum_{j=0}^{m}\phi_{i_j}(y)x^{i_j}+x^n(\overline\alpha+z)
\end{array}
\end{equation}
where $i_m<n$, $\overline\alpha\in{\bf k}$ and  all $\phi_{i_j}(y)$ are nonzero
or $p_1$ is  a 2-point: 
\begin{equation}\label{eqT133}
\begin{array}{ll}
u&=(x^ay^b)^t\\
v&=z\\
w&=\sum_{j=0}^m\phi_{i_j}(z)(x^ay^b)^{i_j}+x^cy^d
\end{array}
\end{equation}
where $ad-bc\ne 0$, all $\phi_{i_j}(y)$ are nonzero, and $i_ja<c$ or $i_jb<d$ for all $j$.
We further have  that  $\phi_{i_j}(0)\ne 0$  for all $j$.

 Let $C$ be the fundamental locus of $\overline f:U\rightarrow V$. There exists $\Phi_1':U_1'\rightarrow U$ which is a product of blow ups of 2-curves (which dominate
an irreducible component of $C$) such that all local forms (\ref{eqT133}) at points $p\in (\overline f\circ\Phi_1')^{-1}(\overline q)$ for $\overline q\in C\cap V$ are such that either $m=-1$ (so that $\sum\phi_{i_j}(z)(x^ay^b)^{i_j}=0)$, or $(x^ay^b)^{i_0}$ divides $x^cy^d$.

Suppose that $\overline q\in C$.
The set of points $p\in (\overline f\circ\Phi_1')^{-1}(\overline q)$ such that  ${\cal I}_C{\cal O}_{U_1'}$ is not invertible is a union of points $p$ such that $p$ has permissible parameters $x,y,z$ of the form 
\begin{equation}\label{eqT146}
u=x^a, v=y, w=x^nz
\end{equation}
with $n<a$ or 
\begin{equation}\label{eqT147}
u=(x^ay^b)^t, v=z, w=x^cy^d
\end{equation}
with $(at-c)(bt-d)<0$.

Let $\overline\Psi_1:V_1\rightarrow V$ be the blow up of $C$ (which has local equations $u=w=0$ at $\overline q\in C$).

We can  blow up curves above $U_1'$ which dominate a component of $C$  to obtain $\overline\Phi_1:U_1\rightarrow U_1'$ such that there exists a factorization 
$\overline f_1:U_1\rightarrow V_1$, and if $\overline q\in C$, $p_1\in(\overline f\circ\Phi_1'\circ\overline\Phi_1)^{-1}(\overline q)$, then an expression (\ref{eqT132}) or (\ref{eqT133}) holds. 

Suppose that $\overline q\in C$, with permissible parameters $u,v,w$ as above.
Let $q_1\in \overline\Psi_1^{-1}(\overline q)$. $q_1$ has permissible parameters $u_1,v_1,w_1$ with $q_1$ a 1-point 
\begin{equation}\label{eqT136}
u=u_1, v=v_1, w=u_1(w_1+\alpha)
\end{equation}
or $q_1$ a 2-point 
\begin{equation}\label{eqT137}
u=u_1v_1, v=w_1, w=u_1.
\end{equation}

Suppose that $p_1\in \overline f_1^{-1}(q_1)$.
\vskip .2truein
\noindent {\bf Case 1} Suppose that $0\ne\alpha$ in (\ref{eqT136}). Then $a=i_0$ and
 $\phi_{i_0}(0)=\alpha$ (or $m=-1$, $a=n$ and $\overline\alpha=\alpha$) if $p_1$ satisfies (\ref{eqT132}), 
 $t=i_0$ and $\phi_{i_0}(0)=\alpha$ if $p_1$ satisfies (\ref{eqT133}).

We thus have that at $p_1$, 
\begin{equation}\label{eqT138}
u_1=x^a, v=y,
w_1=(\phi_{i_0}(y)-\alpha)+\sum_{j=1}^{m}\phi_{i_j}(y)x^{i_j-i_0}+x^{n-i_0}(\overline\alpha+z)
\end{equation}
of the form (\ref{eqT132})
or 
\begin{equation}\label{eqT139}
u_1=(x^ay^b)^t, v=z,
w_1=(\phi_{i_0}(z)-\alpha)+\sum_{j=1}^{m}\phi_{i_j}(z)(x^ay^b)^{i_j-i_0}+x^{c-i_0a}y^{d-i_0b}.
\end{equation}
of the form of (\ref{eqT133}).

If $\overline f_1$ is not toroidal at $p_1$, we have that $u_1=w_1-(\phi_{i_0}(v)-\alpha)=0$ are (formal) local equations of a branch of the fundamental locus of $\overline f_1$ through $q_1$. After possibly replacing $V$ with an open subset of $V$, for $\overline q\in C$,
 $\overline q$ is a general point of a component of $C$, so the fundamental locus of $\overline f_1$ through $q_1$ must be the germ of a nonsingular algebraic curve.

\vskip .2truein
\noindent {\bf Case 2} Suppose that $0=\alpha$ in (\ref{eqT136}).  Then $i_0>a$ 
or $m=-1$ and $n>a$ (or $m=-1$, $n=a$ and $\overline\alpha=0$) in (\ref{eqT132}), or $i_0>t$ or $m=-1$, $(c,d)>t(a,b)$ in (\ref{eqT133}). We have that $u_1,v,w_1$ are permissible parameters at $q_1$ for $\overline f_1$ of the form (\ref{eqT132}) or (\ref{eqT133}).

\vskip .2truein
\noindent {\bf Case 3} (\ref{eqT137}) holds.  Then $a>i_0$ (or $m=-1$, $0\ne\overline\alpha$ and $n<a$) in (\ref{eqT132}) or $t>i_0$ (or $m=-1$ and $(c,d)\le (at,bt)$) in (\ref{eqT133}). 

Suppose that (\ref{eqT132}) holds and $m\ge 0$. We change variables at $p_1$ to get an expression
$$
u=\sum_{i=0}^{\tilde m}\tilde\phi_{\tilde i_j}(y)\tilde x^{\tilde i_j}+\tilde x^{n-i_0+a}(\tilde z+\cdots),
v=y, w=\tilde x^{i_0}
$$
with $\tilde m<n-i_0+a$.
As in the proof of Lemma \ref{LemmaT135}, $\tilde \phi_{\tilde i_j}(0)\ne0$  for all $\tilde i_j$ since $\overline q$ is a general point of a component of $C$.

$q_1$ is a 2-point, and we have: 
\begin{equation}\label{eqT144}
u_1=\tilde x^{i_0}, v_1=\sum_{j=0}^{\tilde m}\tilde\phi_{\tilde i_j}(y)\tilde x^{\tilde i_j-i_0}+\tilde x^{n-2i_0+a}(\tilde z+\cdots),
w_1=y.
\end{equation}

Note that  $u_1=v_1=0$ are local equations of the fundamental locus of $\overline f_1$ at $q_1$
if $\overline f_1$ is not toroidal at $p_1$.

Suppose that (\ref{eqT133}) and (\ref{eqT137}) hold and $m\ge 0$. We change variables at $p_1$ to have an expression
$$
u=\sum_{j=0}^{\tilde m}\tilde\phi_{\tilde i_j}(z)(\tilde x^a\tilde y^b)^{\tilde i_j}+\tilde x^{c-i_0a+ta}\tilde y^{d-i_0b+tb},
v=z, w=(\tilde x^a\tilde y^b)^{i_0}
$$
with $\tilde i_ja<c-i_0a+ta$ or $\tilde i_jb<d-i_0b+tb$ for all $\tilde i_j$.

 As in the proof of Lemma \ref{LemmaT135}, $\tilde\phi_{\tilde i_j}(0)\ne 0$ for all $\tilde i_j$, since $\overline q$ is a general point of a component of $C$.
$q_1$ is  a 2-point and we have 
\begin{equation}\label{eqT145}
u_1=(\tilde x^a\tilde y^b)^{i_0},
v_1=\sum_{j=0}^{\tilde m}\tilde\phi_{\tilde i_j}(z)(\tilde x^a\tilde y^b)^{\tilde i_j-i_0}+\tilde x^{c-2i_0a+ta}\tilde y^{d-2i_0b+tb},
w_1=z.
\end{equation}
$u_1=v_1=0$ are local equations of the fundamental locus of $\overline f_1$ at $q_1$ if $\overline f_1$ is not toroidal at $p_1$.

The fundamental locus $C_1$ of $f_1:U_1\rightarrow V_1$ is a (disjoint) union of nonsingular curves which dominate components of $C$. If $\gamma_1$ is a component of $C_1$ then $\gamma_1$ consists of 1-points or $\gamma_1$ consists of 2-points. We will construct a commutative diagram 
\begin{equation}\label{eqT148}
\begin{array}{rll}
U_2&\stackrel{\overline f_2}{\rightarrow}&V_2\\
\overline\Phi_2\downarrow&&\downarrow\overline\Psi_2\\
U_1&\stackrel{\overline f_1}{\rightarrow}&V_1
\end{array}
\end{equation}
where $\overline \Psi_2:V_2\rightarrow V_1$ is the blow up of $C_1$.

Suppose that $\gamma_1$ is a component of $C_1$ and $q_1\in \gamma_1$.

First suppose that $\gamma_1$ consists of 1-points. Then (\ref{eqT132}) or (\ref{eqT133}) holds at all points $p\in \overline f_1^{-1}(q_1)$. The construction of (\ref{eqT148}) above points of $\gamma_1$ is as in the construction of $\overline f_1$ above.

Suppose that $\gamma_1$  consists of 2-points. Then there exist permissible parameters $u_1,v_1,w_1$ at $q_1$ such that (\ref{eqT144}) or (\ref{eqT145}) holds at all $p\in \overline f_1^{-1}(q_1)$, and $u_1=v_1=0$ are local equations of $\gamma_1$ at $q_1$. 

If $p\in \overline f_1^{-1}(q_1)$ is such that ${\cal I}_{C_1}{\cal O}_{U_2,p}$ is not invertible, then we have permissible parameters $x,y,z$ at $p$ such that
$$
u_1=(x^ay^b)^t, v_1=x^cy^d, w_1=z.
$$
In particular, $\overline f_1$ is toroidal at $p$.

We now blow up curves 2-curves (above $U_1$) which dominate $\gamma_1$ and are supported in the locus where $U_1\rightarrow V_1$ is torodial to obtain the construction of $\overline\Phi_2:U_2\rightarrow U_1$ above $\gamma_1$. $\overline f_2$ is toroidal above the torodial locus of $f$. Let $q_2\in\overline\Psi_2^{-1}(q_1)$. $q_2$ has permissible parameters $u_2,v_2,w_2$ defined by one of the following 3 cases.

$q_2$ is a 1-point 
\begin{equation}\label{eqT171}
u_1=u_2,v_1=u_2(w_2+\alpha), w_1=v_2
\end{equation}
with $\alpha\ne 0$ or
$q_2$ is a 2-point 
\begin{equation}\label{eqT149}
u_1=u_2, v_1=u_2v_2, w_1=w_2
\end{equation}
or
$q_2$ is a 2-point 
\begin{equation}\label{eqT150}
u_1=u_2v_2, v_1=v_2, w_1=w_2.
\end{equation}

As in the case when $q_2$ is a 1-point, we see that if (\ref{eqT171}) holds, then all points above $q_2$ have the form (\ref{eqT132}) or (\ref{eqT133}), and that if (\ref{eqT149}) or (\ref{eqT150}) holds, then all points $p_2$ above $q_2$ have the form (\ref{eqT144}) or (\ref{eqT145}).

We iterate to construct a commutative diagram 

\begin{equation}\label{eqT172}
\begin{array}{rll}
\vdots&&\vdots\\
\downarrow&&\downarrow\\
U_n&\stackrel{\overline f_n}{\rightarrow}&V_n\\
\overline\Phi_n\downarrow&&\downarrow\overline\Psi_n\\
U_{n-1}&\stackrel{\overline f_{n-1}}{\rightarrow}&V_{n-1}\\
\downarrow&&\downarrow\\
\vdots&&\vdots\\
\downarrow&&\downarrow\\
U_1&\stackrel{\overline f_1}{\rightarrow}&V_1\\
\downarrow&&\downarrow\\
U&\stackrel{\overline f}{\rightarrow}&V
\end{array}
\end{equation}
where each $V_r\rightarrow V_{r-1}$ is the blow up of the fundamental locus $C_{r-1}$ of $\overline f_{r-1}$, which is a disjoint union of nonsingular curves which dominate components of $C$.  

All points of $U_n$ have a form (\ref{eqT132}), (\ref{eqT133}), (\ref{eqT144}) or (\ref{eqT145}).
We continue the construction as long as $\overline f_n$ is not toroidal.

Suppose that (\ref{eqT172}) does not converge in a toroidal morphism in a finite number of steps. Then there exists a 0-dimensional valuation $\nu$ of ${\bf k}(X)$ with center on $U$ such that $\overline f_n$ is not toroidal at the center $p_n$ of $\nu$ on $U_n$ for all $n$. Let $q_n$ be the center of $\nu$ on $V_n$.

Suppose that (\ref{eqT132}) holds for $p=p_0$. There exists $r(1)$ such that  $q_{r(1)}$ has permissible parameters $u_{r(1)},v_{r(1)},\overline w_{r(1)}$ defined by
$$
u=u_{r(1)}^e, w=u_{r(1)}^f(\overline w_{r(1)}+\phi_{i_0}(0))
$$
where $\text{gcd}(e,f)=1$ and $\frac{e}{f}=\frac{a}{i_0}$.

The germ of $\overline f$ at $p$ factors through $\overline\Psi_{r(1)}$, so we have

$$
u_{r(1)}=x^{\frac{a}{e}},
\overline w_{r(1)}=(\phi_{i_0}(y)-\phi_{i_0}(0))+\sum_{j=1}^{m}\phi_{i_j}(y)x^{i_j-i_0}+x^{n-i_0}(\overline\alpha+z).
$$
Set $w_{r(1)}=\overline w_{r(1)}-[\phi_{i_0}(v)-\phi_{i_0}(0)]$.

$u_{r(1)},v,w_{r(1)}$ are (formal) regular parameters at $q_{r(1)}$, and $u_{r(1)}=w_{r(1)}=0$ are equations of the fundamental locus at $q_{r(1)}$. We see that at $p=p_{r(1)}$,
 $$
 u_{r(1)}=x^{a(1)}, v=y,
 w_{r(1)}=\sum_{j=0}^{m(1)}\phi(1)_{i(1)_j}(y)x^{i(1)_j}+x^{n(1)}(\overline\alpha+z)
 $$
 where $a(1)=\text{gcd}(a,i_0)$, $m(1)=m-1$, $n(1)=n-i_0$, $i(1)_j=i_{j+1}-i_0$ for $0\le j\le m(1)$,
 $\phi(1)_{i(1)_j}(y)=\phi_{i_{j+1}}(y)$.
 
 We iterate to get for $k\le m+1$, $r(k)$ such that $q_{r(k)}$ has permissible parameters $u_{r(k)},v,\overline w_{r(k)}$ defined by
 
 $$
 u_{r(k-1)}=u_{r(k)}^{e_k},
 w_{r(k-1)}=u_{r(k)}^{f_k}(\overline w_{r(k)}+\phi(k-1)_{i(k-1)_0}(0))
 $$
 where $\text{gcd}(e_k,f_k)=1$.
 
 The rational map $U\rightarrow V_{r(k)}$ is a morphism at $p=p_{r(k)}$.
 Set
 $$
 w_{r(k)}=\overline w_{r(k)}-[\phi(k-1)_{i(k-1)_0}(v)-\phi(k-1)_{i(k-1)_0}(0)].
 $$
 
 We have an expression
 
 $$
 u_{r(k)}=x^{a(k)}, v=y,
 w_{r(k)}=\sum_{j=0}^{m(k)}\phi(k)_{i(k)_j}(y)x^{i(k)_j}+x^{n(k)}(\overline\alpha+z).
 $$
 We have (for $k\le m+1$)
 $a(k)=\text{gcd}(a,i_0,i_1,\ldots,i_{k-1})$, $n(k)=n-i_{k-1}$, $m(k)=m-k$,
 $i(k)_j=i_{j+k}-i_{k-1}$ for $0\le j\le m(k)$.
 
 We further have 
 \begin{equation}\label{eqT263}
 \frac{e_k}{f_k}=\frac{a(k-1)}{i(k-1)_0}=\frac{\text{gcd}(a,i_0,\ldots,i_{k-2})}{i_{k-1}-i_{k-2}}.
 \end{equation}
$u_{r(k)}=w_{r(k)}=0$ are (formal) equations of the fundamental locus at $q_{r(k)}$.
 
 $q_{r(m+1)}$ has permissible parameters $u_{r(m+1)}, v, w_{r(m+1)}$ defined by
 $$
 u_{r(m+1)}=x^{a(m+1)}, v=y, w_{r(m+1)}=x^{n(m+1)}(\overline\alpha +z).
 $$

The rational map $U\rightarrow V_{r(m+1)}$ is a morphism at $p$.
We have 
$$
a(m+1)=\text{gcd}(a,i_0,i_1,\ldots,i_m)
$$
 and $n(m+1)=n-i_m$.

Finally, we see that there exists $r(m+2)$ such that $\overline f_{r(m+2)}$ is toroidal at $p_{r(m+2)}$, a contradiction.

A similar argument holds if (\ref{eqT133}) holds at $p=p_0$.

We conclude that (\ref{eqT172}) converges after a finite number of iterations in a diagram which satisfies the conclusions of Lemma \ref{LemmaT107}.

\end{pf}

\begin{Definition}\label{DefT133}
Suppose that $f:X\rightarrow Y$ is $\tau$-prepared and $q\in G_Y(f,\tau)$ is a 1-point. Then $q$ is {\bf perfect} for $f$ if the fundamental locus of $f$ through $q$ is a (germ of a) nonsingular curve $\gamma$ and if $u,v,w$ are algebraic permissible parameters at $q$ such that $u=w=0$ are local equations of $\gamma$ at $q$ then
 there exist finitely many series $\phi_i(u,v)\in {\bf k}[[u,v]]$ such that
\begin{enumerate}
\item[1.] $u,v,w_i=w-\phi_i(u,v)$ are super parameters at $q$ for all $i$.
\item[2.] For all $p\in f^{-1}(q)$, some $w_i$ is weakly good for $f$ at $p$.
\end{enumerate}
\end{Definition}

\begin{Lemma}\label{LemmaT134}  Suppose that $f:X\rightarrow Y$ is $\tau$-prepared. Let $V\subset G_Y(f,\tau)$ be the set of perfect 1-points. Then $G_Y(f,\tau)-V$ is a finite set.
\end{Lemma}

\begin{pf} Suppose that $\overline q$ is a general point of a curve $C\subset G_Y(f,\tau)$ (so that $\overline q$ is a 1-point). Let $u,v,w$ be algebraic permissible parameters at $\overline q$ such that $u=w=0$ are local equations of $C$.

In a neighborhood of $\overline q$, we construct a diagram (\ref{eqT172}). (\ref{eqT172}) is finite by the conclusions of Lemma \ref{LemmaT107}.

Suppose that $p\in f^{-1}(\overline q)$.
At $p$ we have permissible parameters $x,y,z\in\hat{\cal O}_{X,p}$ such that 
if $p$ is a 1-point: 
\begin{equation}\label{eqT152}
\begin{array}{ll}
u&=x^a\\
v&=y\\
w&=\sum_{j=0}^{m}\phi_{i_j}(y)x^{i_j}+x^n(\beta+z)
\end{array}
\end{equation}
where $\beta\in{\bf k}$, $i_m<n$, all $\phi_{i_j}(y)$ are non zero,
or if $p$ is a 2-point: 
\begin{equation}\label{eqT153}
\begin{array}{ll}
u&=(x^by^c)^a\\
v&=z\\
w&=\sum_{j=0}^m\phi_{i_j}(z)(x^by^c)^{i_j}+x^dy^e
\end{array}
\end{equation}
where $\text{gcd}(b,c)=1$, all $\phi_{i_j}(z)$ are non zero, $i_m(b,c)\not\ge(d,e)$.
In either case, there exists a largest $l\le m$ such that $a\mid i_j$ if $j\le l$.

If at $p$ there is a form (\ref{eqT152}), set
$$
\phi_p(u,v)=\sum_{j\le l}\phi_{i_j}(y)x^{i_j}=\sum_{j\le l}\phi_{i_j}(v)u^{\frac{i_j}{a}}.
$$

If at $p$ there is a form (\ref{eqT153}), set
$$
\phi_p(u,v)=\sum_{j\le l}\phi_{i_j}(z)(x^by^c)^{i_j}=\sum_{j\le l}\phi_{i_j}(v)u^{\frac{i_j}{a}}.
$$
In both cases $w-\phi_p(u,v)$ is weakly good for $f$ at $p$.

By Lemma \ref{LemmaT128}, there exist finitely many points $p_1,\ldots,p_n\in X$
such that if we set $\phi_i(u,v)=\phi_{p_i}(u,v)$, then for all $p\in f^{-1}(\overline q)$, some $w_i=w-\phi_i(u,v)$ is weakly good for $f$ at $p$.

Since the $\phi_{i_j}$ are units in $\hat{\cal O}_{X,p_i}$ by Lemma \ref{LemmaT107}, we see that $u,v,w_i$ satisfy 5 or 6 of Definition \ref{Def357} of super parameters at $p_i$.

Let $p=p_i$ for some $i$, with the notation of (\ref{eqT152}) or (\ref{eqT153}).

Suppose that $\overline p\in f^{-1}(\overline q)$. We must show that $u,v,w_i$ are super parameters at $\overline p$.

At $\overline p$ we have permissible parameters $\overline x,\overline y,\overline z\in\hat{\cal O}_{X,\overline p}$ such that
if $\overline p$ is a 1-point:
\begin{equation}\label{eqT134}
\begin{array}{ll}
u&=\overline x^{\overline a}\\
v&=\overline y\\
w&=\sum_{j=0}^{\overline m}\overline \phi_{\overline i_j}(\overline y)\overline x^{\overline i_j}+\overline x^{\overline n}(\overline\beta+\overline z)
\end{array}
\end{equation}
where $\overline \beta\in{\bf k}$, $\overline i_{\overline m}<\overline n$, all $\overline \phi_{\overline i_j}(\overline y)$ are non zero,

or if $\overline p$ a 2-point: 
\begin{equation}\label{eqT135}
\begin{array}{ll}
u&=(\overline x^{\overline b}\overline y^{\overline c})^{\overline a}\\
v&=\overline z\\
w&=\sum_{j=0}^{\overline m}\overline \phi_{\overline i_j}(\overline z)(\overline x^b\overline y^c)^{\overline i_j}+\overline x^{\overline d}\overline y^{\overline e}
\end{array}
\end{equation}
where $\text{gcd}(\overline b,\overline c)=1$, $\overline i_{\overline m}(\overline b,\overline c)\not\ge(\overline d,\overline e)$, all $\overline \phi_{\overline i_j}(\overline z)$ are non zero.

We know from Lemma \ref{LemmaT135} that all $\phi_{i_j}$ and $\overline\phi_{\overline i_j}$ are units in $\hat{\cal O}_{X,p}$ and $\hat{\cal O}_{X,\overline p}$ respectively.

It will follow that $w_i$ are super parameters at $\overline p$ after we have proven that if 
there exists $t$ with $t\le\text{min}\{l,\overline m\}$ and
$$
\frac{i_j}{a}=\frac{\overline i_j}{\overline a}\text{ and }\phi_{i_j}(0)=\overline\phi_{\overline i_j}(0)
$$
for $0\le j\le t$, then we have equality of  power series in $u$,
$$
\overline\phi_{\overline i_j}(v)=\phi_{i_j}(v)
$$
for $0\le j\le t$, and thus equality of series

$$
\sum_{j=0}^t\phi_{i_j}(v)u^{\frac{i_j}{a}}=\sum_{j=0}^t\overline\phi_{\overline i_j}(v)u^{\frac{\overline i_j}{\overline a}}.
$$

We will prove this in the case when $p=p_i$ satisfies (\ref{eqT152}) and $\overline p$ satisfies (\ref{eqT134}).

The proof of the remaining three cases is similar.

We prove this by induction. First assume that $l\ge 0$ and 
\begin{equation}\label{eqT262}
\frac{i_0}{a}=\frac{\overline i_0}{\overline a}
\end{equation}
and $\phi_{i_0}(0)=\overline \phi_{\overline i_0}(0)$. 

Let $\nu$ be a valuation of ${\bf k}(X)$ such that the center of $\nu$ on $X$ is $p$, and identifying $\nu$ with an extension of $\nu$ to the quotient field of $\hat{\cal O}_{X,p}$ which dominates $\hat{\cal O}_{X,p}$, we have
$$
\nu(w-\sum_{j=0}^k\phi_{i_j}(y)x^{i_j})>\nu(w-\sum_{j=0}^{k-1}\phi_{i_j}(y)x^{i_j})
$$
for $0\le k\le m$.

Let $p_n$ be the center of $\nu$ on $U_n$, $q_n$ be the center of $\nu$ on $V_n$ in
 the commutative diagram (\ref{eqT172}).

With the notation of the proof of Lemma \ref{LemmaT107}, we see that the rational map 
$U\rightarrow V_{r(1)}$ is a morphism at $p=p_{r(1)}$, and $\overline f_{r(1)}(p)=q_{r(1)}$ has permissible parameters $u_{r(1)}, v, w_{r(1)}$ defined by 
\begin{equation}\label{eqT173}
u=u_{r(1)}^e, w=u_{r(1)}^f(\overline w_{r(1)}+\phi_{i_0}(0)),
w_{r(1)}=\overline w_{r(1)}-[\phi_{i_0}(v)-\phi_{i_0}(0)]
\end{equation}
with $\text{gcd}(e,f)=1$ and $\frac{e}{f}=\frac{a}{i_0}$.

$$
u_{r(1)}=w_{r(1)}=0
$$
are local equations of (a branch of) the fundamental locus of $\overline f_{r(1)}:U\rightarrow V_{r(1)}$ at $q_{r(1)}$. We see from (\ref{eqT173}), (\ref{eqT134}) and (\ref{eqT262}) that $U\rightarrow V_{r(1)}$ is a morphism at $\overline p=\overline p_{r(1)}$, that $r(1)=\overline r(1)$, and $\overline f_{r(1)}(\overline p)=q_{r(1)}$. Further,
$$
u_{r(1)}=\overline w_{r(1)}-[\overline\phi_{\overline i_0}(v)-\overline\phi_{\overline i_0}(0)]=0
$$
are also (formal) local equations of (a branch of) the fundamental locus of $U_{r(1)}\rightarrow V_{r(1)}$ at $q_{r(1)}$. Since $\overline q$ is a general point of $C$, the fundamental locus of
$U_{r(1)}\rightarrow V_{r(1)}$ is a nonsingular curve. Thus 
$$
\phi_{i_0}(v)=\overline\phi_{\overline i_0}(v).
$$

(\ref{eqT262}) implies that 
$$
\text{gcd}(a,i_0)=\frac{a}{e}
$$
and
$$
\text{gcd}(\overline a,\overline i_0)=\frac{\overline a}{e}.
$$
Suppose that we further have that $l\ge 1$, $\frac{i_1}{a}=\frac{\overline i_1}{\overline a}$ 
and $\phi_{i_1}(0)=\overline\phi_{\overline i_1}(0)$.
Then from (\ref{eqT262}) we have
$$
\frac{i_1-i_0}{a}=\frac{\overline i_1-\overline i_0}{\overline a}.
$$
Thus
$$
\frac{\text{gcd}(a,i_0)}{i_1-i_0}=\frac{\text{gcd}(\overline a,\overline i_0)}{\overline i_1-\overline i_0}.
$$

We have (with the notation of (\ref{eqT263}) of the proof of Lemma \ref{LemmaT107}) that
$$
\frac{e_2}{f_2}=\frac{\text{gcd}(a,i_0)}{i_1-i_0},
$$
the rational map $U\rightarrow X_{r(2)}$ is a morphism at $p=p_{r(2)}$, and $q_{r(2)}=\overline f_{r(2)}(p)$ has permissible parameters $u_{r(2)},v,w_{r(2)}$ defined by 
$$
u_{r(1)}=u_{r(2)}^{e_2},
v_{r(1)}=v_{r(2)},
 w_{r(1)}=u_{r(2)}^{f_2}(\overline w_{r(2)}+\phi_{i_1}(0)),
 w_{r(2)}=\overline w_{r(2)}-[\phi_{i_1}(v)-\phi_{i_1}(0)].
 $$

 We see that  $\overline f_{r(2)}$ is a morphism at $\overline p$, and 
 $\phi_{i_1}(0)=\overline \phi_{i_1}(0)$ implies 
that
 $\overline f_{r(2)}(\overline p)=q_{r(2)}$. 

We have
 $$
 u_{r(2)}=w_{r(2)}=\overline w_{r(2)}-[\phi_{i_1}(v)-\phi_{i_1}(0)]=0
 $$
 and
 $$
 u_{r(2)}=\overline w_{r(2)}-[\overline \phi_{i_1}(v)-\overline \phi_{i_1}(0)]=0
 $$
 are local equations of (branches of) the fundamental locus of $\overline f_{r(1)}$ at $q_{r(1)}$. Thus, since the fundamental locus of $\overline f_{r(1)}$
 is nonsingular,
$$
\phi_{i_1}(v)=\overline\phi_{\overline i_1}(v).
$$

Assume that $l\ge 2$,
$$
\frac{i_0}{a}=\frac{\overline i_0}{\overline a}, \frac{i_1}{a}=\frac{\overline i_1}{\overline a}\text{ and }
\frac{i_2}{a}=\frac{\overline i_2}{\overline a}
$$
and $\phi_{i_2}(0)=\overline \phi_{\overline i_2}(0)$ (as well as $\phi_{i_0}(0)=\overline\phi_{i_0}(0)$ and
 $\phi_{i_1}(0)=\overline\phi_{i_1}(0)$).

Then
$$
\frac{\text{gcd}(a,i_0,i_1)}{a}=\frac{\text{gcd}(\overline a,\overline i_0,\overline i_1)}{\overline a}.
$$
Now
$$
\frac{i_2-i_1}{a}=\frac{\overline i_2-\overline i_1}{\overline a}
$$
implies
$$
\frac{\text{gcd}(a,i_0,i_1)}{i_2-i_1}=\frac{\text{gcd}(\overline a,\overline i_0,\overline i_1)}{\overline i_2-\overline i_1}.
$$

 The rational map $U\rightarrow V_{r(3)}$ is a morphism at $p$ and $\overline p$, 
  $\overline f_{r(3)}(p)=\overline f_{r(3)}(\overline p)=q_{r(3)}$, 
 and 
 $$
 \phi_{i_2}(v)=\overline \phi_{i_2}(v)
 $$
  since the fundamental locus of $\overline f_{r(3)}$ is nonsingular.

Iterating, we see that if $j\le t$,
$$
\frac{i_j}{a}=\frac{\overline i_j}{\overline a}
$$
 and
$$
\phi_{i_j}(0)=\overline\phi_{i_j}(0)
$$
 for $j\le t$, then

$$
\frac{\text{gcd}(a,i_0,\ldots,i_{j-1})}{i_{j}-i_{j-1}}
=\frac{\text{gcd}(\overline a,\overline i_0,\ldots,\overline i_{j-1})}{\overline i_{j}-\overline i_{j-1}}
$$
for $j\le t$
and 
$$
\phi_{i_{j}}(v)=\overline\phi_{\overline i_{j}}(v)
$$
 for $j\le t$.

We have verified that a general point of every one dimensional component of $G_Y(f,\tau)$ is perfect. Thus the conclusions of the lemma hold.

\end{pf}

\begin{Definition}\label{DefT151} Suppose that $f:X\rightarrow Y$ is $\tau$-prepared. Let
$V$ be the largest open subset of $Y$ on which the conclusions of Lemma \ref{LemmaT107} hold.
Let $\Theta(Y)=\Theta(f,Y)$ be the set of perfect 1-points in $V \cap G_Y(f,\tau)$.

\end{Definition}

\begin{Remark}\label{RemarkT159} Suppose that $f$ is $\tau$-prepared. Then
$G_Y(f,\tau)-\Theta(f,Y)$ is a finite set by Lemma \ref{LemmaT134}, Lemma \ref{LemmaT107} and Lemma \ref{LemmaT128}.
\end{Remark}

\section{Relations}

In this section, we suppose that $Y$ is a nonsingular projective 3-fold with toroidal structure $D_Y$, and
$f:X\rightarrow Y$ is a dominant proper morphism of nonsingular 3-folds, with toroidal structures $D_Y$ and
$D_X=f^{-1}(D_Y)$, such that $D_X$ contains the singular locus of $f$.

\begin{Definition}\label{DefT56}   A  quasi-pre-relation $R$ on $Y$  is an association $U$ from a locally closed subset $U(R)\subset D_Y$, such that $U(R)$ contains no non trivial open subsets of 2-curves or 3-points and $\text{dim }U(R)\le 1$.

If $q\in U(R)$ is a 2-point, 
$$
R(q)=(S_R(q),(E_1)_R(q),(E_2)_R(q),w_{R(q)},u_{R(q)},v_{R(q)},e_R(q),a_R(q),b_R(q),\lambda_R(q))
$$
with $\text{gcd}(a_R(q),b_R(q),e_R(q))=1$, $e_R(q)>1$, $u_{R(q)},v_{R(q)},w_{R(q)}$ are (possibly formal) permissible parameters at $q$ with $u_{R(q)}, v_{R(q)}\in {\cal O}_{Y,q}$, $0\ne\lambda_R(q)\in {\bf k}$.

We will also allow quasi-pre-relations with $a_R(q)=b_R(q)=\infty$, $e_R(q)=1$ and $\lambda_R(q)=1$.

If $q\in U(R)$ is a 1-point,
$$
R(q)=(S_R(q),E_R(q),w_{R(q)},u_{R(q)},v_{R(q)},e_R(q),a_R(q),\lambda_R(q))
$$
with $\text{gcd}(a_R(q),e_R(q))=1$, $e_R(q)>1$, $u_{R(q)},v_{R(q)},w_{R(q)}$ are (possibly formal) permissible parameters at $q$ with $u_{R(q)}, v_{R(q)}\in{\cal O}_{Y,q}$, $0\ne\lambda_R(q)\in {\bf k}$.

We will also allow quasi-pre-relations with $a_R(q)=\infty$, $e_R(q)=1$ and $\lambda_R(q)=1$.

A restriction $R'$ of a quasi-pre-relation $R$ is the association from a locally closed subset $U(R')$ of $U(R)$ such that $R'(q)=R(q)$ for $q\in U(R')$.
\end{Definition}

Suppose that  $R$ is a  quasi-pre-relation and $q\in U(R)$ is a 2-point.
Let 
$$
u=u_{R(q)}, v=v_{R(q)}, w=w_{R(q)},
a=a_R(q), b=b_R(q), e=e_R(q), \lambda=\lambda_R(q).
$$

If $a,b\ne\infty$, then $R(q)$ is determined by the expression 
\begin{equation}\label{eq151}
w^e-\lambda u^av^b.
\end{equation}
Depending on the signs of $a$ and $b$, this expression determines a (formal) germ at $q$ of an (irreducible)  surface
singularity 
\begin{equation}\label{eq152}
F=F_{R(q)}=0
\end{equation}
 of one of the following forms:
$$
F=w^e-\lambda u^av^b=0
$$
if $a,b\ge 0$ and $a+b>0$,
$$
F=w^eu^{-a}-\lambda v^b=0
$$
if $a<0$, $b>0$,
$$
F=w^ev^{-b}-\lambda u^a=0
$$
if $b<0$, $a>0$.

In the remaining case, $a,b\le 0$,
$$
F=w^eu^{-a}v^{-b}-\lambda
$$
is a unit in $\hat{\cal O}_{Y,q}$ and $F(q)\ne 0$.

If $a,b=\infty$, and $q\in U(R)$ is a 2-point, then $R(q)$ is determined by the expression 
\begin{equation}\label{eq359}
F=F_{R(q)}=w_{R(q)}=0.
\end{equation}

Suppose that $R$ is a  quasi-pre-relation, and  $q\in U(R)$ is a 1-point.
Let 
$$
u=u_{R(q)}, v=v_{R(q)}, w=w_{R(q)},
a=a_R(q), e=e_R(q), \lambda=\lambda_R(q).
$$

Observe that if 
 $a\ne\infty$,  then $R(q)$ is determined by the expression 
\begin{equation}\label{eqT107}
w^e-\lambda u^a.
\end{equation}
This expression determines a (formal) germ at $q$ of an (irreducible)  surface
singularity 
\begin{equation}\label{eqT108}
F=F_{R(q)}=0
\end{equation}
 of   the  form
$$
F=w^e-\lambda u^a=0.
$$
if $a> 0$. In the remaining case, $a\le 0$, so that
$$
F=w^eu^{-a}-\lambda
$$
is a unit in $\hat{\cal O}_{Y,q}$ and $F(q)\ne 0$.

If $a=\infty$, and $q\in U(R)$ is a 1-point, then $R(q)$ is determined by the expression 
\begin{equation}\label{eqT109}
F=F_{R(q)}=w_{R(q)}=0.
\end{equation}

A quasi-pre-relation $R$ is resolved if $F_{R(q)}$ is a unit in $\hat{\cal O}_{Y,q}$ for all
$q\in U(R)$ (This includes the case $U(R)=\emptyset$).

\begin{Definition}\label{Def154}
A subvariety $G$ of $Y$ is an admissible center for a quasi-pre-relation $R$ on $Y$ if one of the
following holds:
\begin{enumerate}
\item[1.] $G$ is a 2-point.
\item[2.] $G$ is a 1-point.
\item[3.] $G$ is a 2-curve of $Y$.
\item[4.] $G\subset D_Y$ is a nonsingular curve which contains a 1-point and makes SNCs with $D_Y$. If
$q\in U(R)\cap G$ then  the (formal) germ of $G$ at $q$ is contained in the germ $w_{R(q)}=0$.
\end{enumerate}
\end{Definition}

Observe that admissible centers are possible centers.

Suppose that $R$ is a   quasi-pre-relation on $Y$, $G$ is an admissible center for $R$,
and $\Psi:Y_1\rightarrow Y$ is the blow up of $G$.

Let $W$ be
 the union over $q\in U(R)$ of points $q_1$ in $\Psi^{-1}(q)$ such that
$q_1$ is on the strict transform of $w_{R(q)}=0$. Assume that this is a locally closed subset of $D_{Y_1}$
of dimension $\le 1$ which contains no 2-curves or 3-points (This condition will always be satisfied when $R$ is a pre-relation, Definition \ref{DefT165}).
The transform $R^1$ of
$R$ on $Y_1$ is then the   quasi-pre-relation on $Y_1$ defined by the condition that
$U(R^1)=W$.
For such $q_1$, $R^1(q_1)$ is determined by the following rules:

If $q\in U(R)\cap G$, and 
$$
u=u_{R(q)}, v=v_{R(q)}, w=w_{R(q)},
$$
 then $G$ has local equations of one of the following forms at $q$ (corresponding to the cases of Definition \ref{Def154}):
\begin{enumerate}
\item[1.,2.] $u=v=w=0$,
\item[3.] $u=v=0$,
\item[4. a)] $q$ a 2-point, $u=w=0$ or $v=w=0$ 
\item[4. b)] $q$ a 1-point, $u=w=0$.
\end{enumerate}

If $q_1\in U(R_1)\cap \Psi^{-1}(q)$, then 
$$
u_1=u_{R^1(q_1)}, v_1=v_{R^1(q_1)},  w_1=w_{R^1(q_1)}
$$
are defined, respectively, by
\begin{enumerate}
\item[1., 2.] $u=u_1, v=u_1(v_1+\alpha), w=u_1w_1$ for some $\alpha\in {\bf k}$,
or $u=u_1v_1, v=v_1, w=v_1w_1$,
\item[3.] $u=u_1, v=u_1(v_1+\alpha), w=w_1$ for some $\alpha\in {\bf k}$,
or $u=u_1v_1, v=v_1, w=w_1$,
\item[4. a)] $u=u_1, w=u_1w_1$ or $v=v_1, w=v_1w_1$
\item[4. b)] $u=u_1,v=v_1, w=u_1w_1$.
\end{enumerate}
\vskip .2truein
Suppose that $q_1\in U(R^1)\cap\Psi^{-1}(q)$ and case 1  holds. 
Suppose that
$$
R(q)=(S,E,w,u,v,e,a,b,\lambda).
$$
If $0\ne\alpha$ and $a,b\ne\infty$, we define
$$
a_{R^1}(q_1)=a+b-e, b_{R^1}(q_1)=b, e_{R^1}(q_1)=e, \lambda_{R^1}(q_1)=\lambda\alpha^b.
$$
$R^1(q_1)$ is thus determined by
$$
w_1^e-\lambda\alpha^bu_1^{a+b-e}.
$$
If $a,b\ne\infty$ and
$$
u=u_1, v=u_1v_1, w=u_1w_1,
$$
 we define
$$
a_{R^1}(q_1)=a+b-e, b_{R^1}(q_1)=b, e_{R^1}(q_1)=e,\lambda_{R^1}(q_1)=\lambda.
$$
$R^1(q_1)$ is determined by $w_1^e-\lambda u_1^{a+b-e}v_1^b$.

If $a,b\ne\infty$ and
$$
u=u_1v_1, v=v_1, w=v_1w_1,
$$
we define
$$
a_{R^1}(q_1)=a, b_{R^1}(q_1)=a+b-e, e_{R^1}(q_1)=e,\lambda_{R^1}(q_1)=\lambda.
$$
$R^1(q_1)$ is determined by $w_1^e-\lambda u_1^av_1^{a+b-e}$.

If $a=b=\infty$, we define
$$
a_{R^1}(q_1)=\infty,
b_{R^1}(q_1)=\infty,
e_{R^1}(q_1)=1,
\lambda_{R^1}(q_1)=1,
$$
and $R^1(q_1)$ is determined by $w_{R^1(q_1)}$.

In cases 3 and 4, we define $R^1(q_1)$ in an analogous way.

Suppose that $\Psi_1:Y_1\rightarrow Y$ is a sequence of blow ups of admissible centers for (the transforms of) $R$, $R^1$ is the transform of $R$ on $Y_1$, $q\in U(R)$ and $q_1\in\Psi_1^{-1}(q)\cap U(R^1)$. Let
$$
u=u_{R(q)}, v=v_{R(q)}, w=w_{R(q)}.
$$
$$
u_1=u_{R^1(q_1)}, v_1=v_{R^1(q_1)}, w_1=w_{R^1(q_1)}.
$$
$u$ and $v$ are related to $u_1,v_1$ birationally. That is, $k(u,v)=k(u_1,v_1)$. We  have one of the following expressions: 
\begin{equation}\label{eqT267}
u=u_1^av_1^b,
v=u_1^cv_1^d,
w=u_1^ev_1^fw_1
\end{equation}
with $ad-bc=\pm1$, 
\begin{equation}\label{eqT268}
u=u_1^a\gamma_1(u_1,v_1),
v=u_1^b\gamma_2(u_1,v_1),
w=u_1^c\gamma_3(u_1,v_1)w_1
\end{equation}
where $\gamma_1,\gamma_2,\gamma_3$ are unit series, 
\begin{equation}\label{eqT269}
u=(u_1^av_1^b)^t\gamma_1(u_1,v_1),
v_1=(u_1^av_1^b)^k\gamma_2(u_1,v_1),
w=u_1^ev_1^f\gamma_3(u_1,v_1)w_1
\end{equation}
where $\text{gcd}(a,b)=1$, and $\gamma_1,\gamma_2,\gamma_3$ are unit series.

\begin{Definition}\label{DefT165} A quasi-pre-relation $R$ is a pre-relation if
there exists a nonsingular 3-fold $\overline Y_R$ with toroidal structure $D_{\overline Y_1}$, a pre-relation $R^0$ on $\overline Y_R$ such that $U(R^0)=\{\overline q\}$ is a single point with 
$$
R^0(\overline q)=(\cdots,w_{R^0(\overline q)},u_{R^0(\overline q)}, v_{R^0(\overline q)},\cdots)
$$
 and a sequence of possible blow ups 
$$
\overline\Psi_R:Y=Y_n\rightarrow\cdots\rightarrow Y_1\rightarrow Y_0=\overline Y_R
$$
where each $Y_i$ has a quasi-pre-relation $R^i$ which is the restriction of the transform of $R^{i-1}$, and $Y_{i+1}\rightarrow Y_i$ is an admissible blow up for $R^i$, and $R=R^n$.
\end{Definition}

\begin{Definition}\label{Def199}
A  pre-relation $R$ on $Y$ is algebraic if
there exists an open subset $V$ of $Y$ and a nonsingular irreducible  closed surface $\Omega(R)\subset V$ such that $\Omega(R)$ makes SNCs with $D_Y$, $U(R)\subset\Omega(R)$ and
$S_R(q)$ is the (formal) germ of $\Omega(R)$ at $q$ for all $q\in U(R)$.
Further, if $q\in\Omega(R)\cap D_Y$, then there exist super parameters $u_R,v_R,w_R$ at $q$ such that $w_R=0$ is a local equation of $\Omega(R)$.
\end{Definition}

Suppose that $R$ is algebraic, and $\Psi:Y_1\rightarrow Y$ is an admissible blow up. then after possibly replacing $\Omega(R)$ with an open subset of $\Omega(R)$ (containing $U(R)$), we have that the transform $R^1$ of $R$ is algebraic, where $\Omega(R^1)$ is the
strict transform of $\Omega(R)$ by $\Psi$.

\begin{Definition}\label{DefT57}  Suppose that $f:X\rightarrow Y$ is prepared. 

A primitive  relation $R$ for $f$ is

\begin{enumerate}
\item[1.] A  pre-relation $\overline R$ on $Y$.
\item[2.] A locally closed subset $T(R)\subset  f^{-1}(U(\overline R))$ such that if $p\in T(R)$ and
$f(p)$ is a 2-point with
$$
\overline R(f(p))=(S,E_1,E_2,w,u,v,e,a,b,\lambda),
$$
then $u,v$ are toroidal forms at $p$. If $a,b\ne\infty$, then there exist permissible parameters $x,y,z$ at $p$ for $u,v,w$ such that

\begin{equation}\label{eqT58}
w^e=u^av^b\overline\Lambda(x,y,z)
\end{equation}
where $\overline\Lambda(0,0,0)=\lambda$.

If $a=b=\infty$, then $u,v,w$ have a monomial form (Definition \ref{Def125}) at $p$.

If $f(p)$ is a 1-point with
$$
\overline R(f(p))=(S,E,w,u,v,e,a,\lambda),
$$
then $u,v$ are toroidal forms at $p$.
If $a\ne\infty$, then there exist permissible parameters $x,y,z$ at $p$ for $u,v,w$ such that

\begin{equation}\label{eqT59}
w^e=u^a\overline\Lambda(x,y,z)
\end{equation}
where $\overline\Lambda(0,0,0)=\lambda$.
\end{enumerate}
If $a=\infty$, then $u,v,w$ have a  monomial form (Definition \ref{Def125}) at $p$.

In all cases, we define $R(p)=\overline R(f(p))$.

A relation $R$ for $f$ is a finite set of  pre-relations $\{\overline R_i\}$ on $Y$
with associated primitive  relations $R_i$ for $f$ such that the sets $T(R_i)$ are pairwise disjoint.

We denote $U(R)=\cup_iU(\overline R_i)$ and $T(R)=\cup_iT(R_{i})$, and
define
$$
R(p)=R_{i}(p)
$$
if $p\in T(R_{i})$.

If $U(\overline R_i)\cap U(\overline R_j)\ne\emptyset$, then we further require that $Y_{\overline R_i}=Y_{\overline R_j}$ (with the notation of Definition \ref{DefT165}) and $u_{(\overline R_i)^0(\overline q)}=u_{(\overline R_j)^0(\overline q)}$, $v_{(\overline R_i)^0(\overline q)}=v_{(\overline R_j)^0(\overline q)}$.
This implies  that 
$$
u_{\overline R_i(q)}=u_{\overline R_j(q)}, v_{\overline R_i(q)}=v_{\overline R_j(q)}
$$
if $q\in U(\overline R_i)\cap U(\overline R_j)$.
We will call $\{\overline R_i\}$  the  pre-relations associated to $R$.
We will say that $R$ is algebraic if each $\overline R_i$ is algebraic and 
\begin{equation}\label{eq253}
\Omega(\overline R_i)\cap U(R)=U(\overline R_i)
\end{equation}
for all $i$. We will also denote $\Omega(R_i)=\Omega(\overline R_i)$. For $p\in T(R_i)\subset T(R)$, we denote
$$
R(p)=\left(\begin{array}{l}
S=S_{R}(p),E_1(p), E_2(p),w=w_{ R(p)},u=u_{R(p)},\\
v=v_{R(p)},
e=e_{R}(p),
a=a_{R}(p),b=b_{R}(p),\lambda=\lambda_{R}(p)\end{array}\right)
$$
if $f(p)$ is a 2-point,
$$
R(p)=\left(\begin{array}{l}
S=S_{R}(p),E(p),w=w_{ R(p)},u=u_{R(p)},\\
v=v_{R(p)},e=e_{R}(p),
a=a_{R}(p),\lambda=\lambda_{R}(p)\end{array}\right)
$$
if $f(p)$ is a 1-point.
\end{Definition}

A   relation $R$ is resolved if $T(R)=\emptyset$.

\begin{Definition}\label{Def161}
Suppose that $f:X\rightarrow Y$ is prepared, $R$ is a 
relation for $f$ and
$$
\begin{array}{rll}
X_1&\stackrel{f_1}{\rightarrow}&Y_1\\
\Phi\downarrow&&\downarrow\Psi\\
X&\stackrel{f}{\rightarrow}&Y
\end{array}
$$
is a commutative diagram such that
\begin{enumerate}
\item[1.] $\Psi$ is a product of blow ups which are admissible for all of the pre-relations
$\overline R_i$ associated to $R$ (and their transforms) and $\Phi$ is a
product of blow ups of possible centers
\item[2.] $f_1$ is prepared.
\item[3.] Let $\overline R_i^1$ be the transforms of the $\overline R_i$ on $Y_1$ and let
$$
T_i=\{p\in X_1\mid 
p\in\Phi^{-1}(T(R_{i}))\cap f_1^{-1}(U(\overline R_i^1))\}.
$$

Suppose that  $p\in T_i$ then $u_{\overline R_i^1(f_1(p))}, v_{\overline R_i^1(f_1(p))}$ are toroidal forms at $p$.
\end{enumerate}
 Then the transform $R^1$ of $R$ for $f_1$ is the 
relation for $f_1$
defined by  
$$
T(R^1)=\cup T_i,
$$
$$
R^1(p)=\overline R_i^{1}(f_1(p))
$$
for $p\in T_i$.

It is straightforward to verify that $R^1$ satisfies the conditions of
Definition \ref{DefT57}, substituting from (\ref{eqT58}), (\ref{eqT59}) and (\ref{eqT267}) - (\ref{eqT269}) into (\ref{torf}).
\end{Definition}

\section{Well Prepared Morphisms}

Suppose that $\tau\in{\bf N}$ and $f:X\rightarrow Y$ is a dominant, proper, $\tau$-prepared morphism of nonsingular 3-folds with
toroidal structures $D_Y$ and $D_X=f^{-1}(D_Y)$. Further suppose that the singular locus of $f$ is contained in $D_X$. If $R$ is a  relation for $f$
with associated  pre-relations $\{\overline R_i\}$,
 then for $p\in T(R_i)$ such that $f(p)=q$ is a 2-point, we have
that 
\begin{equation}\label{eq311}
R(p)=\left(\begin{array}{l} 
S_i=S_R(p),E_1=E_{R,1}(p),E_2=E_{R,2}(p),w_i=w_{R(p)},
u=u_{R(p)},\\
v=v_{R(p)},e_i=e_R(p),a_i=a_R(p),b_i=b_R(p),
\overline\lambda_i=\lambda_R(p)
\end{array}\right)
\end{equation}
which we will abbreviate (as in (\ref{eq151})) as 
\begin{equation}\label{eq168}
R(p)=w_i^{e_i}-\overline \lambda_iu^{a_i}v^{b_i},
\end{equation}
with $e_i>1$, if $a_i,b_i\ne\infty$, or (as in (\ref{eq359})) 
\begin{equation}\label{eq254}
R(p)=w_i
\end{equation}
if $a_i,b_i=\infty$. In this case, $u,v,w_i$ have a monomial form (Definition \ref{Def125}) at $p$.

For $p\in T(R_i)$ such that $f(p)=q$ is a 1-point, we have
that 
\begin{equation}\label{eqT104}
R(p)=\left(\begin{array}{l} 
S_i=S_R(p),E_1=E_{R,1}(p),w_i=w_{R(p)},
u=u_{R(p)},\\
v=v_{R(p)},e_i=e_R(p),a_i=a_R(p),
\overline\lambda_i=\lambda_R(p)
\end{array}\right)
\end{equation}
which we will abbreviate (as in (\ref{eq151})) as 
\begin{equation}\label{eqT98}
R(p)=w_i^{e_i}-\overline \lambda_iu^{a_i},
\end{equation}
with $e_i>1$, if $a_i\ne\infty$, or (as in (\ref{eq359})) 
\begin{equation}\label{eqT99}
R(p)=w_i
\end{equation}
if $a_i=\infty$. In this case, $u,v,w_i$ have a monomial form (Definition \ref{Def125}) at $p$.

Recall that if $p'\in T(R)$ is such that $f(p')=f(p)$, then $u_{R(p')}=u_{R(p)}=u$ and $v_{R(p')}=v_{R(p)}=v$.
Let $I$ be an index set for the pre-relations $\{\overline R_i\}$ associated to $R$.

\begin{Definition}\label{DefT60}
Suppose that $\tau\ge 0$. A $\tau$-prepared morphism $f:X\rightarrow Y$ is pre-$\tau$-quasi-well prepared with  relation $R$ if: 
\begin{enumerate}
\item[1.]  $T(R)=G_X(f,\tau)\cap f^{-1}(U(R))$
\item[2.] Suppose that $p\in T(R)$. Then $\tau>0$ implies $R(p)$ has a form (\ref{eq168}) or (\ref{eqT98}), $\tau=0$ implies $R(p)$ has a form (\ref{eq254}) or (\ref{eqT99}).
\item[3.] If $q\in U(\overline R_i)\cap U(\overline R_j)$, then there exists  $\lambda_{ij}(u,v)\in {\bold k}[[u,v]]$, 
with 
$$
u=u_{\overline R_i(q)}=u_{\overline R_j(q)}, v=v_{\overline R_i(q)}=v_{\overline R_j(q)},
w_i=w_{\overline R_i(q)}, w_j=w_{\overline R_j(q)}
$$
 such that 
$$
w_j=w_{i}+\lambda_{ij}(u,v),
$$
 and with the notation of Definition \ref{DefT57}, 
there exists a series
$$
(\lambda_{ij})^0(u_{(\overline R_j)^0(\overline q)},v_{(\overline R_j)^0(\overline q)})
$$
such that
$$
w_{(\overline R_j)^0(\overline q)}-w_{(\overline R_i)^0(\overline q)}=(\lambda_{ij})^0(u_{(\overline R_j)^0(\overline q)},v_{(\overline R_j)^0(\overline q)}),
$$
and $\lambda_{ij}(u,v)$ is obtained from $\lambda_{ij}^0(u_{(\overline R_j)^0(\overline q)},v_{(\overline R_j)^0(\overline q)})$ from the appropriate expression (\ref{eqT267}) - (\ref{eqT269}).
\item[4.] Suppose that $q\in U(\overline R_i)$,
where $\overline R_i$ is a  relation associated to $R$. Then
 $u=u_{\overline R_i(q)},v=v_{\overline R_i(q)},w_i=w_{\overline R_i(q)}$ are super parameters at $q$ (Definition \ref{Def357}).
\end{enumerate}
$f$ is $\tau$-quasi-well prepared with relation $R$ if $f$ is pre-$\tau$-quasi-well prepared with $T(R)=G_X(f,\tau)$.
\end{Definition}

 We will allow a $\tau$-prepared morphism without relation ($U(R)=\emptyset$) as a type of pre-$\tau$-quasi-well prepared morphism.

\begin{Definition}
$f:X\rightarrow Y$ is $\tau$-quasi-well prepared with
 relation $R$ and pre-algebraic structure  if $f$ is
 $\tau$-quasi-well prepared with  relation $R$ and
 $u_{\overline R_i(q)}, v_{\overline R_i(q)}, w_{\overline R_i(q)}\in{\cal O}_{Y,q}$ for all $\overline R_i$ associated to $R$, and $q\in U(\overline R_i)$.
\end{Definition}

\begin{Definition}\label{Def65} $f:X\rightarrow Y$ is $\tau$-well prepared with
 relation $R$   if
\begin{enumerate}
\item[1.] $f$ is $\tau$-quasi-well prepared with   relation
$R$ and pre-algebraic structure.
\item[2.] The primitive pre-relations $\{\overline R_i\}$ associated to $R$ are algebraic, and $R$ is algebraic (Definition \ref{DefT57}).
\item[3.] Suppose that  $q\in U(\overline R_i)\cap U(\overline R_j)$.
 Let $w_i=w_{\overline R_i(q)}$ and
 $w_j=w_{\overline R_j(q)}$,  $u=u_{\overline R_i(q)}=u_{\overline R_j(q)}, v=v_{\overline R_i(q)}=v_{\overline R_j(q)}$. 

Suppose that $q$ is a 2-point. Then there exists a  unit series
$\phi_{ij}\in {\bold k}[[u,v]]$ and $a_{ij}, b_{ij}\in{\bf N}$ (or $\phi_{ij}=0$ with $a_{ij}=b_{ij}=-\infty$)  with 
\begin{equation}\label{eq64}
w_j=w_i+u^{a_{ij}}v^{b_{ij}}\phi_{ij}.
\end{equation}

Suppose that $q$ is a 1-point. Then there exists a  unit series
$\phi_{ij}\in {\bold k}[[u,v]]$ and $c_{ij}\in{\bf N}$ (or $\phi_{ij}=0$ with $c_{ij}=-\infty$) with 
\begin{equation}\label{eqT100}
w_j=w_i+u^{c_{ij}}\phi_{ij}.
\end{equation}

\item[4.] For $q\in U(R)$ a 2-point, set  $I_q=\{i\in I\mid q\in U(\overline R_i)\}$. Then the set 
\begin{equation}\label{eq255}
\left\{(a_{ij},b_{ij})\mid i,j\in I_q\right\}
\end{equation}
 from equation (\ref{eq64}) is totally ordered.
\end{enumerate}
\end{Definition}

\begin{Definition}\label{DefT80} Suppose that $f:X\rightarrow Y$ is pre-$\tau$-quasi-well prepared, $\overline q\in Y$ is a 1-point such that $\overline q\not\in U(R)$. A curve  $C\subset D_Y$  such that $\overline q\in C$ is called a {\bf resolving curve} for $f$ and $R$ at $\overline q$ if
\begin{enumerate}
\item[1.]  $C$ makes SNCs with $D_Y$.
\item[2.] $C\cap G_Y(f,\tau)\subset\{\overline q\}$.
\item[3.] If $q\in C$ is a 2-point, then there exist super parameters $u,v,w$ at $q$ such that $u=w=0$ are local equations of $C$ at $q$.
\item[4.] If $q\in C$ is a 1-point, then there exist permissible parameters $u,v,w$ at $q$ such that
$u,v$ are toroidal forms at $p$ for all $p\in f^{-1}(q)$, and $u=v=0$ are local equations of $C$ at $q$.
\end{enumerate}
\end{Definition}

\begin{Definition}\label{Def66}
Suppose that $f:X\rightarrow Y$ is pre-$\tau$-quasi-well prepared with   relation $R$.
\begin{enumerate}
\item[1.] A 2-point $q\in U(R)$ is prepared for $R$.
\item[2.] A 1-point or 2-point $q\in Y$ such that $q\not\in U(R)$ is prepared for $R$ if there exist
super parameters $u,v,w$ at $q$.
\item[3.] A 2-curve $C\subset Y$ is prepared for $R$.
\item[4.] A 1-point $q\in U(R)$ is prepared.
\item[5.] A resolving curve $C$ for $f$ and $R$ at a 1-point $q\not\in U(R)$  is prepared.
\end{enumerate}
\end{Definition}

If $E$ is a component of $D_Y$, $\overline R_i$ is pre-algebraic, and $q\in U(\overline R_i)$, we will denote
$\overline{E\cdot S_{\overline R_i}(q)}$ for the Zariski closure in $Y$ of the curve germ $u=w_{\overline R_i(q)}=0$ at $q$, where $u=0$ is
a local equation of $E$.

\begin{Definition}\label{Def200}
Suppose that $f:X\rightarrow Y$ is $\tau$-well prepared with  relation $R$ for $f$.
 A nonsingular curve $C\subset D_Y$ which makes SNCs with $D_Y$ is prepared for $R$
of type 6 if
\begin{enumerate}
\item[1.] $C=\overline{E_{\alpha}\cdot S_{\overline R_i}(q_{\beta})}$ for some
component $E_{\alpha}$ of $D_Y$,  pre-relation $\overline R_i$ associated to
$R$ and $q_{\beta}\in U(\overline R_i)$.
\item[2.] $\Omega(\overline R_i)$ contains $C$.
\item[3.] If $C'=\overline{E_{\gamma}\cdot S_{\overline R_j}(q_{\delta})}$ is such that
$C'\subset\Omega(\overline R_j)$, $C\ne C'$,
and $q\in C\cap C'$, then $q\in U(\overline R_i)\cap U(\overline R_j)$ and
$C'=\overline{E_{\gamma}\cdot S_{\overline R_j}(q)}$.
\item[4.] If $j\ne i$ and $C=\overline{E_{\gamma}\cdot S_{\overline R_j}(q_{\delta})}$
then $C$ satisfies 1 and 2 and 3 of this definition (for $\overline R_j$). (In this case we have by
(\ref{eq253}) that
$U(\overline R_j)\cap C=U(\overline R_i)\cap C$).
\item[5.]
Let
$$
I_C=\{j\in I\mid C=\overline{E_{\gamma}\cdot S_{\overline R_j}(q_{\delta})}
\text{ for some }\overline R_j, E_{\gamma}, q_{\delta}\in U(\overline R_j)\}.
$$
 Suppose that $q\in C$ is a 1-point or a 2-point such that $q\not \in U(R)$.
Then there exist $u,v\in{\cal O}_{Y,q}$ such that for $j\in I_C$ there exists  $\tilde w_j\in
{\cal O}_{Y,q}$ such that
\begin{enumerate}
\item $\tilde w_j=0$ is a local equation of $\Omega(\overline R_j)$ and
 $u,v,\tilde w_j$ are permissible parameters at $q$ such that $u=\tilde w_j=0$ are local
equations of $C$ at $q$.
\item   $u,v,\tilde w_j$ are super parameters at $q$.
\item If $i,j\in I_C$ and $q$ is a 1-point, there exist relations
$$
\tilde w_i-\tilde w_j=u^{c_{ij}}\phi_{ij}(u,v)
$$
where $\phi_{ij}$ is a unit series (or $\phi_{ij}=0$ and $c_{ij}=-\infty$).
\item If $i,j\in I_C$ and $q$ is a 2-point (with $q\not\in U(R)$) then
there exist relations
$$
\tilde w_i-\tilde w_j=u^{c_{ij}}v^{d_{ij}}\phi_{ij}(u,v)
$$
where $\phi_{ij}$ is a unit series (or $\phi_{ij}=0$ and $c_{ij}=d_{ij}=-\infty$), and the set $\{(c_{ij},d_{ij})\}$ is totally ordered.

\end{enumerate}
\end{enumerate}
\end{Definition}

If $f:X\rightarrow Y$ is $\tau$-well prepared with  relation $R$, and $\overline R_i$ is a pre-relation
associated to $R$, we will feel free to replace $\Omega(\overline R_i)$ with an open subset of $\Omega(\overline R_i)$ containing $U(\overline R_i)$, and all curves $C=\overline{E\cdot S_{\overline R_i}(q)}$ such that $E$ is a component of
$D_Y$, $q\in U(\overline R_i)$ and $C$ is prepared for $R$ of type 6. This convention will allow some simplification of the statements
of the theorems and proofs.

\begin{Definition}\label{Def130}  $f:X\rightarrow Y$ is
$\tau$-very-well prepared  with  relation $R$ if
\begin{enumerate}
\item[1.] $f$ is $\tau$-well prepared with  relation $R$.
\item[2.] If $E$ is a component of $D_Y$ and $q\in U(\overline R_i)\cap E$,
 then $C=\overline{E\cdot S_{\overline R_i}(q)}$
is prepared for $R$
of type 6 (Definition \ref{Def200}).
\item[3.] For all $\overline R_i$ associated to $R$, let
$$
V_i(Y)=\left\{\gamma=\overline{E_{\alpha}\cdot S_{\overline R_{i}}(q_{\gamma})}\mid
q_{\gamma}\in U(\overline R_i), E_{\alpha}\text{ is a component of }D_Y\right\}.
$$
Then
$$
F_i=\sum_{\gamma\in V_i(Y)}\gamma
$$
is a SNC divisor on $\Omega(\overline R_i)$ whose intersection graph
is a forest (its connected components are trees).
\end{enumerate}
\end{Definition}

If $f:X\rightarrow Y$ is $\tau$-very-well prepared, we will feel free to replace $\Omega(\overline R_i)$ with an open neighborhood of $F_i$ in $\Omega(\overline R_i)$.
This will allow some simplification of the proofs.

\begin{Remark}\label{Remark281}  Suppose that $f:X\rightarrow Y$ is $\tau$-very well
prepared. Then it follows from Definition \ref{Def130} and (\ref{eq253}) that $F_i\cap U(R)=U(\overline R_i)$
for all $\overline R_i$ associated to $R$.
\end{Remark}

\begin{Definition}\label{Def289}
Suppose that $f:X\rightarrow Y$ is pre-$\tau$-quasi-well prepared (or $\tau$-well prepared
or $\tau$-very-well prepared)
with   relation $R$.
Let $\{\overline R_i\}$ be the pre-relations  associated to $R$.
Suppose that $G$ is a
 point or nonsingular curve in $Y$ which is an admissible center for all of the $\overline R_i$. Then $G$ is called
a permissible center for $R$ if there exists a commutative diagram 
 \begin{equation}\label{eq30}
\begin{array}{rll}
X_1&\stackrel{f_1}{\rightarrow}&Y_1\\
\Phi\downarrow&&\downarrow\Psi\\
X&\stackrel{f}{\rightarrow}&Y
\end{array}
\end{equation}
where  $\Psi$ is the blow up of $G$ and $\Phi$ is a sequence of blow ups
$$
X_1=\overline X_n\rightarrow \cdots\rightarrow \overline X_1\rightarrow X
$$
of
nonsingular curves and 3-points $\gamma_i$
 which are possible centers such that
\begin{enumerate}
\item[1.] $f_1$ is $\tau$-prepared and the assumptions of Definition \ref{Def161} hold so that the transform $R^1$ of $R$
for $f_1$ is defined.
\item[2.] 
$$
\tau_{f_1}(p)\le \tau_f(\phi(p))
$$
for $p\in D_{X_1}$.
\item[3.] $f_1:X_1\rightarrow Y_1$ is pre-$\tau$-quasi-well prepared, (or $\tau$-well prepared or
$\tau$-very-well prepared)
with   relation $R^1$.
\end{enumerate}
\end{Definition}

(\ref{eq30}) is called a pre-$\tau$-quasi-well prepared (or $\tau$-well prepared or $\tau$-very-well prepared) diagram of $R$ (and $\Psi$).

\begin{Definition}\label{Def219}
Suppose that  $f:X\rightarrow Y$ is  $\tau$-well prepared (or $\tau$-very-well prepared) with relation
$R$ and $C\subset Y$ is prepared for $R$ of type 6. Then $C$ is a $*$-permissible center
for $R$ if there exists a commutative diagram 
\begin{equation}\label{eq233}
\begin{array}{rll}
X_1&\stackrel{f_1}{\rightarrow}&Y_1\\
\Phi\downarrow&&\downarrow\Psi\\
X&\stackrel{f}{\rightarrow}&Y
\end{array}
\end{equation}
such that
\begin{enumerate}
\item[1.] $f_1$ is $\tau$-prepared and the assumptions of Definition \ref{Def161} hold so that
the transform $R^1$ of $R$ for $f_1$ is defined.
\item[2.]
$$
\tau_{f_1}(p)\le \tau_f(\phi(p))
$$
for $p\in D_{X_1}$.
\item[3.] $f_1:X_1\rightarrow Y_1$ is pre-$\tau$-well prepared (or $\tau$-very-well prepared).
\item[4.] (\ref{eq233}) has a factorization 
\begin{equation}\label{eq220}
\begin{array}{rcl}
X_1=\overline X_m&\stackrel{f_1=\overline f_m}{\rightarrow}&\overline Y_m=Y_1\\
\downarrow&&\downarrow\\
\vdots&&\vdots\\
\downarrow&&\downarrow\\
\overline X_2&\stackrel{\overline f_2}{\rightarrow}&\overline Y_2\\
\overline\Phi_2\downarrow&&\downarrow\overline\Psi_2\\
\overline X_1&\stackrel{\overline f_1}{\rightarrow}&\overline Y_1\\
\overline \Phi_1\downarrow&&\downarrow\overline\Psi_1\\
X&\stackrel{f}{\rightarrow}&Y
\end{array}
\end{equation}
where $\overline\Psi_1$ is the blow up of $C$, 
\begin{equation}\label{eq310}
\begin{array}{rll}
\overline X_1&\stackrel{\overline f_1}{\rightarrow}&\overline Y_1\\
\overline\Phi_1\downarrow&&\downarrow\overline\Psi_1\\
X&\stackrel{f}{\rightarrow}&Y
\end{array}
\end{equation}
is a $\tau$-well prepared diagram of $R$ and $\overline\Psi_1$ of the form (\ref{eq30}), each
$\overline \Psi_{i+1}:\overline Y_{i+1}\rightarrow \overline Y_i$ for $i\ge 1$ is
the blow up of a 2-point $q\in\overline Y_i$ which is prepared for the transform $R^i$
of $R$ on $\overline X_i$ of type 2 of Definition \ref{Def66}, and
$$
\begin{array}{rll}
\overline X_{i+1}&\stackrel{\overline f_{i+1}}{\rightarrow} &\overline Y_{i+1}\\
\overline\Phi_{i+1}\downarrow&&\downarrow\overline\Psi_{i+1}\\
\overline X_{i}&\stackrel{\overline f_i}{\rightarrow}&\overline Y_{i}
\end{array}
$$
is a $\tau$-well prepared diagram of $R^i$ and $\overline\Psi_{i+1}$ of the form of (\ref{eq30}).
\item[5.] Suppose that $E$ is the strict transform of $\overline\Psi_1^{-1}(C)$ on $Y_1$. then $\overline{E\cdot \overline R_i^1(q)}$ is prepared for $R^1$ of type 6 for all primitive relations $\overline R_i^1$ associated to $R^1$, and $q\in U(R_i^1)\cap E$.
\item[6.] Suppose that $\gamma\subset Y$ is a curve which is prepared for $R$ of type 6. Then the strict transform of $\gamma$ on $Y_1$ is prepared for $R^1$ of type 6.

\end{enumerate}
\end{Definition}

\begin{Definition}\label{Def396}
Suppose that  $f:X\rightarrow Y$ is  pre-$\tau$-quasi-well prepared (or $\tau$-well prepared or $\tau$-very-well prepared) with  relation $R$.
Suppose that 
\begin{equation}\label{eq405}
\begin{array}{rll}
X_1&\stackrel{f_1}{\rightarrow}&Y_1\\
\Phi\downarrow&&\downarrow\Psi\\
X&\stackrel{f}{\rightarrow}&Y
\end{array}
\end{equation}
is a commutative diagram 
such that there is a factorization 
\begin{equation}\label{eq412}
\begin{array}{rcl}
X_1=\overline X_m&\stackrel{f_1=\overline f_m}{\rightarrow}&\overline Y_m=Y_1\\
\downarrow&&\downarrow\\
\vdots&&\vdots\\
\downarrow&&\downarrow\\
\overline X_2&\stackrel{\overline f_2}{\rightarrow}&\overline Y_2\\
\overline\Phi_2\downarrow&&\downarrow\overline\Psi_2\\
\overline X_1&\stackrel{\overline f_1}{\rightarrow}&\overline Y_1\\
\overline \Phi_1\downarrow&&\downarrow\overline\Psi_1\\
X&\stackrel{f}{\rightarrow}&Y
\end{array}
\end{equation}

where each commutative diagram
$$
\begin{array}{rll}
\overline X_{i+1}&\rightarrow &\overline Y_{i+1}\\
\overline\Phi_{i+1}\downarrow&&\downarrow\overline\Psi_{i+1}\\
\overline X_{i}&\rightarrow&\overline Y_{i}
\end{array}
$$
is either of the form (\ref{eq30}) or of the form (\ref{eq233}). Then (\ref{eq405}) is called a pre-$\tau$-quasi-well prepared
(or $\tau$-well prepared or $\tau$-very-well prepared) diagram of $R$ (and $\Psi$).
\end{Definition}

\section{Construction of $\tau$-well prepared diagrams}

\begin{Lemma}\label{Lemma31}
Suppose that $f:X\rightarrow Y$ is pre-$\tau$-quasi-well prepared
(or $\tau$-well prepared or $\tau$-very-well prepared) and $C\subset Y$ is a  2-curve. Then $C$ is a permissible center for $R$, and there exists  a pre-$\tau$-quasi-well-prepared
(or $\tau$-well prepared or $\tau$-very-well prepared) diagram
(\ref{eq30}) of $R$ and the blow up $\Psi:Y_1\rightarrow Y$ of $C$
such that $\Phi$ is a product of blow ups of 2-curves. Furthermore,
\begin{enumerate}
\item[1.] $\Phi$ is an isomorphism over $f^{-1}(Y-C)$
\item[2.] Further suppose that $f$ is $\tau$-well prepared. Then
\begin{enumerate}
\item  Let $E$ be the exceptional divisor for $\Psi$. Suppose that
$q\in U(\overline R_i^1)\cap E$ for some $\overline R_i$ associated to $R$. Let $\gamma_i=\overline{S_{\overline R_i^1}(q)\cdot E}$.
Then $\gamma_i=\Psi^{-1}(\Psi(q))$ is a prepared curve for $R^1$ of type 6.
\item If $\gamma$ is a prepared curve for $R$, then the strict transform of $\gamma$ on $Y_1$ is
a prepared curve for $R^1$.
\end{enumerate}
\item[3.] Suppose that $\overline q\in C$, $p\in f^{-1}(\overline q)$, $p'\in \Phi_1^{-1}(p)$, $u,v,w$ are permissible parameters at $\overline q$ such that $u=v=0$ are local equations of $C$, and $w$ is good at $p$. If $w$ is not good at $p'$, then $\tau_{f_1}(p')<\tau_f(p)$.
\end{enumerate}
\end{Lemma}

\begin{pf} This follows from Lemma \ref{Lemma1} and a straight forward extension of the proofs of Lemma 5.2 \cite{C5} and  Lemma 7.11 \cite{C5}.
3 follows from calculations as will be given in detail in the proof of 6 of Lemma \ref{LemmaT79}.
\end{pf}

The proofs of Remarks \ref{RemarkT278} and \ref{Remark424} follow easily from the methods of the proof of Lemma \ref{Lemma31}.

\begin{Remark}\label{RemarkT278} Suppose that $f:X\rightarrow Y$ is pre-$\tau$-quasi-well prepared (or $\tau$-well prepared or $\tau$-very-well prepared) and $C\subset D_X$ is a 2-curve or a 3-point. Let $\Phi_1:X_1\rightarrow X$ be the blow up of $C$, $f_1=f\circ\Phi_1:X_1\rightarrow Y$. Then
\begin{enumerate}
\item[1.] $f_1$ is pre-$\tau$-quasi-well prepared (or $\tau$-well prepared or $\tau$-very well prepared).
\item[2.] Suppose that $p_1\in X_1$, $p=\Phi_1(p_1)$, $q=f_1(p)$, $u,v,w$ are permissible parameters at $q$ such that $w$ is good (weakly good) at $p$ for $f$. If $w$ is not good (weakly good) at $p_1$ for $f$, then $\tau_{f_1}(p_1)<\tau_f(p)$.
\item[3.] $(f\circ\Phi_1)^{-1}(\Theta(f,Y))\subset\Theta(f_1,Y_1)$.
\end{enumerate}
\end{Remark}

\begin{Remark}\label{Remark424}  The proof of Lemma \ref{Lemma31} shows that if $f:X\rightarrow Y$ is pre-$\tau$-quasi-well prepared (or $\tau$-well prepared or
$\tau$-very-well prepared), $C\subset Y$ is a 2-curve, $\Psi:Y_1\rightarrow Y$ is the blow up of $C$ and $\Phi:X_1\rightarrow X$ is a sequence of blow ups of
2-curves and 3-points such that the rational map $f_1:X_1\rightarrow Y_1$ is a morphism, then
$$
\begin{array}{rll}
X_1&\stackrel{f_1}{\rightarrow}&Y_1\\
\Phi\downarrow&&\downarrow\Psi\\
X&\stackrel{f}{\rightarrow}&Y
\end{array}
$$
is pre-$\tau$-quasi-well prepared (or $\tau$-well prepared or $\tau$-very-well prepared) for $R$ and $\Psi$.  In fact, with the above notation, if $f$ satisfies 1 -- 3 of Definition \ref{DefT60}, then $f_1$ satisfies 1 -- 3 of Definition 
\ref{DefT60}. Further, 2 and 3 of the conclusions of Lemma \ref{Lemma31} hold.
\end{Remark}

\begin{Lemma}\label{LemmaT79}  Suppose that $f:X\rightarrow Y$ is pre-$\tau$-quasi-well prepared (or $\tau$-well prepared), $\overline q\in Y$ is a 1-point such that $\overline q\not\in U(R)$, and $C\subset D_Y$ is a resolving curve for $f$ at $\overline q$.
Then there exists a pre-$\tau$-quasi-well prepared ($\tau$-well prepared) diagram
$$
\begin{array}{rll}
X_1&\stackrel{f_1}{\rightarrow}&Y_1\\
\Phi_1\downarrow&&\downarrow\Psi_1\\
X&\stackrel{f}{\rightarrow}&Y
\end{array}
$$
such that 
\begin{enumerate}
\item[1.]
$\Psi_1$ is the blow up of $C$,
\item[2.] $\Phi_1$ is a sequence of blow ups of 2-curves. 
\item[3.]$\tau_{f_1}(p_1)\le \tau_f(\Phi_1(p_1))$ for $p_1\in X_1$. Thus $\tau_{f_1}(p_1)<\tau$ if $p_1\in (\Psi_1\circ f_1)^{-1}(C-\{\overline q\})$.
\item[4.] Suppose that $C_1\subset Y_1$ is a section over $C$,   $q\in C$ is a 1-point, $u,v,w$ are permissible parameters at $q$ such that $u=v=0$ are local equations of $C$, 
 $u,v$ are toroidal forms at $p$ for all $p\in f^{-1}(q)$,
and $q_1\in C_1\cap \Psi_1^{-1}(q)$ is a 1-point. Then there exist permissible parameters $\overline u_1,\overline v_1,w$ at $q_1$ such that $\overline u_1,\overline v_1$ are torodial forms at $p_1$ for all $p_1\in f_1^{-1}(q_1)$, and $\overline u_1=\overline v_1=0$ is a local equation of $C_1$.
\item[5.] Suppose that  $\overline p\in f^{-1}(\overline q)$, $p'\in \Phi_1^{-1}(\overline p)$, $u,v,w$ are permissible parameters at $\overline q$ such that $u=v=0$ are local equations of $C$, and $w$ is good at $\overline p$ for $f$. If $w$ is not good at $p'$ for $f_1$, then $\tau_{f_1}(p')<\tau_f(\overline p)$.
\end{enumerate}
\end{Lemma}

\begin{pf} 

Suppose that $q\in C$ is a 1-point. Then there exist permissible parameters $u,v,w$ at $q$ such that $u=v=0$ are local equations of $C$, and if $p\in f^{-1}(q)$, then $u,v$ are toroidal forms
at $p$. Thus there exist permissible parameters $x,y,z$ at $p$ for $u,v,w$ such that one of the forms (\ref{eqTF01}) or (\ref{eqTF02}) hold.

If $q\in C$ is a 2-point, then there exist super parameters $u,v,w$ at $q$ such that $u=w=0$ are local equations of $C$. Thus if $p\in f^{-1}(q)$, then there exist permissible parameters $x,y,z$ at $p$ for $u,v,w$ such that one of the forms (\ref{eqT221}) - (\ref{eqT224}) hold.

By Lemma 3.13 \cite{C5} and Remark \ref{RemarkT278},
after blowing up 2-curves and 3-points by a morphism $\Phi_0:X_0\rightarrow X$
such that  if $q\in C$ is a 2-point, $u,v,w$ are the above permissible parameters at $q$, and $p\in (f\circ\Phi_0)^{-1}(q)$, then
$(x^ay^b,x^ey^f,x^gy^h)$ is a principal ideal if (\ref{eqT222})  holds at $p$, $((x^ay^b)^k,(x^ay^b)^l,x^cy^d)$ is principal if (\ref{eqT223}) holds at $p$, $(x^ay^bz^c,x^gy^hz^i,z^jy^kz^l)$ is principal if (\ref{eqT224}) holds at $p$.
We further have that $f\circ\Phi_0$ is pre-$\tau$-quasi-well prepared ($\tau$-well prepared). In particular,
$$
\tau_{f\circ\Phi_0}(p)\le \tau_f(\Phi_0(p))
$$
for $p\in X_0$.

We now analyze the points $p\in  X_0$ where ${\cal I}_C{\cal O}_{X_0}$ is not principal.

First suppose that $q\in C$ is a 1-point, $p\in X_0$, and $f\circ\Phi_0(p)=q$.

If $p$ is a 1-point  then we have an expression
$$
u=x^a, v=y
$$
and $u=v=0$ are local equations of $C$. The non principal locus has local equations $x=y=0$

If $p$ is a 2-point  then 
$$
u=(x^ay^b)^k, v=z
$$
and $u=v=0$ are local equations of $C$. The non principal locus has local equations $\{x=z=0\}\cup\{y=z=0\}$

Now suppose that $q\in C$ is a 2-point and $f\circ\Phi_0(p)=q$. 

If $p$ is a 1-point of the form (\ref{eqT221}),  then
$$
\begin{array}{ll}
u&=x^a\\
v&=x^b(\alpha+y)\\
w&=x^c\gamma(x,y)+x^d(z+\beta).
\end{array}
$$
$u=w=0$ are local equations of $C$. If $p$ is in the non principal locus then we have 
(after possibly making a change of variables in $z$) $w=x^dz$ with $d<a$ and $x=z=0$ are local equations of the non principal locus.

If $p$ is a 2-point of the form (\ref{eqT222}), then 
$$
\begin{array}{ll}
u&=x^ay^b\\
v&=x^cy^d\\
w&=x^ey^f\gamma(x,y)+x^gy^h(z+\beta).
\end{array}
$$
$u=w=0$ are local equations of $C$.
If $p$ is in the non principal locus then (after possibly making a change of variables in $z$), we have $w=x^gy^hz$ with $(g,h)<(a,b)$.
Local equations of the non principal locus are
  $x=z=0$ (if $g<a$, $h=b$), $y=z=0$ (if $g=a$, $h<b$) and $\{x=z=0\}\cup\{y=z=0\}$ if $g<a$ and $h<b$.

If $p$ is a 2-point of the form (\ref{eqT223}), then
$$
\begin{array}{ll}
u&=(x^ay^b)^k\\
v&=(x^ay^b)^t(\alpha+z)\\
w&=(x^ay^b)^l\gamma(x^ay^b,z)+x^cy^d
\end{array}
$$
$u=w=0$ are local equations of $C$, and $p$ is in the principal locus.

If $p$ is a 3-point, of the form (\ref{eqT224}), then
$$
\begin{array}{ll}
u&=x^ay^bz^c\\
v&=x^dy^ez^f\\
w&=x^gy^hz^i\gamma+x^jy^kz^l
\end{array}
$$
$u=w=0$ are local equations of $C$, and $p$ is in the principal locus.

We see that the non principal locus of ${\cal I}_C{\cal O}_{X_0}$ is a union of nonsingular curves which are possible centers and are not 2-curves.

Let $U_0\subset X_0$ be the largest open set on which the rational map $\overline f_0=\Psi_1^{-1}\circ f\circ\Phi_0:X_0\rightarrow Y_1$ is a morphism. We will now show that $\overline f_0:U_0\rightarrow Y_1$ is $\tau$-prepared, and $\tau_{\overline f_0}(p)\le\tau_f(\Phi_0(p))$ for $p\in U_0$.

Suppose that $p\in U_0\cap (f\circ\Phi_0)^{-1}(C)$. Then $(f\circ\Phi_0)(p)=q$ is a 2-point, and there exist super parameters $u,v,w$ at $q$ such that $u=w=0$ are local equations of $C$, and one of the forms (\ref{eqT221}) - (\ref{eqT224}) hold for $u,v,w$ at $p$. Let $q_1=\overline f_0(p)$. We have permissible parameters $u_1,v,w_1$ in ${\cal O}_{Y_1,q_1}$ such that either
$$
u=u_1, w=u_1(w_1+\delta)
$$
with $\delta\in{\bf k}$, or
$$
u=u_1w_1, w=w_1.
$$

The most difficult case to analyze is when $p$ is a 3-point, (so that $u,v,w$ satisfy (\ref{eqT224})),  $\tau_{f\circ\Phi_0}(p)\ge 0$,
and $u=u_1w_1, w=w_1$. We will work out this case in detail.
Since $q$ is a 2-point, we will then have that
$$
\tau_{\overline f_0}(p)\le\tau_{f\circ\Phi_0}(p)\le\tau_f(\Phi_0(p))<\tau.
$$

In this case $q_1$ is a 3-point.
$$
u=x^ay^bz^c,
v=x^dy^ez^f,
w=\sum_{i\ge 0}\alpha_ix^{a_i}y^{b_i}z^{c_i}+x^gy^hz^i
$$
with $0\ne\alpha_i$ for all $i$. Set
$$
\gamma=\sum \alpha_ix^{a_i-a_0}y^{b_i-b_0}z^{c_i-c_0}+x^{g-a_0}y^{h-a_0}z^{i-a_0}.
$$
$\gamma\in\hat{\cal O}_{X,p}$ is a unit series, and
$$
w=x^{a_0}y^{b_0}z^{c_0}\gamma,
$$
with $(a_0,b_0,c_0)<(a,b,c)$.
\vskip .2truein
\noindent {\bf Suppose that $(a,b,c)$ and $(a_0,b_0,c_0)$ are linearly independent.}

After possibly interchanging $x,y,z$, we may assume that $a_0b-ab_0\ne 0$. There exist $\lambda_1,\lambda_2\in{\bf Q}$ such that if 
\begin{equation}\label{eqT249}
x=\overline x\gamma^{\lambda_1}, y=\overline y\gamma^{\lambda_2},
\end{equation}
then
$$
\begin{array}{ll}
w&=\overline x^{a_0}\overline y^{b_0}z^{c_0}\\
u&=\overline x^a\overline y^bz^c\\
v&=\overline x^d\overline y^ez^f\gamma^{\lambda}
\end{array}
$$
for some $0\ne\lambda\in{\bf Q}$.

There exists a series $P(x,y,z)$ where the monomials in $x,y,z$ in $P$ with nonzero coefficients are monomials in $x^{a_i-a_0}y^{b_i-b_0}z^{c_i-c_0}$ for $i\ge 0$, such that 
\begin{equation}\label{eqT245}
\gamma^{\lambda}=P
+x^{g-a_0}y^{h-b_0}z^{i-c_0}\Omega.
\end{equation}
Iterating, by substituting (\ref{eqT249}) into successive iterations of (\ref{eqT245}), we see that there exists a series $\overline P(\overline x,\overline y,\overline z)$, where the monomials in $\overline x,\overline y,\overline z$ in $\overline P$ with nonzero coefficients are monomials in $\overline x^{a_i-a_0}\overline y^{b_i-b_0}\overline z^{c_i-c_0}$ for $i\ge 0$ such that
$$
\gamma^{\lambda}
=\overline P
+\overline x^{g-a_0}\overline y^{h-b_0}z^{i-c_0}\overline \Omega.
$$

Comparing the Jacobian determinants
$$
\text{Det}\left(\begin{array}{l}\partial(u,v,w)\\ \partial(x,y,z)\end{array}\right)
$$
and
$$
\text{Det}\left(\begin{array}{l}\partial(w,u,v)\\ \partial(\overline x,\overline y,z)\end{array}\right),
$$
we see that $\overline\Omega$ is a unit series.

There exist rational numbers $\beta_1,\beta_2,\beta_3$ such that we can make a formal change of variables, setting
$$
\tilde x=\overline x\overline\Omega^{\beta_1},
\tilde y=\overline y\overline\Omega^{\beta_2},
\tilde z=z\overline\Omega^{\beta_3}
$$
to get an expression
$$
\begin{array}{ll}
w&=\tilde x^{a_0}\tilde y^{b_0}\tilde c^{c_0}\\
u&=\tilde x^a\tilde y^b\tilde z^c\\
v&=\tilde x^d\tilde y^e\tilde z^f\left(
\overline P(\tilde x,\tilde y,\tilde z)
+\tilde x^{g-a_0}\tilde y^{h-b_0}\tilde z^{i-c_0}\right).
\end{array}
$$

Thus
$$
\begin{array}{ll}
w_1&=\tilde x^{a_0}\tilde y^{b_0}\tilde c^{c_0}\\
u_1&=\tilde x^{a-a_0}\tilde y^{b-b_0}\tilde z^{c-c_0}\\
v&=\tilde x^d\tilde y^e\tilde z^f\left(
\overline P(\tilde x,\tilde y,\tilde z)
+\tilde x^{g-a_0}\tilde y^{h-b_0}\tilde z^{i-c_0}\right).
\end{array}
$$
is an expression of the form of (\ref{eq16}).

We see that $\overline f_0$ is prepared at $p$, and
$$
\begin{array}{ll}
\tau_{\overline f_0}(p)&\le \ell(((a_0,b_0,c_0){\bf Z}+(a-a_0,b-b_0,c-c_0){\bf Z}+(d,e,f){\bf Z}
+\sum_{i\ge1}(a_i-a_0,b_i-b_0,c_i-c_0){\bf Z})\\
&\,\,\, /((a_0,b_0,c_0)+(a-a_0,b-b_0,c-c_0){\bf Z}+(d,e,f){\bf Z}))\\
&=\ell(((a,b,c){\bf Z}+(d,e,f){\bf Z}
+\sum_{i\ge0}(a_i,b_i,c_i){\bf Z})\\
&\,\,\, /((a,b,c){\bf Z}+(d,e,f){\bf Z}+(a_0,b_0,c_0){\bf Z}))\\
&\le\ell(((a,b,c){\bf Z}+(d,e,f){\bf Z}
+\sum_{i\ge0}(a_i,b_i,c_i){\bf Z})\\
&\,\,\, /((a,b,c){\bf Z}+(d,e,f){\bf Z}))\\
&=\tau_{f\circ\Phi_0}(p)\le \tau_f(\Phi_0(p)).
\end{array}
$$
\vskip .2truein
\noindent {\bf Suppose that $(a,b,c)$ and $(a_0,b_0,c_0)$ are linearly dependent.}

Then $(d,e,f)$ and $(a_0,b_0,c_0)$ are linearly independent. 
After possibly interchanging $x,y,z$, we may assume that $db_0-ea_0\ne 0$. There exist $\lambda_1,\lambda_2\in{\bf Q}$ such that if 
$$
x=\overline x\gamma^{\lambda_1}, y=\overline y\gamma^{\lambda_2},
$$
then
$$
\begin{array}{ll}
w&=\overline x^{a_0}\overline y^{b_0}z^{c_0}\\
v&=\overline x^d\overline y^ez^f\\
u&=\overline x^a\overline y^bz^c\gamma^{\lambda}
\end{array}
$$
for some $0\ne\lambda\in{\bf Q}$.

As in the case when $(a,b,c)$ and $(a_0,b_0,c_0)$ are linearly independent, we can make a change of variables to get an expression where the monomials in $\tilde x,\tilde y,\tilde z$ in $\tilde P$ with nonzero coefficients are monomials in $\tilde x^{a_i-a_0}\tilde y^{\tilde b_i-b_0}\tilde z^{c_i-c_0}$ for $i\ge 0$.

$$
\begin{array}{ll}
w&=\tilde x^{a_0}\tilde y^{b_0}\tilde c^{c_0}\\
v&=\tilde x^d\tilde y^e\tilde z^f\\
u&=\tilde x^a\tilde y^b\tilde z^c\left(
\overline P(\tilde x,\tilde y,\tilde z)
+\tilde x^{g-a_0}\tilde y^{h-b_0}\tilde z^{i-c_0}\right).
\end{array}
$$

Thus
$$
\begin{array}{ll}
w_1&=\tilde x^{a_0}\tilde y^{b_0}\tilde c^{c_0}\\
v&=\tilde x^{d}\tilde y^{e}\tilde z^{f}\\
u_1&=\tilde x^{a-a_0}\tilde y^{b-b_0}\tilde z^{c-c_0}\left(
\overline P(\tilde x,\tilde y,\tilde z)
+\tilde x^{g-a_0}\tilde y^{h-b_0}\tilde z^{i-c_0}\right).
\end{array}
$$
is an expression of the form of (\ref{eq16}).

We see that $\overline f_0$ is prepared at $p$, and
$$
\begin{array}{ll}
\tau_{\overline f_0}(p)&\le \ell(((a_0,b_0,c_0){\bf Z}+(d,e,f){\bf Z}+(a-a_0,b-b_0,c-c_0){\bf Z}
+\sum_{i\ge1}(a_i-a_0,b_i-b_0,c_i-c_0){\bf Z})\\
&\,\,\, /((a_0,b_0,c_0)+(d,e,f){\bf Z}+(a-a_0,b-b_0,c-c_0){\bf Z}))\\
&\le \tau_{f\circ\Phi_0}(p)\le\tau_f(\Phi_0(p)).
\end{array}
$$

We conclude that $\overline f_0:U_0\rightarrow Y_1$ is $\tau$-prepared, and  $\tau_{\overline f_0}(p)\le \tau_f(\Phi_0(p))$ for $p\in U_0$.
Thus $\overline f_0$ is pre-$\tau$-quasi-well prepared ($\tau$-well prepared), as $C\cap U(R)=\emptyset$.

Let $Z_0=X_0$.

We now construct a sequence of morphisms $\Lambda_i:Z_i\rightarrow Z_{i-1}$ which are the blow up of a possible curve $C_{i-1}$ contained in the locus where ${\cal I}_{C}{\cal O}_{Z_{i-1}}$ is not invertible.

Let  $h_i:Z_i\rightarrow Y_1$ be the rational map $h_i=\overline f_0\circ \Lambda_1\circ\cdots\circ \Lambda_i$. We will verify by induction that:
\begin{enumerate}
\item[A.] For all $p_1$ in the locus where $h_i$ is a morphism,
$h_i$ is $\tau$-prepared and
$$
\tau_{h_i}(p_1)\le \tau_{f}(\Phi_0\circ\Lambda_1\circ\Lambda_2\circ\cdots\circ\Lambda_i(p_1)).
$$
\item[B.]
Suppose that
$p_1\in Z_i$ is a point where ${\cal I}_C{\cal O}_{Z_i,p_1}$ is not principal.
Let $p=\Lambda_1\circ\cdots\circ\Lambda_i(p_1)$, $q=f\circ\Phi_0(p)$. Then there exist permissible parameters $u,v,w$ at $q$, permissible parameters $x,y,z$ in $\hat{\cal O}_{Z_0,p}$ for $u,v$ and regular parameters $x_1,y_1,z_1$ in $\hat{\cal O}_{Z_i,p_1}$ such that we have one of the following forms:

\begin{enumerate}
\item[1.] $q$ a 1-point, $p$ a 1-point, $p_1$ a 1-point. We have an expression
$$
u=x^a, v=z, w=\sum_{i<n, a_{ij}\ne 0}a_{ij}x^iz^j+x^n(y+\beta)
$$ 
in $\hat{\cal O}_{Z_0,p}$, where $u=v=0$ are local equations of $C$.

We have permissible parameters $x_1,y_1,z_1$ at $p_1$ such that 
$$
x=x_1, y=y_1, z=x_1^bz_1
$$
with $b<a$, and

\begin{equation}\label{eqT225}
\begin{array}{ll}
u&=x_1^a\\
v&=x_1^bz_1\\
w&=\sum a_{ij}x_1^{i+bj}z_1^j+x_1^n(y_1+\beta).
\end{array}
\end{equation}
with $b<a$, where $u=v=0$ are local equations of $C$.

\item[2.]
$q$ a 1-point, $p$ a 2-point. 
We have an expression
$$
u=(x^ay^b)^k, v=z,
w=\sum_{i,l\ge 0}a_{il}z^l(x^ay^b)^i+x^ey^f
$$
in $\hat{\cal O}_{Z_0,p}$,
where the sum is over $i, l$ such that $(ia,ib)\not\ge (e,f)$, $af-eb\ne 0$ and $u=v=0$ are local equations of $C$.

We have permissible parameters $x_1,y_1,z_1$ at $p_1$ such that 
$$
x=x_1, y=y_1, z=x_1^cy_1^dz_1
$$
with $(c,d)<(ak,bk)$, and

\begin{equation}\label{eqT180}
\begin{array}{ll}
u&=(x_1^ay_1^b)^k\\
v&=x_1^cy_1^dz_1\\
w&=\sum_{i,l}a_{il}z_1^lx_1^{ai+cl}y_1^{bi+dl}+x_1^ey_1^f.
\end{array}
\end{equation}
with $(c,d)<(ak,bk)$, $u=v=0$ are local equations of $C$.

\item[3.]
$q$ a 2-point, $p$ a 1-point, $p_1$ a 1-point. We have an expression 
\begin{equation}\label{eqT181}
\begin{array}{ll}
u&=x_1^a\\
v&=x_1^b(\beta+y_1)\\
w&=x_1^dz_1
\end{array}
\end{equation}
in $\hat{\cal O}_{Z_i,p_i}$, with $d<a$, $u=w=0$ a local equation of $C$

\item[4.]
$q$ a 2-point, $p$ a 2-point, $p_1$ a 2-point. We have an expression 
\begin{equation}\label{eqT182}
\begin{array}{ll}
u&=x_1^ay_1^b\\
v&=x_1^cy_1^d\\
w&=x_1^gy_1^hz_1
\end{array}
\end{equation}
with $(g,h)<(a,b)$, $u=w=0$ are local equation of $C$.
\end{enumerate}
\end{enumerate}

We have verified that A and B are true for $h_0=\overline f_0$. Suppose that $A$ and $B$ are true for $h_i$. We will verify that A and B are true for $h_{i+1}$.

We may suppose that $p_1\in C_i$ (recall that $C_i$ is the curve blown up by  $\Lambda_{i+1}$). Then $C_i$ has local equations $x_1=z_1=0$ if $p_1$ satisfies (\ref{eqT225}).
$C_i$ has local equations $x_1=z_1=0$ (or $y_1=z_1=0$) in (\ref{eqT180}). $C_i$ has local equations $x_1=z_1=0$ in (\ref{eqT181}). $C_i$ has local equations $x_1=z_1=0$ (or $y_1=z_1=0$) in (\ref{eqT182}).

After possibly interchanging $x_1$ and $y_1$, we may assume that $x_1=z_1=0$ is a local equation of $C_i$.

Suppose that $p_2\in\Lambda_{i+1}^{-1}(p_1)$. Then 
$\hat{\cal O}_{Z_{i+1},p_2}$ has regular parameters $x_2,y_2,z_2$ such that 

\begin{equation}\label{eqT183}
x_1=x_2,
y_1=y_2,
z_1=x_2(z_2+\alpha)
\end{equation}
with $\alpha\in{\bf k}$,
or 
\begin{equation}\label{eqT184}
x_1=x_2z_2, 
y_1=y_2,
z_1=z_2.
\end{equation}

Suppose that (\ref{eqT225}) holds. Then
$$
\tau_{f\circ\Phi_0}(p)=\ell\left([\sum_{a_{ij}\ne 0}i{\bf Z}+a{\bf Z}]/a{\bf Z}\right).
$$

Suppose that (\ref{eqT225}) and (\ref{eqT183}) hold. Then

\begin{equation}\label{eqT227}
\begin{array}{ll}
u&=x_2^a\\
v&=x_2^{b+1}(z_2+\alpha)\\
w&=\sum a_{ij}x_2^{i+(b+1)j}(z_2+\alpha)^j+x_2^n(y_2+\beta).
\end{array}
\end{equation}

Suppose that $a=b+1$ in (\ref{eqT227}). Then $Z_{i+1}\rightarrow Y_1$ is a morphism at $p_2$.
Let $q_1=h_{i+1}(p_2)$. There exist permissible parameters $u_1,v_1,w$ at $q_1$ defined by
$$
u=u_1,
v=u_1(v_1+\alpha).
$$
We have an expression 
\begin{equation}\label{eqT231}
\begin{array}{ll}
u_1&=x_2^a\\
v_1&=z_2\\
w&=\sum a_{ij}x_2^{i+aj}(z_2+\alpha)^j+x_2^n(y_2+\beta),
\end{array}
\end{equation}
of type (\ref{eqTF01})  so that $h_{i+1}$ is prepared at $p_2$.
We have by (\ref{eqT17}) that
$$
\tau_{h_{i+1}}(p_2)\le\ell((a{\bf Z}+\sum_{a_{ij}\ne 0}[i+a]{\bf Z})/a{\bf Z})=\tau_{ f\circ\Phi_0}(p)\le\tau_f(\Phi_0(p)).
$$

Suppose that $b+1<a$ and $\alpha\ne 0$ in (\ref{eqT227}). Then $Z_{i+1}\rightarrow Y_1$ is a morphism at $p_2$. Let $q_1=h_{i+1}(p_2)$.
There exist permissible parameters $u_1,v_1,w$ at $q_1$ defined by

$$
u=u_1v_1,v=v_1.
$$
There exist regular parameters $\overline x_2,\overline y_2,\overline z_2$ in $\hat{\cal O}_{Z_{i+1},p_2}$ such that 

\begin{equation}\label{eqT232}
\begin{array}{ll}
u_1&=x_2^{a-b-1}(z_2+\alpha)^{-1}=\overline x_2^{a-b-1}\\
v_1&=x_2^{b+1}(z_2+\alpha)=\overline x_2^{b+1}(\overline z_2+\overline\alpha)\\
w&=\sum_sx_2^s(\sum_{i+(b+1)j=s}a_{ij}(z_2+\alpha)^j)+x_2^n(y_2+\beta)\\
&=\sum_s \overline x_2^s(\overline z_2+\alpha)^{\frac{s}{a}}\left[\sum_ja_{ij}(\overline z_2+\alpha)^{j(\frac{a-b-1}{a})}\right]+\overline x_2^n(\overline y_2+\overline\beta).
\end{array}
\end{equation}

of type (\ref{eqTF1}), so that $h_{i+1}$ is prepared at $p_2$.

We observe that 
$$
\sum_{i+(b+1)j=s}a_{ij}(z_2+\alpha)^j\ne 0
$$ 
if and only if some $a_{ij}\ne 0$ with $i+(b+1)j=s$.  To show this,  observe that
$$
x^s[\sum_{i+(b+1)j=s}a_{ij}\left(\frac{y}{x^{b+1}}\right)^j]=\sum_{i+(b+1)j=s}a_{ij}x^iy^j.
$$
We calculate  from (\ref{eqT15}),

$$
\begin{array}{ll}
\tau_{h_{i+1}}(p_2)&=\ell([(a-b-1){\bf Z}+(b+1){\bf Z}\\
&\,\,\,+\sum_{a_{ij}\ne 0, i+(b+1)j\not\ge n}(i+(b+1)j){\bf Z}]/[(a-b-1){\bf Z}+(b+1){\bf Z}])\\
&\le\ell\left([\sum_{a_{ij}\ne 0}i{\bf Z}+a{\bf Z}+(b+1){\bf Z}]/[a{\bf Z}+(b+1){\bf Z}]\right)\\
&\le\tau_{f\circ\Phi_0}(p)\le\tau_f(\Phi_0(p))
\end{array}
$$

Suppose that $b+1<a$, $0=\alpha$ in (\ref{eqT227}). Then we are back in the form (\ref{eqT225}) with a decrease in $a-b$.

Suppose that (\ref{eqT225}) and (\ref{eqT184}) hold. Then

$$
\begin{array}{ll}
u&=x_2^az_2^a\\
 v&=x_2^bz_2^{b+1}\\
w&=\sum a_{ij}x_2^{i+bj}z_2^{i+(b+1)j}+x_2^nz_2^n(y_2+\beta).
\end{array}
$$

$Z_{i+1}\rightarrow Y_1$ is a morphism at $p_2$. Let $q_1=h_{i+1}(p_2)$. 
There exist permissible parameters $u_1,v_1,w$ at $q_1$ defined by

$$
u=u_1v_1,
v=v_1.
$$
We have an expression 
\begin{equation}\label{eqT233}
\begin{array}{ll}
u_1&=x_2^{a-b}z_2^{a-b-1}\\
v _1&=x_2^bz_2^{b+1}\\
w&=\sum a_{ij}x_2^{i+bj}z_2^{i+(b+1)j}+x_2^nz_2^n(y_2+\beta)
\end{array}
\end{equation}
of type (\ref{eqTF21}), so that $h_{i+1}$ is prepared at $p_2$.

From (\ref{eqT10}), we have
$$
\begin{array}{ll}
\tau_{h_{i+1}}(p_2)&=\ell(((a-b,a-b-1){\bf Z}+(b,b+1){\bf Z}\\
&\,\,\,+\sum_{a_{ij}\ne0,
(i+bj,i+(b+1)j)\not\ge(n,n)}(i+bj,i+(b+1)j){\bf Z})\\
&\,\,\, /((a-b,a-b-1){\bf Z}+(b,b+1){\bf Z}))\\
&\le\ell(((a,a){\bf Z}+\sum_{a_{ij}\ne0}(i,i){\bf Z})/(a,a){\bf Z})\\
&=\tau_{f\circ\Phi_0}(p)\le\tau_f(\Phi_0(p)).
\end{array}
$$

Suppose that (\ref{eqT180}) holds. Then
$$
\tau_{f\circ\Phi_0}(p)=\ell\left[(k{\bf Z}+\sum_{a_{il}\ne 0}i{\bf Z})/k{\bf Z}\right].
$$

Suppose that  (\ref{eqT180}) and (\ref{eqT183}) hold. Then 
we have

\begin{equation}\label{eqT226}
\begin{array}{ll}
u&=(x_2^ay_2^b)^k\\
v&=x_2^{c+1}y_2^d(z_2+\alpha)\\
w&=\sum_{i,l}a_{il}(z_2+\alpha)^lx_2^{ai+(c+1)l}y_2^{bi+dl}+x_2^ey_2^f.
\end{array}
\end{equation}
Suppose that $(ak,bk)=(c+1,d)$ in (\ref{eqT226}). Then $Z_{i+1}\rightarrow Y_1$ is a morphism at $p_2$.
Let $q_1=h_{i+1}(p_2)$.
There exist permissible parameters $u_1,v_1,w$ at $q_1$ defined by

$$
\begin{array}{ll}
u&=u_1\\
v&=u_1(v_1+\alpha).
\end{array}
$$
We have an expression 
\begin{equation}\label{eqT237}
\begin{array}{ll}
u_1&=(x_2^ay_2^b)^k\\
v_1&=z_2\\
w&=\sum_{i,l}a_{il}(z_2+\alpha)^lx_2^{ai+(c+1)l}y_2^{bi+dl}+x_2^ey_2^f\\
&=\sum_{i,l}a_{il}(z_2+\alpha)^l(x_2^ay_2^b)^{i+lk}+x_2^ey_2^f
\end{array}
\end{equation}
of type (\ref{eqTF02}), so that $h_{i+1}$ is prepared at $p_2$. From  (\ref{eqT13}),
we have
$$
\tau_{h_{i+1}}(p_2)\le
\ell((k{\bf Z}+\sum_{a_{il}\ne0, (i+lk)(a,b)\not\ge(e,f)}(i+lk){\bf Z}/k{\bf Z})\le\tau_{ f\circ\Phi_0}(p)\le\tau_f(\Phi_0(p)).
$$

Suppose that $(ak,bk)>(c+1,d)$ and $0\ne\alpha$ in (\ref{eqT226}). Then $Z_{i+1}\rightarrow Y_1$ is a morphism at $p_2$. Let $q_1=h_{i+1}(p_2)$.
There exist permissible parameters $u_1,v_1,w$ at $q_1$ defined by

$$
u=u_1v_1, v=v_1.
$$

Suppose that $ad-b(c+1)\ne 0$ in (\ref{eqT226}). 
Then there exist regular parameters $\overline x_2,\overline y_2,\overline z_2$ in $\hat{\cal O}_{Z_{i+1},p_2}$
defined by
$$
\begin{array}{ll}
\overline x_2&=x_2(z_2+\alpha)^{-\frac{bk}{\overline h}}\\
\overline y_2&=y_2(z_2+\alpha)^{\frac{ak}{\overline h}}\\
\overline z_2&=(z_2+\alpha)^{\frac{ebk-fak}{\overline h}}-\overline\alpha
\end{array}
$$
with $\overline h=adk-(c+1)bk$ and $\overline\alpha=\alpha^{\frac{ebk-fak}{\overline h}}$,
such that 
\begin{equation}\label{eqT234}
\begin{array}{ll}
u_1&=x_2^{ak-c-1}y_2^{bk-d}(z_2+\alpha)^{-1}=\overline x_2^{ak-c-1}\overline y_2^{bk-d}\\
v_1&=x_2^{c+1}y_2^d(z_2+\alpha)=\overline x_2^{c+1}\overline y_2^d\\
w&=\sum_{i,l}a_{il}\overline x_2^{ai+(c+1)l}\overline y_2^{bi+dl}+\overline x_2^e\overline y_2^f(\overline z_2+\overline \alpha)
\end{array}
\end{equation}
with $0\ne\overline\alpha\in{\bf k}$,
which is an expression of type (\ref{eqTF21}), so that $h_{i+1}$ is prepared at $p_2$.
From (\ref{eqT10}), we see that

\begin{equation}\label{eqT235}
\begin{array}{ll}
\tau_{h_{i+1}}(p_2)&=\ell(
(ak-c-1,bk-d){\bf Z}+(c+1,d){\bf Z}\\
&\,\,\,+\sum_{a_{il}\ne 0, i(a,b)+l(c+1,d)\not\ge(e,f)}[i(a,b)+l(c+1,d)]{\bf Z}]/\\
&\,\,\,[(ak-c-1,bk-d){\bf Z}+(c+1,d){\bf Z})\\
&\le \ell\left[(k(a,b){\bf Z}+\sum_{a_{il}\ne 0}i(a,b){\bf Z})/k(a,b){\bf Z}\right]\\
&=\ell\left[k{\bf Z}+\sum_{a_{il}\ne 0}i{\bf Z}/k{\bf Z}\right]=\tau_{f\circ\Phi_0}(p)\le\tau_f(\Phi_0(p)).
\end{array}
\end{equation}

Suppose that $ad-b(c+1)=0$ in (\ref{eqT226}), (and $(ak,bk)>(c+1,d)$, $0\ne\alpha$ in (\ref{eqT226})).
There exist positive integers $\overline a,\overline b,\overline t,\overline k$ such that $\text{gcd}(\overline a,\overline b)=1$, and 
$$
(ak-c-1,bk-d)=\overline t(\overline a,\overline b),\,\,\,
(c+1,d)=\overline k(\overline a,\overline b).
$$
From $k(a,b)=(\overline t+\overline k)(\overline a,\overline b)$ and $\text{gcd}(a,b)=1$, $\text{gcd}(\overline a,\overline b)=1$, we see that $(a,b)=(\overline a,\overline b)$ and $\overline t+\overline k=k$.

 There exist regular parameters $\overline x_2,\overline y_2,\overline z_2$ in $\hat{\cal O}_{Z_{i+1},p_2}$, defined by
$$
\begin{array}{ll}
\overline x_2&=x_2(z_2+\alpha)^{-\frac{f}{\overline h\overline t}}\\
\overline y_2&=y_2(z_2+\alpha)^{\frac{e}{\overline h\overline t}}\\
\overline z_2&=(z_2+\alpha)^{1+\frac{\overline k}{\overline t}}-\overline\alpha
\end{array}
$$
where $\overline h=f a-e b$, $\overline\alpha=e^{1+\frac{\overline k}{\overline t}}$, such that 
\begin{equation}\label{eqT238}
\begin{array}{ll}
u_1&=x_2^{ak-c-1}y_2^{bk-d}(z_2+\alpha)^{-1}=(\overline x_2^{a}\overline y_2^{ b})^{\overline t}\\
v_1&=x_2^{c+1}y_2^d(z_2+\alpha)=(\overline x_2^{a}\overline y_2^{b})^{\overline k}(\overline z_2+\overline\alpha)\\
w&=\sum_{i,l}a_{il}(\overline z_2+\overline\alpha)^{(\frac{\overline t}{k}+1)l+\frac{i}{k}}(\overline x_2^a\overline y_2^b)^{i+kl}+\overline x_2^e\overline y_2^f,
\end{array}
\end{equation}
which is an expression of type (\ref{eqTF22}), so that $h_{i+1}$ is prepared at $p_2$.

From (\ref{eqT11}), we have

\begin{equation}\label{eqT239}
\begin{array}{ll}
\tau_{h_{i+1}}(p_2)&= \ell((\overline t{\bf Z}+\overline k{\bf Z}+\sum_{a_{il}\ne0,
(i+kl)(a,b)\not\ge(e,f)}(lk+i){\bf Z})
/(\overline t{\bf Z}+\overline k{\bf Z}))\\
&\le\ell((k{\bf Z}+\overline k{\bf Z}+\sum_{a_{il}\ne0}i{\bf Z})
/(k{\bf Z}+\overline k{\bf Z}))\\
&\le\ell((k{\bf Z}+\sum_{a_{il}\ne0}i{\bf Z})
/k{\bf Z})\\
&=\tau_{f\circ\Phi_0}(p)\le\tau_f(\Phi_0(p))
\end{array}
\end{equation}

Suppose that $(ak,bk)>(c+1,d)$ and $0=\alpha$ in (\ref{eqT226}). Then we are back in the form (\ref{eqT180}), with a decrease in $(ak-c)+(bk-d)$.

Suppose that (\ref{eqT180}) and (\ref{eqT184}) hold. Then 
\begin{equation}\label{eqT236}
\begin{array}{ll}
u&=x_2^{ak}y_2^{bk}z_2^{ak}\\
v&=x_2^{c}y_2^dz_2^{c+1}\\
w&=\sum_{i,l}a_{il}x_2^{ai+cl}y_2^{bi+dl}z_2^{l+ai+cl}+x_2^ey_2^fz_2^e.
\end{array}
\end{equation}

Then $Z_{i+1}\rightarrow Y_1$ is a morphism at $p_2$. Let $q_1=h_{i+1}(p_2)$.
There exist permissible parameters $u_1,v_1,w$ at $q_1$ defined by

$$
u=u_1v_1,
v=v_1.
$$
We have an expression 
\begin{equation}\label{eqT250}
\begin{array}{ll}
u_1&=x_2^{ak-c}y_2^{bk-d}z_2^{ak-c-1}\\
v_1&=x_2^{c}y_2^dz_2^{c+1}\\
w&=\sum_{i,l}a_{il}x_2^{ai+cl}y_2^{bi+dl}z_2^{l+ai+cl}+x_2^ey_2^fz_2^e.
\end{array}
\end{equation}
of type (\ref{eqTF3}), so that $h_{i+1}$ is prepared at $p_2$. We calculate $\tau_{h_{i+1}}(p_2)$ from (\ref{eq16}).

Set
$$
\begin{array}{ll}
H&=
((ak,bk,ak){\bf Z}+(c,d,c+1){\bf Z}\\
&+\sum_{a_{il}\ne 0}(ai+cl,bi+dl,l+ai+cl){\bf Z})/\\
&((ak,bk,ak){\bf Z}+(c,d,c+1){\bf Z})\\
&\cong ((ak,bk,ak){\bf Z}+\sum_{a_{il}\ne0}(ai,bi,ai){\bf Z})/(ak,bk,ak){\bf Z}).
\end{array}
$$
$$
\begin{array}{ll}
\tau_{h_{i+1}}(p_2)&=\ell(
(ak-c,bk-d,ak-c-1){\bf Z}+(c,d,c+1){\bf Z}\\
&+\sum_{a_{il}\ne 0,(ai+cl,bi+dl,l+ai+cl)\not>(e,f,e)}(ai+cl,bi+dl,l+ai+cl){\bf Z})/\\
&((ak-c,bk-d,ak-c-1){\bf Z}+(c,d,c+1){\bf Z}))\\
&\le \ell\left(H\right).
\end{array}
$$

The homomorphism ${\bf Z}\rightarrow {\bf Z}^3$ defined by $x\mapsto (xa,xb,xa)$ induces an isomorphism  
\begin{equation}\label{eqT240}
(k{\bf Z}+\sum_{a_{il}\ne 0}i{\bf Z})/k{\bf Z}\rightarrow H.
\end{equation}
Thus $\tau_{h_{i+1}}(p_2)\le\tau_{f\circ\Phi_0}(p)\le\tau_f(\Phi_0(p))$.

Suppose that (\ref{eqT181}) and (\ref{eqT183}) hold. Then 
\begin{equation}\label{eqT228}
\begin{array}{ll}
u&=x_2^a\\
v&=x_2^b(\beta+y_2)\\
w&=x_2^{d+1}(z_2+\alpha).
\end{array}
\end{equation}
Suppose that $d+1=a$ in (\ref{eqT228}). Then $Z_{i+1}\rightarrow Y_1$ is a morphism at $p_2$.
Let $q_1=h_{i+1}(p_2)$.
There exist permissible parameters $u_1,v,w_1$ at $q_1$ defined by

$$
u=u_1,
w=u_1(w_1+\alpha).
$$
We have an expression
$$
u_1=x_2^a,
v=x_2^b(\beta+y_2),
w_1=z_2
$$
which is  toroidal, so that $h_{i+1}$ is prepared at $p_2$, and
$$
\tau_{h_{i+1}}(p_2)=-\infty\le\tau_{f}(\Phi_0(p)).
$$

Suppose that $d+1<a$ and $0\ne\alpha$ in (\ref{eqT228}). Then $Z_{i+1}\rightarrow Y_1$ is a morphism at $p_2$.
Let $q_1=h_{i+1}(p_1)$.
There exist permissible parameters $u_1,v,w_1$ at $q_1$ defined by

$$
u=u_1w_1,
w=w_1,
$$
so that $q_1$ is a 3-point.

There exist regular parameters $\overline x_2,\overline y_2,\overline z_2$ in $\hat{\cal O}_{Z_{i+1},p_2}$ and $0\ne\overline\alpha,\overline\beta \in{\bf k}$ such that
$$
\begin{array}{ll}
u_1&=x_2^{a-d-1}(z_2+\alpha)^{-1}=\overline x_2^{a-d-1}\\
v&=\overline x_2^b(\overline \beta+\overline y_2)\\
w_1&=\overline x_2^{d+1}(\overline z_2+\overline \alpha)
\end{array}
$$
which is toroidal, so that $h_{i+1}$ is prepared at $p_2$, and
$$
\tau_{h_{i+1}}(p_2)=-\infty\le\tau_{f}(\Phi_0(p)).
$$

Suppose that $d+1<a$ and $0=\alpha$ in (\ref{eqT228}). Then we are back in the form (\ref{eqT181}) with a decrease in $a-d$.

Suppose that (\ref{eqT181}) and (\ref{eqT184}) hold.  Then
$$
\begin{array}{ll}
u&=x_2^az_2^a\\
v&=x_2^bz_2^b(\alpha+y_2)\\
w&=x_2^dz_2^{d+1}.
\end{array}
$$
Thus $Z_{i+1}\rightarrow Y_1$ is a morphism at $p_2$.
Let $q_1=h_{i+1}(p_2)$.
There exist permissible parameters $u_1,v,w_1$ at $q_1$ defined by

$$
u=u_1w_1,w=w_1.
$$
We have a 2-point mapping to a 3-point, and an expression
$$
\begin{array}{ll}
u_1&=x_2^{a-d}z_2^{a-d-1}\\
v&=x_2^bz_2^b(\alpha+y_2)\\
w_1&=x_2^dz_2^{d+1}
\end{array}
$$
which is a toroidal form, so that $h_{i+1}$ is prepared at $p_2$, and
$$
\tau_{h_{i+1}}(p_2)=-\infty\le\tau_{f}(\Phi_0(p)).
$$

Suppose that (\ref{eqT182}) and (\ref{eqT183}) hold. Then 
\begin{equation}\label{eqT229}
\begin{array}{ll}
u&=x_2^ay_2^b\\
v&=x_2^cy_2^d\\
w&=x_2^{g+1}y_2^h(z_2+\alpha).
\end{array}
\end{equation}
Suppose that  $(a,b)=(g+1,h)$ in (\ref{eqT229}). Then $Z_{i+1}\rightarrow Y_1$ is a morphism at $p_2$.
Let $q_1=h_{i+1}(p_2)$.
There exist permissible parameters $u_1,v,w_1$ at $q_1$ defined by

$$
u=u_1,
w=u_1(w_1+\alpha).
$$
We have an expression
$$
\begin{array}{ll}
u_1&=x_2^ay_2^b\\
v&=x_2^cy_2^d\\
w_1&=z_2
\end{array}
$$
which is toroidal, so that $h_{i+1}$ is prepared at $p_2$, and
$$
\tau_{h_{i+1}}(p_2)=-\infty\le\tau_{f}(\Phi_0(p)).
$$

Suppose that $(g+1,h)<(a,b)$ and $0\ne\alpha$ in (\ref{eqT229}). Then $Z_{i+1}\rightarrow Y_1$ is a morphism at $p_2$.
Let $q_1=h_{i+1}(p_2)$.
There exist permissible parameters $u_1,v,w_1$ at $q_1$ defined by

$$
u=u_1w_1, w=w_1.
$$
We have an expression
$$
\begin{array}{ll}
u_1&=x_2^{a-g-1}y_2^{b-h}(z_2+\alpha)^{-1}\\
v&=x_2^cy_2^d\\
w_1&=x_2^{g+1}y_2^h(z_2+\alpha)
\end{array}
$$
which is toroidal after a change of variable in $\hat{\cal O}_{Z_{i+1},p_2}$, so that $h_{i+1}$ is prepared at $p_2$, and
$$
\tau_{h_{i+1}}(p_2)=-\infty\le\tau_{f}(\Phi_0(p)).
$$

Suppose that $(g+1,h)<(a,b)$ and $\alpha=0$ in (\ref{eqT229}). Then we are back in the form (\ref{eqT182}), with a decrease in $(a-g)+(b-h)$.

Suppose that (\ref{eqT182}) and (\ref{eqT184}) hold. Then
$$
\begin{array}{ll}
u&=x_2^ay_2^bz_2^a\\
v&=x_2^cy_2^dz_2^c\\
w&=x_2^gy_2^hz_2^{g+1}.
\end{array}
$$
Thus $Z_{i+1}\rightarrow Y_1$ is a morphism at $p_2$.
Let $q_1=h_{i+1}(p_2)$.
There exist permissible parameters $u_1,v,w_1$ at $q_1$ defined by

$$
u=u_1w_1,
w=w_1.
$$
We have a  3-point mapping   to a 3-point, and an expression
$$
\begin{array}{ll}
u_1&=x_2^{a-g}y_2^{b-h}z_2^{a-g-1}\\
v&=x_2^cy_2^dz_2^c\\
w_1&=x_2^gy_2^hz_2^{g+1}
\end{array}
$$
which is toroidal, so that $h_{i+1}$ is prepared at $p_2$, and
$$
\tau_{h_{i+1}}(p_2)=-\infty\le\tau_{f}(\Phi_0(p)).
$$
We have thus verified A and B for all $h_i$.

Suppose that $p\in Z_i$ is in the locus $\Sigma_i$ where ${\cal I}_C{\cal O}_{Z_i}$ is not locally principal. Define
$$
C(p)=\left\{\begin{array}{ll}
a-b&\text{ if (\ref{eqT225}) holds at $p$}\\
ak-c+bk-d&\text{ if (\ref{eqT180}) holds at $p$}\\
a-d&\text{ if (\ref{eqT181}) holds at $p$}\\
a-g+b-h&\text{ if (\ref{eqT182}) holds at $p$}.
\end{array}\right.
$$
Let
$$
C(Z_i)=\text{max}\{C(p)\mid p\in \Sigma_i\}.
$$
We have shown in the above analysis that 
$$
0\le C(Z_{i+1})< C(Z_i)
$$
if the rational map $h_i$ is not a morphism. Thus there exists a finite $n$ such that 
$f_1=h_{n}$, $\Phi_1=\Phi_0\circ\Lambda_1\circ\cdots\circ\Lambda_n$, 
$X_1=Z_n$ satisfy 1 - 4 of the conclusions of Lemma \ref{LemmaT79}.

We now prove 5. If $p\in f^{-1}(q)$, then we have permissible parameters  $x,y,z$ for $u,v,w$ in $\hat{\cal O}_{X,p}$ such that there is  one of the forms:

$p$  a 1-point 
\begin{equation}\label{eqT96}
u=x^a, v=y, w=g(x,y)+x^n(z+\beta)
\end{equation}

or $p$ a 2-point 
\begin{equation}\label{eqT97}
u=(x^ay^b)^k, v=z, w=h(x^ay^b,z)+x^cy^d.
\end{equation}
We further have that $u=v=0$ are local equations of $C$.

Since $q_1$ is a 1-point, at $q_1$ we have permissible parameters $u_1,v_1,w_1$ defined by
$$
u=u_1, v=u_1(v_1+\alpha), w=w_1
$$
for some $\alpha\in {\bf k}$.
Since $C_1$ is a section over $C$, there exists a series $\lambda(v_1,w)$ such that $u_1=\lambda(v_1,w_1)=0$ are local equations of $C_1$ (in $\hat{\cal O}_{Y_1,q_1}$), and
$$
{\bf k}[[u,v,w]]/(u,v)\rightarrow {\bf k}[[u_1,v_1,w]]/(u_1,\lambda)
$$
is an isomorphism, which implies that $\lambda=(v_1-\phi(w))\gamma$ where $\phi$ is a series and $\gamma$ is a unit series in $\hat{\cal O}_{Y_1,q_1}$. Set
$$
\overline u_1=u_1, \overline v_1=v_1-\phi(w).
$$
$\overline u_1,\overline v_1,w$ are permissible parameters at $q_1$, and $\overline u_1=\overline v_1=0$ are local equations of $C_1$ at $q_1$.

Suppose that $p_1\in f_1^{-1}(q_1)$. First suppose that $p=\Phi_1(p_1)$ has the form (\ref{eqT96}). Then $\hat{\cal O}_{X_1,p_1}$ has regular parameters $x_1,y_1,z_1$ with $x=x_1, y=x_1^a(y_1+\alpha)$, since $X_0\rightarrow X$ is an isomorphism over $p$.
$p_1$ is a 1-point, and substituting into (\ref{eqT96}), we have 
$$
u_1=x_1^a, v_1=y_1, w=g(x_1,x_1^a(y_1+\alpha))+x_1^n(z+\beta)
$$
and thus $\overline u_1,\overline v_1$ are toroidal forms at $p_1$.

Now suppose that $p=\Phi_1(p_1)$ has the form (\ref{eqT97}). 

Let 
$$
p'=\Lambda_1\circ\cdots\circ\Lambda_n(p)\in X_0.
$$
There exist regular parameters $x_1,y_1,z_1$ in $\hat{\cal O}_{X_1,p'}$ such that 
\begin{equation}\label{eqT241}
x=x_1^{\overline a}y_1^{\overline b}, y=x_1^{\overline c}y_1^{\overline d}, z=z_1
\end{equation}
or 
\begin{equation}\label{eqT242}
x=x_1^{\overline a}(y_1+\alpha)^{\overline b},
y=x_1^{\overline c}(y_1+\overline\alpha)^{\overline d}, z=z_1
\end{equation}
with $0\ne\alpha\in{\bf k}$.
If (\ref{eqT241}) holds, then we have 
\begin{equation}\label{eqT243}
\begin{array}{ll}
u&=(x_1^{a'}y_1^{b'})^k\\
v&=z_1\\
w&=h(x_1^{a'}y_1^{b'},z)+x_1^{c'}y_1^{d'}
\end{array}
\end{equation}
with $0\ne\alpha\in{\bf k}$.

If (\ref{eqT242}) holds, then we have $0\ne\alpha'\in{\bf k}$ and permissible parameters $\overline x_1,\overline y_1,z_1$ in $\hat{\cal O}_{X_1,p'}$ such that 
\begin{equation}\label{eqT244}
\begin{array}{ll}
u&=\overline x_1^{a' k}\\
v&=z_1\\
w&=h(\overline x_1^{a'},z_1)+\overline x_1^{b'}(\overline y_1+\alpha').
\end{array}
\end{equation}

Suppose that a form (\ref{eqT243}) holds at $p'$.
Then $\hat{\cal O}_{X_1,p_1}$ has regular parameters $x_2,y_2,z_2$ with
$$
x_1=x_2, y_1=y_2, z_1=x_2^{a'k}y_2^{b'k}(z_2+\alpha).
$$
$p_1$ is a 2-point, and substituting into (\ref{eqT243}), we have
$$
u_1=(x_2^{a'}y_2^{b'})^k,
v_1=z_2,
w=h(x_2^{a'}y_2^{b'},x_2^{a'k}y_2^{b'k}(z_2+\alpha))+x_2^{c'}y_2^{d'}.
$$
Thus $\overline u_1,\overline v_1$ are toroidal forms at $p_1$.

Suppose that a form (\ref{eqT244}) holds at $p'$. Then $\hat{\cal O}_{X_1,p_1}$ has regular parameters $x_2,y_2,z_2$ with 
$$
\overline x_1=x_2,
\overline y_1=y_2,
z_1=x_2^{a'k}(z_2+\alpha).
$$
$p_1$ is a 1-point, and substituting into (\ref{eqT244}), we have
$$
u_1=x_2^{a'k},
v_1=z_2,
w=h(x_2^{a'},x_2^{a'k}(z_2+\alpha))+x_2^{c'}y_2^{d'}.
$$
Thus $\overline u_1,\overline v_1$ are toroidal forms at $p_1$.

We now prove 6. We may assume that $\tau_{f_1}(p')=\tau_f(\overline p)$. There exists a smallest $i$ such that $X_1\rightarrow Z_{i+1}$ is an isomorphism at $p'$. Let $p_2$ be the image of $p'$ in $Z_{i+1}$. Let $p$ be the image of $p_2$ in $Z_0$. 

If $Z_0\rightarrow X$ is not an isomorphism at $p$, then $\overline p$ is a 2-point, and we have permissible parameters $x,y,z$ at $\overline p$ such that 
\begin{equation}\label{eqT246}
\begin{array}{ll}
u&=(x^ay^b)^k\\
v&=z\\
w&=\sum_{i,l\ge 0}a_{il}z^l(x^ay^b)^i+x^ey^f,
\end{array}
\end{equation}
where $\text{gcd}(a,b)=1$ and the sum is over $i$ such that $(ia,ib)\not\ge(e,f)$.
Since we are assuming $w$ is good at $\overline p$ for $f$, we have $a_{il}=0$ if $k$ divides $i$.

Since $Z_0\rightarrow X$ is a sequence of blow ups of 2-curves above $\overline p$, we have regular parameters $x_1,y_1,z$ in $\hat{\cal O}_{X_0,p}$ such that 
\begin{equation}\label{eqT247}
x=x_1^{\overline a}y_1^{\overline b}, y=x_1^{\overline c}y_1^{\overline d}
\end{equation}
with $\overline a\overline d-\overline b\overline c=\pm1$,
(and $\text{gcd}(a\overline a+b\overline c,a\overline b+b\overline d)=1$) or
$$
x=\overline x_1^{\overline a}(\overline y_1+\alpha)^{\overline b}, y=\overline x_1^{\overline c}(\overline y_1+\alpha)^{\overline d}
$$
where $\overline a,\overline c\in{\bf N}$, $\overline b,\overline d\in{\bf Q}$, $0\ne \alpha\in{\bf k}$ are such that 
\begin{equation}\label{eqT248}
x^ay^b=\overline x_1^{\overline aa+\overline c b},
x^ey^f=\overline x_1^{\overline ae+\overline cf}(\overline y_1+\alpha).
\end{equation}

If either (\ref{eqT247}) or (\ref{eqT248}) holds, we see from substitution into (\ref{eqT246}) that $w$ is good at $p$ for $f\circ\Phi_0$.

Since $\overline q$ is a 1-point, $\overline f_0=\Psi_1^{-1}\circ f\circ\Phi_0$ is not a morphism at $p$. Let $p_1$ be the image of $p_2$ in $Z_i$. At $p_1$ we have an expression (\ref{eqT225}) or (\ref{eqT180}). 

First suppose that (\ref{eqT225}) holds at $p_1$. Since $w$ is good at $p$ for $f\circ\Phi_0$, we have $a_{ij}=0$ if $a$ divides $i$. 

Suppose that (\ref{eqT225}) holds and at $p_2$ and there is an expression (\ref{eqT231}).

Suppose that $i+aj\in a{\bf Z}$ for some $i$. Then $i\in a{\bf Z}$, and since $w$ is good at $p$ for $f\circ\Phi_0$, we have that $a_{ij}=0$. Thus $w$ is good at $p_2$ for $h_{i+1}$
(and $w$ is good at $p'$ for $f_1$).

Suppose that (\ref{eqT225}) holds and at $p_2$ there is an expression (\ref{eqT232}).

If $w$ is not good at $p_2$ for $h_{i+1}$, there exists $a_{ij}\ne 0$ with 
$i+(b+1)j<n$ and $i+(b+1)j\in a{\bf Z}+(b+1){\bf Z}$. Thus $i\in a{\bf Z}+(b+1){\bf Z}$ which implies that $i$ is in the kernel of the surjective projection homomorphism
$$
[\sum_{a_{ij}\ne 0}i{\bf Z}+a{\bf Z}]/a{\bf Z}]\rightarrow 
[\sum_{a_{ij}\ne 0}i{\bf Z}+a{\bf Z}+(b+1){\bf Z}]/[a{\bf Z}+(b+1){\bf Z}].
$$
thus 
$$
\tau_{f_1}(p')=\tau_{h_{i+1}}(p_2)<\tau_{f\circ\Phi_0}(p)\le\tau_f(\overline p),
$$
 a contradiction.

Suppose that (\ref{eqT225}) holds and at $p_2$ there is an expression (\ref{eqT233}).

Suppose that $w$ is not good at $p_2$ for $h_{i+1}$. Then there exists $a_{ij}\ne0$ with
$$
(i+bj,i+(b+1)j)\not\ge (n,n)
$$ 
and
$$
(i+bj,i+(b+1)j)\in
(a-b,a-b-1){\bf Z}+(b,b+1){\bf Z}=
(a,a){\bf Z}+(b,b+1){\bf Z},
$$
which implies that $(i,i)\in (a,a){\bf Z}+(b,b+1){\bf Z}$.
Thus $i\in a{\bf Z}$, a contradiction, and we conclude that $w$ is good at $p_2$ for $h_{i+1}$.

Now suppose that (\ref{eqT180}) holds at $p_1$. Since $w$ is good at $p$ for $f\circ\Phi_0$, we have $a_{il}=0$ if $k$ divides $i$.

Suppose that (\ref{eqT180}) holds and at $p_2$ there is an expression (\ref{eqT237}).
Suppose that 
$(i+lk)(a,b)\not\ge(e,f)$ and
$i+lk\in k{\bf Z}$. then $i\in k{\bf Z}$, and since $w$ is good at $p$ for $f\circ\Phi_0$, we have that $a_{il}=0$. Thus $w$ is good at $p_2$ for $h_{i+1}$.

Suppose that (\ref{eqT180}) holds and at $p_2$ there is an expression (\ref{eqT234}). Then by a similar calculation to the above, (\ref{eqT235}) is a strict inequality if $w$ is not good at $p_2$ for $h_{i+1}$.

Suppose that (\ref{eqT180}) holds and at $p_2$ there is an expression (\ref{eqT238}). Then (\ref{eqT239}) is a strict inequality if $w$ is not good at $p_2$ for $h_{i+1}$.

Suppose that (\ref{eqT180}) holds and $p_2$ has an expression (\ref{eqT250}). The  homomorphism (\ref{eqT240})
has a nontrivial kernel if $w$ is not good at $p_2$ for $h_{i+1}$.

In all of these cases, we have a contradiction to our assumption that $\tau_{f_1}(p')=\tau_f(\overline p)$.

\end{pf}

\begin{Lemma}\label{LemmaT55} Suppose that $f:X\rightarrow Y$ is pre-$\tau$-quasi-well prepared with relation $R$, $q\in Y-U(R)$ is a 2-point such that $q\in G_Y(f,\tau)$, and $u,v,w$ are permissible parameters at $q$. Then there exists a pre-$\tau$-quasi-well diagram
$$
\begin{array}{rll}
X_1&\stackrel{f_1}{\rightarrow}&Y_1\\
\Phi\downarrow&&\downarrow\Psi\\
X&\stackrel{f}{\rightarrow}&Y
\end{array}
$$
such that $\Phi$ and $\Psi$ are products of blow ups of 2-curves, $f_1$ is pre-$\tau$-quasi-well prepared, $\Phi$ is an isomorphism over $f^{-1}(Y-\Sigma(Y))$, and if $q_1\in \Psi^{-1}(q)$ (with permissible parameters $u_1,v_1,w$), then if $p_1\in f^{-1}(q)$
is such that $\tau_{f_1}(p_1)=\tau$, then there exists a series $\phi_{p_1}(u_1,v_1)$ such that
$w-\phi_{p_1}(u_1,v_1)$ is good at $p_1$.
\end{Lemma}

\begin{pf} 
Let $u,v,w$ be permissible parameters at $q$.
By Lemma \ref{LemmaT47} and Lemma \ref{Lemma31}, there exist sequences of blow ups of 2-curves
$\Psi_0:Y_0\rightarrow Y$ and $\Phi_0:X_0\rightarrow X$ making a pre-$\tau$-quasi-well diagram
$$
\begin{array}{rll}
X_0&\stackrel{f_0}{\rightarrow}&Y_0\\
\Phi_0\downarrow&&\downarrow\Psi_0\\
X&\stackrel{f}{\rightarrow}&Y
\end{array}
$$
such that  if $q_0\in\Psi_0^{-1}(q)$ is a 2-point, then there are permissible parameters $(u_0,v_0,w_0)$ at $q_0$ such that 
$$
u=u_0^av_0^b, v=u_0^c,v_0^d
$$
with $ad-bc=\pm1$,
and there are no points of the form 2 (b) or 2 (c) of Definition \ref{Def31} for $f_0$ above $q_0$. 

We construct an infinite commutative diagram of morphisms 
\begin{equation}\label{eqT175}
\begin{array}{rll}
\vdots&&\vdots\\
\downarrow&&\downarrow\\
X_n&\stackrel{f_n}{\rightarrow}&Y_n\\
\Phi_n\downarrow&&\downarrow\Psi_n\\
\vdots&&\vdots\\
\Phi_2\downarrow&&\downarrow\Psi_2\\
X_1&\stackrel{f_1}{\rightarrow}&Y_1\\
\Phi_1\downarrow&&\downarrow\Psi_1\\
X_0&\stackrel{f_0}{\rightarrow}&Y_0
\end{array}
\end{equation}
as follows. Order the 2-curves of $Y_0$ which intersect $\Psi_0^{-1}(q)$, and let $\Psi_1:Y_1\rightarrow Y_0$ be the blow up of the 2-curve $C$ of smallest order. Then construct
(by  Lemma \ref{Lemma31}) a pre-$\tau$-quasi-well diagram 
\begin{equation}\label{eqT112}
\begin{array}{rll}
X_1&\stackrel{f_1}{\rightarrow}&Y_1\\
\Phi_1\downarrow&&\downarrow\Psi_1\\
X&\stackrel{f}{\rightarrow}&Y
\end{array}
\end{equation}
where  $\Psi_1$ is a product of blow up of 2-curves and $\Phi_1$ is an isomorphism above $f^{-1}(Y-C)$. Order the 2-curves of $Y_1$ which intersect $(\Psi_0\circ\Psi_1)^{-1}(q)$
so that the 2-curves contained in the exceptional divisor of $\Psi_1$ have larger order than the order of the (strict transforms of the) 2-curves of $Y_0$.

Let $\Psi_2:Y_2\rightarrow Y_1$ be the blow up of the 2-curve $C_1$ on $Y_1$ of smallest order, and construct a pre-$\tau$-quasi-well diagram
$$
\begin{array}{rll}
X_2&\stackrel{f_2}{\rightarrow}&Y_2\\
\Phi_2\downarrow&&\downarrow\Psi_2\\
X_1&\stackrel{f_1}{\rightarrow}&Y_1
\end{array}
$$
as in (\ref{eqT112}). We now iterate to construct (\ref{eqT175}). Let 
$$
\overline\Psi_n=\Psi_0\circ\Psi_1\circ\cdot\circ\Psi_n:Y_n\rightarrow Y,
$$
$$
\overline\Phi_n=\Phi_0\circ\Phi_1\circ\cdots\circ\Phi_n:X_n\rightarrow X.
$$

For all $q_n\in \overline\Psi_n^{-1}(q)$ there exist permissible parameters $u_n,v_n,w$ at $q_n$ such that
$$
u=u_n^a(v_n+\alpha)^b,
v=u_n^c(v_n+\alpha)^d
$$
with $ad-bc\ne 0$ and $\alpha\in {\bf k}$. $q_n$ is a 2-point if $\alpha=0$, and a 1-point if $\alpha\ne 0$

Suppose that $p_n\in (\overline\Psi_n\circ f_n)^{-1}(q)$. We will say that $p_n$ is good for $q_n=f_n(p_n)$ if  $\tau_{f_n}(p_n)<\tau$ or if $\tau_{f_n}(p_n)=\tau$ and there exists a series $\phi_{p_n}(u_n,v_n)$ such that $w-\phi_{p_n}(u_n,v_n)$ is good for $f_n$ at $p_n$.

We first observe that if $p_n$ is good for $q_n$ at $p_{n+1}\in\Phi_{n+1}^{-1}(p_n)$, then $p_{n+1}$ is good for $q_{n+1}$. (This follows from 3 of Lemma \ref{Lemma31}).

Let $\nu$ be a zero-dimensional valuation of ${\bf k}(X)$ whose center on $Y$ is $q$, and let $p_n$ be the center of $\nu$ on $X_n$, $q_n=f_n(p_n)$.

We will show that there exists $n_0$ such that $n\ge n_0$ implies $p_n$ is good for $q_n$.

Once we have established this, it will follow from compactness of the Zariski-Riemann manifold of $X$ \cite{Z} that there exists $n'$ such that all $p\in (\overline\Psi_{n'}\circ f_{n'})^{-1}(q)$ are good for $q'=f_{n'}(p)$, so that the conclusions of the theorem hold.

We may identify $\nu$ with an extension of $\nu$ to the quotient field of $\hat{\cal O}_{Y,q}$ which dominates $\hat{\cal O}_{Y,q}$.

If $\nu(u)$ and $\nu(v)$ are rationally dependent, then there exists $n_0$ such that $q_{n_0}$ is a 1-point, which implies that $p_{n_0}$ is good (by Remark \ref{RemarkT174}), and thus $p_n$ is good for all $n\ge n_0$.

So we may assume that  $\nu(u)$ and $\nu(v)$ are rationally independent. We then have  that
$$
u=u_n^{a_n}v_n^{b_n},
v=u_n^{c_n}v_n^{d_n}
$$
with $a_nd_n-b_nc_n=\pm 1$ for all $n$. We thus have  (for $n>>0$) that $q_n$ is a 2-point, and $p_n$ has one of the forms (\ref{eqT105}) or (\ref{eqT106}) below:

$p_n$ is a 2-point 
\begin{equation}\label{eqT105}
\begin{array}{ll}
u_n&=x_n^ay_n^b\\
v_n&=x_n^cy_n^d\\
w&=\gamma_n+x_n^ey_n^f(z_n+\beta)
\end{array}
\end{equation}
with
$$
\gamma_n=\sum \alpha_iM_i,
$$
$\beta, \alpha_i\in{\bf k}$,
$$
M_i=x_1^{a_i}y_1^{b_i}\text{ and }M_i^{e_i}=u_1^{k_i}v_1^{l_i}
$$
with $k_i,l_i,e_i\in{\bf Z}$, $e_i>0$ and $a_i,b_i\in{\bf N}$, or

$p_n$ is a 3-point 
\begin{equation}\label{eqT106}
\begin{array}{ll}
u_n&=x_n^ay_n^bz_n^c\\
v_n&=x_n^dy_n^ez_n^f\\
w&=\gamma_n+N
\end{array}
\end{equation}
with $N=x_n^{g_n}y_n^{h_n}z_n^{i_n}$, $\text{rank}(u_n,v_n,N)=3$,
$$
\gamma_n=\sum \alpha_iM_i,
$$
$\alpha_i\in{\bf k}$,
$$
M_i=x_1^{a_i}y_1^{b_i}z_1^{c_i}\text{ and }M_i^{e_i}=u_1^{k_i}v_1^{l_i}
$$
with $k_i,l_i,e_i\in{\bf Z}$, $e_i>0$ and $a_i,b_i,c_i\in{\bf N}$.

Further, for $n>>0$, either all points $p_n$ have the form (\ref{eqT105}) or all $p_n$ have the form (\ref{eqT106}).

It is shown in the proof of Lemma 5.4 \cite{C5} that if (\ref{eqT106}) holds, then for $n>>0$, we have $k_i,l_i\in{\bf N}$ for all $i$, which implies there exists a good form for $p_n$ at $q_n$. Essentially the same argument shows that the same statement  holds if (\ref{eqT105}) holds for $n>>0$.
\end{pf}

\begin{Lemma}\label{Lemma32}
Suppose that $f:X\rightarrow Y$ is pre-$\tau$-quasi-well prepared
(or $\tau$-well prepared or $\tau$-very-well prepared) and $q\in U(R)$ is a  2-point
(prepared of type 1 in Definition \ref{Def66}).
 Then $q$ is a permissible center for $R$, and  there exists
 a pre-$\tau$-quasi-well prepared (or $\tau$-well prepared or $\tau$-very-well prepared)
diagram (\ref{eq30}) of $R$ and the blow up $\Psi:Y_1\rightarrow Y$ of $q$ such that:
\begin{enumerate}
\item[1.] $\Phi$ is an isomorphism over $f^{-1}(Y-\Sigma(Y))$.
\item[2.] Suppose that $f$ is $\tau$-well prepared. Then
\begin{enumerate}
\item
Let $E$ be the exceptional divisor of $\Psi$.
Suppose that $q_1\in U(\overline R_i^1)\cap E$. Let $\gamma_i=\overline{S_{\overline R_i^1}(q_1)
\cdot E}$.
 Then $\gamma_i$ is  a prepared
curve for $R^1$ of type 6. Suppose that $q'\in U(\overline R_j^1)\cap E$.
 Let $\gamma_j=\overline{S_{\overline R_j^1}(q')\cdot E}$. Then either
\begin{enumerate}
\item $\gamma_i=\gamma_j$ or
\item $\gamma_i,\gamma_j$ intersect transversally at a 2 point on $E$ (their tangent spaces have distinct directions at this point and are otherwise disjoint). 
\end{enumerate}
\item If $\gamma$ is a prepared curve on $Y$ then the strict transform of $\gamma$ is
a prepared curve on $Y_1$.
\end{enumerate}
\end{enumerate}
\end{Lemma}

The proof of Lemma \ref{Lemma32} is a straightforward generalization of Lemma 7.13 \cite{C5}.

\begin{Lemma}\label{Lemma171} Suppose that $f:X\rightarrow Y$ is pre-$\tau$-quasi-well
prepared (or $\tau$-well prepared or $\tau$-very-well prepared)
with   relation $R$. Suppose that $q\in Y$ is
a 1-point or a 2-point such that $q\not\in U(R)$ and $q$ is prepared of type 2 of Definition \ref{Def66} for $R$. Then $q$ is a permissible center for $R$ and there exists  a pre-$\tau$-quasi-well prepared (or $\tau$-well prepared
or $\tau$-very-well prepared)
diagram (\ref{eq30}) of $R$ and the blow up $\Psi:Y_1\rightarrow Y$ of $q$ such that:
\begin{enumerate}
\item[1.] $\Phi$ is an isomorphism over $f^{-1}(Y-\Sigma(Y))$
\item[2.] Suppose that $f$ is $\tau$-well prepared.
If $\gamma$ is a prepared curve on $Y$ then the strict transform of $\gamma$ is a prepared
curve on $Y_1$.
\end{enumerate}
\end{Lemma}

The proof of Lemma \ref{Lemma171} is a straight forward generalization of Lemma 7.14 \cite{C5}

\begin{Lemma}\label{LemmaT63} Suppose that $f:X\rightarrow Y$ is 
pre-$\tau$-quasi-well prepared (or $\tau$-well prepared or $\tau$-very-well prepared), $q\in U(R)$ is a 1-point (which is prepared of type 4 of Definition \ref{Def66}).
 Then there exists a 
 pre-$\tau$-quasi-well prepared (or $\tau$-well prepared or $\tau$-very-well prepared)
 diagram 
\begin{equation}\label{eqT293}
\begin{array}{rll}
X_1&\stackrel{f_1}{\rightarrow}&Y_1\\
\Phi\downarrow&&\downarrow\Psi\\
X&\stackrel{f}{\rightarrow}&Y
\end{array}
\end{equation}
where $\Psi$ is the blow up of $q$
such that 
\begin{enumerate}
\item[1.] $\Phi$ is a sequence of blow ups of 2-curves over $f^{-1}(Y-\{q\})$.
\item[2.]
Let $E$ be the exceptional divisor of $\Psi$.
Suppose that $q_1\in U(\overline R_i^1)\cap E$. Let $\gamma_i=\overline{S_{\overline R_i^1}(q_1)
\cdot E}$.
 Then $\gamma_i$ is  a prepared
curve for $R^1$ of type 6. Suppose that $q'\in U(\overline R_j^1)\cap E$.
 Let $\gamma_j=\overline{S_{\overline R_j^1}(q')\cdot E}$. Then either
\begin{enumerate}
\item $\gamma_i=\gamma_j$ or
\item $\gamma_i,\gamma_j$ intersect transversally at a 2 point on $E$ (their tangent spaces have distinct directions at this point and are otherwise disjoint). 
\end{enumerate}
\item[3.] If $\gamma$ is a prepared curve on $Y$ then the strict transform of $\gamma$ is
a prepared curve on $Y_1$.
\end{enumerate}
\end{Lemma}

The proof of Lemma \ref{LemmaT63} is a variant of the proof of Lemma 7.13 \cite{C5}, keeping in
mind the simpler forms 5 and 6 of Definition \ref{Def357} of super parameters above a 1-point, and the simpler from (\ref{eqT100}) of 3 of Definition \ref{Def65} at $q$.

\begin{Lemma}\label{LemmaT78}
Suppose that $f:X\rightarrow Y$ is pre-$\tau$-quasi-well prepared, $\overline q\in Y$ is a 1-point such that $\overline q\not\in U(R)$, and $C\subset D_Y$ is an integral curve such that $\overline q\in C$. Suppose that $C$ satisfies 2 and 4 of Definition \ref{DefT80} of a resolving curve. Then there exists a pre-$\tau$-quasi-well prepared diagram
$$
\begin{array}{rll}
X_1&\stackrel{f_1}{\rightarrow}&Y_1\\
\Phi_1\downarrow&&\downarrow\Psi_1\\
X&\stackrel{f}{\rightarrow}&Y
\end{array}
$$
such that 
\begin{enumerate}
\item[1.] $\Psi_1$ is a product of blow ups of 2-curves and 2-points, $\Phi_1$ is a product of possible blow ups, $\Phi_1$ is an isomorphism over $f^{-1}(Y-\Sigma(Y))$.
\item[2.]  Let $\overline C$ be the strict transform of $C$ on $Y_1$. Then
$\overline C$ is a resolving curve for $f_1$ and $R^1$ at $\overline q$.
\end{enumerate}
\end{Lemma}

\begin{pf} There exists a sequence of blow ups of 2-curves $\Psi_1:Y_1\rightarrow Y$ such that the strict transform $C_1$ of $C$ on $Y_1$ contains no 3-points. Let 
\begin{equation}\label{eqT74}
\begin{array}{rll}
X_1&\stackrel{f_1}{\rightarrow}&Y_1\\
\Phi_1\downarrow&&\downarrow\Psi_1\\
X&\stackrel{f}{\rightarrow}&Y
\end{array}
\end{equation}
be the pre-$\tau$-quasi-well prepared diagram obtained by iterating the construction  of Lemma \ref{Lemma31}.

Now by Lemma \ref{Lemma31} and Lemma 5.6 \cite{C5}, there exists a pre-$\tau$-quasi-well prepared diagram

\begin{equation}\label{eqT75}
\begin{array}{rll}
X_2&\stackrel{f_2}{\rightarrow}&Y_2\\
\Phi_2\downarrow&&\downarrow\Psi_2\\
X_1&\stackrel{f_1}{\rightarrow}&Y_1
\end{array}
\end{equation}
obtained by iterating the construction of Lemma \ref{Lemma31} such that $\Psi_2$ and $\Psi_2$ are products of blow ups of 2-curves, and for all 2-points $q$ on the strict transform $C_2$ of $C$ on $Y_2$, there exist super parameters $u,v,w$ at $q$.

Let $\Psi_3:Y_3\rightarrow Y_2$ be the blow up of the 2-points on $C_2$. By Lemma \ref{Lemma171}, there exists a pre-$\tau$-quasi-well prepared diagram 
\begin{equation}\label{eqT76}
\begin{array}{rll}
X_3&\stackrel{f_3}{\rightarrow}&Y_3\\
\Phi_3\downarrow&&\downarrow\Psi_3\\
X_2&\stackrel{f_2}{\rightarrow}&Y_2.
\end{array}
\end{equation}

By iterating the above construction, (and by embedded resolution of plane curve singularities, c.f. Section 3.4, Exercise 3.13 \cite{C6}) we eventually construct a pre-$\tau$-quasi-well prepared diagram
$$
\begin{array}{rll}
X'&\stackrel{f'}{\rightarrow}&Y'\\
\Phi'\downarrow&&\downarrow\Psi'\\
X&\stackrel{f}{\rightarrow}&Y
\end{array}
$$
such that the strict transform $C'$ of $C$ on $Y'$ is nonsingular and makes SNCs with $D_{Y'}$. If $q\in C'$ is a 2-point, let $u,v,w$ be permissible parameters at $q$ such that $u=w=0$ are local equations of $C'$. 

By Lemma \ref{Lemma31} and Lemma 5.6 \cite{C5}, there exists a pre-$\tau$-quasi-well prepared diagram
$$
\begin{array}{rll}
X''&\stackrel{f''}{\rightarrow}&Y''\\
\Phi''\downarrow&&\downarrow\Psi''\\
X'&\stackrel{f'}{\rightarrow}&Y'
\end{array}
$$
such that 3 of Definition \ref{DefT80} holds for the strict transform of $C$ on $Y''$. Thus
 the conclusions of the lemma hold.
\end{pf}

\begin{Lemma}\label{Lemma67}
Suppose that $f:X\rightarrow Y$ is $\tau$-very-well prepared
  and $C\subset Y$ is a prepared curve of type 6
(of Definition \ref{Def200}). Further suppose that
$q_{\delta}\in C\cap U(\overline R_j)$ for some $\overline R_j$ associated to $R$ implies $C=\overline{E\cdot S_{\overline R_j}(q_{\delta})}$ for some component $E$ of $D_X$. 
Then $C$ is a *-permissible center for $R$, and   there exists a $\tau$-very-well prepared diagram
$$
\begin{array}{rll}
X_1&\stackrel{f_1}{\rightarrow}&Y_1\\
\Phi_1\downarrow&&\downarrow\Psi_1\\
X&\stackrel{f}{\rightarrow}&Y
\end{array}
$$
of $R$  of the form of (\ref{eq233}). 
\end{Lemma}

The proof of Lemma \ref{Lemma67} is a straight forward generalization of Lemma 7.15 \cite{C5}

\begin{Remark}\label{Remark293}
\begin{enumerate}
\item[1.] Suppose that  $f:X\rightarrow Y$ is pre-$\tau$-quasi-well prepared and $C\subset D_Y$ is a nonsingular (integral) curve which makes SNCs with $D_Y$ and contains a 1-point such that 
\begin{enumerate}
\item $q\in C\cap U(\overline R_i)$ for some  pre-relation $\overline R_i$ associated to $R$ implies the (formal) germ of $C$ at $q$ is contained in $S_{\overline R_i}(q)$, and 
\item  $q\in C-U(R)$ implies there exist super parameters
$u,v,w$ at $q$ such that $u=w=0$ are local equations of $C$
at $q$. 
\end{enumerate}
Then there exists a pre-$\tau$-quasi-well prepared diagram
$$
\begin{array}{rll}
\overline X_1&\stackrel{\overline f_1}{\rightarrow}&\overline Y_1\\
\overline\Phi_1\downarrow&&\downarrow\overline\Psi_1\\
X&\stackrel{f}{\rightarrow} &Y
\end{array}
$$
where $\overline\Psi_1$ is the blow up of $C$. 
\item[2.] Further suppose that 
$f:X\rightarrow Y$ is $\tau$-well prepared, and if $\gamma=\overline{E\cdot R_k(q_{\alpha})}$ is prepared for
$R$ of type 6, then either $C=\gamma$ or $q\in C\cap \gamma$ implies $q\in U(\overline R_k)$ and the germ of $C$ at $q$ is contained in $S_{R_k}(q)$. 
Then
there exists a $\tau$-well prepared diagram
$$
\begin{array}{rll}
X_1&\stackrel{f_1}{\rightarrow}&Y_1\\
\Phi\downarrow&&\downarrow\Psi\\
X&\stackrel{f}{\rightarrow}&Y
\end{array}
$$
 where $\Psi$ is the blow up $\overline\Psi_1$ of
$C$, possibly followed by  blow ups of 2-points which are prepared for
the transform of $R$ (of type 2 of Definition \ref{Def66}) if $C$ is prepared of type 6 for $R$, such that
\begin{enumerate}
\item If $\gamma\subset Y$ is prepared for $f$, (and $\gamma\ne C$) then the strict transform of $\gamma$
is prepared for $f_1$.
\item If $C\subset Y$ is prepared for $f$ (of type 6) and $q\in U(\overline R_i)\cap C$ for some $i$, then $\overline{E\cdot S_{\overline R_i^1}(q')}$ is prepared for $f_1$ (of type 6) for all $q'\in E\cap U(\overline R_i^1)$, where $E$ is the component of $D_{Y_1}$ dominating $C$.
\end{enumerate}
\end{enumerate}
\end{Remark}

The proof of Remark \ref{Remark293} is a straight forward generalization of Lemma 7.15 \cite{C5} and Remark 7.16 \cite{C5}

\section{Construction of a $\tau$-very well prepared morphism}

Suppose that $f:X\rightarrow Y$ is a dominant, proper morphism of nonsingular 3-folds, with toroidal structures
$D_Y$ and $D_X=f^{-1}(D_Y)$.

\begin{Theorem}\label{TheoremT21} Suppose that $\tau\ge 0$, $f:X\rightarrow Y$ is $\tau$-prepared, $q\in Y$ is a 1-point or a 2-point, and $u,v,w$ are permissible parameters at $q$ such that $u,v$ are toroidal forms at $p$ for all $p\in f^{-1}(q)$ and $u,v\in{\cal O}_{Y,q}$. Then there exists an affine neighborhood $V$ of $q$ in $Y$ and a commutative diagram 
\begin{equation}\label{eqT81}
\begin{array}{rll}
X_1&\stackrel{f_1}{\rightarrow}&Y_1\\
\Phi_1\downarrow&&\downarrow\Psi_1\\
f^{-1}(V)&\stackrel{f}{\rightarrow}&V
\end{array}
\end{equation}
such that $u=v=0$ are local equations of a nonsingular curve $\gamma_0$ in $V$, $\Psi_1$ is a product of possible blow ups of 2-curves and resolving curves which are sections over $\gamma_0$, $\Phi_1$ is a product of possible blow ups, $f_1$ is $\tau$-prepared, $\tau_{f_1}(p_1)\le\tau_f(\Phi_1(p_1))$ for $p_1\in X_1$, and there exist permissible parameters $u_{q'},v_{q'},w$ at all $q'\in\Psi_1^{-1}(q)$ such that $u,v$ are related to $u_{q'}, v_{q'}$ birationally  and $u_{q'},v_{q'},w$ are super parameters at $q'$. 
\end{Theorem}

\begin{pf} Let $V$ be an affine neighborhood of $q$ on which $u,v,w$ are uniformizing parameters
and the intersection of the fundamental locus of $f$ with the curve $u=v=0$ on $V$ is $\{q\}$. Let $W=f^{-1}(V)$ and $\overline f=f\mid W$.  We have a smooth morphism $\pi_0:V\rightarrow S_0=\text{spec}({\bf k}[u,v])$. Give $S_0$ the toroidal structure 
$$
D_{S_0}=\left\{\begin{array}{ll}
uv=0&\text{ if $q$ is a 2-point}\\
u=0&\text{ if $q$ is a 1-point.}
\end{array}\right.
$$ 
Let $\overline q=\pi_0(q)$. The morphism $\pi_0\circ \overline f:W\rightarrow S_0$ is toroidal with respect to $D_{S_0}$ and $(\pi_0\circ\overline f)^{-1}(D_{S_0})$. We construct a commutative diagram 

\begin{equation}\label{eqT21}
\begin{array}{lllll}
\vdots&&\vdots&&\vdots\\
W_n&\stackrel{\overline f_n}{\rightarrow}&V_n&\stackrel{\pi_n}{\rightarrow}&S_n\\
\overline \Phi_n\downarrow&&\overline \Psi_n\downarrow&&\Lambda_n\downarrow\\
\vdots&&\vdots&&\vdots\\
\overline \Phi_2\downarrow&&\overline \Psi_2\downarrow&&\Lambda_2\downarrow\\
W_1&\stackrel{\overline f_1}{\rightarrow}&V_1&\stackrel{\pi_1}{\rightarrow}&S_1\\
\overline\Phi_1\downarrow&&\overline\Psi_1\downarrow&&\Lambda_1\downarrow\\
W&\stackrel{\overline f}{\rightarrow}&V&\stackrel{\pi_0}{\rightarrow}&S_0.
\end{array}
\end{equation}

Here $\Lambda_1:S_1\rightarrow S_0$ is the blow up of $\overline q$, $\overline \Psi_1$ is the blow up of $\gamma_0=\pi_0^{-1}(\overline q)$, which is the curve with equations $u=v=0$ in $V$. If $q$ is a 1-point then $\gamma_0$ is a resolving curve for $\overline f$ at $q$.  $\overline\Phi_1$ is the morphism of Lemma \ref{Lemma31} (if $q$ is a 2-point) or Lemma \ref{LemmaT79} (if $q$ is a 1-point). $\pi_1\circ \overline f_1:W_1\rightarrow S_1$ is toroidal.

Suppose that $\overline q_1\in\Lambda_1^{-1}(q)$.
 Then ${\cal O}_{S_1,\overline q_1}$ has regular parameters $u_1,v_1$ defined by 
 \begin{equation}\label{eqT277}
 u=u_1, v=u_1(v_1+\alpha)
 \end{equation}
 for some $\alpha\in{\bf k}$, or 
 \begin{equation}\label{eqT278}
 u=u_1v_1, v=v_1.
 \end{equation}

$\Lambda_2:S_2\rightarrow S_1$ is the blow up of all points  $\overline q_1$ above $\overline q$ such that there exists a point $q_1$ in $\pi_1^{-1}(\overline q_1)$ such that $u_1,v_1,w$ are not super parameters at $q_1$. $\overline \Psi_2:V_2\rightarrow V_1$ is the blow up of the (disjoint) curves $\gamma_1=\pi_1^{-1}(\overline q_1)$, and $\overline \Phi_2:W_2\rightarrow W_1$ is the morphism of Lemma \ref{Lemma31} or Lemma \ref{LemmaT79}. $\pi_2\circ \overline f_2:W_2\rightarrow S_2$ is toroidal.

We continue in this way to construct (\ref{eqT21}) as long as $\overline f_n:W_n\rightarrow V_n$ does not satisfy the conclusions of the theorem.

Suppose that the algorithm never ends.

Let $\nu$ be  a 0-dimensional valuation of ${\bf k}(X)$ whose center on $Y$ is $q$,  and let $p_i$ be the center of $\nu$ on $W_i$, $q_i$ the center of $\nu$ on $V_i$. Let $\overline q_i=\pi_i(q_i)$. We may suppose that $\overline\Psi_i$ is not an isomorphism at the center of $\nu$ for all $i$. There exist permissible parameters $u_i,v_i,w_i$ at $q_i$ for all $i$ such that $u_i,v_i$ are regular parameters in ${\cal O}_{S_i,\overline q_i}$, obtained by iteration of (\ref{eqT277}) and (\ref{eqT278}), as determined by $\nu$. We will show that there exists $j_0$ such that $u_j,v_j,w$ are super parameters at $p_j$ for all $j\ge j_0$. If $u_i,v_i,w$ are super parameters at $p_i$ for some $i$, we have that $u_j,v_j,w$ are super parameters at $p_j$ for all $j\ge i$. Suppose that $u_i,v_i,w$ are not super parameters at $p_i$ for all $i$.

We may identify $\nu$ with an extension of $\nu$ to the quotient field of $\hat{\cal O}_{X,p}$ which dominates $\hat{\cal O}_{X,p}$.

There exist permissible parameters $x_i,y_i,z_i$ in $\hat{\cal O}_{W_i,p_i}$ such that 
 we have one of the following forms:

 $p_i$ a 1-point, $q_i$ (and $\overline q_i$) a 1-point
\begin{equation}\label{eqT24}
u_i=x_i^{a_i},v_i=y_i
\end{equation}
or $p_i$ a 2-point, $q_i$ (and $\overline q_i$) a 1-point 
\begin{equation}\label{eqT25}
u_i=(x_i^{a_i}y_i^{b_i})^k, v_i=z_i
\end{equation}

$p_i$ a 1-point, $q_i$ (and $\overline q_i$) a 2-point 
\begin{equation}\label{eqT26}
u_i=x_i^{a_i},v_i=x_i^{b_i}(y_i+\alpha_i)
\end{equation}
with $0\ne\alpha_i$
or $p_i$ a 2-point, $q_i$ (and $\overline q_i$) a 2-point 
\begin{equation}\label{eqT27}
u_i=x_i^{a_i}y_i^{b_i}, v_i=x_i^{c_i}y_i^{d_i}
\end{equation}

$p_i$ a 2-point, $q_i$ (and $\overline q_i$) a 2-point 
\begin{equation}\label{eqT28}
u_i=(x_i^{a_i}y_i^{b_i})^{k_i},v_i=(x_i^{a_i}y_i^{b_i})^{t_i}(z_i+\alpha_i)
\end{equation}
with $\alpha_i\ne 0$
or $p_i$ a 3-point, $q_i$ (and $\overline q_i$) a 2-point 
\begin{equation}\label{eqT29}
u_i=x_i^{a_i}y_i^{b_i}z_i^{ci}, v_i=x_i^{d_i}y_i^{e_i}z_i^{f_i}
\end{equation}

Suppose that $f\in \hat{\cal O}_{S_0,\overline q}$. Then (by embedded resolution of plane curve singularities, c.f. Section 3.4, Exercise 3.13 \cite{C6}) there exists an $n_0$ such that $uvf=0$ is a SNC divisor in $\hat{\cal O}_{S_n,\overline q_n}$ for all $n\ge n_0$.

We further have that if $g\in \hat{\cal O}_{W,p}$ is a fractional series in $u$ and $v$, then $g\in \hat{\cal O}_{W_n,p_n}$ is a fractional series in $u_n$ and $v_n$.

First suppose that $q_i$ is a 2-point for all $i$. Then $\nu(u)$ and $\nu(v)$ are rationally independent, and we have that (\ref{eqT27}) or (\ref{eqT29}) hold for all $i$.

We either have that $p_i$ is a 3-point for all $i$ (and a form (\ref{eqT29}) holds for all $i$) or some $p_i$ is a 2-point, and thus (\ref{eqT27}) holds for all $i$ sufficiently large.

If $p_i$ is a 3-point for all $i$, it follows from Lemma 5.6 \cite{C5} that $u_i,v_i,w$ are super parameters at $p_i$ for $i>>0$, a contradiction.

Suppose that $p_i$ is a 2-point for some $i$. We may then assume that $p_i$ is a 2-point for all $i$, and (\ref{eqT27}) holds for all $i$.

We have an expression
$$
\begin{array}{ll}
u&=x^ay^b\\
v&=x^cy^d\\
w&=f_1(x,y)+x^ly^mz
\end{array}
$$
in $\hat{\cal O}_{W,p}$.
There exists $n\in{\bf N}$ and $\overline a,\overline b\in{\bf Z}$ such that 
$$
x^n=u^{\overline a}v^{\overline b},
y^n=u^{\overline c}v^{\overline d}.
$$

$\nu(u^{\overline a}v^{\overline b})=\nu(x^n)>0$, $\nu(u^{\overline c}v^{\overline d})=\nu(y^n)>0$ imply that for $i>>0$, 
$$
u^{\overline a}v^{\overline b}=u_i^{\overline a_i}v_i^{\overline b_i},
u^{\overline c}v^{\overline d}=u_i^{\overline c_i}v_i^{\overline d_i}
$$
with $\overline a_i,\overline b_i,\overline c_i,\overline d_i\in{\bf N}$.

We have
$$
\begin{array}{ll}
u_i&=x_i^{a_i}y_i^{b_i}\\
v_i&=x_i^{c_i}y_i^{d_i}\\
w&=f_1(x,y)+(x_i^{e_i}y_i^{f_i})^l(x_i^{g_i}y_i^{h_i})^mz.
\end{array}
$$
$$
f_1(x,y)=f_1((u^{\overline a}v^{\overline b})^{\frac{1}{n}},(u^{\overline c}v^{\overline d})^{\frac{1}{n}})=g_1(u_i^{\frac{1}{n}},v_i^{\frac{1}{n}})\in{\bf k}[[u_i^{\frac{1}{n}},v_i^{\frac{1}{n}}]].
$$
Let $\omega\in{\bf k}$ be a primitive $n$-th of unity. Set
$$
f=\prod_{i,j=1}^ng_1(\omega^iu_i^{\frac{1}{n}},\omega^jv_i^{\frac{1}{n}})\in {\bf k}[[u_i,v_i]].
$$
Recall that for $i\ge n_0$, $fuv=0$ is a SNC divisor in  $\hat{\cal O}_{S_i,\overline q_i}$. Since $\overline q_i$ is a 2-point for all $i$,
$$
fuv=u_i^{\overline a}v_i^{\overline b}\gamma
$$
where $\gamma$ is a unit series in $\hat{\cal O}_{S_i,\overline q_i}$. As $g_1\mid f$ in $\hat{\cal O}_{W_i,p_i}$, we have that $u_i,v_i,w$ are super parameters at $p_i$.

The final case is when $q_i$ is a 1-point for some $i$.

Suppose that $q_i$ is a 1-point for some $i$. Without loss of generality, we may assume that $q_i=q$ and $p_i=p$.

If $p$ is a 1-point (and $q$ is a 1-point), we have permissible parameters $x,y,z$ in $\hat{\cal O}_{W,p}$ such that
$$
u=x^a,v=y,
w=f_1(x,y)+x^rz
$$
of the form (\ref{eqT24}).

Set
$$
f=\prod_{i=1}^af_1(\omega^iu^{\frac{1}{a}},v)\in {\bf k}[[u,v]]
$$
where $\omega\in{\bf k}$ is a primitive $a$-th root of unity.

If $p$ is a 2-point (and $q$ is a 1-point), we have permissible parameters $x,y,z\in\hat{\cal O}_{W,p}$ such that
$$
u=(x^ay^b)^k,v=z,
w=f_1(x^ay^b,z)+x^cy^d
$$
of the form (\ref{eqT25}).
Set
$$
f=\prod_{i=1}^af_1(\omega^iu^{\frac{1}{k}},v)\in {\bf k}[[u,v]]
$$
where $\omega\in {\bf k}$ is a primitive $k$-th root of unity.

Recall that $fuv=0$ is a SNC divisor in $\hat{\cal O}_{S_i,\overline q_i}$ for $i\ge n_0$.

We will now show that if $i\ge n_0$ and $q_i$ is a 2-point, then $u_i,v_i,w$ are super parameters at $q_i$.
In these cases $\overline q_i$  is a 2-point,  $u_iv_i=0$ is a 
local equation of $D_{V_i}$ and $u=0$ is a local equation of $D_{V_i}$.
Thus $fuv=u_i^mv_i^n\gamma$, where $m,n>0$ and $\gamma\in\hat{\cal O}_{S_i,\overline q_i}$ is a unit series.
 Since $f_1\mid f$ in $\hat{\cal O}_{X_i,p_i}$, and $f_1$ is a fractional series in $u_i$ and $v_i$, we have the desired conclusion.

We have reduced to the case where $q_i$ is a 1-point for all $i$.

Since $fuv=0$ is a SNC divisor for $i\ge n_0$, the only cases where $u_i,v_i,w$ are not super parameters at $p_i$ are if $u_i, v_i$ satisfy (\ref{eqT24}) at $p_i$, and 
\begin{equation}\label{eqT30}
u_i=x_i^a, v_i=y_i, w=x_i^s(y_i-\phi(x_i^a))^r
\gamma(x_i,y_i)+x_i^c(z_i+\alpha)
\end{equation}
where $\gamma$ is a unit series, $\alpha\in{\bf k}$, or $p_i$ satisfies (\ref{eqT25}), and 
\begin{equation}\label{eqT32}
u_i=(x_i^ay_i^b)^k, v_i=z_i, w=(x_i^ay_i^b)^s(z_1-\phi((x_i^ay_i^b)^k))^r
\gamma(x_i^ay_i^b,z_i)+x_i^cy_i^d
\end{equation}
where $\text{gcd}(a,b)=1$, $\text{ord }\phi>0$ and $\gamma$ is a unit series.

We consider the case when (\ref{eqT32}) holds at $p_i$. The case (\ref{eqT30}) is similar.

By the constructions of Lemma \ref{Lemma31} and Lemma \ref{LemmaT79}, we have a factorization of $W_{i+1}\rightarrow W_i$, by monoidal transforms
$$
W_{i+1}=Z_m\stackrel{\Omega_m}{\rightarrow} Z_{m-1}\stackrel{\Omega_{m-1}}{\rightarrow} \cdots \stackrel{\Omega_2}{\rightarrow} Z_1 \stackrel{\Omega_1}{\rightarrow} Z_0=W_i
$$
Let $a_j$ be the center of $\nu$ on $Z_j$. We may assume that each morphism $Z_l\rightarrow Z_{l-1}$ is not an isomorphism at the center of $\nu$.

If $a_j<m$, we have permissible parameters $x_{ij},y_{ij},z_{ij}$ in $\hat{\cal O}_{Z_j,a_j}$ such that
$a_j$ is  a 2-point and

\begin{equation}\label{eqT33}
\begin{array}{ll}
u_i&=(x_{ij}^{a}y_{ij}^{b})^{k},\\
v_i&=z_{ij}x_{ij}^{e_i}y_{ij}^{f_i},\\
w&=(x_{ij}^{a}y_{ij}^{b})^{s}(z_{ij}x_{ij}^{e_i}y_{ij}^{f_i}-\phi((x_{ij}^{a}y_{ij}^{b})^{k})^{r}
\gamma(x_{ij}^{a}y_{ij}^{b},z_{ij}x_{ij}^{e_i}y_{ij}^{f_i})+x_{ij}^{c}y_{ij}^{d}
\end{array}
\end{equation}
with $(e_i,f_i)<(ak,bk)$.  

After possibly interchanging $x_{ij}$ and $y_{ij}$, we have that $e_i<ak$ and there are regular parameters $x_{i,j+1}, y_{i,j+1}, z_{i,j+1}$ in $\hat{\cal O}_{Z_{j+1},a_{j+1}}$ defined by 
\begin{equation}\label{eqT34}
x_{ij}=x_{i,j+1},
z_{ij}=x_{i,j+1}(z_{i,j+1}+\alpha).
\end{equation}
or 
\begin{equation}\label{eqT35}
x_{ij}=x_{i,j+1}z_{i,j+1}, z_{ij}=z_{i,j+1}.
\end{equation}

Suppose that (\ref{eqT34}) holds. Then 
\begin{equation}\label{eqT82}
\begin{array}{ll}
u_i&=(x_{i,j+1}^{a}y_{i,j+1}^{b})^k\\
v_i&=(z_{i,j+1}+\alpha)x_{i,j+1}^{e_i+1}y_{i,j+1}^{f_i}\\
w&=(x_{i,j+1}^{a}y_{i,j+1}^{b})^{s}[x_{i,j+1}^{e_i+1}y_{i,j+1}^{f_i}(z_{i,j+1}+\alpha)-\phi((x_{i,j+1}^{a}y_{i,j+1}^{b})^{k})]^{r}\\
&\gamma(x_{i,j+1}^{a}y_{i,j+1}^{b},x_{i,j+1}^{e_i+1}y_{i,j+1}^{f_i}(z_{i,j+1}+\alpha))+x_{i,j+1}^{c}y_{i,j+1}^{d}.
\end{array}
\end{equation}

Suppose that
$0\ne\alpha$, $(e_i+1,f_i)<(ak,bk)$ and $af_i-b(e_i+1)\ne 0$ in (\ref{eqT82}).
Then the rational map $Z_{j+1}\rightarrow V_{i+1}$ is a morphism near $a_{j+1}=p_{i+1}$.

We  make a change of variables, to get permissible parameters $\overline x,\overline y,\overline z\in\hat{\cal O}_{W_{i+1},p_{i+1}}$ and $0\ne\beta\in{\bf k}$ such that 
$$
\begin{array}{ll}
u_{i+1}&=\frac{u_i}{v_i}=x_{i,j+1}^{ak-(e_i+1)}y_{i,j+1}^{bk-f_i}(z_{i,j+1}+\alpha)^{-1}=\overline x^{ak-(e_i+1)}\overline y^{bk-f_i}\\
v_{i+1}&=v_i=(z_{i,j+1}+\alpha)x_{i,j+1}^{e_i+1}y_{i,j+1}^{f_i}=\overline x^{e_i+1}\overline y^{f_i}\\
w&=(\overline x^{a}\overline y^{b})^{s}[\overline x^{e_i+1}\overline y^{f_i}-\phi((\overline x^{a}\overline y^{b})^{k})]^{r}
\gamma(\overline x^{a}\overline y^{b},\overline x^{e_i+1}\overline y^{f_i})+\overline x^{c}\overline y^{d}(\overline z+\beta)\\
&=\overline x^{as+r(e_i+1)}\overline y^{bs+rf_i}(1-\frac{\phi((\overline x^{a}\overline y^{b})^{k})}{\overline x^{e_i+1}\overline y^{f_i}})^{r}\gamma(\overline x^{a}\overline y^{b},\overline x^{e_i+1}\overline y^{f_i})+\overline x^{c}\overline y^{d}(\overline z+\beta)
\end{array}
$$
and $u_{i+1},v_{i+1},w$ are super parameters at $p_{i+1}$.

Suppose that (\ref{eqT34}) holds, $0\ne\alpha$, $(e_i+1,f_i)<(ak,bk)$ and $(e_i+1,f_i)=t_i(a,b)$ for some integer $t_i$. Then we make a change of variable  to get 
permissible parameters $\overline x,\overline y,\overline z\in\hat{\cal O}_{W_{i+1},p_{i+1}}$
such that
$$
\begin{array}{ll}
u_{i+1}&=\frac{u_i}{v_i}=(\overline x^{a}\overline y^{b})^{k-t_i}\\
v_{i+1}&=v_i=(\overline x^{a}\overline y^{b})^{t_i}(\overline z+\alpha)\\
w&=(\overline x^{a}\overline y^{b})^{s}(\overline z+\alpha)^{\frac{s}{k}}[(\overline x^{a}\overline y^{b})^{t_i}(\alpha+\overline z)-\phi((\overline x^{a}\overline y^{b})^{k}(\overline z+\alpha))]^{r}\\
&\,\,\,\gamma(\overline x^{a}\overline y^{b}(\overline z+\alpha)^{\frac{1}{k}},(\overline x^{a}\overline y^{b})^{t_i})(\overline z+\alpha))+\overline x^{c}\overline y^{d}\\
&=(\overline x^{a}\overline y^{b})^{s+t_ir}(\overline z+\alpha)^{\frac{s}{k}}(\alpha+\overline z-\frac{\phi((\overline x^{a}\overline y^{b})^{k}(\overline z+\alpha))}{(\overline x^{a}\overline y^{b})^{t_i}})^{r}\\
&\,\,\,\gamma(\overline x^{a}\overline y^{b}(\overline z+\alpha)^{\frac{1}{k}},(\overline x^{a}\overline y^{b})^{t_i}(\overline z+\alpha))+\overline x^{c}\overline y^{d}
\end{array}
$$
as $t_i<k$, and $u_{i+1},v_{i+1},w$ are thus super parameters at $p_{i+1}$.

Suppose that (\ref{eqT34}) holds, and $(e_i+1,f_i)=(ak,bk)$.  Then 
the rational map $Z_{j+1}\rightarrow V_{i+1}$ is a morphism near $a_{j+1}=p_{i+1}$, and the 1-point $q_{i+1}$ has permissible parameters defined by
$$
u_{i}=u_{i+1}, v_i=u_{i+1}(v_{i+1}+\alpha).
$$
Substituting into (\ref{eqT82}), we see that
$$
\begin{array}{ll}
u_{i+1}&=(x_{i,j+1}^{a}y_{i,j+1}^{b})^{k}\\
v_{i+1}&=z_{i,j+1}\\
w&=(x_{i,j+1}^{a}y_{i,j+1}^{b})^{s+kr}[z_{i,j+1}+\alpha
-\frac{\phi((x_{i,j+1}^{a}y_{i,j+1}^{b})^{k})}{(x_{i,j+1}^{a}y_{i,j+1}^{b})^{k}}]^{r}\gamma
+x_{i,j+1}^{c}y_{i,j+1}^{d}.
\end{array}
$$
Thus $u_{i+1},v_{i+1},w$ are super parameters at $p_{i+1}$, or we have a form (\ref{eqT32}), with
\begin{equation}\label{eqT166}
(c-as,d-bs)>(c-(a(s+kr),d-b(s+kr)).
\end{equation}
In fact, $(c-(a(s+kr),d-b(s+kr))$ must decrease from $(c-as,d-bs)$ by at least $(1,1)$.

If we have (\ref{eqT34}) with $(e_i+1,f_i)<(ak,bk)$ with $\alpha=0$, then (\ref{eqT82}) is back in the form (\ref{eqT33}) (with a decrease in $(ak-e_i)+(bk-f_i)$.

Suppose that (\ref{eqT35}) holds.
The rational map $Z_{j+1}\rightarrow V_{i+1}$ is a morphism near $a_{j+1}=p_{i+1}$, and the 2-point $q_{i+1}$ has permissible parameters defined by
 $u_i=u_{i+1}v_{i+1}, v_i=v_{i+1}$. Substituting into (\ref{eqT33}), we have

$$
\begin{array}{ll}
u_{i+1}&=x_{i,j+1}^{ak-e_i}y_{i,j+1}^{bk-f_i}z_{i,j+1}^{ak-e_i-1}\\
v_{i+1}&=x_{i,j+1}^{e_i}y_{i,j+1}^{f_i}z_{i,j+1}^{e_i+1}\\
w&=(x_{i,j+1}^{a}y_{i,j+1}^{b}z_{i,j+1}^{a})^{s}(x_{i,j+1}^{e_i}y_{i,j+1}^{f_i}z_{i,j+1}^{e_i+1}-\phi((x_{i,j+1}^{a}y_{i,j+1}^{b}z_{i,j+1}^{a})^{k}))^{r}
\gamma\\
&\,\,\,+x_{i,j+1}^{c}y_{i,j+1}^{d}z_{i,j+1}^{c}\\
&=(x_{i,j+1}^{a}y_{i,j+1}^{b}z_{i,j+1}^{a})^{s}(x_{i,j+1}^{e_i}y_{i,j+1}^{f_i}z_{i,j+1}^{e_i+1})^{r}
(1-\frac{\phi((x_{i,j+1}^{a}y_{i,j+1}^{b}z_{i,j+1}^{a})^{k})}{x_{i,j+1}^{e_i}y_{i,j+1}^{f_i}z_{i,j+1}^{e_i+1}})^{r}
\gamma\\
&\,\,\,+x_{i,j+1}^{c}y_{i,j+1}^{d}z_{i,j+1}^{c}
\end{array}
$$
Thus $u_{i+1},v_{i+1},w$ are super parameters at the 3-point $p_{i+1}$.

We thus have that $u_i,v_i,w$ are super parameters at $p_i$ for $i>>0$ unless each $p_i$ has a form (\ref{eqT32}), and by (\ref{eqT166}) we eventually get
$$
(c-as,d-bs)<(0,0).
$$
We thus have that $w=x_i^{c}y_i^{d}\gamma$ where $\gamma$ is a unit series. Thus $u_i,v_i,w$ are super parameters at $p_i$.

We have shown that for any 0-dimensional valuation $\nu$ of ${\bf k}(X)$ whose center is $q$ on $Y$, there exists $j_1$ such that $u_j,v_j,w$ are super parameters at $p_j$ for $j\ge j_1$. By compactness of the Zariski Riemann manifold \cite{Z}, it follows that
 the sequence (\ref{eqT21}) must terminate after a finite number of steps, in $W_n\rightarrow V_n$ satisfying the conclusions of the theorem.

\end{pf}

\begin{Theorem}\label{Theorem6}  
 Suppose that $\tau\ge0$, $f:X\rightarrow Y$
is pre-$\tau$-quasi-well  prepared (or $\tau$-well prepared) with
 relation $R$, and $C\subset Y$ is a reduced (but possibly not irreducible) curve
consisting of components of the fundamental locus of $f$ which contain a 1-point of $Y$. 
 Then there exist  sequences of possible  blow ups  $\Phi_1:X_1\rightarrow X$ and $\Psi_1:Y_1\rightarrow Y$ such that
 there is a commutative diagram
 
\begin{equation}\label{eq140}
\begin{array}{rll}
X_1&\stackrel{f_1}{\rightarrow}&Y_1\\
\Phi_1\downarrow&&\downarrow\Psi_1\\
X&\stackrel{f}{\rightarrow}&Y
\end{array}
\end{equation}
satisfying
\begin{enumerate}
\item[1.] 
$$
\tau_{f_1}(p_1)\le\tau_f(\Phi_1(p_1))
$$
for $p_1\in D_{X_1}$.
\item[2.] $f_1$ is pre-$\tau$-quasi-well prepared  with respect to a relation $R^1$.
\item[3.] $\Phi_1^{-1}(T(R))\cap G_{X_1}(f_1,\tau)\subset T(R^1)$.
\item[4.] If $f$ is $\tau$-quasi-well prepared ($\tau$-well prepared) then (\ref{eq140}) is a $\tau$-quasi-well prepared ($\tau$-well prepared) diagram.
\item[5.]  The strict transform $\overline C$ of $C$
on $Y_1$ is nonsingular and
makes SNCs with $D_{Y_1}$.
\item[6.] If $C_j$ is an irreducible component of
 $\overline C$ and $q\in U(\overline R_i^1)$ for some $\overline R_i^1$ associated to $R^1$ is such that $q\in C_j$ then the germ of $C_j$ at $q$ is contained in $S_{\overline R_i^1}(q)$.
\item[7.] If $f$ is $\tau$-well prepared, $C_j$ is an irreducible  component of $\overline C$,
$E$ is a component of $D_{Y_1}$, $q_{\alpha}\in U(R_k^1)\cap E$,
and $\gamma=\overline{E\cdot R_k^1(q_{\alpha})}$ is prepared for $R^1$ of type 6 (of Definition \ref{Def200}), then either $C_j=\gamma$ or
$q\in C_j\cap \gamma$ implies $q\in U(\overline R_k^1)$ (and thus the germ of $C_j$ at $q$ is contained in $S_{R_k}(q)$ by 5).
\item[8.] The components $C_j$ of $\overline C$  are
permissible centers (or *-permissible centers if $f$ is $\tau$-well prepared and $C_j$ is prepared of type 6) for $R^1$.
\item[9.] $\psi_1$ is a sequence of blow ups of prepared 1-points, prepared 2-points, 2-curves and resolving curves. 
\end{enumerate}

There is a  pre-$\tau$-quasi-well prepared (or $\tau$-well prepared) diagram 
$$
\begin{array}{rll}
X_2&\stackrel{f_2}{\rightarrow}&Y_2\\
\Phi_2\downarrow&&\downarrow\Psi_2\\
X_1&\stackrel{f_1}{\rightarrow}&Y_1
\end{array}
$$
where $\Psi_2$ is the blow up of $\overline C$, possibly followed by blow ups of 2-points which are
prepared of type 2 of Definition \ref{Def66} for the transform of $R$ if $f$ is $\tau$-well prepared, and $C$ contains a component which is prepared of type 6 for $R$. 
\end{Theorem}

\begin{pf} Let $\{q_1,\ldots,q_m\}$ be the 1-points of $C$ which are not contained in $U(R)$, and for which there do not exist super parameters $u,v,w\in{\cal O}_{Y,q}$  such that $u=w=0$ are local equations of $C$ at $q$. This set is finite by Lemma \ref{LemmaT135}. 

\vskip .2truein
\noindent{\bf Step 1.}  Let $\gamma_1$ be a general curve through $q_1$ on $D_Y$.

If $q\in\gamma_1$ is a 1-point, then there exist permissible parameters $u,v,w\in{\cal O}_{Y,q}$, such that $u,v$ have a toroidal form at $p$ for all $p\in f^{-1}(q)$, $u=v=0$ are local equations of $\gamma_1$ at $q$, and
either $u,v,w$ are toroidal forms at $p$ for all $p\in f^{-1}(q)$, or 
$u=w=0$ is a local equation of $C$ at $q$. This follows from (the proof of) Lemma \ref{LemmaT1} since $\gamma_1$ intersects the fundamental locus of $f$ transversally at general points of one dimensional components of the fundamental locus.

Since $\gamma_1$ is a general curve through $q_1$, $(\gamma_1-\{q_1\})\cap G_Y(f,\tau)\subset\Theta(f,Y)$ by Remark \ref{RemarkT159}. Suppose that $q\in(\gamma_1-\{q_1\})\cap(G_Y(f,\tau)-U(R))$.

Since $q$ is perfect for $f$ and by Lemma \ref{LemmaT128}, there exist locally closed subsets $\overline V_1,\ldots, \overline V_n$ of $X$ which are a partition of $G_X(f,\tau)\cap f^{-1}(q)$ and series 
$$
\phi_1(u,v),\ldots,\phi_n(u,v)
$$
 such that $u,v,w_i=w-\phi_i$ are super parameters for $f$ at $q$ and $w_i$ is weakly good at $p$ for $p\in \overline V_i$. Here $u,v,w\in{\cal O}_{Y,q}$ are permissible parameters at the 1-point $q$ such that $u,v$ have a toroidal form at $p$, $u=v=0$ are local equations of $\gamma_1$ at $q$, and $u=w=0$ are local
equations of $C$ at $q$.

If $p\in f^{-1}(q)$ is a 1-point and $\tau_f(p)>0$, then $w_i=0$ is supported on $D_X$ at $p$, since $u,v,w_i$ are super parameters at $q$.

 By blowing up 2-curves above $X$, by a map $\Phi_1:X_1\rightarrow X$, with induced map $f_1=f\circ \Phi_1:X_1\rightarrow Y$,
 and  substitution of local forms of $\Phi_1$ in (\ref{eqT279}), we obtain that $w_i=0$ is supported on $D_{X_1}$  at $p_1$ for all $p_1\in f_1^{-1}(q)$ such that $\tau_{f_1}(p_1)>0$. $u,v,w_p$ are super parameters for $f_1$ at $q$.

By Lemma \ref{LemmaT127} and Remark \ref{RemarkT278}, there exist locally closed subsets $V_1,\ldots, V_n$ of $X_1$ which are a partition of $G_{X_1}(f_1,\tau)\cap f_1^{-1}(q)$  such that
\begin{enumerate}
\item[1.] $u,v,w_i=w-\phi_i(u,v)$ are super parameters for $f_1$ at $q$
\item[2.] $w_i$ is weakly good at $p$ for $p\in V_i$
\item[3.] $w_i=0$ is supported on $D_{X_1}$ at $p$ for $p\in V_i$ if $\tau>0$.
\end{enumerate}

Suppose that $\tau>0$. 
Then there exists a relation $w_i^{e_i}-u^{a_i}\Lambda_i=0$ for $p \in V_i\cap D_{X_1}$, where $\Lambda_i$ is a unit on $V_i$, $e_i,a_i\in{\bf N}$, $e_i>1$ and $\text{gcd}(e_i,a_i)=1$.
Since $u,v,w_i$ are super parameters at the 1-point $q$, we see from (\ref{eqT280}) and (\ref{eqT279}) that there exists $\lambda_q^i\in{\bf k}$ such that $\Lambda_i(p)=\lambda_q^i$ for $p\in V_i$.

Suppose that $\tau=0$. Then $w_i=0$ is a monomial form at $p$ for $p\in V_i$ by Remark \ref{RemarkT284}.

We may now define relations $R_{q,i}$ for $f_1$ by $T(R_{q,i})=V_i$, $U(R_{q,i})=\{q\}$, and define $R_{q,i}(p)$ for $p\in T(R_{q,i})$ by $R_{q,i}=w_i^{e_i}-\lambda_q^iu^{a_i}$ if $\tau>0$ and
by $R_{q,i}=w_i$ if $\tau=0$.

We extend the transform $R^1$ of $R$ on $X_1$ by adding in the new relations $R_{q,i}$
for all $q\in (\gamma_1-\{q\})\cap(G_Y(f,\tau)-U(R))$. $f_1:X_1\rightarrow Y$ is then pre-$\tau$-quasi-well prepared for $R^1$.  We have that
$$
(\gamma_1-\{q_1\})\cap G_Y(f_1,\tau)\subset U(R^1).
$$

By Lemma \ref{Lemma32} and Lemma \ref{LemmaT63}, and since for $q\in \gamma_1\cap U(R^1_i)$, the germ of $\gamma_1$ at $q$ is not contained in $S_{R^1_i}(q)$ for any relation $R_i^1$ associated to $R^1$ ($\gamma_1$ is general), there exists a pre-$\tau$-quasi-well prepared ($\tau$-well prepared) diagram
$$
\begin{array}{rll}
X_2&\stackrel{f_2}{\rightarrow}&Y_2\\
\Phi_2\downarrow&&\downarrow\Psi_2\\
X_1&\stackrel{f_1}{\rightarrow}&Y
\end{array}
$$
such that $\Psi_2$ is a product of blow ups of prepared 1-points of type 4 (of Definition \ref{Def66}) and prepared 2-points of type 1 (of Definition \ref{Def66})  such that  the strict transform $\gamma_1^2$ of $\gamma_1$ on $Y_2$ satisfies 1,2 and 4 of Definition \ref{DefT80} of
 a resolving curve for $f_2$ at $q_1$.
Further, if
 $q\in \gamma_1^2$ is a 2-point, we have that $\tau_{f_2}(p)<\tau$ for $p\in f_2^{-1}(q)$. By Lemma 5.6 \cite{C5} and Remark \ref{Remark424}  there exists a pre-$\tau$-quasi-well prepared 
($\tau$-well prepared) diagram
$$
\begin{array}{rll}
X_3&\stackrel{f_3}{\rightarrow}&Y_3\\
\Phi_3\downarrow&&\downarrow\Psi_3\\
X_2&\stackrel{f_2}{\rightarrow}&Y_2
\end{array}
$$
where $\Phi_3$ is a product of blow ups of 2-curves and 3-points and $\Psi_3$ is a  product of blow ups of 2-curves such that if $q\in \gamma_1^3$ is a 2-point, where $\gamma_1^3$ is the strict transform  of $\gamma_1$ on $Y_3$, then there exist super parameters $u,v,w$ for $f_3$ at $q$ such that $u=w=0$ are local equations of $\gamma_1^3$ at $q$. Thus $\gamma_1^3$ is a resolving curve for $f_3$ at $\overline q$.

Let 
$$
\begin{array}{rll}
X_4&\stackrel{f_4}{\rightarrow}&Y_4\\
\Phi_4\downarrow&&\downarrow\Psi_4\\
X_3&\stackrel{f_1}{\rightarrow}&Y_3
\end{array}
$$
be the pre-$\tau$-quasi-well prepared ($\tau$-well prepared) diagram of Lemma \ref{LemmaT79} where $\Psi_4$ is the blow up of $\gamma_1^3$. The strict transform $C^4$ of $C$ on $Y_4$ intersects $\Psi_4^{-1}(q_1)$ in a 2-point. $Y_4\rightarrow Y$ is an isomorphism over a neighborhood of $\{q_2,\ldots,q_m\}$.
\vskip .2truein
\noindent{\bf Step 2.}

Iterate the  construction of Step 1 for the points $q_2,\ldots,q_m$. We obtain a commutative diagram
$$
\begin{array}{rll}
\tilde X&\stackrel{\tilde f}{\rightarrow}&\tilde Y\\
\tilde \Phi\downarrow&&\downarrow\tilde \Psi\\
X&\stackrel{f}{\rightarrow}&Y
\end{array}
$$
such that $\tilde f$ is pre-$\tau$-quasi-well prepared ($\tau$-well prepared) with respect to a relation $\tilde R$, and if $\overline C$ is the strict transform of $C$ on $\tilde Y$, and if $q\in \overline C-U(\tilde R)$ is a 1-point, then there exist super parameters $u,v,w$ at $q$ such that $u=w=0$ are local equations of $\overline C$ at $q$.
\vskip .2truein
\noindent{\bf Step 3.}
By Lemma 5.6 \cite{C5} and Remark \ref{Remark424}, there exists a pre-$\tau$-quasi-well prepared ($\tau$-well prepared) diagram
$$
\begin{array}{rll}
\tilde X_1&\stackrel{\tilde f_1}{\rightarrow}&\tilde Y_1\\
\tilde \Phi_1\downarrow&&\downarrow\tilde \Psi_1\\
\tilde X&\stackrel{\tilde f}{\rightarrow}&\tilde Y
\end{array}
$$
where $\tilde\Phi_1$ is a product of blow ups of 2-curves and 3-points, and $\tilde\Psi_1$ is a product of
 blow ups of 2-curves such that if $\overline C_1$ is the strict transform of $\overline C$ on 
 $\tilde Y_1$, and $q\in\overline C_1$ is a 2-point not in $U(\tilde R_1)$, then there exist 
 super parameters $u,v,w$ at $q$, so that $q$ is prepared of type 2 of Definition \ref{Def66}. Let
$$
\begin{array}{rll}
\tilde X_2&\stackrel{\tilde f_2}{\rightarrow}&\tilde Y_2\\
\tilde \Phi_2\downarrow&&\downarrow\tilde \Psi_2\\
\tilde X_1&\stackrel{\tilde f_1}{\rightarrow}&\tilde Y_1
\end{array}
$$
be the pre-$\tau$-quasi-well prepared ($\tau$-well prepared) diagram obtained by blowing up the 2-points on $\overline C_1-U(\tilde R^1)$ (by Lemma \ref{Lemma171}).
\vskip .2truein
\noindent {\bf Step 4.}
By embedded resolution of plane curve singularities (c.f. Section 3.4 and Exercise 3.13 \cite{C6}), we can iterate Step 3 to construct a pre-$\tau$-quasi-well prepared ($\tau$-well prepared) diagram
$$
\begin{array}{rll}
\tilde X_3&\stackrel{\tilde f_3}{\rightarrow}&\tilde Y_3\\
\tilde \Phi_3\downarrow&&\downarrow\tilde \Psi_3\\
\tilde X_2&\stackrel{\tilde f_2}{\rightarrow}&\tilde Y_2
\end{array}
$$
such that if $\overline C_3$ is the strict transform of $C$ on $\tilde Y_3$ and $q\in\overline C_3-U(\tilde R^3)$ is a 2-point, then there exist super parameters $u,v,w$ at $q$ such that $u=w=0$ are local equations of $\overline C_3$.
\vskip .2truein
\noindent{\bf Step 5.} By embedded resolution of plane curve singularities (c.f. Section 3.4 and Exercise 3.13 \cite{C6}), Lemma \ref{Lemma31} (for 3-points in $\overline C_3$), Lemma \ref{Lemma32} and Lemma \ref{LemmaT63}, there exists a pre-$\tau$-quasi-well prepared ($\tau$-well prepared) diagram
$$
\begin{array}{rll}
\tilde X_4&\stackrel{\tilde f_4}{\rightarrow}&\tilde Y_4\\
\tilde \Phi_4\downarrow&&\downarrow\tilde \Psi_4\\
\tilde X_3&\stackrel{\tilde f_3}{\rightarrow}&\tilde Y_3
\end{array}
$$
such that the strict transform $\overline C_4$ of $C$ on $\tilde Y_4$ satisfies the hypotheses of 1 of Remark \ref{Remark293}.
 \vskip .2truein
 \noindent {\bf Step 6.}
 If $f$ is $\tau$-well-prepared, then we must perform a final sequence of blowups.
 Assume that $f$ (and thus $\tilde f_4$) is $\tau$-well prepared (for the transform $\tilde R^4$ of $R$).
Let
$$
\Sigma=\left\{\begin{array}{ll}
q\in\overline C_4& \text{such that }q\in \gamma\text{ and }q\not\in  U(\tilde R^4)\\
&\text{ where $\gamma$ is a curve which is prepared of type 6 for $\tilde R^4$.}
\end{array}\right\} .
$$

$\Sigma$ is a finite set of points which are prepared of type 2 of Definition \ref{Def66} (by 5 b) of Definition \ref{Def200}).

By Lemma \ref{Lemma171}, and embedded resolution of plane curve singularities, we can construct a $\tau$-well prepared diagram
$$
\begin{array}{rll}
\tilde X_5&\stackrel{\tilde f_5}{\rightarrow}&\tilde Y_5\\
\tilde\Phi_5\downarrow&&\downarrow\tilde\Psi_5\\
\tilde X_4&\stackrel{\tilde f_4}{\rightarrow}&\tilde Y_4
\end{array}
$$
such that the hypotheses of 2 of Remark \ref{Remark293} hold, and thus the conclusions of the theorem hold.

\end{pf}

\begin{Theorem}\label{CorollaryT126} Suppose that $f:X\rightarrow Y$ is pre-$\tau$-quasi-well-prepared ($\tau$-well prepared) with relation $R$ and $R'$ is a restriction of $R$ with $U(R')=G_Y(f,\tau)-\Theta(f,Y)$. Then there exists a pre-$\tau$-quasi-well prepared ($\tau$-well prepared) diagram of $R$
$$
\begin{array}{rll}
X_1&\stackrel{f_1}{\rightarrow}&Y_1\\
\Phi_1\downarrow&&\downarrow\Psi_1\\
X&\stackrel{f}{\rightarrow}&Y
\end{array}
$$
such that  $\Phi_1$, $\Psi_1$ are products of blow ups of possible centers and $f_1$ is $\tau$-quasi-well prepared ($\tau$-well prepared) with respect to the transform $(R^1)'$ of $R'$. Further, there exists a finite set of points $\Omega=\{q_1,\ldots,q_r\}\subset Y$ such that $f_1$ is toroidal on $f_1^{-1}(Y_1-\Psi_1^{-1}(\Omega))$, and $\Psi_1\circ f_1(T((R^1)'))\subset U(R')$.
\end{Theorem}

\begin{pf} Let $C\subset Y$ be the union of the one dimensional components in the fundamental locus of $f$ which contain a 1-point.
By Lemma 
\ref{LemmaT107}, there exists a Zariski open subset $\overline Y$ of $Y$ such that $\overline Y\cap G_Y(f,\tau)=\Theta(f,Y)$, $\overline Y$ contains a generic point of each component of $C$, and there exists a commutative diagram
$$
\begin{array}{rll}
\overline X_n&\stackrel{\overline f_n}{\rightarrow}&\overline Y_n\\
\overline\Phi_n\downarrow&&\downarrow\overline\Psi_n\\
\overline X_{n-1}&\stackrel{\overline f_{n-1}}{\rightarrow}&\overline Y_{n-1}\\
\overline\Phi_{n-1}\downarrow&&\downarrow\overline\Psi_{n-1}\\
\vdots&&\vdots\\
\downarrow&&\downarrow\\
\overline X_1&\stackrel{\overline f_1}{\rightarrow}&\overline Y_1\\
\overline\Phi_1\downarrow&&\downarrow\overline\Psi_1\\
\overline f^{-1}(\overline Y)=\overline X&\stackrel{\overline f}{\rightarrow}&\overline Y\\
\end{array}
$$
where $\overline f_n$ is toroidal, each $\overline\Psi_i$ is the blow up of a nonsingular curve $\gamma_i$ (in the fundamental locus of $\overline f_{i-1}$) dominating a component of $C$, and $\overline\Phi_i$ is a sequence of blow ups of nonsingular curves dominating $\gamma_i$.
Further, we can choose $\overline Y$ so that $\overline f_n$ is $\tau$-quasi-well prepared ($\tau$-well prepared) for the transform of the restriction of $R$ to $\overline Y$.

We have that $q\in (C-\Theta(f,Y))\cap G_Y(f,\tau)$ implies $q\in U(R')$.

We apply Theorem \ref{Theorem6} to $C$.  We construct a pre-$\tau$-quasi-well prepared 
($\tau$-well prepared) diagram of $R$. Let $R^1$ be the transform of $R$ on $X_1$, and let $(R^1)'$ be the transform of $R'$ on $X_1$ (which is a restriction of $R^1$). Let the diagram be
$$
\begin{array}{rll}
X_1&\stackrel{ f_1}{\rightarrow}& Y_1\\
\Phi_1\downarrow&&\downarrow\Psi_1\\
 X&\stackrel{ f}{\rightarrow}& Y,\\
\end{array}
$$
which, after possibly replacing $\overline Y$ with a proper open subset of $\overline Y$, restricts to

$$
\begin{array}{rll}
\overline X_1&\stackrel{\overline f_1}{\rightarrow}&\overline Y_1\\
\overline\Phi_1\downarrow&&\downarrow\overline\Psi_1\\
\overline X&\stackrel{\overline f}{\rightarrow}&\overline Y\\
\end{array}
$$

over $\overline Y$ and such that $G_{Y_1}(f_1,\tau)\subset U((R^1)')\cup \Phi_1^{-1}(f^{-1}(\Theta(f,Y)))$.
We iterate this construction for the Zariski closure of $\gamma_i$ in $Y_i$, using Theorem \ref{Theorem6} if $\gamma_i$ contains a 1-point, and Lemma \ref{Lemma31} if $\gamma_i$ is a 2-curve, for $2\le i\le n$, to achieve the conclusions of the theorem.
\end{pf}

\begin{Theorem}\label{Theorem169} Suppose that $f:X\rightarrow Y$ is $\tau$-prepared.
 Then there exists a commutative diagram
$$
\begin{array}{rll}
X_1&\stackrel{f_1}{\rightarrow}&Y_1\\
\Phi\downarrow&&\downarrow\Psi\\
X&\stackrel{f}{\rightarrow}&Y
\end{array}
$$
such that $\Phi$ and $\Psi$ are  products of possible blow ups, 
and there exists a  relation $R^1$ for $f_1$
such that $f_1$ is $\tau$-quasi-well
prepared with  relation $R^1$. 
\end{Theorem}

\begin{pf}  We will give an inductive construction.
The set $G_Y(f,\tau)-\Theta(f,Y)$ is finite
 by Remark \ref{RemarkT159}. Suppose that there exists a relation $R$ on $X$ such that $\Omega=U(R)\subset G_Y(f,\tau)-\Theta(f,Y)$. (Initially, $\Omega=\emptyset$). Let $\Lambda(Y)=G_Y(f,\tau)-\Theta(f,Y)-\Omega$.

Remark \ref{RemarkT174}, Lemma \ref{LemmaT55} and Lemma \ref{Lemma31} imply there exists a pre-$\tau$-quasi-well prepared diagram (for $R$)
$$
\begin{array}{rll}
X_1&\stackrel{f_1}{\rightarrow}&Y_1\\
\Phi_1\downarrow&&\downarrow\Psi_1\\
X&\stackrel{f}{\rightarrow}&Y
\end{array}
$$
obtained by blowing up 2-curves such that if $q_1\in\Psi_1^{-1}(\Lambda(Y))$, 
then there exist algebraic permissible parameters $u_1,v_1,w_1$ at $q_1$ such that if 
$p_1 \in f_1^{-1}(q_1)\cap G_{X_1}(f_1,\tau)$, then there exist
good parameters $u_1,v_1,w_1-\phi_{p_1}(u_1,v_1)$ at $p_1$ for $f_1$.

We have $\Psi_1^{-1}(\Theta(f,Y))\cap G_{Y_1}(f_1,\tau)\subset\Theta(f_1,Y_1)$. Let $R^1$ be the transform of $R$ on $X_1$. We  restrict $R^1$ by removing from $U(R^1)$ the points of $\Theta(f_1,Y_1)$, so that $U(R^1)\cap\Theta(f_1,Y_1)=\emptyset$. Let $\Omega_1=U(R^1)$, a finite set of points (by Remark \ref{RemarkT159}).
We also have that a general point of each curve contained in $G_{Y_1}(f_1,\tau)$ is in $\Theta(f_1,Y_1)$, since $f_1$ is $\tau$-prepared. Let

$$
\Lambda(Y_1)=G_{Y_1}(f_1,\tau)-\Theta(f_1,Y_1)-\Omega_1\subset\Psi_1^{-1}(\Lambda(Y)).
$$
$\Lambda(Y_1)$ is a finite set of points. By Lemma \ref{LemmaT128},
and our construction of $f_1$, for each $q\in\Lambda(Y_1)$, we can associate 
algebraic permissible parameters 
\begin{equation}\label{eqT161}
u_q,v_q,w_q
\end{equation}
 at $q$,
 locally closed subsets $A_1,\ldots, A_{n(q)}$ of $f_1^{-1}(q)\cap G_{X_1}(f_1,\tau)$ such that
 $\{A_1,\ldots,A_{n(q)}\}$ is a partition of
 $f_1^{-1}(q)\cap G_{X_1}(f_1,\tau)$,
and
permissible parameters 
$$
u_q,v_q,w_{qi}=w_q-\phi_i(u_q,v_q)
$$
at $q$ for $1\le i\le n(q)$ such that for $p_1\in A_i$, $w_{qi}$ is good  at $p_1$ for $f_1$. 
We can assume that $u_q=v_q=0$ are local equations of a general curve through $q$ on $D_Y$ if $q$ is a 1-point (by Bertini's Theorem and Remark \ref{RemarkT174}).

Now fix $q\in\Lambda(Y_1)$. By Theorem \ref{TheoremT21}, applied to $u_q,v_q,w_{qi}$ for  $1\le i\le n(q)$, there exists an affine neighborhood $V_q$ of $q$ and a commutative diagram 
\begin{equation}\label{eqT142}
\begin{array}{rll}
\overline X&\stackrel{\overline f}{\rightarrow} &\overline Y\\
\overline\Phi\downarrow&&\downarrow\overline\Psi\\
f_1^{-1}(V_q)&\rightarrow&V_q
\end{array}
\end{equation}
satisfying the conclusions of Theorem \ref{TheoremT21} for all $w_{qi}$.
 If $\overline q\in\overline\Psi^{-1}(q)$, then there exist
$\overline u,\overline v\in{\cal O}_{\overline Y,q}$  such that $\overline u,\overline v,w_i$ are permissible parameters at $\overline q$,  $\overline u,\overline v,w_i$ are super parameters at $\overline q$ for all $i$, and by 3 of Lemma \ref{Lemma31} and 5 of Lemma \ref{LemmaT79}, for 
$$
p\in\overline f^{-1}(\overline q)\cap G_{\overline X}(\overline f,\tau)\cap\overline\Phi^{-1}(A_i),
$$
 $w_i$ is good at $p$.

Observe that we can modify the construction of (\ref{eqT81}) in Theorem \ref{TheoremT21},
by a diagram (\ref{eqT21}), by performing an arbitrary sequence of blow ups of 2-curves above each $V_i$ and $W_i$ before constructing $V_{i+1}$ and $W_{i+1}$.

Suppose that our fixed $q\in\Lambda(Y_1)$ is 1-point.
Recall that we have chosen $u_q, v_q, w_{qi}$ so that $u_q=v_q=0$ are local equations of a general curve $C$ on $D_{Y_1}$ through $q$, so that $C$ makes SNCs with $D_{Y_1}$, $C$ intersects 2-curves of $Y_1$ at general points, and $C-\{q\}$ intersects the fundamental locus $\gamma$ of $f_1$ transversally at general points of irreducible 1-dimensional components of $\gamma$.
Thus we have $G_{Y_1}(f_1,\tau)\cap (C-\{q\})\subset \Theta(f_1,Y_1)$ 
and $C$ intersects $\Theta(f_1,Y_1)$ transversally at general points of curves in $G_{Y_1}(f_1,\tau)$. For $\overline q\in G_{Y_1}(f_1,\tau)\cap (C-\{q\})$, there exist algebraic permissible parameters $$
u_{\overline q}, v_{\overline q}, w_{\overline q}
$$
 at the 1-point $\overline q$ such that $u_{\overline q}=v_{\overline q}=0$ are local equations of $C$, and since $\overline q$ is perfect for $f$, there exist series $\phi_i(u_{\overline q},v_{\overline q})\in {\bf k}[[u_{\overline q},v_{\overline q}]]$ such that $u_{\overline q},v_{\overline q},w_{\overline qi}=w_{\overline q}-\phi_i$ are super parameters at $\overline q$ for all $i$, and for $p\in f_1^{-1}(\overline q)$, some $w_{\overline qi}$ is weakly good at $p$.

  Now by Lemma \ref{LemmaT128}, for  $\overline q\in G_{Y_1}( f_1,\tau)\cap(C-\{q\})$, there exist locally closed subsets $V_i\subset X_1$ for $1\le i\le n(\overline q)$ such that $f_1^{-1}(\overline q)\cap G_{X_1}(f_1,\tau)$ is the disjoint union of $V_1,\ldots, V_{n(\overline q)}$, and for $p\in V_i$, $w_{\overline qi}$ is weakly good for $f_1$ at $p$.

By Lemma \ref{LemmaT127} and Remark \ref{RemarkT278}, there exists a sequence of blow ups of 2-curves $\hat\Phi:\hat X\rightarrow X_1$ with induced pre-$\tau$-quasi-well prepared morphism
$\hat f=f\circ\hat\Phi:\hat X\rightarrow Y_1$ such that for all $\overline q\in G_{Y_1}(f_1,\tau)\cap( C-\{q\})$, if $p\in\hat f^{-1}(\overline q)$, then for all $i$,
$u_{\overline q}, v_{\overline q}, w_{\overline qi}$ are super parameters at $p$ and if $\tau_{\hat f}(p)>0$, then $w_{\overline qi}=0$ is a local equation of a divisor supported on $D_{\hat X}$ at $p$. Further, if
$p\in\hat\Phi^{-1}(V_i)\cap G_{\hat X}(\hat f,\tau)$, then $w_{\overline qi}$ is weakly good for $\hat f$ at $p$.

We define new primitive relations $R^1_{\overline q,i}$ for the finitely many
$$
\overline q\in G_{Y_1}(\hat f,\tau)\cap (C-\{q\})
\subset G_{Y_1}(f_1,\tau)\cap (C-\{q\})
$$
and $1\le i\le n(\overline q)$.

Let $U(R^1_{\overline q,i})=\{\overline q\}$,
$T(R_{\overline q,i}^1)=\hat\Phi^{-1}(V_i)\cap G_{\hat X}(\hat f,\tau)$. 

Suppose that $\tau>0$.

If $p\in V_i$ is a 1-point, then we have an expression
$$
u_{\overline q}=x^a, v_{\overline q}=y, w_{\overline q i}=x^c\gamma
$$
where $x,y,z$ are regular parameters in $\hat{\cal O}_{\hat X,p}$, $\gamma\in\hat{\cal O}_{\hat X,p}$ is a unit series and $a\not\,\mid c$. Let $d=\text{gcd}(a,c)<a$, $e=\frac{a}{d}>1$ and $\overline c=\frac{c}{d}$.

We define $R^1_{\overline q,i}(p) = w_{\overline qi}^e-\gamma(0,0,0)^eu_{\overline q}^{\overline c}$.

If $p\in V_i$ is a 2-point, then we have an expression
$$
u_{\overline q}=(x^ay^b)^k, v_{\overline q}=z, w_{\overline qi}=(x^ay^b)^l\gamma
$$
where $x,y,z$ are regular parameters in $\hat{\cal O}_{\hat X,p}$, $\gamma\in\hat{\cal O}_{\hat X,p}$ is a unit series and $k\not\,\mid l$. 
Let $d=\text{gcd}(k,l)<k$, $e=\frac{k}{d}>1$ and $\overline c=\frac{l}{d}$.
We define
$R^1_{\overline q,i}(p)=w_{\overline qi}^e-\gamma(0,0,0)^eu_{\overline q}^{\overline c}$.

Suppose that $\tau=0$. We then define $R^1_{\overline q,i}(p)$ by the relation $w_{\overline qi}=0$ for $p\in V_i$ (by Remark \ref{RemarkT284}).

We extend the transform $\hat R^1$ of $R^1$ on $\hat X$ to include the primitive relations $R^1_{\overline q,i}$ for $\overline q\in G_{Y_1}(\hat f,\tau)\cap (C-\{q\})$ which we have just defined. Observe that we now have
$$
(C-\{q\})\cap G_{Y_1}(\hat f,\tau)\subset U(\hat R^1).
$$

By embedded resolution of plane curve singularities (c.f. Section 3.4 and Exercise 3.13 \cite{C6}), and Lemmas \ref{Lemma32} and \ref{LemmaT63}, there exists a pre-$\tau$-quasi-well prepared diagram
$$
\begin{array}{rll}
\tilde X&\stackrel{\tilde f}{\rightarrow}&\tilde Y\\
\tilde\Phi\downarrow&&\downarrow\tilde\Psi\\
\hat X&\stackrel{\hat f}{\rightarrow}&Y_1
\end{array}
$$
such that $\tilde\Psi$ is a sequence of blow ups of prepared 1-points and 2-points (of types 4 and 1 of Definition \ref{Def66}) such that  if $\tilde C$ is the strict transform of $C$ on $\tilde Y$, then $\tilde C\cap G_{\tilde Y}(f,\tau)=\{q\}$ (since the germ of $C$ at $\overline q\in G_{Y_1}(f_1,\tau)\cap(C-\{q\})$ is not contained in the surface germ $w_{\overline qi}=0$ for any $i$). Thus $\tilde C$ satisfies 1 and 2 of Definition \ref{DefT80} of a resolving curve. 4 of the definition holds since $f_1$ is $\tau$-prepared and $C$ intersects the fundamental curve of $Y_1$ transversally at general points.

By Lemma 5.6 \cite{C5} and Lemma \ref{Lemma31}, there exists a pre-$\tau$-quasi-well prepared diagram
$$
\begin{array}{rll}
\tilde X_2&\stackrel{\tilde f_2}{\rightarrow}&\tilde Y_2\\
\tilde \Phi_2\downarrow&&\downarrow\tilde \Psi_2\\
\tilde X&\stackrel{\tilde f}{\rightarrow}&\tilde Y
\end{array}
$$
where $\tilde \Psi_2$ and $\tilde \Phi_2$ are products of blow ups of 2-curves
such that the strict transform $\tilde C_2$ of $C$ satisfies 3 of Definition \ref{DefT80}.
Observe that $\tilde\Psi\circ\Psi_2$ is an isomorphism over  $q$, and
$\tilde X_2\rightarrow X_1$ is a sequence of blow ups of 2-curves over $f_1^{-1}(q)$
Thus $\tilde C_2$ is a resolving curve for $\tilde f_2$ at $q$.

We may identify $q$ with a point of $\tilde Y_2$. Let $\Phi':\tilde X_2\rightarrow X_1$ be our morphism $\Phi'=\hat\Phi\circ\tilde\Phi\circ\tilde\Phi_2$. Let $\tilde A_i=(\Phi')^{-1}(A_i)\cap G_{\tilde X_2}(\tilde f_2,\tau)$ for $1\le i\le n(q)$. We have that $u_q,v_q,w_{qi}$ are permissible parameters at $q$ such that $w_{qi}$ is good for $\tilde f_2$ at all $p\in\tilde A_i$ (by Remark \ref{RemarkT278}).

After possibly replacing the neighborhood $V_q$ of $q$ in (\ref{eqT142}) with a smaller neighborhood of $q$, we may identify $V_q$ with a neighborhood of $q$ in $\tilde Y_2$.

Let
$$
\begin{array}{rll}
X_2&\stackrel{f_2}{\rightarrow}&Y_2\\
\Phi_2\downarrow&&\downarrow\Psi_2\\
\tilde X_2&\stackrel{\tilde f_2}{\rightarrow}&\tilde Y_2
\end{array}
$$

be the pre-$\tau$-quasi-well prepared diagram of the conclusions of  Lemma \ref{LemmaT79}, where $\Psi_2$ is the blow up of $\tilde C_2$.

As $(f_1\circ\Phi')^{-1}(V_q)\rightarrow f_1^{-1}(V_q)$ is a sequence of blow ups of 2-curves, as commented after the construction of (\ref{eqT142}), we can assume that if we restrict the diagram
$$
\begin{array}{rll}
X_2&\rightarrow&Y_2\\
\Phi''\downarrow&&\downarrow\Psi''\\
X_1&\rightarrow&Y_1
\end{array}
$$
that we have constructed to $V_q$, we obtain
 the diagram
$$
\begin{array}{rll}
W_1&\rightarrow&V_1\\
\overline\Phi_1\downarrow&&\downarrow\overline\Psi_1\\
f_1^{-1}(V_q)=W&\rightarrow&V=V_q
\end{array}
$$
of (\ref{eqT21})
constructed in the proof of Theorem \ref{TheoremT21}. 

Now suppose that our fixed $q\in\Lambda(Y_1)$ is a 2-point. Then, by Lemma \ref{Lemma31} we construct a 
$\tau$-quasi-well prepared diagram for $R^1$ and the blow up $\Psi''$ of $C$
$$
\begin{array}{rll}
X_2&\stackrel{f_2}{\rightarrow}&Y_2\\
\Phi''\downarrow&&\downarrow\Psi''\\
X_1&\rightarrow Y_1
\end{array}
$$
where $\Phi''$ is a sequence of blow ups of 2-curves.

We can assume that if we restrict the diagram
$$
\begin{array}{rll}
X_2&\rightarrow&Y_2\\
\downarrow&&\downarrow\\
X_1&\rightarrow&Y_1
\end{array}
$$
to $V_q$, we obtain the diagram
$$
\begin{array}{rll}
W_1&\rightarrow&V_1\\
\overline\Phi_1\downarrow&&\downarrow\overline\Psi_1\\
f_1^{-1}(V_q)=W&\rightarrow&V=V_q
\end{array}
$$
of (\ref{eqT21}) constructed in the proof of Theorem \ref{TheoremT21}.

Let $R^2$ be our relation  on $X_2$. We have (by Remark \ref{RemarkT278}) 
$$
(\Psi'')^{-1}(\Theta(f_1,Y_1)-C)\subset\Theta(f_2,Y_2).
$$
We restrict  $R^2$ if necessary, so that $U(R^2)\cap\Theta(f_2,Y_2)=\emptyset$. Let $\Omega_2=U(R^2)$

With the notation of the proof of Theorem \ref{TheoremT21},
let $C_2$ be the Zariski closure of $\gamma_1$, the curve blown up in $V_2\rightarrow V_1$,  in $Y_2$. $C_2$ is a section over $C$. Either $C_2$ is a 2-curve or $C_2$ contains a 1-point. 

Suppose that $C_2$ contains a 1-point. 
Then $q_1=(\Psi'')^{-1}(q)\cap C_2$ is a 1-point.

In this case,
by 4 of the conclusions of Lemma \ref{LemmaT79}, at all 1-points $\overline q_1\in C_2$, there exist permissible parameters $u_1,v_1,w$ such that $u_1=v_1=0$ are local equations of $C_2$ at $\overline q_1$,
and $u_1,v_1$ are toroidal forms at all $p\in f_2^{-1}(\overline q_1)$.
Thus $C_2$ satisfies 2 and 4 of Definition \ref{DefT80} of a resolving curve for $f_2$ at $q_1$.

We can apply Lemma \ref{LemmaT78}  to construct a pre-$\tau$-quasi-well prepared diagram 
\begin{equation}\label{eqT313}
\begin{array}{rll}
X_3&\stackrel{f_3}{\rightarrow}&Y_3\\
\downarrow&&\downarrow\\
X_2&\stackrel{f_2}{\rightarrow}&Y_2
\end{array}
\end{equation}
satisfying the conclusions of Lemma \ref{LemmaT78},
so that the strict transform of $C_2$ is a resolving curve  for $f_3$ at $q_1$.

The vertical arrows of (\ref{eqT313}) are products of blow ups of 2-curves  above $V_q$.

Let 
$$
\begin{array}{rll}
X_4&\rightarrow&Y_4\\
\Phi_4\downarrow&&\downarrow\Psi_4\\
X_3&\rightarrow&Y_3
\end{array}
$$
be the pre-$\tau$-quasi-well prepared diagram of the conclusions of  Lemma \ref{LemmaT79}, where $\Psi_4$ 
is the blow up of  the strict transform of $C_2$.

Suppose that $C_2$ is a 2-curve. Then we construct (from Lemma \ref{Lemma31}) 
\begin{equation}\label{eqT335}
\begin{array}{rll}
X_4&\rightarrow&Y_4\\
\downarrow&&\downarrow\\
X_2&\rightarrow&Y_2
\end{array}
\end{equation}
as a $\tau$-quasi-well prepared diagram for $R^2$ and the blow ups of $C_2$.

As commented after the construction of (\ref{eqT142}), we may assume that the diagram (\ref{eqT335})
restricts to the diagram
$$
\begin{array}{rll}
W_2&\rightarrow&V_2\\
\downarrow&&\downarrow\\
W_1&\rightarrow&V_1
\end{array}
$$
of (\ref{eqT21}) above $V_q$.

Since (\ref{eqT21}) is finite, after finitely many iterations, we achieve a pre-$\tau$-quasi-well prepared diagram
$$
\begin{array}{rll}
X_5&\stackrel{f_5}{\rightarrow}&Y_5\\
\Phi_5\downarrow&&\downarrow\Psi_5\\
X_1&\rightarrow &Y_1
\end{array}
$$
which restricts to  the diagram (\ref{eqT142}) above $V_q$,  $\Psi_5$ is an isomorphism over $\Lambda(Y_1)-\{q\}$ and $\Phi_5$ is a sequence of blow ups of 2-curves over $\Lambda(Y_1)-\{q\}$.

By Lemma \ref{LemmaT127}, and Remark \ref{RemarkT278}, 3 of Lemma \ref{Lemma31} and 5 of Lemma \ref{LemmaT79}, there exists a sequence of blow ups of 2-curves $\Phi_6:X_6\rightarrow X_5$ such that $f_6:X_6\rightarrow Y_5$ is pre-$\tau$-quasi-well prepared, and (with the notation introduced with (\ref{eqT142})), there exist algebraic permissible parameters 
$u_q,v_q,w_q\in{\cal O}_{Y_1,q}$ 
and $w_{qi}=w_q-\phi_i(u_q,v_q)$ such that for $\overline q$ in the finite set
$$
\Sigma=(\Psi_5\circ\Psi_6)^{-1}(q)\cap G_{Y_5}(f_6,\tau)-\Theta(f_6,Y_5),
$$
 we have algebraic permissible parameters $u_{\overline q},v_{\overline q},w_q$ at $\overline q$, and series 
$$
\overline \phi_i(u_{\overline q},v_{\overline q})=\phi_i(u_q,v_q)
$$
 such that $u_{\overline q},v_{\overline q}, w_{qi}=w_q-\overline \phi_i(u_{\overline q},v_{\overline q})$ are super parameters for all $i$, and for
$$
p\in V_{\overline qi}=f_6^{-1}(\overline q)\cap (\Phi_5\circ\Phi_6)^{-1}(A_i)\cap G_{X_6}(f_6,\tau),
$$
$w_{\overline qi}=w_q-\overline\phi_i(u_{\overline q},v_{\overline q})$ is good  for $f_6$ at $p$, and $w_{\overline qi}=0$ is supported on $D_{X_6}$ at $p$ if $\tau>0$.

We define new primitive relations $R^6_{\overline q,i}$ for $\overline q\in\Sigma$ and
 $1\le i\le n(q)$.

Let $U(R^6_{\overline q,i})=\{\overline q\}$,
$T(R_{\overline q,i}^6)=V_{\overline qi}$.

 Suppose that $\tau>0$.

If $p\in V_{\overline qi}$ is a 1-point, then we have an expression
$$
u_{\overline q}=x^a, v_{\overline q}=y, w_{\overline qi}=x^c\gamma
$$
where $x,y,z$ are regular parameters in $\hat{\cal O}_{X_6,p}$, $\gamma\in\hat{\cal O}_{X_6,p}$ is a unit series and $a\not\,\mid c$.

Let
$$
a'=\frac{a}{\text{gcd}(a,c)}>1,\,\,\,c'=\frac{c}{\text{gcd}(a,c)}.
$$

We define $R^6_{\overline q,i}(p)=w_{\overline qi}^{a'}-\gamma(0,0,0)^{a'}u_{\overline q}^{c'}$.

If $p\in V_i$ is a 2-point, then we have an expression
$$
u_{\overline q}=(x^ay^b)^k, v_{\overline q}=z, w_{\overline qi}=(x^ay^b)^l\gamma
$$
where $x,y,z$ are regular parameters in $\hat{\cal O}_{X_6,p}$, $\gamma\in\hat{\cal O}_{X_6,p}$ is a unit series and $k\not\,\mid\, l$.

Let
$$
k'=\frac{k}{\text{gcd}(k,l)}>1,\,\,\,l'=\frac{l}{\text{gcd}(k,l)}.
$$

 We define
$R^6_{\overline q,i}(p)=w_{\overline qi}^{k'}-\gamma(0,0,0)^{k'}u_{\overline q}^{l'}$.

Suppose that $\tau=0$. We then define $R^6_{\overline q,i}(p)=w_{\overline qi}$.
By Remark \ref{RemarkT284}, $w_{\overline qi}$ is a monomial form at $p$.

 We can extend the transform $R^6$ of $R^5$ on $X_6$ to include these new primitive relations $R_{\overline q,i}^6$.

 We now restrict $R^6$ to remove the points of $\Theta(f_6,Y_5)$ from $U(R^6)$. 
  We have a natural identification 
  $$
  \Lambda(Y_6)=G_{Y_5}(f_6,\tau)-\Theta(f_6,Y_5)-U(R^6)
  $$
   with $\Lambda(Y_1)-\{q\}$, by our construction, and since $Y_6\rightarrow Y_1$ is an isomorphism over
  $\Lambda(Y_1)-\{q\}$.

By induction on $|\Lambda(Y_1)|$, we may iterate the above procedure to construct
a commutative diagram
$$
\begin{array}{rll}
X_7&\stackrel{f_7}{\rightarrow}&Y_7\\
\downarrow&&\downarrow\\
X&\rightarrow&Y
\end{array}
$$
such that $f_7$ is pre-$\tau$-quasi-well prepared with $U(R^7)=G_{Y_7}(f_7,\tau)-\Theta(f_7,Y_7)$.

Now by Theorem \ref{CorollaryT126}  there exists a pre-$\tau$-quasi-well diagram for $R^7$
$$
\begin{array}{rll}
X_8&\stackrel{f_8}{\rightarrow}&Y_8\\
\downarrow&&\downarrow\\
X_7&\rightarrow&Y_7
\end{array}
$$
such that $f_8$ is $\tau$-quasi-well prepared for the transform of $R^7$.
\end{pf}

\begin{Lemma}\label{Lemma348}  Suppose that $\tau\ge0$, $f:X\rightarrow Y$
is $\tau$-quasi-well prepared with   relation $R$.
Further suppose there exists a $\tau$-quasi-well prepared diagram for $R$
\begin{equation}\label{eq336}
\begin{array}{rll}
\tilde X&\stackrel{\tilde f}{\rightarrow}&\tilde Y\\
\tilde\Phi\downarrow&&\downarrow\tilde\Psi\\
X&\stackrel{f}{\rightarrow}&Y,
\end{array}
\end{equation}
where $\tilde R$ is the transform of $R$ on $\tilde X$, such that if $q_1\in U(\tilde R)$ is on a component $E$ of $D_{\tilde Y}$ such that $\tilde\Psi(E)$ is not a point, then
$T(\tilde R)\cap \tilde f^{-1}(q_1)=\emptyset$. 
Then there exists a commutative diagram
$$
\begin{array}{rll}
X_1&\stackrel{f_1}{\rightarrow}&Y_1\\
\Phi\downarrow&&\downarrow\Psi\\
\tilde X&\stackrel{\tilde f}{\rightarrow}&\tilde Y
\end{array}
$$
such that $\Phi$, $\Psi$ are products of blow ups of possible centers, and
$f_1$ is $\tau$-quasi-well prepared with   relation
$R^1$ and pre-algebraic structure. Further, $R^1$ is algebraic. ($R^1$ will in general not be the transform of $\tilde R$.)
\end{Lemma}

\begin{pf}
Given a diagram (\ref{eq336}), we will define a new  relation $\tilde R'$
on $\tilde X$ for $\tilde f$.  This is accomplished as follows. We have that $\mid U(\tilde R)-\Theta(\tilde f,\tilde Y)\mid<\infty$ by Remark \ref{RemarkT159}. Suppose that
$q_1\in U(\tilde R)-\Theta(\tilde f,\tilde Y)$ is such that $\tilde f^{-1}(q_1)\cap T(\tilde R)\ne\emptyset$. Let
$$
J_{q_1}=\{i\mid T(\tilde R_i)\cap\tilde f^{-1}(q_1)\ne \emptyset\}.
$$
Let $q=\tilde\Psi(q_1)$. For $j\in J_{q_1}$, let 
$$
u=u_{\overline R_j(q)}, v=v_{\overline R_j(q)}, w_j=w_{\overline R_j(q)}.
$$
 Let
$$
u_1=u_{\overline{\tilde R_j}(q_1)},
v_1=v_{\overline{\tilde R_j}(q_1)},
w_{j,1}=w_{\overline{\tilde R_j}(q_1)}.
$$
Since $\tilde\Psi$ is a composition of admissible blow ups for the transforms of the pre-relations $\overline R_i$ on $Y$, by the description of admissible blow ups (\ref{eqT267}) - (\ref{eqT269}) following Definition \ref{Def154},
we have in $\hat{\cal O}_{Y_1,q}$, one of the following  relations: $q_1$ a 2-point 
\begin{equation}\label{eq381}
u=u_1^{\tilde a}v_1^{\tilde b},
v=u_1^{\tilde c}v_1^{\tilde d},
w_j=u_1^{\tilde e}v_1^{\tilde f}w_{j,1}
\end{equation}
with $\tilde a\tilde d-\tilde b\tilde c=\pm 1$, or $q_1$ a 1-point, 
\begin{equation}\label{eqT294}
u=u_1,
v=u_1^{\tilde c}v_1,
w_j=u_1^{\tilde e}w_{j,1}
\end{equation}
or $q_1$ a 1-point 
\begin{equation}\label{eqT114}
u=u_1^{\tilde a}\gamma_1(u_1,v_1),
v=u_1^{\tilde b}\gamma_2(u_1,v_1),
w_j=u_1^{\tilde c}\gamma_3(u_1,v_1)w_{j,1}
\end{equation}
where $\gamma_1,\gamma_2,\gamma_3$ are unit series and $\gamma_1,\gamma_2,\gamma_3\in{\cal O}_{Y_1,q}$
 or $q_1$ a 2-point 
\begin{equation}\label{eqT275}
u=(u_1^{\tilde a}v_1^{\tilde b})^{\tilde t}\gamma_1(u_1,v_1),
v_1=(u_1^{\tilde c}v_1^{\tilde d})^{\tilde k}\gamma_2(u_1,v_1),
w_j=u_1^{\tilde e}v_1^{\tilde f}\gamma_3(u_1,v_1)w_{j,1}
\end{equation}
where $\gamma_1,\gamma_2,\gamma_3$ are unit series, and $\gamma_1,\gamma_2,\gamma_3\in{\cal O}_{Y_1,q}$, $\tilde a,\tilde b>0$.

By assumption, if $q_1$ is on a component $E$ of $D_{\tilde Y}$ we must have $\tilde\Psi(E)$ is a point.  Thus in (\ref{eq381}) we have $\tilde a,\tilde b,\tilde c,\tilde d,\tilde e,\tilde f$ all nonzero. If (\ref{eqT294}) holds, we have $\tilde c,\tilde e>0$.
In (\ref{eqT114}) we have $\tilde a,\tilde b,\tilde c>0$. In (\ref{eqT275}) we have $\tilde t,\tilde k,\tilde e,\tilde f>0$.

Suppose that $p_1\in T(\tilde R_j)\cap \tilde f^{-1}(q_1)$ and $\tau>0$.
Let $p=\tilde\Phi(p_1)\in T(R_j)$. On $X_1$, $w_j=0$ is a divisor supported on $D_{X_1}$ at $p_1$. In (\ref{eq381}), (\ref{eqT114}), (\ref{eqT275}) we see that $u=0$ is a   local equation of $D_{X_1}$ at $p_1$,
and $v=0$ is a local equation of $D_{X_1}$ at $p_1$. 
In (\ref{eqT294}), we have that $u=0$ is a local equation of $D_{X_1}$ at $p_1$, $w_j=0$ is a local equation of $D_{X_1}$ at $p_1$, and $v=0$ is a local equation of a divisor that contains $D_{X_1}$ at $p_1$.
Thus in all cases,
 there exists a natural number $r$ such that $w_j$ divides $u^r$ and $v^r$
in $\hat{\cal O}_{\tilde X,p_1}$.

Define 
$$
\eta=\eta(q_1)=\left\{\begin{array}{ll}
\text{max}\{2r^2,\tilde e,\tilde f\}&\text{if (\ref{eq381}) or (\ref{eqT275}) holds}\\
\text{max}\{2r^2,\tilde c\}&\text{if (\ref{eqT114})  holds}\\
\text{max}\{2r^2,\tilde e\}&\text{if (\ref{eqT294}) holds}
\end{array}\right.
$$
where the maximum is over $p_1\in T(\tilde R_j)\cap\tilde f^{-1}(q_1)$.

Fix $j\in J_{q_1}$. There exists $\sigma(u,v,w_j)\in {\bold k}[[u,v,w_j]]=\hat{\cal O}_{Y,q}$ such that
the order of the series $\sigma$ is greater than $\eta$ and $w_j+\sigma\in{\cal O}_{Y,q}$. Let 
$$
w_{j}^*=w_j+\sigma(u,v,w_j).
$$
For $p_1\in T(\tilde R_j)\cap\tilde f^{-1}(q_1)$, we have
$$
w_{j}^*=w_j\gamma_{p_1j}
$$
where $\gamma_{p_1j}\in\hat{\cal O}_{X_1,p_1}$ is a  unit series. If (\ref{eq381}) holds at $q_1$, set 
\begin{equation}\label{eqT164}
w_{q_1,j}=\frac{w_{j}^*}{u_1^{\tilde e}v_1^{\tilde f}}=w_{j,1}+\frac{\sigma(u_1^{\tilde a}v_1^{\tilde b},u_1^{\tilde c}v_1^{\tilde d},
u_1^{\tilde e}v_1^{\tilde f}w_{j,1})}{u_1^{\tilde e}v_1^{\tilde f}}
\in\hat{\cal O}_{\tilde Y,q_1}\cap {\bold k}(Y)={\cal O}_{\tilde Y,q_1},
\end{equation}
since $u_1,v_1\in{\cal O}_{\tilde Y,q_1}$ (this is part of the definition of a pre-$\tau$-quasi-well prepared morphism).

We further have 
$$
w_{q_1,j}=w_{j,1}\gamma_{p_1,j}.
$$

There is a similar argument if (\ref{eqT294}), (\ref{eqT114}) or (\ref{eqT275}) holds.

For $k\in J_{q_1}$, there exists $\lambda_{jk}(u,v)\in {\bold k}[[u,v]]$ such that $w_k=w_j+\lambda_{jk}(u,v)$.
 Write
$$
\lambda_{jk}(u,v)=\alpha_{k}(u,v)+h_{k}(u,v)
$$
 where $\alpha_{k}(u,v)$ is a polynomial, and $h_{k}(u,v)$  is a series of order greater than $\eta$. Set
$$
w_{k}^*=w_j+\sigma(u,v,w_j)+\lambda_{jk}(u,v)-h_{k}(u,v)\in{\cal O}_{Y,q}.
$$
$$
\begin{array}{ll}
w_{k}^*&=w_k+\sigma(u,v,w_k-\lambda_{jk}(u,v))-h_{k}(u,v)\\
&=w_k+\overline\sigma_{k}(u,v,w_k)
\end{array}
$$
where $\overline\sigma_{k}$ is a series of order greater than $\eta$. 

Suppose that (\ref{eq381}) holds at $q_1$.
Set
$$
w_{q_1,k}=\frac{w_{k}^*}{u_1^{\tilde e}v_1^{\tilde f}}.
$$
From (\ref{eq381}) we see that $u_1,v_1,w_{q_1,k}$ are permissible parameters at $q_1$ and 
$$
w_{k}=u_1^{\tilde e}v_1^{\tilde f}w_{q_1,k},
$$
 so that $w_{q_1,k}\in{\cal O}_{\tilde Y,q_1}$.
We further have that  
\begin{equation}\label{eqT329}
w_{q_1,k}=w_{k,1}\gamma_{p_1k}
\end{equation}
 for some unit series $\gamma_{p_1k}\in \hat{\cal O}_{\tilde X_1,p_1}$.

We have that for $k\in J_{q_1}$, 
\begin{equation}\label{eqT264}
w_{q_1,k}=w_{k,1}+\frac{\overline\sigma_k(u_1^{\tilde a}v_1^{\tilde b},u_1^{\tilde c}v_1^{\tilde d},u_1^{\tilde e}v_1^{\tilde f}w_{k,1})}{u_1^{\tilde e}v_1^{\tilde f}}
\end{equation}
with 
$$
\frac{\overline\sigma_k}{u_1^{\tilde e}v_1^{\tilde f}}\in \hat{\cal O}_{\tilde Y,q_1}.
$$

We further have
$$
w_{q_1,k}-w_{q_1,j}=\frac{w^*_{k}-w^*_{j}}{u_1^{\tilde e}v_1^{\tilde f}}
=\frac{\lambda_{jk}(u,v)-h_{k}(u,v)}{u_1^{\tilde e}v_1^{\tilde f}}
\in {\bold k}((u_1,v_1))\cap {\bold k}[[u_1,v_1,w_{j,1}]]={\bold k}[[u_1,v_1]].
$$

There is a similar argument if (\ref{eqT294}), (\ref{eqT114}) or (\ref{eqT275}) holds.
In these cases (\ref{eqT164}) becomes: 

$$
w_{q_1j}=\frac{w_j^*}{u_1^{\tilde e}},
w_{q_1j}=\frac{w_j^*}{u_1^{\tilde c}},
w_{q_1j}=\frac{w_j^*}{u_1^{\tilde e}v_1^{\tilde f}}
$$
respectively.

In all these cases, an equation (\ref{eqT329}) holds.

In case (\ref{eqT294}), a variant of  equation (\ref{eqT264}) holds. 

Suppose that (\ref{eqT275}) holds and $k\in J_{q_1}$. Then there exist series $\overline\sigma_k(u,v,w_k)$ such that $\text{ord }\overline\sigma_k>\eta$ and such that 
$$
w_k^*=w_k+\overline\sigma_k(u,v,w_k).
$$
Thus 
\begin{equation}\label{eqT296}
w_{q_1,k}=\gamma_3(u_1,v_1)w_{k,1}
+\frac{\overline\sigma_k((u_1^{\tilde a}v_1^{\tilde b})^{\tilde t}\gamma_1(u_1,v_1),(u_1^{\tilde c}v_1^{\tilde d})^{\tilde k}\gamma_2(u_1,v_1),u_1^{\tilde e}v_1^{\tilde f}\gamma_3(u_1,v_1)w_{j,1})}
{u_1^{\tilde e}v_1^{\tilde f}}
\end{equation}
with
$$
\frac{\overline\sigma_k}{u_1^{\tilde e}v_1^{\tilde f}}\in\hat{\cal O}_{\tilde Y,q_1}.
$$
There is a similar expression if (\ref{eqT114}) holds.

Now suppose that $\tau=0$. We make the same argument as the $\tau>0$ case if $w_{j,1}$ satisfies a form of Definition \ref{Def125} with $\alpha\ne 0$. A similar, but slightly different argument is required if $w_{j1}$ satisfies 1, 2 or 5 of Definition \ref{Def125}, with $\alpha=0$.

We now define the new  relation $\tilde R'$ on $\tilde X$ for $\tilde f$.
Set $T(\tilde R')=T(\tilde R)-\tilde f^{-1}(\Theta(\tilde f,\tilde Y))$, $U(\tilde R')=\tilde f(T(\tilde R))-\Theta(\tilde f,\tilde Y)$.
$U(\tilde R')$ is a finite set by Remark \ref{RemarkT159}. For $q_1\in U(\tilde R')$, we define primitive relations $R_{q_1,k}$ as follows. Set
$U(R_{q_1,k})=\{q_1\}$, $T(R_{q_1,k})=T(\tilde R_k)\cap \tilde f^{-1}(q_1)$.

For $p_1\in T(R_{q_1,k})$ define

$$
u_{R_{q_1,k}(p_1)}=u_{\tilde R_k(p_1)},
v_{R_{q_1,k}(p_1)}=v_{\tilde R_k(p_1)},
w_{R_{q_1,k}(p_1)}=w_{q_1,k}.
$$
If $\tau>0$ and $q_1$ is a 2-point (from \ref{eqT329}), we define $R_{q_1,k}(p_1)$ by
$$
a_{R_{q_1,k}}(p_1)=a_{\tilde R_k}(p_1), b_{R_{q_1,k}}(p_1)=b_{\tilde R_k}(p_1), e_{R_{q_1,k}}(p_1)=e_{\tilde R_k}(p_1)
$$
and
$$
\lambda_{R_{q_1,k}}(p_1)=\lambda_{\tilde R_k}(p_1)\gamma_{p_1k}(0,0,0)^{e_{\tilde R_k}(p_1)}.
$$
If $\tau=0$, and $q_1$ is a 2-point, we define
$$
a_{R_{q_1,k}}(p_1)=b_{R_{q_1,k}}(p_1)=\infty.
$$

If $q_1$ is a 1-point, we define $R_{q_1,k}(p_1)$ in an analogous way.

From the above calculations, we see that $f:\tilde X\rightarrow \tilde Y$ with the relation $\tilde R'$ satisfies 1 - 3 of the
conditions of Definition \ref{DefT60} of a pre-$\tau$-quasi-well prepared
morphism. 

Recall that all exponents are positive in (\ref{eq381}) - (\ref{eqT275}), and thus in (\ref{eqT264}) and (\ref{eqT296}). Thus
we can choose a possibly larger $\eta(q_1)$ so that $\text{ord}(\sigma_k)$ is sufficiently large in (\ref{eqT264}) and (\ref{eqT296}) that $\tilde R'$
satisfies 4 of the conditions of Definition \ref{DefT60}
(as well as 1 -- 3). Thus $\tilde f$ is pre-$\tau$-quasi-well prepared with respect to $\tilde R'$.

For $q_i\in U(R_{q_i,k})$, let $\Omega(R_{q_i,k})$ be an affine neighborhood of $q_i$ on the surface with local equation $w_{q_i,k}=0$ at $q_i$, such that
\begin{enumerate}
\item[1.] $\Omega(R_{q_i,k})$ is nonsingular and makes SNCs with $D_{\tilde Y}$.
\item[2.] $\Omega(R_{q_i,k})\cap U(\tilde R')=\{q_i\}$
\end{enumerate}

We now restrict $\tilde R'$ so that $U(\tilde R')=G_{\tilde Y}(\tilde f,\tau)-\Theta(\tilde f,\tilde Y)$, $\tilde f$ is pre-$\tau$-quasi-well prepared with relation $\tilde R'$ and  $\tilde R'$ is algebraic.

By Theorem \ref{CorollaryT126} there exists a pre-$\tau$-quasi-well prepared diagram for $\tilde R'$
$$
\begin{array}{rll}
X'&\stackrel{f'}{\rightarrow}&Y'\\
\Phi'\downarrow&&\downarrow\Psi'\\
\tilde X&\rightarrow&\tilde Y
\end{array}
$$
such that $f'$ is $\tau$-quasi-well prepared for the transform $R'$ of $\tilde R'$.
By our construction, $R'$ has
pre-algebraic structure. Further, $R'$ is algebraic. Thus the conclusions of Lemma \ref{Lemma348} hold.
\end{pf}

\begin{Remark}\label{RemarkT119} Suppose that $f:X\rightarrow Y$ is $\tau$-quasi-well prepared and $q\in U(\overline R_i)$ is a 1-point. Let
$$
u=u_{\overline R_i(q)}, v=v_{\overline R_i(q)}, w_i=w_{\overline R_i(q)}.
$$
Then
$f$ is not toroidal above $q$ if and only if  $q$ is contained in the fundamental locus of $f$.

Suppose that $f$ is not toroidal above $q$. Then the germ of the fundamental locus of $f$ at $q$ is a nonsingular (algebraic) curve, and
$u=w_i=0$ are (formal) local equations of  the fundamental locus of $f$ at $q$. 
\end{Remark}

\begin{pf} 
Suppose that $f$ is not toroidal above $q$. $u,v,w_i$ are super parameters at $p$ for all $p\in f^{-1}(q)$. From consideration of the local forms 5 and 6 of Definition \ref{Def357} of super parameters at a 1-point $q$, we see that $u=w_i=0$ are local equations of (a formal branch of) the fundamental locus of $f$ at $q$. Since the fundamental locus of $f$ is algebraic, we obtain the conclusions of the remark.
\end{pf}

\begin{Theorem}\label{TheoremT124} Suppose that $\tau\ge 0$, $f:X\rightarrow Y$ is $\tau$-quasi-well prepared  with relation $R$. Then there exists a $\tau$-quasi-well prepared diagram
$$
\begin{array}{rll}
\tilde X&\stackrel{\tilde f}{\rightarrow}&\tilde Y\\
\tilde\Phi\downarrow&&\downarrow\tilde\Psi\\
X&\stackrel{f}{\rightarrow}&Y
\end{array}
$$
where $\tilde R$ is the transform of $R$ on $\tilde X$ such that if $q_1\in U(\tilde R)$ is on a component $E$ of $D_{\tilde Y}$ such that $\tilde\Psi(E)$ is not a point, then $T(\tilde R)\cap \tilde f^{-1}(q_1)=\emptyset$.
\end{Theorem}

\begin{pf}

\noindent{\bf Step 1.}
Let $A_0$ be the set of 2-points  $q\in Y$ such that
$q\in U(\overline R_i)$ for some $\overline R_i$ associated to $R$ and  $f^{-1}(q)\cap T(R_i)\ne\emptyset$.  $A_0$ is a finite set since $f$ is $\tau$-prepared.

For $q\in A_0\cap U(\overline R_i)$, 
set 
\begin{equation}\label{eq345}
u=u_{\overline R_i(q)}, v=v_{\overline R_i(q)},
w_i = w_{\overline R_i(q)}.
\end{equation}
Let $\Psi_1:Y_1\rightarrow Y$ be the blowup of all $q\in A_0$, and let
$$
\begin{array}{rll}
X_1&\stackrel{f_1}{\rightarrow}&Y_1\\
\Phi_1\downarrow&&\downarrow\Psi_1\\
X&\stackrel{f}{\rightarrow}&Y
\end{array}
$$
be a $\tau$-quasi-well prepared  diagram of $R$ and $\Psi_1$. Such a diagram exits by
Lemma \ref{Lemma32}. Suppose that $q\in A_0\cap U(\overline R_i)$ and $q_1\in \Psi_1^{-1}(q)$ is a 1-point lying on a component $E$ of $D_{Y_1}$ and $f_1^{-1}(q_1)\cap T(\overline R^1_i)\ne\emptyset$. Then $q_1$ has regular parameters $u_1,v_1,w_{i,1}$ with 
$$
u=u_1,
v=u_1(v_1+\alpha),
w_i=u_1w_{i,1}
$$
with $0\ne\alpha\in{\bf k}$,
which implies that $E$ has local equation $u_1=0$, so that $\Psi_1(E)$ is a point.

Let $A_1$ be the set of all 2-points
$q_1\in Y_1$ such that for some $i$, $q_1\in U(\overline R_i^1)$, $f_1^{-1}(q_1)\cap T(R_i^1)\ne\emptyset$ and $q_1$
 is on a component $E$ of $D_{Y_1}$ such that $\Psi_1(E)$ is not a point. We have $A_1\subset \Psi_1^{-1}(A_0)$.
Let $\Psi_2:Y_2\rightarrow Y_1$ be the blowup of all
$q_1\in A_1$, and let (by Lemma \ref{Lemma32})
$$
\begin{array}{rll}
X_2&\stackrel{f_2}{\rightarrow}&Y_2\\
\Phi_2\downarrow&&\downarrow\Psi_2\\
X_1&\stackrel{f_1}{\rightarrow}&Y_1
\end{array}
$$
be a $\tau$-quasi-well prepared  diagram of $R^1$ and $\Psi_2$.  Continue in this way to construct (for arbitrary $n$) a sequence of $n$ blow ups of sets of 2-points
$\Psi_{k+1}:Y_{k+1}\rightarrow Y_{k}$ for $0\le k\le n-1$ with $\tau$-quasi-well prepared diagrams
$$
\begin{array}{rll}
X_{k+1}&\stackrel{f_{k+1}}{\rightarrow}&Y_{k+1}\\
\Phi_{k+1}\downarrow&&\downarrow\Psi_{k+1}\\
X_{k}&\stackrel{f_k}{\rightarrow}&Y_{k}
\end{array}
$$
of $R^k$ and $\Psi_{k+1}$.
We have a resulting $\tau$-quasi-well prepared  diagram of $R$ 

\begin{equation}\label{eq98}
\begin{array}{rll}
X_{n}&\stackrel{f_{n}}{\rightarrow}&Y_{n}\\
\Phi\downarrow&&\downarrow\Psi\\
X&\stackrel{f}{\rightarrow}&Y.
\end{array}
\end{equation}
Suppose that $q_{n}\in Y_{n}$ is a 2-point such that  $q_n$ is on
a component $E$ of $D_{Y_n}$ such that $\Psi(E)$ is  not a point, and
$q_{n}\in U(\overline R_i^{n})$, $f_{n}^{-1}(q_{n})\cap T(R_i^{n})\ne\emptyset$ for some $i$. We have permissible parameters 
\begin{equation}\label{eq414}
u_1=u_{\overline R_i^{n}(q_{n})},v_1=v_{\overline R_i^{n}(q_{n})},
w_{i,1}=w_{\overline R_i^{n}(q_{n})}
\end{equation}
 at $q_{n}$ such that for $\Psi(q_n)=q$ and with notation of (\ref{eq345}), 

\begin{equation}\label{eq96}
\begin{array}{ll}
u&=u_1\\
v&=u_1^nv_1\\
w_i&=u_1^nw_{i,1}
\end{array}
\end{equation}
or

$$
\begin{array}{ll}
u&=u_1v_1^n\\
v&=v_1\\
w_i&=v_1^nw_{i,1}.
\end{array}
$$

Suppose that $p\in  T(R_i)\cap f^{-1}(q)$ is a 1-point. First suppose that
$\tau>0$. There are permissible parameters $x,y,z$ at $p$ and $0\ne\alpha\in{\bf k}$, where $\gamma\in\hat{\cal O}_{X,p}$ is a unit, such that we have
$$
u=x^a, v=x^b(\alpha+y),
w=x^c\gamma.
$$

From the construction of $(\ref{eq98})$ and the algorithm of Lemma \ref{Lemma32}, we have that $\Phi$ is an isomorphism at points of $\Phi^{-1}(p)\cap f_n^{-1}(q_n)$.
Thus for $n\ge \text{max}\{\frac{a}{b},\frac{b}{a}\}$, $f_n^{-1}(q_n)\cap \Phi^{-1}(p)=\emptyset$.

Suppose that $\tau=0$ (and $p\in T(R_i)\cap f^{-1}(q)$ is a 1-point). Then there exist permissible parameters $x,y,z$ at $p$, $\beta\in {\bf k}$, and $0\ne\alpha\in{\bf k}$ such that
$$
u=x^a,
v=x^b(\alpha+y),
w=x^c(\beta+z).
$$

If $\beta\ne 0$, then $\Phi$ is an isomorphism at points of $\Phi^{-1}(p)\cap f_n^{-1}(q_n)$.
For $n\ge\text{max}\{\frac{a}{b},\frac{b}{a}\}$, we have $f_n^{-1}(q_n)\cap \Phi^{-1}(p)=\emptyset$.

If $\beta=0$, then $\Phi$ is a product of blow ups of sections over the curve $x=z=0$
at points of $f_n^{-1}(q_n)$, and in this case also, $\Phi^{-1}(p)\cap f_n^{-1}(q_n)=\emptyset$ for $n\ge\text{max}\{\frac{a}{b},\frac{b}{a}\}$.

Suppose that $p\in T(R_i)\cap f^{-1}(q)$ is a 2-point. First suppose that
$\tau>0$. There are permissible parameters $x,y,z$ at $p$ and $0\ne\alpha\in{\bf k}$ such that we have one of the forms 
\begin{equation}\label{eqT265}
u=(x^ay^b)^k, v=(x^ay^b)^t(\alpha+z), w_i=(x^ay^b)^l\gamma,
\end{equation}
where  $0\ne\alpha\in{\bf k}$, $\text{gcd}(a,b)=1$, $\gamma\in\hat{\cal O}_{X,p}$ is a unit,

or we have 
\begin{equation}\label{eqT163}
u=x^ay^b, v=x^cy^d, w_i=x^ey^f\gamma
\end{equation}
where  $ad-bc\ne 0$, $\gamma\in\hat{\cal O}_{X,p}$ is a unit.

Suppose that (\ref{eqT265}) holds. From the construction of $(\ref{eq98})$ and the algorithm of Lemma \ref{Lemma32},
we have that $\Phi$ is a sequence of blow ups of 2-curves at points of $\Phi^{-1}(p)\cap f_n^{-1}(q_n)$,
so that for $n\ge \text{max}\{\frac{k}{t},\frac{t}{k}\}$, $f_n^{-1}(q_n)\cap \Phi^{-1}(p)=\emptyset$,

Suppose that (\ref{eqT163}) holds.
If $f_n^{-1}(q_n)\cap \Phi^{-1}(p)\cap T(R_i^n)\ne\emptyset$ for all $n$, we will show that, after possibly interchanging $u,v$ and $x,y$,
(\ref{eqT163}) must be 
\begin{equation}\label{eqT117}
u=x^a, v=x^cy^d, w_i=x^ey^f\gamma,
\end{equation}
with $d,f>0$, and (\ref{eq96}) holds.  

We see this as follows. 

Suppose that $p_n\in f_n^{-1}(q_n)\cap T(R^n_i)$ and $\Phi(p_n)=p$.

Let $\nu$ be a valuation of ${\bf k}(X)$ whose center on $X_n$ is $p_n$. 
We identify $\nu$ with an extension of $\nu$ to the quotient field of $\hat{\cal O}_{X_n,p_n}$ which dominates $\hat{\cal O}_{X_n,p_n}$.
$q_n$ has permissible parameters  (\ref{eq414}).
After possibly interchanging $u$ and $v$, we have a relation (\ref{eq96}), so that $\nu(v)>n\nu(u)$, and $\nu(w_i)>n\nu(u)$. We can reindex $x,y,z$ so that $0<\nu(x)\le\nu(y)$. Then
$$
(c+d-nb)\nu(y)\ge (d-nb)\nu(y)+(c-na)\nu(x)>0,
$$
and 
\begin{equation}\label{eqT266}
(e+f-nb)\nu(y)\ge (f-nb)\nu(y)+(e-na)\nu(x)>0.
\end{equation}
Thus if $b\ne 0$, and $n>c+d$, we have a contradiction.

Taking $n>c+d$ for all $c,d$ in local forms (\ref{eqT163}) for 2-points $p\in T(R)$, we achieve that $b=0$ in all local forms (\ref{eqT163}) which are the images of 2-points $p_n\in T(R^n)$ which map to a point $q_n$ of $Y_n$ which is on a component $E$ of $D_{Y_n}$ such that $\Psi(E)$ is not a point. $d>0$ since $ad-bc\ne0$.

We have $f>0$ if $n>>0$. In fact, if $f=0$ in (\ref{eqT163}), we then have $e>0$, and for $n>\frac{a}{e}$, we have a contradiction to (\ref{eqT266}).

Suppose that $\tau=0$ (and $p\in T(R_i)\cap f^{-1}(q)$ is a 2-point). Then there exist
permissible parameters $x,y,z$ at $p$ and $0\ne\alpha\in{\bf k}$ such that we have one 
of the forms: 
\begin{equation}\label{eqT285}
u=(x^ay^b)^k,
v=(x^ay^b)^t(\alpha+z),
w_i=x^cy^d
\end{equation}
with $0\ne\alpha\in{\bf k}$, $ad-bc\ne 0$, $\text{gcd}(a,b)=1$, or we have 
\begin{equation}\label{eqT286}
u=x^ay^b,
v=x^cy^d,
w_i=x^ey^f(\beta+z)
\end{equation}
with $ad-bc\ne 0$, $\beta\in{\bf k}$.

Suppose that (\ref{eqT285}) holds. From the construction of $(\ref{eq98})$, we have that $\Phi$ is a sequence of blow ups of 2-curves at points of $\Phi^{-1}(p)\cap f_n^{-1}(q_n)$, so that for
$n\ge\text{max}\{\frac{k}{t},\frac{t}{k}\}$, $f_n^{-1}(q_n)\cap\Phi^{-1}(p)=\emptyset$.

Suppose that (\ref{eqT286}) holds. If $\beta\ne0$, the analysis of (\ref{eqT163}) shows that if $f_n^{-1}(q_n)\cap \Phi^{-1}(p)\cap T(R_i^n)\ne\emptyset$ for all $n$, then after possibly interchanging $u,v$ and interchanging $x,y$ (\ref{eqT286}) has the form (\ref{eqT117}), with $\gamma=\beta+z$, $d,f>0$ and  (\ref{eq96}) holds.

Suppose that $\beta=0$ (in (\ref{eqT286})), and $f_n^{-1}(q_n)\cap\Phi^{-1}(p)\cap T(R_i^n)\ne\emptyset$ for all $n$. 
After possibly interchanging $u$ and $v$, we may assume that (\ref{eq96}) holds.
Let $\nu$ be a valuation of ${\bf k}(X)$ whose center on $X_n$ is $p_n$. We identify $\nu$ with an extension of $\nu$ to the quotient field of $\hat{\cal O}_{X_n,p_n}$ which dominates $\hat{\cal O}_{X_n,p_n}$. Then we see that 
\begin{equation}\label{eqT306}
\nu(v)>n\nu(u)
\text{ and }
\nu(w_i)>n\nu(u)
\end{equation}
for all $n$.

Thus for $n$ sufficiently large, after possibly interchanging $x$ and $y$, (\ref{eqT286}) must be
\begin{equation}\label{eqT287}
u=x^a,
v=x^cy^d,
w_i=x^ey^fz
\end{equation}
with $d>0$.

We will now show that we can take $n$ sufficiently large that $f>0$ in (\ref{eqT287}).

Suppose that $f=0$ in (\ref{eqT287}). 
Then we have 
\begin{equation}\label{eqT288}
\nu(y)>n\nu(x)\text{ and }
\nu(z)>n\nu(x)
\end{equation}
for all $n$.

In the algorithm of Lemma \ref{Lemma32} (explicitly worked out in Lemma 7.13 \cite{C5}), we see that $\Phi_1:X_1\rightarrow X$ can be factored by morphisms
$$
X_1=Z_m\rightarrow \cdots\rightarrow Z_2\rightarrow Z_1\rightarrow X
$$
where $Z_1\rightarrow X$ is a sequence of blow ups of 2-curves and 3-points, and each $Z_{i+1}\rightarrow Z_i$ is the blow up of a possible curve containing a 1-point, which is in
the locus where ${\cal I}_q{\cal O}_{Z_i}$ is not invertible.

By (\ref{eqT288}), we see that there exist permissible parameters $\tilde x_1,\tilde y_1,\tilde z_1$ at the center $\tilde p_1$ of $\nu$ on $Z_1$ such that
$$
x=\tilde x_1,
y=\tilde x_1^g\tilde y_1,
z=\tilde z_1.
$$
We have
$$
u=\tilde x_1^a,
v=\tilde x_1^{c+gd}\tilde y_1^d,
w_i=\tilde x^e\tilde z_1
$$
with $a\le c+gd$, since from the construction of $Z_1\rightarrow X_1$, we have that $(u,v){\cal O}_{Z_1,\tilde p_1}$ is invertible.

If $e\ge a$, then ${\cal I}_q{\cal O}_{Z_1,\tilde p_1}$ is invertible, and $X_1\rightarrow Z_1$
is an isomorphism above $\tilde p_1$.
Then at $p_1\in X_1$, there are permissible parameters $\hat x_1,\hat y_1,\hat z_1$ such that
$$
u_{\overline R_i^1(q_1)}=\hat x_1^a,
v_{\overline R_i^1(q_1)}=\hat x_1^{c+gd-a}\hat y_1^d,
w_{\overline R_i^1(q_1)}=\hat x_1^{e-a}\hat z_1.
$$
If $e-a=0$, we have that $f_1$ is toroidal at $p_1$ (so that $\tau_{f_1}(p_1)=-\infty$).

If $e<a$, we have that $X_1\rightarrow Z_1$ is not an isomorphism above $\tilde p_1$. Without loss of generality, we may assume that each $Z_{i+1}\rightarrow Z_i$ is not an isomorphism at the center of $\nu$.

We have that $\tilde x_1=\tilde z_1=0$ are (formal) local equations at $\tilde p_1$ of the curve blown up in $Z_2\rightarrow Z_1$.

Let $\tilde p_2$ be the center of $\nu$ on $\tilde Z_2$. By (\ref{eqT288}), we have that there are permissible parameters $\tilde x_2,\tilde y_2,\tilde z_2$ at $\tilde p_2$ such that
$\tilde x_1=\tilde x_2$, $\tilde y_1=\tilde y_2$, $\tilde z_1=\tilde x_2\tilde z_2$.

We have
$$
u=\tilde x_2^a,
v=\tilde x_2^{c+gd}\tilde y_2^d,
w_i=\tilde x_2^{e+1}\tilde z_2.
$$
We see that in $X_1=Z_m$, there are regular parameters $\hat x_1,\hat y_1,\hat z_1$ such that
$$
u_{\overline R_i^1(q_1)}=\hat x_1^a,
v_{\overline R_i^1(q_1)}=\hat x_1^{c+gd-a}\hat y_1^d,
w_{\overline R_i^1(q_1)}=\hat z_2.
$$
thus $f_1$ is toroidal at $p_1$.

Iterating this analysis for the morphisms $\Psi_2,\ldots,\Psi_n$, we see that for $n>>0$, $f_n$ is toroidal at $p_n$.
In fact, if we take $n\ge\frac{e}{a}$ in (\ref{eq98}), we see that $f_n^{-1}(q_n)\cap T(R_i^n)\cap\Phi^{-1}(p)=\emptyset$.

We can thus take $n$ sufficiently large so that $f>0$ in all local forms (\ref{eqT287}), which are the images of 2-points $p_n\in T(R^n)$ which map to a point $q_n$ of $Y_n$ which is on a component $E$ of $D_{Y_n}$ such that $\Psi(E)$ is not a point.

Suppose that $p\in X$ is a 3-point such that
$p\in T(R_i)\cap f^{-1}(q)$. Then
 there are
 permissible parameters $x,y,z$ for $u,v,w_i$ at $p$ such that 
\begin{equation}\label{eq97}
\begin{array}{ll}
u&=x^ay^bz^c\\
v&=x^dy^ez^f\\
w_i&=x^gy^hz^i\gamma
\end{array}
\end{equation}
where $\gamma$ is a unit series.

We will show that we can choose $n$ sufficiently large in the diagram (\ref{eq98}), so that if $p_{n}\in X_{n}$ is a 3-point such that
$p_{n}\in\Phi^{-1}(p)\cap T(R_i^{n})$ and $q_n=f_{n}(p_{n})$ is on a component $E$ of $D_{Y_n}$ such that $\Psi(E)$ is  not a point,
  then (\ref{eq97}) must have one of the following
forms (after possibly interchanging $u,v$ and $x,y,z$): 
\begin{equation}\label{eq99}
\begin{array}{ll}
u&=x^ay^b\\
v&=x^dy^ez^f\\
w_i&=x^gy^hz^i\gamma
\end{array}
\end{equation}
where $b\ne 0$, $f\ne 0$,  $i\ne 0$ and (\ref{eq96}) holds, or 
\begin{equation}\label{eq100}
\begin{array}{ll}
u&=x^a\\
v&=x^dy^ez^f\\
w_i&=x^gy^hz^i\gamma
\end{array}
\end{equation}
with $e$ and $f\ne 0$,  $h$ or $i\ne 0$ and (\ref{eq96}) holds. 

We will now prove this statement.

Let $\nu$ be any valuation of ${\bold k}(X)$ which has center
 $p_{n}$ on $X_n$. 
 We identify $\nu$ with an extension of $\nu$ to the quotient field of $\hat{\cal O}_{X_n,p_n}$ which dominates $\hat{\cal O}_{X_n,p_n}$.
 
$q_n$
has permissible parameters (\ref{eq414}).

After possibly interchanging $u$ and $v$,  we have a relation (\ref{eq96}), so that
$\nu(v)>n\nu(u)$. We can reindex $x,y,z$ so that 
$$
0<\nu(x)\le\nu(y)\le\nu(z).
$$
 Then
$$
(f+e+d-nc)\nu(z)\ge (f-nc)\nu(z)+(e-nb)\nu(y)+(d-na)\nu(x)>0.
$$
If $c\ne 0$, and $n>f+e+d$, we have a contradiction.
Thus taking  $n>f+e+d$ for all $d,e,f$ in local forms (\ref{eq97}) for 3-points 
$p\in T(R)$,
we achieve that $c=0$ in all local forms (\ref{eq97})   which are the images of 3-points  $p_n\in T(R^{n})$ which
map to a point $q_{n}$ of $Y_{n}$ which is on a component $E$ of $D_{Y_n}$ such that $\Psi(E)$ is not a point. 

If $i=0$ (and $c=0$) in (\ref{eq97}) we have
$$
(h+g-nb)\nu(y)\ge (h-nb)\nu(y)+(g-na)\nu(x)>0
$$

so that if $b\ne 0$ and $n>h+g$ we have a contradiction.  Thus, by taking $n\gg 0$ in (\ref{eq98}), we see that if $b\ne 0$, then a form (\ref{eq99}) must hold at $p$ (since $uv=0$ is a local equation of $D_X$ at $p$ implies $f\ne 0$).
If $b=c=0$  in (\ref{eq97}), then a similar calculation shows that a form (\ref{eq100}) must hold at $p$ (for $n \gg0$).

We observe that in (\ref{eq99})   we have 
\begin{equation}\label{eq101}
(z)\cap\hat{\cal O}_{Y,q}=(v,w_i).
\end{equation}

Suppose that (\ref{eq100}) holds. If $i\ne 0$ then 
\begin{equation}\label{eq399}
(z)\cap\hat{\cal O}_{Y,q}=(v,w_i).
\end{equation}
If $h\ne 0$, then 
\begin{equation}\label{eq400}
(y)\cap \hat{\cal O}_{Y,q}=(v,w_i).
\end{equation}

Suppose that (\ref{eqT117}) or (\ref{eqT287}) hold. Then 
\begin{equation}\label{eqT297}
(y)\cap\hat{\cal O}_{Y,q}=(v,w_i).
\end{equation}

We will show that in (\ref{eq99}), $v=w_i=0$ is a formal branch of an algebraic curve $C$ in the fundamental locus of
$f:X\rightarrow Y$.
Let $R={\cal O}_{Y,q}$, $S={\cal O}_{X,p}$.

Since $p$ is a 3-point, there exist regular parameters $\overline x,\overline y,\overline z$ in ${\cal O}_{X,p}$
and units $\lambda_1,\lambda_2,\lambda_3\in\hat{\cal O}_{X,p}$ such that $\overline x=x\lambda_1$, $\overline y=y\lambda_2$, $\overline z=z\lambda_3$.

 $\overline z=0$ is a  local equation for
a component of $D_X$.  We have that $v\in (\overline z)\cap R$ and $u\not\in(\overline z)\cap R$
so that $(\overline z)\cap R=(v)$ or $(\overline z)\cap R=a$ where $a$ is a height two prime
containing $v$.
We have $(z\hat S)\cap\hat R=(v,w_i)$. Suppose that $(\overline z)\cap R=(v)$. We then have an induced morphism
$$
\hat R/(v)\rightarrow \hat S/(z)
$$
which is an inclusion by the Zariski Subspace Theorem (Theorem 10.14 \cite{Ab3}). This is impossible, so that $a$ is a height 2 prime in $R$, and defines a curve $C$, which is necessarily in the fundamental locus of $f$
since $\overline z=0$ is a local equation at $p$ of a component of $D_{X}$ which dominates $C$.
A similar argument shows that in (\ref{eq100}), (\ref{eqT117}) and (\ref{eqT287}), $v=w_i=0$ is a formal branch of an algebraic curve $C$ in the fundamental
locus of $f$.

\vskip .2truein
\noindent{\bf Step 2.} Let $C$ be the reduced curve in $Y$ whose components are the curves
in the fundamental locus of $f$ which are not 2-curves.
Let $\overline C$ be the reduced curve in $Y_{n}$ which is the strict transform of $C$.
The components of $\overline C$ are then in the fundamental locus of $f_n$.
 By Theorem \ref{Theorem6}, we can perform a sequence of blow ups
of prepared 1-points, prepared 2-points, 2-curves and resolving curves $\Psi':Y'\rightarrow Y_{n}$ so that we can construct a
$\tau$-quasi-well prepared  diagram of $\Psi'$ and $R^n$ 
\begin{equation}\label{eq105}
\begin{array}{rll}
X'&\stackrel{f'}{\rightarrow}&Y'\\
\Phi'\downarrow&&\downarrow\Psi'\\
X_{n}&\rightarrow&Y_{n}
\end{array}
\end{equation}
where $R'$ is the transform of $R^n$ on $X'$,
such that  the strict transform $\tilde C$ of $\overline C$ on $Y'$ is  nonsingular,  and makes SNCs with $D_{Y'}$.
If $q'\in U(R_i')\cap \tilde C$ for some $i$ then
 the germ at $q'$ of $\tilde C$ is contained in $S_{R_i'}(q')$, and
  the (disjoint) components of $\tilde C$ are permissible centers for $R'$.

Let $\Psi(1):Y(1)\rightarrow Y'$ be the blow up of  
$\tilde C$.  By Theorem \ref{Theorem6}, 
we have a $\tau$-quasi-well prepared  diagram of $\Psi(1)$ and $R'$ 
\begin{equation}\label{eq113}
\begin{array}{rll}
X(1)&\stackrel{f(1)}{\rightarrow}& Y(1)\\
\Phi(1)\downarrow&& \downarrow\Psi(1)\\
X'&\stackrel{f'}{\rightarrow}&Y'.
\end{array}
\end{equation}

Let $R(1)$ be the transform of $R'$ on $X(1)$.

Suppose that $q\in U(R)\subset Y$ is a 2-point. Suppose that $\tilde q\in(\Psi\circ\Psi')^{-1}(q)$ and $(f')^{-1}(\tilde q)\cap T(R_i')\ne\emptyset$. Then $q\in A_0$.

Let
$$
\tilde u=u_{\overline R'_i(\tilde q)}, \tilde v=v_{\overline R_i'(\tilde q)}, \tilde w_i=w_{\overline R_i'(\tilde q)},
$$
$$
u=u_{\overline R(q)}, v=v_{\overline R(q)}, w_i=w_{\overline R_i(q)}.
$$

We have (since $q\in A_0$ and $\Psi\circ\Psi'$ is a sequence of blow ups of points and 2-curves at $\tilde q$) one of the following forms (by (\ref{eqT267}) - (\ref{eqT269}): 
\begin{equation}\label{eqT281}
u=\tilde u^a\tilde v^b,
v=\tilde u^c\tilde v^d,
w_i=\tilde u^l\tilde v^m\tilde w_i
\end{equation}
with $ad-bc=\pm 1$ and $l>0$ if $c>0$, $m>0$ if $b>0$, or

\begin{equation}\label{eqT282}
u=(\tilde u^a\tilde v^b)^t\gamma_1(\tilde u,\tilde v),
v=(\tilde u^a\tilde v^b)^k\gamma_2(\tilde u,\tilde v),
w_i=\tilde u^e\tilde v^f\tilde w_i\gamma_3(\tilde u,\tilde v)
\end{equation}
where $\gamma_1,\gamma_2,\gamma_3$ are unit series, $a,b,t,k,e,f>0$, $\text{gcd}(a,b)=1$, or 
\begin{equation}\label{eqT283}
u=\tilde u^a\gamma_1(\tilde u,\tilde v),
v=\tilde u^b\gamma_2(\tilde u,\tilde v),
w_i=\tilde u^c\tilde w_i\gamma_3(\tilde u,\tilde v)
\end{equation}
where $\gamma_1,\gamma_2,\gamma_3$ are unit series, $a,b,c>0$.

We have that $(\Psi\circ\Psi')(E)=q$ for all components $E$ of $D_{Y'}$ containing $\tilde q$ 
(which implies $\tilde q\not\in\tilde C$) unless $\tilde q$ is a 2-point, and we have an expression 
\begin{equation}\label{eqT270}
\begin{array}{ll}
u&=\tilde u\\
v&=\tilde u^c\tilde v\\
w_i&=\tilde u^l\tilde v^m\tilde w_i
\end{array}
\end{equation}
with $c,l>0$ or
$$
\begin{array}{ll}
u&=\tilde u\tilde v^b\\
v&=\tilde v\\
w_i&=\tilde u^l\tilde v^m\tilde w_i
\end{array}
$$
with $b,m>0$.

Let $q^*=\Psi'(\tilde q)$. $q^*$ is a 2-point and $f^{-1}(q^*)\cap T(R_i')\ne\emptyset$.
After possibly interchanging $u$ and $v$ we have that a form (\ref{eq96}) holds at $q$, and thus by (\ref{eq101}) - (\ref{eqT297}) that $v=w_i=0$ are local equations of a formal component of $C$ at $q$. (\ref{eqT270}) thus holds at $\tilde q$, and since $\Psi'$ is a sequence of blow ups of points and 2-curves at $\tilde q$,
 $m=0$ in (\ref{eqT270}). We thus have an expression

\begin{equation}\label{eq342}
u=\tilde u, v=\tilde u^e\tilde v, w_i=\tilde u^f\tilde w_i
\end{equation}
for some $e,f>0$. $\tilde v=0$ is a local equation of the strict transform of $D_Y$ at $\tilde q$, and
$\tilde v=\tilde w_i=0$ are local equations of  $\tilde C$ at $\tilde q$ since $\tilde C$ is nonsingular. 

Suppose that $\overline q\in (\Psi\circ\Psi'\circ\Psi(1))^{-1}(q)$ and $f(1)^{-1}(\overline q)\cap T(R_i(1))\ne\emptyset$. Let $\tilde q=\Psi(1)(\overline q)$. then $(f')^{-1}(\tilde q)\cap T(R_i')\ne\emptyset$.
If $\tilde q\not\in\tilde C$ (so that $\overline q=\tilde q$) then we have seen that $(\Psi\circ\Psi')(E)=q$ for all components $E$ of $D_{Y'}$ containing $\tilde q$.
Suppose that $\tilde q\in\tilde C$. Then an expression (\ref{eq342}) holds at $\tilde q$.

$\Psi(1)$ is the blow up of $\tilde v=\tilde w_i=0$
above $\tilde q$. 
Since $\overline q\in U(R_i(1))$, we must have
$$
\tilde u=u_{\overline R_i(1)(\overline q)}, \tilde v=v_{\overline R_i(1)(\overline q)},
\tilde w_i=v_{\overline R_i(1)(\overline q)}w_{\overline R_i(1)(\overline q)}.
$$

Substituting into (\ref{eq342}), we have 
\begin{equation}\label{a*}
u=u_{\overline R_i(1)(\overline q)},
v=u_{\overline R_i(1)(\overline q)}^ev_{\overline R_i(1)(\overline q)},
w_i=u_{\overline R_i(1)(\overline q)}^fv_{\overline R_i(1)(\overline q)}w_{\overline R_i(1)(\overline q)}
\end{equation}
with $e,f\ge 1$.

Suppose that $q\in U(R)\subset Y$ is a 1-point.
Then $\Psi$ is an isomorphism over $q$.  Suppose that $\tilde q\in(\Psi\circ\Psi')^{-1}(q)$ and $(f')^{-1}(\tilde q)\cap T(R_i')\ne\emptyset$.
Let
$$
\tilde u=u_{\overline R_i'(\tilde q)},
\tilde v=v_{\overline R_i'(\tilde q)},
\tilde w_i=w_{\overline R_i'(\tilde q)},
$$
$$
u=u_{\overline R_i(q)},
v=v_{\overline R_i(q)},
w_i=w_{\overline R_i(q)}.
$$
We have that $u=w_i=0$ are local equations of $C$ 
since $T(R_i)\cap f^{-1}(q)\ne\emptyset$, by Remark \ref{RemarkT119}.
Then $\Psi'$ is either an isomorphism at $\tilde q$, or factors at $\tilde q$ as the blow up of $q$, followed by a sequence of blow ups of 2-points and 2-curves.

First suppose that $\Psi'$ is not an isomorphism at $\tilde q$.
We have one of the following forms (by (\ref{eqT267}) - (\ref{eqT269}): 
\begin{equation}\label{eqT154}
u=\tilde u^a\tilde v^b,
v=\tilde u^c\tilde v^d,
w_i=\tilde u^l\tilde v^m\tilde w_i
\end{equation}
with $ad-bc=\pm 1$ and $l>0$ if $c>0$, $m>0$ if $b>0$, or

\begin{equation}\label{eqT333}
u=(\tilde u^a\tilde v^b)^t\gamma_1(\tilde u,\tilde v),
v=(\tilde u^a\tilde v^b)^k\gamma_2(\tilde u,\tilde v),
w_i=\tilde u^e\tilde v^f\tilde w_i\gamma_3(\tilde u,\tilde v)
\end{equation}
where $\gamma_1,\gamma_2,\gamma_3$ are unit series, $a,b,t,k,e,f>0$, $\text{gcd}(a,b)=1$, or 
\begin{equation}\label{eqT334}
u=\tilde u^a\gamma_1(\tilde u,\tilde v),
v=\tilde u^b\gamma_2(\tilde u,\tilde v),
w_i=\tilde u^c\tilde w_i\gamma_3(\tilde u,\tilde v)
\end{equation}
where $\gamma_1,\gamma_2,\gamma_3$ are unit series, $a,b,c>0$.

We have that $\Psi'(E)=q$ for all components $E$ of $D_{Y'}$ containing $\tilde q$, 
(which implies $\tilde q\not\in \tilde C$)
unless
we have an expression 
\begin{equation}\label{eqT276}
\begin{array}{ll}
u&=\tilde u\tilde v^b\\
v&=\tilde v\\
w_i&=\tilde u^l\tilde v^m\tilde w_i
\end{array}
\end{equation}
with $b,m>0$.

Since $\Psi'$ is an isomorphism over a generic point of every component of $C$, we have that 
 $\tilde q$ is a 2-point and

\begin{equation}\label{eqT118}
u=\tilde u\tilde v^m,
v=\tilde v, w_i=\tilde v^n\tilde w_i
\end{equation}
with $m,n\ge 1$. $\tilde u=\tilde w_i=0$ are local equations of $\tilde C$ at $\tilde q$.

Suppose that $\overline q\in (\Psi\circ\Psi'\circ\Psi(1))^{-1}(q)$ and $f(1)^{-1}(\overline q)\cap T(R_i(1))\ne\emptyset$. Let $\tilde q=\Psi(1)(\overline q)$. then $(f')^{-1}(\tilde q)\cap T(R_i')\ne\emptyset$.
If $\tilde q\not\in\tilde C$ (so that $\overline q=\tilde q$) then we have seen that $(\Psi\circ\Psi')(E)=q$ for all components $E$ of $D_{Y'}$ containing $\tilde q$.
Suppose that $\tilde q\in\tilde C$. Then an expression (\ref{eqT118}) holds at $\tilde q$.
Then $\overline q$ is a 2-point with
$$
\tilde u=u_{\overline R_i(1)(\overline q)},
\tilde v=v_{\overline R_i(1)(\overline q)},
\tilde w_i=u_{\overline R_i(1)(\overline q)}w_{\overline R_i(1)(\overline q)},
$$
and thus 
\begin{equation}\label{eqT120}
u=u_{\overline R_i(1)(\overline q)}v_{\overline R_i(1)(\overline q)}^m,
v=v_{\overline R_i(1)(\overline q)},
w_i=u_{\overline R_i(1)(\overline q)}v_{\overline R_i(1)(\overline q)}^nw_{\overline R_i(1)(\overline q)}
\end{equation}
with $m,n\ge 1$.

Now suppose that $\Psi'$ is an isomorphism over $q$ ($q\in U(R)$ is a 1-point), and $\overline q\in\Psi(1)^{-1}(q)$ is such that
$f(1)^{-1}(q)\cap T(R_i(1))\ne\emptyset$. then $\overline q$ is a 1-point, and (by Remark \ref{RemarkT119}) 
\begin{equation}\label{eqT121}
u=u_{\overline R_i(1)(\overline q)},
v=v_{\overline R_i(1)(\overline q)},
w_i=u_{\overline R_i(1)(\overline q)}w_{\overline R_i(1)(\overline q)}.
\end{equation}

\vskip .2truein
\noindent{\bf Step 3.}
We now apply steps 1 and 2 of the proof to $f(1):X(1)\rightarrow Y(1)$ and $R(1)$. We construct a $\tau$-quasi-well prepared  diagram
$$
\begin{array}{rll}
X(2)&\stackrel{f(2)}{\rightarrow}&Y(2)\\
\overline \Phi(2)\downarrow&&\downarrow\overline\Psi(2)\\
X(1)&\stackrel{f(1)}{\rightarrow} &Y(1),
\end{array}
$$
where $R(2)$ is the transform of $R(1)$ on $X(2)$. Suppose  that $q_2\in U(R_i(2))\subset Y(2)$ is on a component $E_2$ of $D_{Y(2)}$ such that $\Psi\circ\Psi'\circ\Psi(1)\circ\overline\Psi(2)(E_2)$ is not a point of $Y$ and there exists a point $p_2\in f(2)^{-1}(q_2)\cap T(R_i(2))\subset X(2)$. Let $q_1=\overline\Psi(2)(q_2)$. $\overline\Psi_2(2)(E_2)$ is necessarily not a point of $Y(1)$.

Suppose that $q_1$ is a 2-point. Then we have an expression (analogous to (\ref{a*})): 
\begin{equation}\label{b*}
\begin{array}{ll}
u_{\overline R_i(1)(q_1)}&=u_{\overline R_i(2)(q_2)}\\
v_{\overline R_i(1)(q_1)}&=u_{\overline R_i(2)(q_2)}^ev_{\overline R_i(2)(q_2)}\\
w_{\overline R_i(1)(q_1)}&=u_{\overline R_i(2)(q_2)}^fv_{\overline R_i(2)(q_2)}w_{\overline R_i(2)(q_2)}
\end{array}
\end{equation}
or
\begin{equation}\label{c*}
\begin{array}{ll}
u_{\overline R_i(1)(q_1)}&=u_{\overline R_i(2)(q_2)}v_{\overline R_i(2)(q_2)}^e\\
v_{\overline R_i(1)(q_1)}&=v_{\overline R_i(2)(q_2)}\\
w_{\overline R_i(1)(q_1)}&=u_{\overline R_i(2)(q_2)}v_{\overline R_i(2)(q_2)}^fw_{\overline R_i(2)(q_2)}
\end{array}
\end{equation}
with $e,f\ge 1$.

Let $q=(\Psi\circ\Psi'\circ\Psi(1))(q_1)$. Since $q_1$ lies on a component $E_1$ of $D_{Y(1)}$ which does not contract to a point of $Y$, and $f(1)^{-1}(q_1)\cap T(R_i(1))\ne\emptyset$, a form (\ref{a*}) or (\ref{eqT120}) holds for
$$
u=u_{\overline R_i(q)},
v=v_{\overline R_i(q)},
w_i=w_{\overline R_i(q)}
$$
and
$$
u_{\overline R_i(1)(q_1)},
v_{\overline R_i(1)(q_1)},
w_{\overline R_i(1)(q_1)}.
$$

Suppose that (\ref{a*}) holds at $q_1$.
Substituting (\ref{b*}) and (\ref{c*}) into (\ref{a*}), we see that since $q_2$ is on a component $E_2$ of $D_{Y(2)}$ which does not contract to $q$, then we have that (\ref{b*}) holds, and an expression 
\begin{equation}\label{d*}
\begin{array}{ll}
u=u_{R_i(q)}&=u_{\overline R_i(2)(q_2)}\\
v=v_{R_i(q)}&=u_{\overline R_i(2)(q_2)}^{e_2}v_{\overline R_i(2)(q_2)}\\
w=w_{R_i(q)}&=u_{\overline R_i(2)(q_2)}^{f_2}v_{\overline R_i(2)(q_2)}^2w_{\overline R_i(2)(q_2)}
\end{array}
\end{equation}
with $e_2,f_2\ge 2$.

Suppose that  (\ref{eqT120}) holds at $q_1$. Substituting (\ref{b*}) or (\ref{c*}) into (\ref{eqT120}), we see that since $q_2$ is on a component  of $D_{Y(2)}$ which does not contract to $q$,  we have that (\ref{c*}) holds we have and an expression 
\begin{equation}\label{eqT155}
u=u_{\overline R_i(2)(q_2)}v_{\overline R_i(2)(q_2)}^{e_2},
v=v_{\overline R_i(2)(q_2)},
w_i=u_{\overline R_i(2)(q_2)}^2v_{\overline R_i(2)(q_2)}^{f_2}w_{\overline R_i(2)(q_2)}
\end{equation}
with $e_2,f_2\ge 2$.

Suppose that $q_1$ is a 1-point. Then we have an expression (analogous to (\ref{eqT120})) 
\begin{equation}\label{eqT156}
\begin{array}{ll}
u_{\overline R_i(1)(q_1)}&=u_{\overline R_i(2)(q_2)}v_{\overline R_i(2)(q_2)}^m,\\
v_{\overline R_i(1)(q_1)}&=v_{\overline R_i(2)(q_2)},\\
w_{\overline R_i(1)(q_1)}&=u_{\overline R_i(2)(q_2)}v_{\overline R_i(2)(q_2)}^nw_{\overline R_i(2)(q_2)}
\end{array}
\end{equation}
with $m,n\ge 1$, or (analogous to (\ref{eqT121})) 
\begin{equation}\label{eqT157}
u_{\overline R_i(1)(q_1)}=u_{\overline R_i(2)(q_2)},
v_{\overline R_i(1)(q_1)}=v_{\overline R_i(2)(q_2)},
w_{\overline R_i(1)(q_1)}=u_{\overline R_i(2)(q_2)}w_{\overline R_i(2)(q_2)}.
\end{equation}

Since $q_1$  lies on a component $E_1$ of $D_{Y(1)}$ which does not contract to a point of $Y$, and $f(1)^{-1}(q_1)\cap T(R_i(1))\ne\emptyset$,
a
form (\ref{eqT121}) holds at $q_1$ for
$$
u=u_{\overline R_i(q)}, v=v_{\overline R_i(q)},
w_i=w_{\overline R_i}(q)
$$
and
$$
u_{\overline R_1(1)(q_1)}, v_{\overline R_i(1)(q_1)}, w_{\overline R_i(1)(q_1)}.
$$

Substituting (\ref{eqT156}) and (\ref{eqT157}) into (\ref{eqT121}), we see that since $q_1$ is on a component $E_1$ of $D_{Y(1)}$ which does not contract to $q$,  we have 
\begin{equation}\label{eqT271}
u=u_{\overline R_i(2)(q_2)}v_{\overline R_i(2)(q_2)}^m,
v=v_{\overline R_i(2)(q_2)},
w_i=u_{\overline R_i(2)(q_2)}^2v_{\overline R_i(2)(q_2)}^nw_{\overline R_i(2)(q_2)}
\end{equation}
with $m,n\ge 0$.

Iterating steps 1 and 2, we construct a sequence of $\tau$-quasi-well prepared diagrams
\begin{equation}\label{g*}
\begin{array}{rll}
\vdots&&\vdots\\
\downarrow&&\downarrow\\
X(n)&\stackrel{f(n)}{\rightarrow}&Y(n)\\
\overline \Phi(n)\downarrow&&\downarrow\overline\Psi(n)\\
X(n-1)&
\stackrel{f(n-1)}{\rightarrow} &Y(n-1)\\
\downarrow&&\downarrow\\
\vdots&&\vdots\\
\downarrow&&\downarrow\\
X(1)&\stackrel{f(1)}{\rightarrow}&Y(1)\\
\overline \Phi(1)\downarrow&&\downarrow\overline\Psi(1)\\
X&\stackrel{f}{\rightarrow} &Y.
\end{array}
\end{equation}

We continue this algorithm  as long as there exists $q_n\in U(R_i(n))$ for some $i$ such that $q_n$
is on a component $E$ of $D_{Y(n)}$ which does not contract to a point of $Y$, and $f(n)^{-1}(q_n)\cap T(R_i(n))\ne\emptyset$.
\vskip .2truein
\noindent{\bf Step 4.}
Suppose that the algorithm never terminates. 

Let $\nu$ be a 0-dimensional valuation of ${\bf k}(X)$. We will say that $\nu$ is resolved on $X(n)$ if the center of $\nu$ on $X(n)$ is at a point $p_n$ of $X(n)$ such that either $p_n\not\in T(R(n))$ or $p_n\in T(R(n))$ and all components $E$ of $D_{Y(n)}$ containing $q_n=f(n)(p_n)$ contract to a point of $Y$.

By our construction, if  $\nu$ is resolved on $X(n)$, then $\nu$ is resolved on $X(m)$ for all $m\ge n$. Further, the set of $\nu$ in the Zariski-Riemann manifold $\Omega(X)$ of $X$ \cite{Z} which are resolved on $X$ is an open subset of $\Omega(X)$.

Suppose that  $\nu$ is a 0-dimensional valuation of ${\bf k}(X)$ such that $\nu$ is not resolved on $X(n)$ for all $n$. 
Let $p_n$ be the center of $\nu$ on $X(n)$, $q_n$ be the center of $\nu$ on $Y(n)$.

For all $n$, we identify $\nu$ with an extension of $\nu$ to the quotient field of $\hat{\cal O}_{X_n,p_n}$ which dominates $\hat{\cal O}_{X_n,p_n}$.

First suppose that the center of $\nu$ on $Y$ is a 2-point.

 There exists an $i$ such that for all $n$,  $q_n\in U(R_i(n))$ and $p_n\in f(n)^{-1}(q_n)\cap T(R_i(n))$. 
We have expressions (after possibly interchanging $u$ and $v$) 
\begin{equation}\label{e*}
\begin{array}{ll}
u=u_{\overline R_i(q)}&=u_{\overline R_i(n)(q_n)}\\
v=v_{\overline R_i(q)}&=u_{\overline R_i(n)(q_n)}^{e_n}v_{\overline R_i(n)(q_n)}\\
w=w_{\overline R_i(q)}&=u_{\overline R_i(n)(q_n)}^{f_n}v_{\overline R_i(n)(q_n)}^nw_{\overline R_i(n)(q_n)}
\end{array}
\end{equation}
with $e_n,f_n\ge n$ for all $n$.

 From (\ref{e*}), we see that 
 
 \begin{equation}\label{eqT289}
 \nu(v)>n\nu(u)\text{ for all }n\in{\bf N}.
 \end{equation}
  Thus $\nu$ is a composite valuation, and there exists a prime ideal $P$ of the valuation ring $V$ of
$\nu$ such that $v\in P$, $u\not\in P$. Let $\nu_1$ be a valuation whose valuation ring is $V_P$.
We have $\nu_1(u)=0$, $\nu_1(v)>0$. From (\ref{e*}) we see that 
\begin{equation}\label{f*}
\nu_1(w)>n\nu_1(v)>0
\end{equation}
for all $n\in{\bf N}$.

At $p=p_0\in X$, we have a form (\ref{eqT117}), (\ref{eqT287}), (\ref{eq99}) or (\ref{eq100}).
In (\ref{eqT117}) we have $\nu_1(x)=0$ and $d>0$, a contradiction to (\ref{f*}).
In (\ref{eq99}) we have $\nu_1(y)=0$.    $uv=0$ is a local equation of $D_X$ at $p$. Thus either $a>0$ or $d>0$. If $a>0$ then $\nu_1(x)=0$ and $\nu_1(z)>0$, a contradiction to (\ref{f*}) since $f\ne 0$. If $d>0$, we again have a contradiction to (\ref{f*}).
In (\ref{eq100}) we have   a contradiction to (\ref{f*}), since $\nu_1(x)=0$ and $e,f>0$.

Suppose that (\ref{eqT287}) holds. In this case a more detailed analysis is required.
By our construction and with the notation of steps 1 and 2, there exists a factorization
$$
\begin{array}{rll}
X(1)&\stackrel{f(1)}{\rightarrow}&Y(1)\\
\Phi(1)\downarrow&&\downarrow\Psi(1)\\
X'&\stackrel{f'}{\rightarrow}&Y'\\
\Phi'\downarrow&&\downarrow\Psi'\\
X_n&\stackrel{f_n}{\rightarrow}&Y_n\\
\Phi\downarrow&&\downarrow\Psi\\
X&\stackrel{f}{\rightarrow}&Y.
\end{array}
$$

Recall that
$p_1$ is the center of $\nu$ on $X(1)$, $p$ is the center of $\nu$ on $X$,
$q_1$ is the center of $\nu$ on $Y(1)$, and  $q$ is the center of $\nu$ on $Y$.

Let $p'$ be the center of $\nu$ on $X'$, $q'$ be the center of $\nu$ on $Y'$.

At $q'$, $Y'\rightarrow Y$ is a sequence of blow ups of prepared 2-points of type 1, 
and 2-curves, and at $p'$,
$$
\begin{array}{rll}
X'&\rightarrow &Y'\\
\downarrow&&\downarrow\\
X&\rightarrow&Y
\end{array}
$$
is obtained by iterating the constructions of Remark \ref{Remark424} and Lemma \ref{Lemma32}.

Let
$$
u'=u_{\overline R_i'(q')},
v'=v_{\overline R_i'(q')},
w'_i=w_{\overline R_i'(q')}.
$$

By equations (\ref{eqT289}), (\ref{f*}) and (\ref{eqT287}) we see that 
\begin{equation}\label{eqT298}
\nu(y)>n\nu(x) \text{ and }\nu_1(x)=0,\,\,\, \nu_1(z)>n\nu_1(y)
\end{equation}
  for all $n\in{\bf N}$.

We see (by a variant of the analysis of (\ref{eqT287}), using (\ref{eqT298})), that there exist permissible parameters $x',y',z'$ at $p'$ such that
$$
u'=(x')^a,
v'=(x')^{c'}(y')^d,
w_i'=(x')^{e'}(y')^fz'
$$
where $a,d,f$ are the constants of (\ref{eqT287}) and $c',e'\in{\bf N}$.
$\Psi(1)$ is the blow up of the curve $\tilde C$ which has local equations $v'=w_i'=0$ at $q'$.
The construction of $X(1)\rightarrow X'$ (Remark \ref{Remark293}) is analogous to
the analysis of Lemma \ref{LemmaT79}. There exists a factorization
$$
X(1)=W_m\rightarrow \cdots\rightarrow W_2\rightarrow W_1\rightarrow X'
$$
where $W_1\rightarrow X'$ is a sequence of blow ups of 2-curves and 3-points, and each $W_{i+1}\rightarrow W_i$ is a curve containing a 1-point in the locus where 
${\cal I}_{\tilde C}{\cal O}_{W_i}$ is not invertible. 

Let $\tilde p_i$ be the center of $\nu$ on $W_i$. By the analysis of the proof of Lemma \ref{LemmaT79}, and (\ref{eqT298}), we see that there exist permissible parameters $\tilde x_1,\tilde y_1,\tilde z_1$ at $\tilde p_1$ such that 
$$
u'=\tilde x_1^a,
v'=\tilde x_1^{c_1}\tilde y_1^d,w_i'=\tilde x_1^{e_1}\tilde y_1^f\tilde z_1
$$
with $(c_1,d)\le(e_1,f)$, or $d\ge f$ and $(c_1,d)\not\le (e_1,f)$.

If $(c_1,d)\le(e_1,f)$, then $X(1)\rightarrow W_1$ is an isomorphism at $p_1$, and we have 
$$
u'=u_{\overline R_i(1)(q_1)},
v'=v_{\overline R_i(1)(q_1)},
w_i'=w_{\overline R_i(1)(q_1)}v_{\overline R_i(1)(q_1)}.
$$
We have 
\begin{equation}\label{eqT299}
u_{\overline R_i(1)(q_1)}=\tilde x_1^a,
v_{\overline R_i(1)(q_1)}=\tilde x_1^{c_1}\tilde y_1^d,
w_{\overline R_i(1)(q_1)}=\tilde x_1^{e_1-c_1}\tilde y_1^{f-d}\tilde z_1.
\end{equation}

Suppose that $d\ge f$ and $(c,d)\not\le(e_1,f)$. We may assume without loss of generality that each $W_{i+1}\rightarrow W_i$ is not an isomorphism at the center of $\nu$.

$W_2\rightarrow W_1$ is  the blow up of a curve which either has local equations 
\begin{equation}\label{eqT290}
\tilde x_1=\tilde z_1=0
\end{equation}
(and $c_1>e_1$), or 
\begin{equation}\label{eqT291}
\tilde y_1=\tilde z_1=0
\end{equation}
(and $d>f$).

By (\ref{eqT298}), there exist permissible parameters $\tilde x_2,\tilde y_2,\tilde z_2$ at $\tilde p_2$ such that
$$
\tilde x_1=\tilde x_2, \tilde y_1=\tilde y_2, \tilde z_1=\tilde x_2\tilde z_2
$$
if (\ref{eqT290}) holds,
$$
\tilde x_1=\tilde x_2, \tilde y_1=\tilde y_2, \tilde z_1=\tilde y_2\tilde z_2
$$
if (\ref{eqT291}) holds.

We then have that
$$
u'=\tilde x_2^a, v'=\tilde x_2^{c_1}\tilde y_2^d,
w_i'=\tilde x_2^{e_1+1}\tilde y_2^f\tilde z_1
$$
or
$$
u'=\tilde x_2^a, v'=\tilde x_2^{c_1}\tilde y_2^d, w_i'=\tilde x_2^{e_1}\tilde y_2^{f+1}\tilde z_1.
$$
By iteration of this analysis for local equations of $W_{i+1}\rightarrow W_i$, we see that at $p_1=\tilde p_m$, we have permissible parameters $\tilde x_m,\tilde y_m,\tilde z_m$ such that 
$$
u'=\tilde x_m^a, v'=\tilde x_m^{c_1}\tilde y_m^d,
w_i'=\tilde x_m^{e_1}\tilde y_m^d\tilde z_m
$$
with $e_1\ge c_1$.

We have

\begin{equation}\label{eqT292}
u_{\overline R_i(1)(q_1)}=\tilde x_m^a,
v_{\overline R_i(1)(q_1)}=\tilde x_m^{c_1}\tilde y_m^d,
w_{\overline R_i(1)(q_1)}=\tilde x_m^{e_1-c_1}\tilde z_m.
\end{equation}

We see by the analysis of (\ref{eqT287}) that $\nu$ is resolved on $X(2)$ if (\ref{eqT292}) holds, so we must have that (\ref{eqT299}) holds at $\tilde p_1$. Observe that we must have a reduction $f_1=f-d<f$ in (\ref{eqT299}) from (\ref{eqT287}).

iterating this analysis, we see that we must reach the case (\ref{eqT292}) after a finite number of iterations of step 2. This is a contradiction to the assumption that $\nu$ is never resolved on $X(n)$.

Now suppose that the center of $\nu$ on $Y$ is a 1-point. Then there exists an $i$ such that for all $n$, $q_n\in U(R_i(n))$ and $p_n\in f(n)^{-1}(q_n)\cap T(R_i(n))$, and either $q_n$ is a 2-point for $n>>0$ or $q_n$ is a 1-point for all $n$.

$q=q_0\in Y$ is a 1-point and $p=p_0\in f^{-1}(q)\cap T(\overline R_i)$.  Let
$$
u=u_{\overline R_i(q)},
v=v_{\overline R_i(q)},
w=w_{\overline R_i(q)}.
$$
If $\tau>0$, there exist permissible parameters $x,y,z$ at $p$ such that one of the following forms hold:
$p$ a 1-point 
\begin{equation}\label{eqT122}
u=x^a,v=y,
w_i=x^b(\gamma(x,y)+x^{c-b}z)
\end{equation}
where $\gamma$ is a unit series
or $p$ a 2-point 
\begin{equation}\label{eqT123}
u=(x^ay^b)^k,
v=z,
w_i=(x^ay^b)^t(\gamma(x^ay^b,z)+x^cy^d)
\end{equation}
where $\gamma$ is a unit series.

If $\tau=0$, then there exist permissible parameters $x,y,z$ at $p$ such that either $p$ is a 1-point and 
\begin{equation}\label{eqT300}
u=x^a, v=y, w_i=x^b(\beta+z)
\end{equation}
with $\beta\in{\bf k}$, $b>0$, or $p$ is a 2-point and 
\begin{equation}\label{eqT301}
u=(x^ay^b)^k, v=z, w_i=x^cy^d
\end{equation}
with $ad-bc\ne0$.

Suppose that there exists $n_0$ such that $q_n$ is a 2-point for all $n\ge n_0$.  Suppose that $n\ge n_0$. Then by comparing (\ref{eqT121}), (\ref{eqT120}) and (\ref{b*}) or (\ref{c*}), we see that (\ref{b*}) cannot occur (since we assume some component of $D_{Y(n)}$ containing $q_n$ does not contract to $q$).
We have an expression:
$$
\begin{array}{ll}
u&=u_{\overline R_i(n)}v_{\overline R_i(n)}^{e_n},\\
v&=v_{\overline R_i(n)},\\
w_i&=u_{\overline R_i(n)}^{f_n}v_{\overline R_i(n)}^{g_n}w_{\overline R_i(n)},
\end{array}
$$
and  $e_n,f_n,g_n\ge n-n_0$.

We have $\nu(u)>n\nu(v)>0$ for all $n\in{\bf N}$. Thus $\nu$ is a composite valuation, and there exists a prime ideal $P$ of the valuation ring $V$ of $\nu$ such that $u\in P$, $v\not\in P$. Let $\nu_1$ be a valuation whose valuation ring is $V_P$.
We have $\nu_1(v)=0$, $\nu_1(u)>0$. We further have 
\begin{equation}\label{eqT273}
\nu_1(w_i)>n\nu_1(u)>0
\end{equation}
 for all $n\in{\bf N}$.

Suppose that $\tau>0$ (and $q_n$ is a 2-point for $n\ge n_0$).
At $p=p_0\in X$, we have a form (\ref{eqT122}) or (\ref{eqT123}). In (\ref{eqT122}) we have
$$
\nu_1(w_i)=\frac{b}{a}\nu_1(u)
$$
and in (\ref{eqT123}) we have
$$
\nu_1(w_i)=\frac{t}{k}\nu_1(u),
$$
a contradiction (to (\ref{eqT273}).

Suppose that $\tau=0$. If $\beta\ne 0$ in (\ref{eqT300}), then the analysis is the same as for the $\tau>0$ case, so we may assume that $\beta=0$ if (\ref{eqT300}) holds.

If (\ref{eqT300}) holds with $\beta=0$, or if (\ref{eqT301}) holds, we finish the analysis in a similar way to the proof when $q$ is a 2-point, and $\tau=0$, given above.

The final case is when $q_n$ is a 1-point for all $n$. From (\ref{eqT121}) (and Remark \ref{RemarkT119}) we see that
$$
u=u_{\overline R_i(\overline q)}=u_{\overline R_i(n)(q_n)},
v=v_{\overline R_i(\overline q)}=v_{\overline R_i(n)(q_n)},
w_i=w_{\overline R_i(\overline q)}=u_{\overline R_i(n)(q_n)}^{n}w_{\overline R_i(n)(q_n)},
$$
so that 
\begin{equation}\label{eqT272}
\nu(w_i)>n\nu(u)>0
\end{equation} for all $n\in{\bf N}$.

Suppose that $\tau>0$ (and $q_n$ is a 1-point for all $n$).
At $p=p_0\in X$ we have a form (\ref{eqT122}) or (\ref{eqT123}), which implies
$$
\nu(w_i)=\frac{b}{a}\nu(u)
$$
or
$$
\nu(w_i)=\frac{t}{k}\nu(u),
$$
a contradiction to (\ref{eqT272}).

When $\tau=0$, the proof follows in a similar way to the proof when $q$ is a 2-point and $\tau=0$.

 We have shown that for all 0-dimensional valuations $\nu$ of ${\bf k}(X)$, there exists $n$ such that $\nu$ is resolved on $X(n)$.

By compactness of the Zariski-Riemann manifold \cite{Z} there exists $N$ such that  all $\nu\in \Omega(X)$ are resolved on $X(N)$, a contradiction to our assumption that (\ref{g*}) is of infinite length. The diagram
$$
\begin{array}{rll}
X(N)&\rightarrow &Y(N)\\
\downarrow&&\downarrow\\
X&\rightarrow&Y
\end{array}
$$
thus satisfies the  conclusions of Theorem \ref{TheoremT124}.

\end{pf}

\begin{Theorem}\label{Theorem80} Suppose that $\tau\ge0$, $f:X\rightarrow Y$
is $\tau$-quasi-well prepared with  relation $R$. 
 Then there exists a commutative diagram
$$
\begin{array}{rll}
X_1&\stackrel{f_1}{\rightarrow}&Y_1\\
\Phi\downarrow&&\downarrow\Psi\\
X&\stackrel{f}{\rightarrow}&Y
\end{array}
$$
such that $\Phi$ and $\Psi$ are products of blow ups of possible centers and 
$f_1$ is $\tau$-quasi-well prepared with   relation
$R^1$ and pre-algebraic structure. Further, $R$ has an algebraic structure. ($R^1$ will in general not be the transform of $R$.)
\end{Theorem}

\begin{pf}

By Theorem \ref{TheoremT124} there exists a $\tau$-quasi-well prepared diagram (\ref{eq336})
as in the hypothesis of Lemma \ref{Lemma348}. Then Lemma \ref{Lemma348} implies that the conclusions of Theorem \ref{Theorem80} hold.
\end{pf}

\begin{Lemma}\label{Lemma33} Suppose that $\tau\ge 0$, $f:X\rightarrow Y$ is $\tau$-quasi-well prepared 
with  relation $R$ and pre-algebraic structure (or $\tau$-well prepared with  relation $R$),
$q\in U(R)$ is a 2-point  and $p\in f^{-1}(q)\cap T(R_i)$ for some $i$.  Suppose that $E$ is a component of $D_Y$
containing $q$. Let $C=\overline{E\cdot S_{\overline R_i}(q)}$.

Let $\Psi_n:Y_n\rightarrow Y$ be obtained by blowing up $q$,  then blowing up the point $q_1$ which is the intersection of
the exceptional divisor over $q$ and  the strict transform of $C$ on $Y_1$, and iterating this procedure $n$ times, blowing up the intersection point of the last exceptional divisor with the strict transform of $C$. Let
$$
\begin{array}{rll}
X_n&\stackrel{f_n}{\rightarrow}&Y_n\\
\Phi_n\downarrow&&\downarrow\Psi_n\\
X&\stackrel{f}{\rightarrow}&Y
\end{array}
$$
be a $\tau$-quasi-well prepared (or $\tau$-well prepared) diagram of  $R$ and $\Psi_n$ obtained from Lemma \ref{Lemma32}
(so that $\Phi_n$ is an isomorphism above $f^{-1}(Y-\Sigma(Y)))$.

Suppose that for all $n>0$ there exists a point $p_n\in\Phi_n^{-1}(p)\cap
T(R_{i}^n)$ such that
$f_n(p_n)=q_n\in \Psi_n^{-1}(q)\cap C_n$, where $C_n$ is the strict transform of $C$ on $Y_n$. Then
$C$ is a component of the fundamental locus of $f$.
 \end{Lemma}

\begin{pf} The proof follows from Step 1 of the proof of Theorem \ref{TheoremT124}.
Let
$$
u=u_{\overline R_i(q)},
v=v_{\overline R_i(q)},
w_i=w_{\overline R_i(q)}.
$$
After possibly interchanging $u$ and $v$,  $v=w_i=0$ are local equations of $C$ at $q$.

 By our construction of $\Psi_n$,  we  have that 
$$
u_1=u_{\overline R_i^n(q_n)},
v_1=v_{\overline R_i^n(q_n)},
w_{i,1}=w_{\overline R_i^n(q_n)}
$$
are defined by
$$
u=u_1,
v=u_1^nv_1,
w_{i}=u_1^nw_{i,1}.
$$

We must have an expression (\ref{eqT117}), (\ref{eqT287}), (\ref{eq99}) or (\ref{eq100}). It follows,
as in the analysis of step 1 of Theorem \ref{TheoremT124} that $C$ is in the fundamental locus of $f$.

\end{pf}

\begin{Theorem}\label{TheoremT325} Suppose that $\tau\ge 0$ and  $f:X\rightarrow Y$ is
$\tau$-quasi-well
prepared with   relation $R$ and pre-algebraic structure.  Further suppose that $R$ has an algebraic structure.
 Then there exists a commutative
diagram
$$
\begin{array}{rll}
X_1&\stackrel{f_1}{\rightarrow}&Y_1\\
\Phi\downarrow&&\downarrow\Psi\\
X&\stackrel{f}{\rightarrow}&Y
\end{array}
$$
such that $\Phi$,$\Psi$ are products of blow ups of possible centers and $f_1$ is
$\tau$-well prepared with   relation $R^1$.
(In general, $R^1$ is not the transform of $R$). 
\end{Theorem}

\begin{pf}
1 and 2 of Definition \ref{Def65} of a $\tau$-well prepared relation $R$ hold by assumption. 

For $i\ne j$, Let $H_{ij}$ be the set of points $q\in U(R_i)\cap U(R_j)$ such that neither an expression (\ref{eq64}) nor (\ref{eqT100}) of 3 of Definition \ref{Def65} holds between $w_{\overline R_i(q)}$ and $w_{\overline R_j(q)}$. We see by 3 of Definition \ref{DefT60} (and (\ref{eqT267}) - (\ref{eqT269}))  that $H_{ij}$ is a finite set.

Let $\tilde\Psi_1:\tilde Y_1\rightarrow Y$ be the blow up of the union of the sets
$$
\{q\in H_{ij}\mid f^{-1}(H_{ij})\cap[T(R_i)\cup T(R_j)]\ne\emptyset\}.
$$
 By Lemmas \ref{Lemma32} and \ref{LemmaT63}, there exists a $\tau$-quasi-well prepared diagram
$$
\begin{array}{rll}
\tilde X_1&\stackrel{\tilde f_1}{\rightarrow}&\tilde Y_1\\
\tilde\Phi_1\downarrow&&\downarrow\tilde\Psi_1\\
X&\stackrel{f}{\rightarrow}&Y.
\end{array}
$$
We define  finite sets $H_{ij}^1$ in the same way for the transform $\tilde R^1$ of $R$, and iterate to construct a $\tau$-quasi-well prepared diagram 
\begin{equation}\label{eqT318}
\begin{array}{rll}
\tilde X_n&\stackrel{\tilde f_n}{\rightarrow}&\tilde Y_n\\
\tilde\Phi_n\downarrow&&\downarrow\tilde\Psi_n\\
\tilde X_{n-1}&\stackrel{\tilde f_{n-1}}{\rightarrow}&\tilde Y_{n-1}\\
\vdots&&\vdots\\
\downarrow&&\downarrow\\
\tilde X_1&\stackrel{\tilde f_1}{\rightarrow}&\tilde Y_1\\
\tilde\Phi_1\downarrow&&\downarrow\tilde\Psi_1\\
X&\stackrel{f}{\rightarrow}&Y,
\end{array}
\end{equation}
continuing as long as $\tilde f_n^{-1}(H_{ij}^n)\cap[T(R^n_i)\cup T(R^n_j)]\ne \emptyset$ for some $i\ne j$.

Suppose that (\ref{eqT318}) doesn't terminate after a finite number of blow ups $n$. Then there exist $i\ne j$ and a valuation $\nu$ of the function field ${\bf k}(X)$ of $X$ such that the center $p_n$ of $\nu$ on $\tilde X_n$ and the center $q_n$ of $\nu$ on $\tilde Y_n$ satisfy $q_n\in H_{ij}^n$, $p_n\in \tilde f_n^{-1}(H_{ij}^n)\cap T(\tilde R_i^n)$ for all $n$.

For all $n$, we may identify $\nu$ with an extension of $\nu$ to the quotient field of $\hat{\cal O}_{\tilde X_n,p_n}$ which dominates $\hat{\cal O}_{\tilde X_n,p_n}$.

Let
$$
u_n=u_{\overline R_i^n(q_n)},
v_n=v_{\overline R_i^n(q_n)},
w_{ni}=w_{\overline R_i^n(q_n)},
w_{nj}=w_{\overline R_j^n(q_n)}.
$$
We have 
\begin{equation}\label{eqT319}
u_n=u_{n+1}, v_n=u_{n+1}(v_{n+1}+\alpha_{n+1}),
w_{ni}=u_{n+1}w_{n+1,i},
w_{nj}=u_{n+1}w_{n+1,j}
\end{equation}
with $\alpha_{n+1}\in{\bf k}$, or 
\begin{equation}\label{eqT320}
u_n=u_{n+1}v_{n+1},
v_n=v_{n+1},
w_{ni}=v_{n+1}w_{n+1,i},
w_{nj}=v_{n+1}w_{n+1,j}
\end{equation}
for all $n$.

The relation
$$
w_{nj}=w_{ni}+\lambda_{ij}^n(u_n,v_n)
$$
of 3 of Definition \ref{DefT60} transforms as
$$
\lambda_{ij}^{n+1}(u_{n+1},v_{n+1})=\frac{\lambda_{ij}^n(u_{n+1},u_{n+1}(v_{n+1}+\alpha_{n+1}))}{u_{n+1}}
$$
if (\ref{eqT319}) holds, and transforms as
$$
\lambda_{ij}^{n+1}(u_{n+1},v_{n+1})=\frac{\lambda_{ij}^n(u_{n+1}v_{n+1},v_{n+1})}{v_{n+1}}
$$
if (\ref{eqT320}) holds.

Let ${\cal F}_n$ be the germ at $q_n$ of the divisor $\lambda_{ij}^n=0$. By embedded resolution of plane curve singularities, there exists $n_0$ such that $D_{\tilde X_n}+ {\cal F}_n$ is a SNC divisor at $q_n$ for all $n\ge n_0$.  If $q_n$ is a 2-point for some $n\ge n_0$, we have that $q_n\not\in H_{ij}^n$, a contradiction, so we have that $q_n$ is a 1-point for all $n\ge n_0$. Thus a form (\ref{eqT319}) holds for all $n\ge n_0$. We have
$$
w_{n_0i}=u_{n_0}^{n-n_0}w_{ni}.
$$
Thus 
\begin{equation}\label{eqT321}
\nu(w_{n_0i})>n\nu(u_{n_0})>0
\end{equation}
for all positive $n$.  

Suppose that $\tau>0$.
Since $q_{n_0}\in U(\tilde R_i^{n_0})$ is a 1-point, we have permissible parameters $x,y,z$ at $p_{n_0}$ such that 
\begin{equation}\label{eqT322}
u_{n_0}=x^a,
v_{n_0}=y,
w_{n_0i}=x^c\gamma
\end{equation}
where $\gamma\in\hat{\cal O}_{\tilde X_{n_0},p_{n_0}}$ is a unit series, or 
\begin{equation}\label{eqT323}
u_{n_0}=(x^ay^b)^k,
v_{n_0}=z,
w_{n_0,i}=(x^ay^b)^l\gamma
\end{equation}
where $\gamma\in\hat{\cal O}_{\tilde X_{n_0},p_{n_0}}$ is a unit series.
In either case, we have a contradiction to (\ref{eqT321}).

If $\tau=0$, then we have one of the following forms at $p_{n_0}$.
$$
u_{n_0}=x^a,
v_{n_0}=y,
w_{n_0i}=x^c(z+\beta)
$$
with $\beta\in{\bf k}$, or
$$
u_{n_0}=(x^ay^b)^k,
v_{n_0}=z,
w_{n_0,i}=x^cy^d
$$
with $a,b>0$ and $ad-bc\ne0$.

We obtain a contradiction to (\ref{eqT321}), as in the $\tau=0$ case, unless we have the form
$$
u_{n_0}=x^a,
v_{n_0}=y,
w_{n_0i}=x^cz.
$$

Now by an argument similar to the case $\tau=0$ of Theorem \ref{TheoremT124}, we obtain $$
\tau_{\tilde f_n}(p_n)=-\infty
$$
 for $n>>0$, a contradiction.

Thus (\ref{eqT318}) terminates in a finite number of steps $n$. Let $(R^n)'$ be the restriction of $\tilde R^n$, defined by
$$
U((R^n)'_i)=\tilde f_n(T(\tilde R_i^n)).
$$
Let 
$$
\Omega((R^n)_i')=\Omega(\tilde R_i^n)-[U(\tilde R_i^n)-U((R^n)')].
$$

We have that $\tilde f_n:\tilde X_n\rightarrow\tilde Y_n$ with relation $(R^n)'$ satisfies 1, 2 and 3 of Definition \ref{Def65} of a $\tau$-well prepared morphism.

There exists a sequence of blow ups of 2-curves $\Psi_1:Y_1\rightarrow \tilde Y_n$ such that
 4 (as well as 3) of Definition \ref{Def65} hold for the transforms $\{\overline R_i^1\}$ of the
$\{(\overline R^n)'_i\}$ on $Y_1$ at all 2-points of $U((\overline R^n)'_i)$, by  Lemma 5.14 \cite{C6}).  By Lemma \ref{Lemma31}, there exists a $\tau$-quasi-well
prepared diagram
$$
\begin{array}{rll}
X_1&\stackrel{f_1}{\rightarrow}&Y_1\\
\Phi_1\downarrow&&\downarrow\Psi_1\\
\tilde X_n&\stackrel{\tilde f_n}{\rightarrow}&\tilde Y_n
\end{array}
$$
of $(R^n)'_i$ and $\Psi_1$
where $\Phi_1$ is a product of blow ups of 2-curves. 
Let $R^1$ be the transform of $(R^n)'$ on $X_1$.
$f_1$ is $\tau$-well prepared with relation $R^1$.
\end{pf}

\begin{Theorem}\label{TheoremT324} Suppose that $\tau\ge 0$ and  $f:X\rightarrow Y$ is
$\tau$-well
prepared with   relation $R$.  
 Then there exists a commutative
diagram
$$
\begin{array}{rll}
X_1&\stackrel{f_1}{\rightarrow}&Y_1\\
\Phi\downarrow&&\downarrow\Psi\\
X&\stackrel{f}{\rightarrow}&Y
\end{array}
$$
such that $\Phi$,$\Psi$ are products of blow ups of possible centers and $f_1$ is
$\tau$-very-well prepared with   relation $R^1$.
(In general, $R^1$ is not the transform of $R$). 
\end{Theorem}

 \begin{pf}
 Let $\Omega=G_{Y}(f,\tau)-\Theta(f,Y)$. 
$\Omega$ is a finite set by Remark \ref{RemarkT159}. 
 
 Let $R'$ be the restriction of $R$ to $\Omega$. 
 $R'$ has algebraic structure determined by the algebraic structure of $R$.

 By Theorem \ref{CorollaryT126}, there exists a $\tau$-well prepared
 diagram
 
 $$
 \begin{array}{rll}
 X_2&\stackrel{f_2}{\rightarrow}&Y_2\\
 \Phi_2\downarrow&&\downarrow\Psi_2\\
 X&\stackrel{f}{\rightarrow}&Y
 \end{array}
 $$
 for $R$
 where $\Phi_2$ and $\Psi_2$ are products of blow ups of possible centers
 such that $f_2$ is $\tau$-well prepared  for the transform $R^2$ of $R'$, and
 there exists an open subset $V\subset Y$ such that $Y-V$ is a finite set of points and $f_2$ is toroidal over $\Psi_2^{-1}(V)$. Further, $V\cap G_{Y}(f,\tau)=\Theta(f,Y)$.
Since $\Omega$ is finite, we may modify the $R'$, so that $U(R_i')$ is a single point $\{q_i\}$ for all primitive relations $R_i'$ associated to $R'$.

 Let $R_i^2$ be a primitive relation associated to $R^2$. Then $U(R'_{i})=\{q_i\}$
 for some $q_i\in\Omega$.
 $\Omega(R'_{i})$ is a  neighborhood of $q_i$ on a surface in $Y$, and
 $\Omega(R_{i}^2)\rightarrow\Omega(R'_{i})$ is a projective birational map.
 Suppose that $E$ is a component of $D_{Y_2}$ such that $\gamma=E\cdot\Omega(R_{i}^2)$ dominates a curve (containing $q_i$) of 
 $\Omega(R'_{i})$. Then a general point $\eta$ of $\gamma$ is a 1-point over which $f_2$ is toroidal, since $\Psi_2(\eta)\in V$.
 Hence $f_2$ is finite over a general point of $\gamma$, and $\gamma$ is not in the fundamental locus of $f_2$. Further, by Remark \ref{RemarkT119}, if $q\in\gamma\cap U(R_i^2)$, and $f_2^{-1}(q)\cap T(R_i^2)\ne\emptyset$, then $q$ is a 2-point.

 By Lemma \ref{Lemma33}, there exists a $\tau$-well prepared diagram (for $R^2$)
 
 $$
 \begin{array}{rll}
 X_3&\stackrel{f_3}{\rightarrow}&Y_3\\
 \Phi_3\downarrow&&\downarrow\Psi_3\\
 X_2&\stackrel{f_2}{\rightarrow}&Y_2
 \end{array}
 $$
 where $\Psi_3$ is a product of blow ups of prepared 2-points (of type 1 in Definition \ref{Def66}) such that if $R_i^3$ is a primitive relation
 associated to the transform $R^3$ of $R^2$, $E$ is a component of $D_{Y_3}$, and $\gamma=E\cdot\Omega(R_i^3)$ is not exceptional for
 $\Omega(R_i^3)\rightarrow\Omega(R_i')$, then $\gamma\cap f_3(T(R_i^3))=\emptyset$.
 We further have that $f_3(T(R^3))\cap\gamma=\emptyset$. We see this as follows. Suppose that $p\in f_3(T(R^3))\cap\gamma$.
 $\Psi_3(E)$ is a component of $D_{Y_2}$. Thus $p$ is on the strict transform $\gamma'$ of $\overline{\Psi_3(E)\cdot\Omega(R_j^2)}$ on $Y_3$,
 which is not exceptional for $\Psi_2$. Since $\gamma'=\overline{E\cdot\Omega(R_j^3)}$, we have a contradiction.

Let
$$
W_3=\left\{\begin{array}{l} \gamma=\overline{S_{R_i^3}(q)\cdot E} 
\text{ such that $E$ is a component of $D_{Y_3}$},\\
\text{$R_i^3$ is associated to $R^3$ and $q\in f_3(T(R^3))\cap U(\overline R_i^3)$}\end{array}\right\}.
$$
We have that all $\gamma\in W_3$ contract to a point on $Y_1$.
Let
$$
Z_3=\left\{
\begin{array}{l}
q\in U(R^3)-f_3(T(R^3))\text{ such that there exist }\gamma_i,\gamma_j\in W_3\\
\text{ such that }\gamma_i\ne\gamma_j\text{ and }q\in \gamma_i\cap\gamma_j.
\end{array}\right\}
$$

Suppose that $q\in Z_3$. Then there exist $\gamma_i=\overline{S_{R_i^3}(p_i)\cdot E_1}\in W_3$ and $\gamma_j=\overline{S_{R_j^3}(p_j)\cdot E_2}\in W_3$ such that $q\in\gamma_i\cap \gamma_j$ and $\gamma_i\ne\gamma_j$.

$\gamma_i$ and $\gamma_j$ are exceptional, so they contract to the common point $q_i=q_j\in U(R')$.
thus $\gamma_i$ and $\gamma_j$ are contained in $\Omega(R_i^3)$ and $\Omega(R_j^3)$ respectively.

The points of $Z_3$ are prepared 2-points for $R^3$ (of type 1 of Definition \ref{Def66}).
 Let $\Psi_4:Y_4\rightarrow Y_3$ be the
blow up of $Z_3$. By Lemma \ref{Lemma32}, there exists a $\tau$-well prepared diagram
$$
\begin{array}{rll}
X_4&\stackrel{f_4}{\rightarrow}&Y_4\\
\Phi_4\downarrow&&\downarrow\Psi_4\\
X_3&\stackrel{f_3}{\rightarrow}&Y_3
\end{array}
$$
of  $R^3$ and $\Psi_4$. Let $R^4$ be the transform of $R^3$ on $X_4$.

Define
$$
W_4=\left\{\begin{array}{l} \gamma=\overline{S_{\overline R_i^4}(q)\cdot E}\text{ such that $E$ is a component of }D_{Y_4},\\
\text{$R_i^4$ is associated to $R^4$ and $q\in \Psi_4^{-1}(f_3(T(R^3)))\cap U(\overline R_i^4)$}\end{array}\right\},
$$
$$
Z_4=\left\{
\begin{array}{l}
q\in U(R^4)-\Psi_4^{-1}(f_3(T(R^3)))\text{ such that there exist }\gamma_i,\gamma_j\in W_4\\
\text{ such that }\gamma_i\ne\gamma_j\text{ and }q\in \gamma_i\cap\gamma_j.
\end{array}\right\}
$$

  We necessarily
have that the curves in $W_4$ are strict transforms of curves in $W_3$. We can
iterate, blowing up $Z_4$, and constructing a $\tau$-well prepared diagram, and repeating
until we eventually construct a $\tau$-well prepared diagram of $R^4$
$$
\begin{array}{rll}
X_5&\stackrel{f_5}{\rightarrow}&Y_5\\
\Phi_5\downarrow&&\downarrow\Psi_5\\
X_4&\stackrel{f_4}{\rightarrow}&Y_4
\end{array}
$$
such that $\Psi_5$ is a sequence of blow ups of prepared 2-points (of type 1 of Definition \ref{Def66}) and if 
$\gamma_1=\overline{S_{\overline R_i^5}(q_i)\cdot E_i}$,
$\gamma_2=\overline{S_{\overline R_j^5}(q_j)\cdot E_j}$,
for $q_i\in U(\overline R_i^5)\cap (\Psi_4\circ\Psi_5)^{-1}(f_3(T(R^3)))$
(where $R^5$ is the transform of $R^4$) and $E_1, E_2$
components of $D_{Y_5}$, are such that $\gamma_1\ne\gamma_2$, then $\gamma_1\cap\gamma_2
\subset U(R^5)\cap (\Psi_4\circ\Psi_5)^{-1}(f_3(T(R^3))$.

We now construct pre-relations $\overline R_{i}^*$ on $Y_5$ with associated
primitive relations $R_{i}^*$ for $f_5$.

Let $T(R_{i}^*)= T(R_i^5)$
and let 
\begin{equation}\label{eq382}
U(R_{i}^*)=U(R_i^5)\cap(\Psi_4\circ\Psi_5)^{-1}(f_3(T(R^3))).
\end{equation}

For $q'\in U(\overline R_{i}^*)$, define $\overline R_{i}^*(q')=\overline R_i^5(q')$.
For $p\in T(R_{i}^*)$ define $R_{i}^*(p)=\overline R_{i}^5(f_5(p))$.
Let $R^*$ be the  relation for $f_5$ defined by the $R_{i}^*$.
Let $\Omega(\overline R_i^*)=\Omega(\overline R_i^5)$.

 For all $\overline R_i^*$, let 
$$
V_i(Y_5)=\left\{
\gamma=\overline{E_{\alpha}\cdot S_{\overline R_i^*}(q)}\text{ such that $q\in U(\overline R_i^*)$, $E_{\alpha}$
is a component of $D_{Y_5}$}\right\}.
$$
Recall that these curves are all exceptional for $\Psi_2\circ\Psi_3\circ\Psi_4\circ\Psi_5$. 

By our construction, Lemmas \ref{LemmaT79}, \ref{Lemma31}, \ref{Lemma32}, \ref{Lemma171} and \ref{LemmaT63},
and Remark \ref{RemarkT278} and 2 of Remark \ref{Remark293}, every curve $\gamma\in V_i(Y_5)$ is prepared for $R^5$ of type 6. By (\ref{eq382}), we
now conclude that every curve $\gamma\in V_i(Y_5)$ is prepared for $R^*$ of type 6. 1 and 2 of Definition \ref{Def130}
thus  hold for $f_5$ and $R^*$.
3 of Definition \ref{Def130} holds for $f_5$ and $R^*$ since for all $\overline R_i^*$, $V_i(Y_5)$ consists of exceptional
curves of $\Omega(\overline R_{i}^*)$ contracting to a nonsingular point $q_i\in\Omega(\overline R_i')$.
Thus $f_5$ is $\tau$-very-well prepared with  relation
 $R^*$.

\end{pf}

\begin{Theorem}\label{Theorem268} Suppose that $f:X\rightarrow Y$ is prepared and
$\tau=\tau_f(X)\ge 0$.
 Then there exists a commutative diagram
$$
\begin{array}{rll}
X_1&\stackrel{f_1}{\rightarrow}&Y_1\\
\Phi\downarrow&&\downarrow\Psi\\
X&\stackrel{f}{\rightarrow}&Y
\end{array}
$$
such that $f_1$ is prepared, $\Phi$, $\Psi$ are products of blowups of 2-curves,
$\tau_{f_1}(X_1)\le\tau$,  and 
$G_{Y_1}(f_1,\tau)$ contains no 3-points and no 2-curves, so that $f_1$ is $\tau$-prepared.
\end{Theorem}

\begin{pf} This is immediate from Lemma 4.1 \cite{C7} and Lemma \ref{Lemma1}. \end{pf}

\begin{Theorem}\label{Theorem269} Suppose that $f:X\rightarrow Y$ is prepared, and
$\tau=\tau_{f}(X)\ge0$.   Then there exists a commutative diagram
$$
\begin{array}{rll}
X_1&\stackrel{f_1}{\rightarrow}&Y_1\\
\Phi_1\downarrow&&\downarrow\Psi_1\\
X&\stackrel{f_1}{\rightarrow}&Y
\end{array}
$$
such that $\Phi_1$ and $\Psi_1$ are products of blow ups of possible centers and $f_1$ is $\tau$-very-well prepared with a relation $R^1$.
\end{Theorem}

\begin{pf} 
By Theorem \ref{Theorem268}, there exists a commutative diagram  
$$
\begin{array}{rll}
X_1&\stackrel{f_1}{\rightarrow}&Y_1\\
\Phi\downarrow&&\downarrow \Psi\\
X&\stackrel{f}{\rightarrow}&Y
\end{array}
$$
such that $\Phi$ and $\Psi$ are products of 2-curves,  and  $f_1$ is $\tau$-prepared. Now by Theorems \ref{Theorem169},
 \ref{Theorem80}, \ref{TheoremT325} and \ref{TheoremT324}, there exists a commutative diagram
$$
\begin{array}{rll}
X_2&\stackrel{f_2}{\rightarrow}&Y_2\\
\downarrow&&\downarrow\\
X_1&\stackrel{f_1}{\rightarrow}&Y_1
\end{array}
$$
where the vertical arrows are products of blow ups of possible centers  such that $f_2$ is $\tau$-very-well prepared.

\end{pf}

\section{Toroidalization}
Suppose that $f:X\rightarrow Y$ is a proper, birational morphism of nonsingular  3-folds with toroidal structures $D_Y$ and $D_X=f^{-1}(D_Y)$, such that $D_X$ contains the singular locus of $f$.

\begin{Theorem}\label{Theorem270} Suppose that $\tau\ge 0$ and $f:X\rightarrow Y$  is $\tau$-very-well prepared with  relation $R$.
Then there exists a $\tau$-very-well prepared diagram
$$
\begin{array}{rll}
X_1&\stackrel{f_1}{\rightarrow}&Y_1\\
\downarrow&&\downarrow\\
X&\stackrel{f}{\rightarrow}&Y
\end{array}
$$
such that the transform $R^1$ of $R$ is resolved ($T(R')=\emptyset$).
In particular, $f_1$ is prepared and $\tau_{f_1}(X_1)<\tau$. 
\end{Theorem}

\begin{pf} 
Fix a  pre-relation $\overline R_t$ associated to $R$  on $Y$, with associated primitive relation $R_t$.
By induction on the number of pre-relations associated to $R$, it suffices to resolve $R_t$
by a $\tau$-very-well prepared diagram (of $R$).

Recall (Definition \ref{Def130})
$$
V_t(Y)=\left\{\begin{array}{l}\overline{E\cdot S}\text{ such that } E\text{ is a component of }D_Y,\\
S=S_{\overline R_t}(q)\text{ for some }q\in U(\overline R_t)
\end{array}\right\}.
$$

$F_t=\sum_{\gamma\in V_t(Y)}\gamma$ is a SNC divisor on $\Omega(\overline R_t)$ whose intersection graph is a forest.

If $\gamma_1=\overline{E_1\cdot S_{\overline R_t}(q_1)}\in V_t(Y)$ and $q\in \gamma_1$, we will say that $\gamma_1$ is good at $q$ if
whenever $q\in U(\overline R_i)$ for some $i$,
 then  $S_{\overline R_i}(q)$ contains the germ of  $\gamma_1$ at $q$
(so that $\gamma_1=\overline{E_1\cdot S_{\overline R_i}(q)}\subset \Omega(\overline R_i)$).
 Otherwise, say that $\gamma_1$ is bad at $q$. Say that $\gamma_1$ is good
if $\gamma_1$ is good at $q$ for all $q\in\gamma_1$.

Let $Y_0= Y$, $X_0=X$, $f_0=f$. We will show that there exists a sequence of
$\tau$-very-well prepared diagrams of the transform of $R$, 
\begin{equation}\label{eq287}
\begin{array}{rll}
X_{i+1}&\stackrel{f_{i+1}}{\rightarrow}&Y_{i+1}\\
\Phi_{i+1}\downarrow&&\Psi_{i+1}\downarrow\\
X_i&\stackrel{f_i}{\rightarrow}&Y_i
\end{array}
\end{equation}
for $0\le i\le m-1$ such that   the transform $R_t^m$ of $ R_t$ on
$X_m$ is resolved. 

Suppose that  $\gamma_1\in V_t(Y)$ and $q\in \gamma_1$
is a bad point. By Remark \ref{Remark281}, we have that $q\in U(\overline R_t)$.
By (\ref{eqT100}) of Definition \ref{Def65}, we have that $q$ is a 2-point.
Suppose that  $E_1,E_2$ are the two components of $D_{Y}$
containing $q$, and $\gamma_1=\overline{E_1\cdot S_{\overline R_t}(q)}$. Let
$\gamma_2=\overline{E_2\cdot S_{\overline R_t}(q)}$.

We will show that $q$ is a good point of $\gamma_2$.

$\gamma_1$ not good at $q$  implies  there exists
$j\ne t$ such that $q\in U(\overline R_j)$ and the germ of $\gamma_1$ at $q$ is
not contained in $S_{\overline R_j}(q)$.
Let
$$
u=u_{\overline R_t}(q),
v=v_{\overline R_t}(q),
w_t=w_{\overline R_t}(q).
$$
After possibly interchanging $u$ and $v$ we have that $u=w_t=0$ are local equations of
$\gamma_1$, $v=w_t=0$ are local equations of $\gamma_2$ at $q$. Let $w_j= w_{\overline R_j}(q)$.
In the equation
$$
w_j=w_t+u^{a_{tj}}v^{b_{tj}}\phi_{tj}
$$
of (\ref{eq64}) of Definition \ref{Def65} we thus have $a_{tj}=0$.

If $q$ is not a good point for $\gamma_2$ then there exists $k\ne t$ such that
$q\in U(\overline R_k)$ and the germ of $\gamma_2$ at $q$ is not contained in $S_{\overline R_k}(q)$. Let $w_k=w_{\overline R_k}(q)$. In the equation
$$
w_k=w_t+u^{a_{tk}}v^{b_{tk}}\phi_{tk}
$$
of (\ref{eq64}) we thus have $b_{tk}=0$. But we must have
$$
(0,b_{tj})\le (a_{tk},0)\text{ or }(a_{tk},0)\le (0,b_{tk})
$$
by 4 of Definition \ref{Def65}, which is impossible. Thus $q$ is a good point for $\gamma_2$.

Suppose that  all $\gamma\in V_t(Y)$ are bad.
Pick $\gamma_1\in V_t(Y)$. Since $\gamma_1$ is bad there exists
 $\gamma_2\in V_t(Y)
-\{\gamma_1\}$ such that $\gamma_2$ is good at $q_1=\gamma_1\cap\gamma_2$
(as shown above).
$\gamma_1\cap\gamma_2$ is a single point since $V_t(Y)$ is a forest.
Since $\gamma_2$ is bad and $V_t(Y)$ is a forest, there exists $\gamma_3\in V_t(Y)$ which intersects $\gamma_2$
at a single point $q_2$ and is disjoint from $\gamma_1$ such that $\gamma_3$ is good at
$q_2$. Since $V_t(Y)$ is a finite set, and the intersection graph of $V_t(Y)$ is a forest, we must eventually find a curve which is good,
a contradiction.

Let $\gamma\in V_t(Y)$ be a good curve, so that it is prepared for $R$
of type 6, and is a *-permissible center (Lemma \ref{Lemma67}) and let $\Psi_1':
Y_1'\rightarrow Y$ be the blow up of $\gamma$.

By Lemma \ref{Lemma67} we can construct a $\tau$-very-well prepared diagram
of the form of (\ref{eq233}) of Definition \ref{Def219} 
\begin{equation}\label{eq321}
\begin{array}{rll}
X_1&\stackrel{f_1}{\rightarrow}&Y_1\\
\downarrow&&\downarrow\\
\downarrow&&Y'_1\\
\downarrow&&\downarrow\Psi_1'\\
X&\rightarrow&Y.
\end{array}
\end{equation}
where $Y_1\rightarrow Y'_1$ is a sequence of blow ups of 2-points which are prepared
for the transform of $R$ of type 2 of Definition \ref{Def66}.
Observe that if $\gamma_1\in V_t(Y)$ is a good curve, with $\gamma_1\ne\gamma$, then the strict transform of $\gamma_1$ is a good curve in $V_t(Y_1)$.

 We now
iterate this process. 
We order the curves in $V_t(Y)$,  and choose $\gamma=\overline{E\cdot S_{R_t}(q)}\in V_t(Y)$ in the construction of the diagram (\ref{eq321}) so that it is the minimum
good curve in $V_t(Y)$.

We inductively define a sequence of $\tau$-very well prepared diagrams (\ref{eq287}) by blowing up the good curve in $V_t(Y_i)$
with smallest order, and then constructing a very well prepared diagram (\ref{eq287}) of the form of (\ref{eq321}).
Then we define the total ordering on $V_t(Y_{i+1})$ so that the ordering of strict transforms
of elements of $V_t(Y_{i})$ is preserved, and these strict transforms have smaller
order than the element of $V_t(Y_{i+1})$ which is not a strict transform of an element of $V_t(Y_{i})$.
We repeat, as long as $R_t^i$ is not resolved ($T(R_t^i)\ne\emptyset$).

Suppose that the algorithm does not converge in the construction of $f_m:X_m\rightarrow Y_m$ such that the transform  $R_t^m$ of $R_t$ is resolved. Then there exists a diagram 
\begin{equation}\label{eq232}
\begin{array}{rll}
\vdots&&\vdots\\
\downarrow&&\downarrow\\
X_n&\stackrel{f_n}{\rightarrow}&Y_n\\
\Phi_n\downarrow&&\downarrow\Psi_n\\
X_{n-1}&\stackrel{f_{n-1}}{\rightarrow}&Y_{n-1}\\
\downarrow&&\downarrow\\
\vdots&&\vdots\\
\downarrow&&\downarrow\\
X_0=X&\stackrel{f_0=f}{\rightarrow}&Y_0=Y
\end{array}
\end{equation}
constructed by infinitely many iterations of the algorithm such that $T(R_t^n)\ne\emptyset$ for all $n$.

Suppose that $q_n\in U(\overline R_t^n)$ is an infinite sequence of points such that $\Psi_n(q_n)=q_{n-1}$ for all $n$ and $\Psi_n$
is not an isomorphism for infinitely many $n$.
$q_n$ is either a 2-point or a 1-point for all $n$. 

First suppose that $q_n$ is a 2-point for all $n$.

By construction, the restriction of
 $\Psi_n$ to $S_{\overline R_t^{n}}(q_{n})$ is an isomorphism onto  $S_{\overline R_t^{n-1}}(q_{n-1})$ for all $n$. Thus
the restriction
$$
\overline\Psi_n=\Psi_1\circ\cdots\circ\Psi_n:S_{\overline R_t^n}(q_n)
\rightarrow S_{\overline R_t}(q)
$$
is an isomorphism, where $q=q_0=\Psi_1\circ\cdots\circ\Psi_n(q_n)$. Without loss of generality, we may assume that no $\Psi_n$ is an
isomorphism (on $Y_n$) at $q_n$. We have permissible parameters
$u_i=u_{\overline R_t^i}(q_i),v_i=v_{\overline R_t^i}(q_i),w_{t,i}=w_{\overline R_t^i}(q_i)$
at $q_i$ for all $i$ such that either 
\begin{equation}\label{eq230}
u_i=u_{i+1},
v_i=v_{i+1},
w_{t,i}=u_{i+1}w_{t,i+1}
\end{equation}
or 
\begin{equation}\label{eq231}
u_i=u_{i+1},
v_i=v_{i+1},
w_{t,i}=v_{i+1}w_{t,i+1}.
\end{equation}

Suppose there exists $k\ne t$  such that $q_n\in U(\overline R_k^n)$ for all $n$.

Let $w_{k,i}=w_{\overline R_k^i}(q_i)$ for $i\ge 0$.

The relations
$$
w_{k,i}-w_{t,i}=u_i^{a_{tk}}v_i^{b_{tk}}\phi_{t,k}
$$
of (\ref{eq64}) of Definition \ref{Def65}
transform to
$$
w_{k,i+1}-w_{t,i+1}=u_{i+1}^{a_{tk}-1}v_{i+1}^{b_{tk}}\phi_{t,k}
$$
under (\ref{eq230}), and transform to
$$
w_{k,i+1}-w_{t,i+1}=u_{i+1}^{a_{tk}}v_{i+1}^{b_{tk}-1}\phi_{t,k}
$$
under (\ref{eq231}). But we see that after a finite number of iterations $q_n\not\in U(\overline R_k^n)$, unless $a_{tk}=b_{tk}=-\infty$.
Thus there exists $n_0$, such that whenever $n\ge n_0$, $q_n\not\in U(\overline R_k^n)$ if $k\ne t$ and $a_{kt},b_{kt}\ne-\infty$ .

Now suppose that $q_n$ is a 1-point for all $n$. We have permissible parameters
$$
u_i=u_{\overline R_t^i}(q_i),
v_i=v_{\overline R_t^i}(q_i),
w_{t,i}=w_{\overline R_t^i}(q_i)
$$
at $q_i$ for all $i$ such that
$$
u_i=u_{i+1}, v_i=v_{i+1},
w_{t,i}=u_{i+1}w_{t,i+1}
$$
for all $i$.

Suppose that there exists $k\ne t$ such that $q_n\in U(\overline R_k^n)$ for all $n$. Let
$$
w_{k,i}=w_{\overline R_k^i}(q_i)
$$
for $i\ge 0$. The relation
$$
w_{k,i}-w_{t,i}=u_i^{c_{t,k}}\phi_{t,k}
$$
of (\ref{eqT100}) of Definition \ref{Def65} transforms to
$$
w_{k,i+1}-w_{t,i+1}=u_{i+1}^{c_{t,k}-1}\phi_{t,k}.
$$
Thus after a finite number of iterations, $q_n\not\in U(\overline R_k^n)$ unless $c_{t,k}=-\infty$. Thus there exists $n_0$ such that when $n\ge n_0$, $q_n\not\in U(\overline R_k^n)$ if $k\ne t$ and $c_{t,k}\ne-\infty$.

By our ordering, we have that there exists an $n_0$ such that if $n\ge n_0$, $\gamma\in V_t(Y_n)$ is good and if $k$ is such that $\gamma\cap U(\overline R_k^n)\ne\emptyset$ then
$a_{kt},b_{kt}=-\infty$ (or $c_{rk}=-\infty$), so that the Zariski closures of $\Omega(\overline R_k^n)$ and $\Omega(\overline R_t^n)$ are the same.
Thus all elements of $V_t(Y_n)$ are good for $n\ge n_0$, since otherwise, there would be a bad curve $\gamma_1\in V(Y_n)$ which intersects
a good curve $\gamma_2$ at a point $q'$ at which $\gamma_1$ is not good. But then we must have that there exists $k\ne t$ such that 
the Zariski closure of $\Omega(\overline R_k^n)$ is not equal to the Zariski closure of $\Omega(\overline R_t^n)$,
 and 
$q'\in U(\overline R_k^n)$, so that $\gamma_2\cap U(\overline R_k^n)\ne\emptyset$, a contradiction.

Our birational morphism  of $\Omega(\overline R_t^n)$ to $\Omega(\overline R_t)$ is an isomorphism in a neighborhood of $U(\overline R_t^n)$. Thus we have a
natural identification of $V_t(Y_n)$ and $V_t(Y)$, and we see that for $n\ge n_0$, the $\Psi_n$ cyclically blow up the different curves of $V_t(Y)$.

Since $R_t^n$ is (by assumption) not resolved for all $n$, there are points
$p_n\in T(R_t^n)\subset X_n$ such that $\Phi_n(p_n)=p_{n-1}$, and $f_n(p_n)=q_n\in U(\overline R_t^n)$  for all $n$.
Without loss of generality, we may assume that no $\Psi_n$ is an
isomorphism  at $q_n$.

We have that all $q_n$ are 1-points or all $q_n$ are 2-points.

First suppose that all $q_n$ are 2-points.

With the above notation at $q_n=f_n(p_n)$, we have that  (\ref{eq230}) and (\ref{eq231}) must alternate in the diagram (\ref{eq232}) for $n\ge n_0$, by
our ordering  of $V_t(Y_n)$. Let $p=p_0\in X=X_0$, $q=q_0=f(p)$.

We have 
\begin{equation}\label{eq322}
u=u_n,
v=v_n,
w_t=u_n^{a_n}v_n^{b_n}w_{t,n}
\end{equation}
where 
$$
u=u_{\overline R_t}(q), v=v_{\overline R_t}(q), w_t=w_{\overline R_t}(q)
$$
and
 $a_n$, $b_n$ are positive integers which both go to infinity as $n$ goes to infinity.

There exists (by Theorem 4 of  Section 4, Chapter VI \cite{ZS}) a  valuation $\nu$ of ${\bold k}(X)$ which dominates the  
(non-Noetherian) local ring $\cup_{n\ge 0}{\cal O}_{X_n,p_n}$, and thus dominates the local rings
${\cal O}_{X_n,p_n}$ for all $n$.  
Without loss of generality, we may identify $\nu$ with an extension of $\nu$ to the quotient field of $\hat{\cal O}_{X_n,p_n}$ which dominates $\hat{\cal O}_{X_n,p_n}$ for all $n$.

Let $x,y,z$ be permissible parameters for $u,v,w_t$ at $p$. Suppose that $p$ is a 3-point. Write (in $\hat{\cal O}_{X,p}$) 
\begin{equation}\label{eq385}
\begin{array}{ll}
u&=x^ay^bz^c\\
v&=x^dy^ez^f\\
w_t&=x^gy^hz^i\gamma
\end{array}
\end{equation}
where $xyz=0$ is a local equation of $D_X$ at $p$ and $\gamma$ is a unit series.

 We may permute $x,y,z$ so that $0<\nu(x)\le\nu(y)\le \nu(z)$. We have (from (\ref{eq322}))
$$
\nu(w_t)-n\nu(u)-n\nu(v)>0
$$
for all $n\in {\bf N}$. Thus
$$
0<(g-na-nd)\nu(\overline x)+(h-ne-nb)\nu(\overline y)+(i-nf-nc)\nu(\overline z)
\le ((g+h+i)-nf-nc)\nu(\overline z)
$$
for all $n$. Thus $f=c=0$, but this is impossible, since $uv=0$ is a local equation of $D_X$ at $p$.

There is a similar but simpler algorithm if $p$ is a 1-point or a 2-point and $\tau>0$, since $w_t=0$ is supported on $D_X$ at $p$.

Suppose that $\tau=0$ (and all $q_n$ are 2-points). We have that $w_t=0$ is a divisor supported on $D_X$ at $p$, so that the argument for $\tau>0$ works in this case also, or we have one of the following two special forms:

$p$ a 1-point 
\begin{equation}\label{eqT302}
u=x^a,
v=x^b(\alpha+y),
w_t=x^cz
\end{equation}
or $p$ a 2-point 
\begin{equation}\label{eqT303}
u=x^ay^b,
v=x^cy^d,
w_t=x^ey^fz
\end{equation}
with $ad-bc\ne0$.

In the diagram (\ref{eq321}) we have that $q_1\in U(R_i^1)$ implies $Y_1\rightarrow Y_1'$ is an isomorphism near $q_1$.  $X_1\rightarrow X$ factors in a neighborhood of $q_1$ as a diagram
$$
X_1=W_m\stackrel{\Lambda_m}{\rightarrow}W_{m-1}\rightarrow\cdots\rightarrow W_1\stackrel{\Lambda_1}{\rightarrow}X
$$
where $\Lambda_1:W_1\rightarrow X$ is a sequence of blow ups of 2-curves and 3-points, and each $\Lambda_{j+1}$ is the blow up of a possible curve $\Sigma_j$ containing a 1-point and the center
$\overline p_j$ of $\nu$ on $W_j$, such that ${\cal I}_{\gamma}{\cal O}_W$ is not invertible.

Without loss of generality, we may assume that $u=w_t=0$ are local equations of $\gamma$ at $q$.

We have a form (\ref{eqT302}) or (\ref{eqT303}) at the center $\overline p_1$ of $\nu$ on $W_1$, where $(a,b)\le(e,f)$ or $(a,b)>(e,f)$ if (\ref{eqT303}) holds.

Suppose that $a\le c$ in (\ref{eqT302}), or $(a,b)\le(e,f)$ in (\ref{eqT303}).
Then ${\cal I}_{\gamma}{\cal O}_{W_1,\overline p_1}$ is invertible. Thus $W_m=W_1$, and in a neighborhood of $q_1$, $X_1\rightarrow X$ is a sequence of blow ups of 2-curves. Further, we have an expression of the form (\ref{eqT302}) or (\ref{eqT303}) at $p_1$.

Suppose that $a>c$ in (\ref{eqT302}). Then $\Lambda_2$ is the blow up of a curve with local equations $x=z=0$ at $\overline p_1$.   Since $\nu(z)>n\nu(x)$ for all $n\in{\bf N}$, at $\overline p_2$ we have regular parameters $x_2,y_2,z_2$ defined by
$$
x=x_2, y=y_2, z=x_2z_2
$$
and we thus have
$$
u=x_2^a,
v=x_2^b(\alpha+y_2),
w_t=x_2^{c+1}z_2.
$$

Iterating, we see that $\overline p_m$ has permissible parameters $x_m,y_m,z_m$ such that
$$
u=x_m^a,
v=x_m^b(\alpha+y_m),
w_t=x_m^az_m.
$$
The permissible parameters $u_1, v_1, w_{t,1}$ at $q_1$ are defined by
$$
u=u_1, v=v_1, w_t=u_1w_{t1}.
$$
 Thus
$$
u_1=x_m^a, v_1=x_m^b(\alpha+y_m), w_{t1}=z_m
$$
and we have $\tau_{f_1}(p_1)=-\infty$, a contradiction.

We have a similar analysis if $(a,b)>(e,f)$ in (\ref{eqT303}), leading to the conclusions that $\tau_{f_1}(p_1)=-\infty$, a contradiction.

We thus see that for all $n$ in (\ref{eq232}), $X_n\rightarrow X_{n-1}$ factors as a sequence of blow ups of 2-curves at $q_n$.

Suppose that there exists $n_0$ such that $p_{n_0}$ is a 1-point, and thus $p_n$ is a 1-point for all $n\ge n_0$. Then we see that $X_{n+1}\rightarrow X_n$ is an isomorphism at $p_n$ for all $n\ge n_0$.  At $p_{n_0}$ there are regular parameters $x,y,z$ such that 
\begin{equation}\label{eqT304}
u_{n_0}=x^a,
v_{n_0}=x^b(\alpha+y),
w_{tn_0}=x^cz,
\end{equation}
and $q_n$ has regular parameters
$$
u_{n_0}=u_n,
v_{n_0}=v_n,
w_{tn_0}=u_n^{n-n_0}w_{tn}.
$$
Substituting into (\ref{eqT304}), we have a contradiction as soon as
$$
n>\frac{c}{a}+n_0.
$$

Now suppose that $p_n$ is a 2-point for all $n$. Then a form (\ref{eqT303}) holds at $p_1$, and at $p_n$, we have regular parameters $x_n,y_n,z_n$ defined by 
\begin{equation}\label{eqT305}
x=x_n^{r_{11}^n}y_n^{r_{12}^n},
y=x_n^{r_{21}^n}y_n^{r_{22}^n},
z=z_n
\end{equation}
such that $r_{11}^nr_{22}^n-r_{12}^nr_{21}^n=\pm1$.

Now from (\ref{eq322}), (\ref{eqT303}) and (\ref{eqT305}), we see that 
$$
\nu(x^ey^f)-n\nu(x^ay^b)-n\nu(x^cy^d)>0
$$
for all $n\in{\bf N}$, a contradiction, since $ad-bc\ne 0$.

The argument is simpler in the case when $q_n$ is a 1-point for all $n$. (\ref{eq322}) becomes 
\begin{equation}\label{eqT179}
u=u_n, v=v_n, w_t=u_n^{a_n}w_{t,n}
\end{equation}
where $a_n$ goes to infinity as $n$ goes to infinity. 

Suppose that $\tau>0$. If $p\in f^{-1}(q)$, $u=0$ is a local equation of $D_X$ at $p$, and $w_t=0$ is supported on $D_X$ at $p$.  Thus (\ref{eqT179}) leads to a contradiction.

If $\tau=0$, there is a similar argument to the above case of $q_n$ a 2-point for all $n$ and $\tau=0$.

Thus the algorithm converges in a morphism $f_m:X_m\rightarrow Y_m$ such that $T(R_t^m)= \emptyset$, and after iterating for each primitive relation
associated to $R$, we obtain the construction of $f_1:X_1\rightarrow Y_1$, as in the conclusions of the theorem, such  that
$f_1$ is prepared and $\tau_{f_1}(X_1)<\tau$.

\end{pf}

\noindent{\bf Proof of Theorem \ref{Theorem1}} First suppose that $X$ and $Y$ are proper over ${\bf k}$. By resolution of singularities and resolution of indeterminacy \cite{H} (cf. Section 6.8 \cite{C6}),
and by \cite{M}, there exists a commutative diagram
$$
\begin{array}{rll}
X_1&\stackrel{f_1}{\rightarrow}&Y_1\\
\Phi_1\downarrow&&\downarrow\Psi_1\\
X&\stackrel{f}{\rightarrow}&Y
\end{array}
$$
where $\Phi_1$, $\Psi_1$ are products of possible blow ups of points and nonsingular curves supported above $D_X$ and $D_Y$, such that $X_1$ and $Y_1$
are nonsingular and projective. Further, $D_{Y_1}=\Psi_1^{-1}(D_Y)$ and $D_{X_1}=\Phi_1^{-1}(D_X)=f_1^{-1}(D_{Y_1})$ are SNC divisors,
and $D_{X_1}$ contains the locus where $f_1$ is not smooth. By 
Theorem \ref{Theorem2},  we can construct a commutative diagram
$$
\begin{array}{rll}
X_2&\stackrel{f_2}{\rightarrow}&Y_2\\
\Phi_2\downarrow&&\downarrow\Psi_2\\
X_1&\stackrel{f_1}{\rightarrow}&Y_1
\end{array}
$$
such that $\Phi_2$ and $\Psi_2$ are products of possible blow ups of points and nonsingular curves, such that $f_2$ is prepared for $D_{Y_2}=\Psi_2^{-1}(D_{Y_1})$ and $D_{X_2}=\Phi_2^{-1}(D_{X_1})$.

Now by descending induction on $\tau=\tau_{f_2}(X_2)$ and Theorems \ref{Theorem269} and \ref{Theorem270}, there
exists a commutative diagram
$$
\begin{array}{rll}
X_3&\stackrel{f_3}{\rightarrow}&Y_3\\
\Phi_3\downarrow&&\downarrow\Psi_3\\
X_2&\stackrel{f_2}{\rightarrow}&Y_2
\end{array}
$$
such that $\Phi_2$ and $\Psi_3$ are products of blow ups of possible centers, $f_3$ is prepared, and $\tau_{f_3}(X_3)=-\infty$.

Thus $f_3$ is toroidal, and the conclusions of the theorem follow.

Now suppose that $X$ and $Y$ are (not necessarily proper) abstract varieties.  There exist (by \cite{N}) proper ${\bf k}$-varieties $\overline X$ and $\overline Y$ such that $X$ is an open subset of $\overline X$, and $Y$ is an open subset of $\overline Y$. After possibly modifying $\overline X$ and $\overline Y$ by blowing up $\overline X-X$ and $\overline Y-Y$, we may assume that $F_1=\overline X-X$ and $F_2=\overline Y-Y$ are closed subsets of pure codimension 1 in $\overline X$, $\overline Y$ respectively.

Let $\overline D_X$ be the Zariski closure of $D_X$ in $\overline X$, $\overline D_Y$ be the Zariski closure of $D_Y$ in $\overline Y$.

Let $D_{\overline X}=\overline D_X+F_1$, $D_{\overline Y}=\overline D_Y+F_2$. By resolution of indeterminancy, after possibly modifying $\overline X$ by blowing up subvarieties of $\overline X$ supported above $D_{\overline X}$, we have that the rational map $\overline f:\overline X\rightarrow \overline Y$ which extends $f:X\rightarrow Y$ is a morphism.

The hypotheses of Theorem \ref{Theorem1} are satisfied for $\overline f:\overline X\rightarrow \overline Y$, so by the first part of this proof, there exists a commutative diagram
$$
\begin{array}{rll}
\overline X_1&\stackrel{\overline f_1}{\rightarrow}&\overline Y_1\\
\overline\Phi\downarrow&&\downarrow\overline\Psi\\
\overline X&\stackrel{\overline f}{\rightarrow}&\overline Y
\end{array}
$$
satisfying the conclusions of Theorem \ref{Theorem1}.

Let $X_1=\Phi^{-1}(X)$, $Y_1=\Psi^{-1}(Y)$, $f_1=\overline f_1\mid X_1$, $\Phi=\overline\Phi\mid D_1$, $\Psi=\overline\Psi\mid Y_1$. Then
$$
\begin{array}{rll}
X_1&\stackrel{f_1}{\rightarrow}&Y_1\\
\Phi\downarrow&&\downarrow\Psi\\
X&\stackrel{f}{\rightarrow}&Y
\end{array}
$$
satisfying the conclusions of Theorem \ref{Theorem1}.

\section{List of technical terms}

(See also Section \ref{SectionNotation}, Notation)
\vskip .2truein

admissible center: Definition \ref{Def154}

fundamental locus: Remark \ref{Remark1}

$G_X(f,\tau)$: Definition \ref{Def326}

$G_Y(f,\tau)$: Definition \ref{Def326}

good at $p$ for $f$: Definition \ref{DefT40}

monomial form: Definition \ref{Def125}

perfect for $f$: Definition \ref{DefT133}.

permissible center: Definition \ref{Def289}

*-permissible center: Definition \ref{Def219}

permissible parameters: before and after Definition \ref{torf}

possible center:  Section \ref{SectionNotation}, Notation

prepared point or curve of type 1-5: Definition \ref{Def66}

prepared curve of type 6: Definition \ref{Def200}

prepared morphism

\hskip .2truein prepared morphism of 3-folds: Definition \ref{Def31}

\hskip .2truein $\tau$-prepared morphism: Definition \ref{Def326}

\hskip .2truein pre-$\tau$-quasi-well prepared: Definition \ref{DefT60}

\hskip .2truein $\tau$-quasi-well prepared: Definition \ref{DefT60}

\hskip .2truein $\tau$-well prepared: Definition \ref{Def65}

\hskip .2truein $\tau$-very-well prepared: Definition \ref{Def130}

relation

\hskip .2truein quasi-pre-relation: Definition \ref{DefT56}

\hskip .2truein pre-relation: Definition \ref{DefT165}

\hskip .2truein algebraic pre-relation: Definition \ref{Def199}

\hskip .2truein primitive relation: Definition \ref{DefT57}

\hskip .2truein relation: Definition \ref{DefT57}

\hskip .2truein algebraic relation: Definition \ref{DefT57}

resolved quasi-pre-relation: after Definition \ref{DefT56}

resolved relation: after Definition \ref{DefT57}

resolving curve: Definition \ref{DefT80}

super parameters: Definition \ref{Def357}

$\tau_f(p)$: Definition \ref{Def221}

$\tau_f(X)$: after Definition \ref{Def221}

$\tau$-quasi-well prepared diagram: after Definitions \ref{Def289} and \ref{Def396}

$\Theta(f,Y)$: Definition \ref{DefT151}

toroidal forms for $u,v,w$: after Definition \ref{Def247}

toroidal forms for $u,v$: Definition \ref{torf}

torodial morphism: Definition \ref{Def247}

transform of a pre-relation: after Definition \ref{Def154}

transform of a relation: after Definition \ref{Def161}

weakly good at $p$ for $f$:  Definition \ref{DefT327}

\vskip.5truein
\noindent
Department of Mathematics

\noindent University of Missouri

\noindent Columbia, MO  65211

\end{document}